\documentclass{memo-l}

\usepackage[mathscr]{euscript}
\usepackage{mathtools}
\usepackage{amssymb}
\usepackage{stackrel}
\usepackage{pb-diagram,pb-xy} 
\usepackage[all,cmtip,arrow]{xy} 
\usepackage{quiver}
\usepackage[bbgreekl]{mathbbol}

\DeclareSymbolFontAlphabet{\mathbb}{AMSb}
\DeclareSymbolFontAlphabet{\mathbbl}{bbold}

\usepackage{hyperref} 
\hypersetup{colorlinks=true,linkcolor=blue,citecolor=magenta}

\makeindex

\numberwithin{equation}{chapter}

\theoremstyle{plain}
\newtheorem{theorem}[equation]{Theorem}
\newtheorem{lemma}[equation]{Lemma}
\newtheorem{corollary}[equation]{Corollary}
\newtheorem{proposition}[equation]{Proposition}

\theoremstyle{definition}
\newtheorem{definition}[equation]{Definition}

\newtheorem{remark}[equation]{Remark}

\newtheorem{example}[equation]{Example}
\newtheorem{examples}[equation]{Examples}
\newtheorem{warning}[equation]{Warning}

\newcommand{\isom}{\cong}
\newcommand{\homeq}{\simeq}
\newcommand{\smsh}{\wedge}
\newcommand{\wdge}{\vee}
\newcommand{\Wdge}{\bigvee}

\newcommand{\N}{\mathbb{N}}
\newcommand{\R}{\mathbb{R}}

\newcommand{\weq}{\; \tilde{\longrightarrow} \;}
\newcommand{\lweq}{\; \tilde{\longleftarrow} \;}
\newcommand{\into}{\hookrightarrow}

\newcommand{\colim}{\operatorname{colim} }
\newcommand{\holim}{\operatorname{holim} }
\newcommand{\hofib}{\operatorname{hofib} }
\newcommand{\hocofib}{\operatorname{hocofib} }
\newcommand{\Fun}{\operatorname{Fun} }
\newcommand{\End}{\operatorname{End} }
\newcommand{\creff}{\operatorname{cr} }
\newcommand{\Hom}{\operatorname{Hom} }
\newcommand{\Jet}{\operatorname{Jet} }
\newcommand{\Exc}{\operatorname{Exc} }
\newcommand{\Alg}{\operatorname{Alg} }

\newcommand{\Diff}{\mathbb{D}\mathrm{iff}}

\newcommand{\un}[1]{\underline{#1}}

\newcommand{\cat}[1]{\mathscr{#1}}
\newcommand{\bcat}[1]{\mathbb{#1}}

\newcommand{\Weil}{\mathbb{W}\mathrm{eil}}
\newcommand{\Weilinfty}{\Weil_\infty}

\newcommand{\uWI}{\un{\Weil}_\infty}
\newcommand{\E}{\mathbb{E}_\infty}
\newcommand{\Mon}{\mathbb{M}\mathrm{on}}
\newcommand{\EMon}{\mathbb{E}_\infty\Mon}
\newcommand{\ERig}{\mathbb{E}_\infty\mathbb{R}\mathrm{ig}}
\newcommand{\Mfld}{\mathbb{M}\mathrm{fld}}
\newcommand{\dMfld}{\mathbb{DM}\mathrm{fld}}
\newcommand{\Sch}{\mathbb{S}\mathrm{ch}}
\newcommand{\CRing}{\mathbb{CR}\mathrm{ing}}
\newcommand{\CRig}{\mathbb{CR}\mathrm{ig}}

\newcommand{\Fin}{\mathsf{Fin}}
\newcommand{\finsset}{\mathsf{Kan}_{\mathrm{fin},*}}

\newcommand{\spectra}{{\mathscr{S}p}}
\newcommand{\finbased}{\mathscr{S}_{\mathrm{fin},*}}
\newcommand{\spaces}{\mathscr{S}}

\newcommand{\bDelta}{\mathbbl{\Delta}}
\newcommand{\msset}{\mathsf{sSet^+}}
\newcommand{\msmon}{\mathsf{sMon^+}}

\newcommand{\sset}{\mathsf{sSet}}
\newcommand{\scset}{\mathsf{sSet^{sc}}}
\newcommand{\Mod}{\mathbf{Mod}}
\newcommand{\Tan}{\mathbf{Tan}}

\newcommand{\ModWeil}{\mathsf{Mod^+_{\Weilinfty}}}
\newcommand{\ModscWeil}{\mathsf{Mod^{sc}_{\Weilinfty}}}
\newcommand{\qCat}{\mathsf{Cat_{\infty}}}
\newcommand{\FPCM}{\mathsf{FPCM}}

\newcommand{\Catinf}{\mathbb{C}\mathrm{at}_\infty}
\newcommand{\bCatinf}{\mathbf{Cat}_\infty}
\newcommand{\Catdiff}{{\mathbb{C}\mathrm{at}_\infty^{\mathrm{diff}}}}
\newcommand{\RelzCatdiff}{{\mathbb{R}\mathrm{el}_0\mathbb{C}\mathrm{at}^{\mathrm{diff}}_{\infty}}}
\newcommand{\RelCatdiff}{\mathbb{R}\mathrm{el}\mathbb{C}\mathrm{at}^{\mathrm{diff}}_{\infty}}
\newcommand{\Catst}{\mathbb{C}\mathrm{at}_\infty^{\mathrm{diff,st}}}
\newcommand{\Catab}{{\mathbb{C}\mathrm{at}^{\mathrm{ab}}}}

\newcommand{\CATdiff}{\mathbb{C}\mathrm{AT}_\infty^{\mathrm{diff}}}
\newcommand{\Catpr}{\mathbb{C}\mathrm{at}_\infty^{\mathrm{pres}}}
\newcommand{\RelzCATdiff}{{\mathbb{R}\mathrm{el}_0\mathbb{C}\mathrm{AT}^{\mathrm{diff}}_{\infty}}}
\newcommand{\ReloCATdiff}{{\mathbb{R}\mathrm{el}_1\mathbb{C}\mathrm{AT}^{\mathrm{diff}}_{\infty}}}
\newcommand{\RelCATdiff}{{\mathbb{R}\mathrm{el}\mathbb{C}\mathrm{AT}_{\infty}^{\mathrm{diff}}}}
\newcommand{\Topinf}{\mathbb{T}\mathrm{opos}_\infty}
\newcommand{\MonCat}{\mathbf{MonCat}_\infty}
\newcommand{\SpanFin}{\mathrm{Span}(\Fin)}

\newcommand{\Exinf}{\mathrm{Ex}^{\infty}}

\begin{document}

\frontmatter

\title{Tangent $\infty$-Categories and {G}oodwillie Calculus}

\author{Kristine Bauer}
\address{Department of Mathematics and Statistics, University of Calgary, 2500 University Drive NW, Calgary, AB, 	T2N 1N4, Canada}
\email{bauerk@ucalgary.ca}
\thanks{The first author was supported by the Natural Sciences and Engineering Research Council of Canada.}

\author{Matthew Burke}
\address{Department of Mathematics and Statistics, University of Calgary, 2500 University Drive NW, Calgary, AB, T2N 1N4, Canada}
\email{matthew.burke@cantab.net}

\author{Michael Ching}
\address{Department of Mathematics, Amherst College, Amherst, MA 01002, United States of America}
\email{mching@amherst.edu}
\thanks{The third author was supported by the National Science Foundation under grant DMS-1709032.}


\begin{abstract}
We make precise the analogy between Goodwillie's calculus of functors in homotopy theory and the differential calculus of smooth manifolds by introducing a higher-categorical framework of which both theories are examples. That framework is an extension to $\infty$-categories of the \emph{tangent categories} of Cockett and Cruttwell, introduced originally by Rosick\'{y}. The basic data of a tangent $\infty$-category consist of an endofunctor that plays the role of the tangent bundle construction, together with various natural transformations that mimic the structure possessed by the ordinary tangent bundles of smooth manifolds.

The role of the tangent bundle functor in Goodwillie calculus is played by Lurie's tangent bundle for $\infty$-categories, introduced to generalize the cotangent complexes of Andr\'{e}, Quillen, and Illusie. We show that Lurie's construction admits the additional structure maps and satisfies the conditions needed to form a tangent $\infty$-category, which we refer to as the \emph{Goodwillie tangent structure}.

Cockett and Cruttwell (and others) have started to develop various aspects of differential geometry in the abstract context of tangent categories, and we begin to apply those ideas to Goodwillie calculus. For example, we show that the role of Euclidean spaces in the calculus of manifolds is played in Goodwillie calculus by the stable $\infty$-categories. We also show that Goodwillie's $n$-excisive functors are the direct analogues of \emph{$n$-jets} of smooth maps between manifolds; to state that connection precisely, we develop a notion of tangent $(\infty,2)$-category and show that Goodwillie calculus is best understood in that context.
\end{abstract}

\maketitle

\setcounter{tocdepth}{0}

\tableofcontents

\mainmatter

\thispagestyle{plain}

\chapter*{Introduction}

Goodwillie's calculus of functors, developed in the series of papers~\cite{goodwillie:1990,goodwillie:1991,goodwillie:2003}, provides a systematic way to apply ideas from ordinary calculus to homotopy theory. For example, central to this functor calculus is the `Taylor tower', an analogue of the Taylor series, which comprises a sequence of `polynomial' approximations to a functor between categories of topological spaces or spectra.

The purpose of this paper is to explore a slightly different analogy, also proposed by Goodwillie, where we view the categories of topological spaces or spectra (or, to be more precise, the corresponding $\infty$-categories) as analogues of smooth manifolds and the functors between those categories as playing the role of smooth maps. Our main goal is to make this analogy precise by introducing a common framework within which both Goodwillie calculus and the calculus of smooth manifolds exist as examples. That framework is the theory of \emph{tangent categories}\footnote{The phrase `tangent category' is unfortunately used in the literature to describe two different objects, both of which feature extensively in this paper. We follow Cockett and Cruttwell's terminology, and use `tangent category' to refer to the broader notion: a category whose objects admit tangent bundles. The notion referred to by, for example, Harpaz-Nuiten-Prasma in the title of~\cite{harpaz/nuiten/prasma:2019} will for us be called the `tangent bundle of an ($\infty$-)category'.} of Cockett and Cruttwell~\cite{cockett/cruttwell:2014}, extended to $\infty$-categories.

A central object in the calculus of manifolds is the tangent bundle construction: associated to each smooth manifold $M$ is another smooth manifold $TM$ along with a projection $p_M: TM \to M$ as well as various other structure maps that, among other things, make $p_M$ into a vector bundle for each $M$. The first step in merging functor calculus with differential geometry is to describe the `tangent bundle' for an $\infty$-category, such as that of spaces or spectra.

Fortunately for us, the analogue of the tangent bundle construction in homotopy theory is well-known. In the generality we want in this paper that construction is given by Lurie in~\cite[7.3.1]{lurie:2017}, but the ideas go back at least to work of Andr\'{e}~\cite{andre:1974} and Quillen~\cite{quillen:1970} on cohomology theories for commutative rings, and Illusie~\cite{illusie:1971} on cotangent complexes in algebraic geometry.

For an object $X$ in an ordinary category $\cat{C}$, we can define the `tangent space' to $\cat{C}$ at $X$ to be the category of abelian group objects in the slice category $\cat{C}_{/X}$ of objects over $X$. These abelian group objects are also called \emph{Beck modules} after their introduction by Beck in a 1967 PhD thesis~\cite{beck:2003}, and they were used by Quillen~\cite{quillen:1967} as the coefficients for the definition of cohomology in an arbitrary model category.

Quillen's approach was refined by Basterra and Mandell~\cite{basterra/mandell:2005}, building on previous work of Basterra~\cite{basterra:1999} on ring spectra, to what is now known as \emph{topological} Quillen homology. In that refinement, the abelian group objects are replaced by cohomology theories, or spectra in the sense of stable homotopy theory. In particular, for an object $X$ in an $\infty$-category $\cat{C}$, the \emph{tangent space} to $\cat{C}$ at $X$:
\begin{equation} \label{eq:TXC} T_X\cat{C} := \spectra(\cat{C}_{/X}) \end{equation}
is the $\infty$-category of spectra in the corresponding slice $\infty$-category of objects over $X$. This $\infty$-category $T_X\cat{C}$ is the `stabilization' of $\cat{C}_{/X}$: its best approximation by an $\infty$-category that is stable.\footnote{An $\infty$-category is \emph{stable} when it admits a null object, it has finite limits and colimits, and when pushout and pullback squares coincide; see~\cite[1.1]{lurie:2017} for an extended introduction. By~\cite[1.4.2.23]{lurie:2017}, any $\infty$-category $\cat{D}$ with finite limits admits a stabilization $\spectra(\cat{D})$ which is universal among stable $\infty$-categories with a finite-limit-preserving functor to $\cat{D}$.}

For example, when $\cat{C}$ is the $\infty$-category of topological spaces, $T_X\cat{C}$ can be identified with the $\infty$-category of spectra parameterized\footnote{See~\cite{may/sigurdsson:2006} for parameterized spectra.} over the topological space $X$. Basterra and Mandell~\cite{basterra/mandell:2005} prove that when $\cat{C}$ is the $\infty$-category of commutative ring spectra, $T_R\cat{C}$ is the $\infty$-category of $R$-module spectra. Various other calculations of these tangent spaces have been done recently in a series of papers by Harpaz, Nuiten and Prasma, including~\cite{harpaz/nuiten/prasma:2019} (for algebras over operads of spectra),~\cite{harpaz/nuiten/prasma:2018} (for $\infty$-categories themselves), and~\cite{harpaz/nuiten/prasma:2019b} (for $(\infty,2)$-categories).

Lurie's construction~\cite[7.3.1.10]{lurie:2017} collects all of these individual tangent spaces together into a `tangent bundle': for each suitably nice $\infty$-category $\cat{C}$, there is a functor
\[ p_{\cat{C}}: T\cat{C} \to \cat{C} \]
whose fibre over $X$ is precisely $T_X\cat{C}$. One of the aims of this paper is to give substance to the claim that this $p_{\cat{C}}$ is analogous to the ordinary tangent bundle $p_M: TM \to M$ for a smooth manifold $M$, by presenting both as examples of the notion of `tangent category'.

\subsection*{Tangent categories}
A categorical framework for the tangent bundle construction on smooth manifolds was first provided by Rosick\'{y} in~\cite{rosicky:1984}. In that work, he describes various structures that can be built on the tangent bundle functor
\[ T: \Mfld \to \Mfld \]
where $\Mfld$ denotes the category of smooth manifolds and smooth maps. Of course, there is a natural transformation $p: T \to \mathrm{Id}$ that provides the tangent bundle projection maps. There are also natural transformations $0: \mathrm{Id} \to T$, given by the zero section for each tangent bundle, and $+: T \times_{\mathrm{Id}} T \to T$, capturing the additive structure of those vector bundles.

Rosick\'{y}'s work was largely unused until resurrected in 2014 by Cockett and Cruttwell~\cite{cockett/cruttwell:2014} in order to describe connections between calculus on manifolds and structures appearing in logic and computer science, such as the differential $\lambda$-calculus of Ehrhard and Regnier~\cite{ehrhard/regnier:2003}. Cockett and Cruttwell define a tangent structure on a category $\bcat{X}$ to consist of an endofunctor $T: \bcat{X} \to \bcat{X}$ together with five natural transformations:
\begin{itemize}
	\item the \emph{projection} $p: T \to \mathrm{Id}$
	\item the \emph{zero section} $0: \mathrm{Id} \to T$
	\item the \emph{addition} $+: T \times_{\mathrm{Id}} T \to T$
	\item the \emph{flip} $c: T^2 \to T^2$
	\item the \emph{vertical lift} $\ell: T \to T^2$
\end{itemize}
for which there is a large collection of diagrams that are required to commute. As indicated above, the first three of these natural transformations make $TM$ into a bundle of commutative monoids over $M$, for any object $M \in \bcat{X}$.\footnote{In Rosick\'{y}'s original work, the addition map was required to have fibrewise negatives so that $TM$ is a bundle of abelian groups. Cockett and Cruttwell relaxed that condition, and we take full advantage of this relaxation since Lurie's tangent bundle does not admit negatives in that sense.} The maps $c$ and $\ell$ express aspects of the `double' tangent bundle $T^2M = T(TM)$ that are inspired by the case of smooth manifolds: its symmetry in the two tangent directions and a canonical way of `lifting' tangent vectors to the double tangent bundle.

In addition to the required commutative diagrams, there is one additional condition, referred to by Cockett and Cruttwell~\cite[2.1]{cockett/cruttwell:2018} as the `universality of the vertical lift'. This axiom states that certain squares of the form
\begin{equation} \label{eq:vertical-lift} \begin{diagram}
		\node{TM \times_M TM} \arrow{s} \arrow{e} \node{T^2M} \arrow{s,r}{T(p)} \\
		\node{M} \arrow{e,t}{0} \node{TM}
\end{diagram} \end{equation}
are required to be pullback diagrams in $\bcat{X}$. (See~\ref{rem:T} for a precise statement.) When $\bcat{X} = \Mfld$ and $M$ is a smooth manifold, this vertical lift axiom can be expressed as a collection of diffeomorphisms
\[ T(T_xM) \isom T_xM \times T_xM \]
varying smoothly with $x \in M$.

Thus, the vertical lift axiom tells us a familiar fact: that the tangent bundle of a tangent space $T_xM$ (indeed of any Euclidean space) is trivial. Cockett and Cruttwell's insight was that this axiom is also the key to translating many of the constructions of ordinary differential geometry into the context of abstract tangent categories. Since the publication of~\cite{cockett/cruttwell:2014}, a small industry has developed around this task, providing, for example:
\begin{itemize}
	\item a \emph{Lie bracket for vector fields} (described by Rosick\'{y} in his original work on tangent categories~\cite{rosicky:1984} and refined in~\cite{cockett/cruttwell:2014});
	\item an analogue for \emph{smooth vector bundles}: the `differential bundles' introduced in~\cite{cockett/cruttwell:2018};
	\item notions of \emph{connection}, \emph{torsion} and \emph{curvature}, developed by Cockett and Cruttwell~\cite{cockett/cruttwell:2017}, also studied by Lucyshyn-Wright~\cite{lucyshyn-wright:2018};
	\item analogues of \emph{affine spaces}, described by Blute, Cruttwell and Lucyshyn-Wright~\cite{blute/cruttwell/lucyshyn-wright:2019};
	\item versions of \emph{differential forms} and \emph{cohomology}, studied by Cruttwell and Lucyshyn-Wright~\cite{cruttwell/lucyshyn-wright:2018}; and
	\item \emph{Lie algebroids}, studied by the second author and MacAdam~\cite{burke/macadam:2019}.
\end{itemize}
In Chapter~\ref{sec:jet} of this paper, we add to this list the notion of `$n$-jet' of smooth maps in order to make precise the connection with Goodwillie's $n$-excisive functors. We do not address any of these other topics here, though we expect each of them has an analogue for tangent $\infty$-categories, which might be worth studying. Note, however, that some of these concepts, such as the Lie bracket, depend on the existence of additive inverses in the tangent bundle, and hence are not applicable in the main example of this paper where those inverses do not exist.

\subsection*{Tangent $\infty$-categories}

Our paper is split into two parts, and the first part focuses on extending the theory of tangent categories to the $\infty$-categorical context. The principal goal of this part is Definition~\ref{def:tangent-structure-infty}, which introduces the notion of tangent structure on an $\infty$-category $\bcat{X}$, and which recovers Cockett and Cruttwell's definition when $\bcat{X}$ is an ordinary category.

A challenge in making this definition is that the large array of commutative diagrams listed by Cockett and Cruttwell in their definition of tangent category would require an even larger array of higher cohering homotopies if translated directly to the $\infty$-categorical framework. Fortunately, work of Leung~\cite{leung:2017} provides an alternative definition of tangent category based on a category of `Weil-algebras' already used in algebraic geometry and synthetic differential geometry.

In Definition~\ref{def:weil1} we describe a symmetric monoidal category $\Weil$ whose objects are certain commutative semirings of the form
\[ \N[x_1,\dots,x_n]/(x_ix_j) \]
where the relations consist of quadratic monomials that include the squares $x_i^2$. The monoidal structure on $\Weil$ is given by the tensor product. These objects form only a subset of the more general Weil-algebras introduced by Weil in~\cite{weil:1953} to study tangent vectors on manifolds in terms of infinitesimals, but they are sufficient to capture the notion of tangent category.

Leung's main result~\cite[14.1]{leung:2017} is that the structure of a tangent category $\bcat{X}$ is precisely captured by a monoidal functor
\begin{equation} \label{eq:weil-tan} T^{\bullet}: \Weil \to \Fun(\bcat{X},\bcat{X}) \end{equation}
from $\Weil$ to the category of endofunctors on $\bcat{X}$ under composition, or equivalently, to an \emph{action} of $\Weil$ on $\bcat{X}$. This action is subject to the additional condition that certain pullbacks be preserved by $T^{\bullet}$; in particular, we have the `universality of vertical lift' axiom described in (\ref{eq:vertical-lift}).

The simplest non-trivial example of a Weil-algebra is the `dual numbers' object $W = \N[x]/(x^2)$ which, under Leung's formulation, corresponds to the ordinary tangent bundle functor $T = T^W: \bcat{X} \to \bcat{X}$. Evaluating the functor $T^{\bullet}$ on morphisms in $\Weil$, i.e.\ the homomorphisms between Weil-algebras, provides the various natural transformations that make up the tangent structure. For example, $\Weil$ encodes the structure that makes each projection map $TM \to M$ into a bundle of commutative monoids over $M$.

In order to generalize tangent structures to $\infty$-categories, we construct in Chapter~\ref{sec:weil} an $\infty$-categorical version of $\Weil$, which we denote $\Weilinfty$. The monoidal $\infty$-category $\Weilinfty$ can be described in a very similar manner to $\Weil$, but based on $E_\infty$-semirings (introduced by May~\cite{may:1972} as `$E_\infty$-ring spaces') instead of ordinary commutative semirings. An $E_\infty$-semiring is a topological space with two binary operations, addition and multiplication, that satisfy the usual axioms of a semiring, but only up to homotopy. There are also higher coherence conditions which have no analogue in the discrete case.

Since $E_\infty$-semirings are complicated topological objects, we give an alternative construction of the $\infty$-category $\Weilinfty$ as the nerve of a certain bicategory in which all the $2$-morphisms are invertible. This construction highlights the close relationship between $\Weil$ and $\Weilinfty$, and makes it easy to see, for example, that the homotopy category of $\Weilinfty$ is isomorphic to $\Weil$.\footnote{In the original version of this paper, we used the $1$-category $\Weil$ in our definition of tangent $\infty$-structure. We are grateful to Thomas Nikolaus and Maxime Ramzi for pointing out that, with this definition, our construction of the Goodwillie tangent structure entails strictly commutative structures on the $\infty$-category of spectra which cannot exist. There was a mistake in our original construction; specifically, the natural transformation $\alpha$ in the previous version of Definition~\ref{def:alpha} was not sufficiently well-defined. That observation led us to the introduction of the monoidal $\infty$-category $\Weilinfty$ defined in Chapter~\ref{sec:weil}.}

Our Definition~\ref{def:tangent-structure-infty} of tangent structure on an $\infty$-category $\bcat{X}$ then follows the same format as Leung's characterization of tangent categories; it is a functor of monoidal $\infty$-categories of the form (\ref{eq:weil-tan}), but with $\Weil$ replaced by $\Weilinfty$, and that preserves those same pullbacks (though now with pullback understood in the $\infty$-categorical sense). Under this definition, any ordinary tangent category, such as $\Mfld$ with its usual tangent bundle construction, is also a tangent $\infty$-category with structure map given by composing the relevant map (\ref{eq:weil-tan}) with the projection $\Weilinfty \to \Weil$ from $\Weilinfty$ to (the nerve of) its homotopy category.

There are, however, tangent $\infty$-categories that do not arise from an ordinary tangent category. One example is given by the `derived smooth manifolds' of~\cite{spivak:2010}. There, Spivak defines an $\infty$-category $\dMfld$ that contains $\Mfld$ as a full subcategory, but which admits all pullbacks, not only those along transverse pairs of smooth maps. We show in Proposition~\ref{prop:dMfld} that the tangent structure on $\Mfld$ extends naturally to $\dMfld$ using a universal property described by Carchedi and Steffens~\cite{carchedi/steffens:2019}.

Another example is the $\infty$-category $\E\spectra$ of $E_\infty$-ring spectra, a central concept in stable homotopy theory. Just as the category $\CRing$ of ordinary commutative rings admits a tangent structure (in fact, two tangent structures, one on $\CRing$ and one on $\CRing^{op}$), we show in Examples~\ref{ex:ring-spec} and~\ref{ex:ring-spec-op} that $\E\spectra$ and $\E\spectra^{op}$ are tangent $\infty$-categories in an analogous way.

\subsection*{Tangent $(\infty,2)$-categories and other objects}

The definition of tangent structure on an $\infty$-category can easily be extended to a wide range of other types of objects, and we examine this generalization in Chapter~\ref{sec:infty2}. Let $\bcat{X}$ be an object in an $(\infty,2)$-category $\mathbf{C}$ which, for the purposes of this introduction, one may view simply as a category enriched in $\infty$-categories. Then $\bcat{X}$ admits a monoidal $\infty$-category $\Hom_{\mathbf{C}}(\bcat{X},\bcat{X})$ of endomorphisms. We define (in~\ref{def:tangent-object}) a tangent structure on the object $\bcat{X}$ to be a monoidal functor
\[ T: \Weilinfty \to \Hom_{\mathbf{C}}(\bcat{X},\bcat{X}) \]
that preserves the appropriate pullbacks. This definition recovers our notion of tangent $\infty$-category (and hence ordinary tangent categories) when $\mathbf{C}$ is the $(\infty,2)$-category of $\infty$-categories. Taking $\mathbf{C}$ to be a suitable $(\infty,2)$-category of $(\infty,2)$-categories, we also obtain a notion of tangent $(\infty,2)$-category. In particular, our definition entails a notion of tangent bicategory, in which the Cockett-Cruttwell axioms hold up to coherent $2$-isomorphisms.

\subsection*{The Goodwillie tangent structure on differentiable $\infty$-categories}

In the second part of this paper, we construct a specific tangent $\infty$-category for which the tangent bundle functor is equivalent to that defined by Lurie, and which encodes the theory of Goodwillie calculus. The existence of this tangent structure (which we refer to as the \emph{Goodwillie tangent structure}) justifies the analogy between functor calculus and the calculus of smooth manifolds.

The underlying $\infty$-category for this tangent structure is $\Catdiff$, a subcategory of Lurie's $\infty$-category of $\infty$-categories~\cite[3.0.0.1]{lurie:2009}. The objects in $\Catdiff$ are those $\infty$-categories $\cat{C}$ that are \emph{differentiable} in the sense introduced in~\cite[6.1.1.6]{lurie:2017}: those $\cat{C}$ that admit finite limits and sequential colimits, which commute. This condition is satisfied by many $\infty$-categories of interest, including any compactly generated $\infty$-category and any $\infty$-topos. The morphisms in the $\infty$-category $\Catdiff$ are those functors between differentiable $\infty$-categories that preserve sequential colimits.

The tangent bundle construction from~\cite[7.3.1.10]{lurie:2017} can be described explicitly as the functor $T: \Catdiff \to \Catdiff$ given by
\[ T\cat{C} := \Exc(\finbased,\cat{C}) \]
where $\finbased$ denotes the $\infty$-category of finite pointed spaces, and $\Exc(-,-)$ denotes the $\infty$-category of functors that are \emph{excisive} in the sense of Goodwillie (those that map pushouts in $\finbased$ to pullbacks in $\cat{C}$). With this definition in mind, the Goodwillie tangent structure on $\Catdiff$ consists of the following natural transformations:
\begin{itemize}
	\item the \emph{projection} map $p: T\cat{C} \to \cat{C}$ is evaluation at the null object:
	\[ L \mapsto L(*); \]
	\item the \emph{zero section} $0: \cat{C} \to T\cat{C}$ maps an object $X$ of $\cat{C}$ to the constant functor with value $X$;
	\item \emph{addition} $+: T\cat{C} \times_{\cat{C}} T\cat{C} \to T\cat{C}$ is the fibrewise product:
	\[ (L_1,L_2) \mapsto L_1(-) \times_{L_1(*) = L_2(*)} L_2(-); \]
	\item identifying $T^2\cat{C}$ with the $\infty$-category of functors $\finbased \times \finbased \to \cat{C}$ that are excisive in each variable individually, the \emph{flip} $c: T^2\cat{C} \to T^2\cat{C}$ is the symmetry in those two variables:
	\[ L \mapsto [(X,Y) \mapsto L(Y,X)]; \]
	\item the \emph{vertical lift} map $\ell: T\cat{C} \to T^2\cat{C}$ is precomposition with the smash product:
	\[ L \mapsto [(X,Y) \mapsto L(X \smsh Y)]. \]
\end{itemize}
A complete definition of the Goodwillie tangent structure requires the construction of a monoidal functor
\[ T: \Weilinfty \to \Fun(\Catdiff,\Catdiff). \]
For a Weil-algebra of the form
\[ A = \N[x_1,\dots,x_n]/(x_ix_j), \]
with $n$ generators, we define the corresponding endofunctor $T^A: \Catdiff \to \Catdiff$ by
\[ T^A(\cat{C}) := \Exc^A(\finbased^n,\cat{C}) \subseteq \Fun(\finbased^n,\cat{C}), \]
the full subcategory consisting of those functors $\finbased^n \to \cat{C}$ that satisfy the property of being `$A$-excisive' (see Definition~\ref{def:A-excisive}). When $A = \N[x]/(x^2)$, $A$-excisive is excisive, and we recover Lurie's definition of the tangent bundle $T\cat{C}$. We define the functor $T$ on a morphism $\phi$ in $\Weilinfty$ via precomposition with a certain functor
\[ \tilde{\phi}: \finbased^{n'} \to \finbased^{n} \]
whose construction mimics the algebra homomorphism underlying $\phi$; see Definition~\ref{def:phi} for details.

In Chapter~\ref{sec:underlying} we prove the key homotopy-theoretic results needed to check that the definitions above do indeed form a tangent structure. Principal among those is the vertical lift axiom analogous to that in (\ref{eq:vertical-lift}) above; this axiom is verified in Proposition~\ref{prop:vertical-lift}, and involves in detail the classification of multilinear functors in Goodwillie calculus~\cite[Sec.\ 5]{goodwillie:2003} and splitting results for functors with values in a stable $\infty$-category. It is there that the technical heart of the construction of the Goodwillie tangent structure lies.

The full definition of $T$ as a monoidal functor between monoidal $\infty$-categories is rather involved and relies on a model for $\Catdiff$ based on `relative' $\infty$-categories; see~\cite{mazel-gee:2019}. The specifics of this construction are in Chapter~\ref{sec:tangent-structure-catdiff}. The reader interested in understanding the central ideas of our construction, rather than the intricate details, should focus on Chapter~\ref{sec:underlying} where the most important of the underlying definitions are provided.

\subsection*{Tangent functors, differential objects and jets}

Much of this paper is concerned with the definition of tangent $\infty$-category and the construction of the Goodwillie tangent structure on the $\infty$-category $\Catdiff$. However, we also begin the task of extending the current tangent category literature to the $\infty$-categorical context. In this paper, we focus on three specific aspects of that theory: functors between tangent categories, differential objects, and jets.

As with any categorical structure, it is crucial to identify the appropriate morphisms. Cockett and Cruttwell introduced in~\cite[2.7]{cockett/cruttwell:2014} `morphisms of tangent structure', and in Chapter~\ref{sec:tangent-functor} we extend the strong version of their notion to tangent $\infty$-categories. Briefly, a tangent functor between tangent $\infty$-categories is a functor that commutes (up to higher coherent equivalences) with the corresponding actions of $\Weilinfty$. Our definition, using work of Garner~\cite{garner:2018}, reduces to Cockett and Cruttwell's definition in the case of ordinary tangent categories.

Differential objects were introduced by Cockett and Cruttwell~\cite[4.8]{cockett/cruttwell:2014} in order to axiomatize the role of Euclidean spaces in the theory of smooth manifolds: these are the objects that are tangent \emph{spaces}. In Chapter~\ref{sec:differential} we describe our extension of the theory of differential objects to tangent $\infty$-categories. Our description provides a new perspective on differential objects even in the setting of ordinary tangent categories; see Proposition~\ref{prop:diff}.

Another role for differential objects is in making the connection between tangent categories and the `cartesian differential categories' of Blute, Cockett and Seely~\cite{blute/cockett/seely:2009}. Roughly speaking, a cartesian differential category is a tangent category in which every object has a canonical differential structure; for example, the subcategory of $\Mfld$ whose objects are the Euclidean spaces $\R^n$. We show that a similar relationship (Theorem~\ref{thm:cartesian-differential-structure}) also holds in the tangent $\infty$-category setting after passage to homotopy categories.

In Chapter~\ref{sec:catdiff-diff} we analyse the notion of differential object in the specific context of the Goodwillie tangent structure on $\Catdiff$. It is not hard to see that the differential objects in $\Catdiff$ are precisely the \emph{stable} $\infty$-categories. This result is not at all surprising; the role of stabilization is built into our tangent structure via the tangent spaces described in (\ref{eq:TXC}). It does, however, confirm Goodwillie's intuition that categories of spectra should be viewed, from the point of view of functor calculus, as analogues of Euclidean spaces.

We also deduce the existence of a cartesian differential category whose objects are the stable $\infty$-categories and whose morphisms are natural equivalence classes of functors. This result extends work of the first author and Johnson, Osborne, Riehl, and Tebbe~\cite{bauer/johnson/osborne/riehl/tebbe:2018}, which describes a similar construction for (chain complexes of) abelian categories in the context of Johnson and McCarthy's `abelian functor calculus' variant of Goodwillie's theory~\cite{johnson/mccarthy:2004}. In fact, the paper~\cite{bauer/johnson/osborne/riehl/tebbe:2018} provided much of the inspiration for our development of tangent $\infty$-categories and for the construction of the Goodwillie tangent structure.

In Chapter~\ref{sec:jet} we turn to the notion of `$n$-jet' of a morphism which does not appear explicitly in the tangent category literature, though it has been studied in the context of synthetic differential geometry, e.g.\ see~\cite[2.7]{kock:2010}. In our case, the importance of $n$-jets is that they correspond in the Goodwillie tangent structure on $\Catdiff$ to the $n$-excisive functors, i.e.\ Goodwillie's analogues of degree $n$ polynomials.

In an arbitrary tangent $\infty$-category $\bcat{X}$, we say that two morphisms $F,G: \cat{C} \to \cat{D}$ determine the same $n$-jet at a (generalized) point $x \in \cat{C}$ if they induce equivalent maps on the $n$-fold tangent spaces $T^n_x\cat{C}$ at $x$. For smooth manifolds, this definition recovers the ordinary notion of $n$-jet, the equivalence class of smooth maps that agree to order $n$ in a neighbourhood of the point $x$. For $\infty$-categories, we then prove an analogous result (Theorem~\ref{thm:jet}) which says that a natural transformation $\alpha: F \to G$ between two functors $\cat{C} \to \cat{D}$ induces an equivalence $P_n^XF \weq P_n^XG$ between Goodwillie's $n$-excisive approximations at $X$ in $\cat{C}$ if and only if $\alpha$ induces an equivalence on $T^n_X\cat{C}$.

The significance of the previous result is that it shows that the notion of $n$-excisive equivalence, and hence $n$-excisive functor, can be recovered directly from the Goodwillie tangent structure. To make proper sense of this claim, though, we need to be able to talk about non-invertible natural transformations between functors of $\infty$-categories. This observation reveals that Goodwillie calculus is better understood in the context of tangent structures on an $(\infty,2)$-category.

In Theorem~\ref{thm:CATdiff-T} we show that there is an $(\infty,2)$-category $\CATdiff$ of differentiable $\infty$-categories which admits a Goodwillie tangent structure extending that on $\Catdiff$. This tangent $(\infty,2)$-category completely encodes Goodwillie calculus and the notion of Taylor tower.

\subsection*{Conjectures and connections} \label{intro:con}

Lurie introduced the tangent bundle $T\cat{C}$ not directly in relation to Goodwillie calculus but as part of the development of deformation theory in an $\infty$-category $\cat{C}$; see~\cite[7.4]{lurie:2017}. That theory is controlled by the cotangent complex functor, a certain section
\[ L: \cat{C} \to T\cat{C} \]
of the tangent bundle projection map; see~\cite[7.3.2.14]{lurie:2017}. The functor $L$ does not appear to have an analogue in the general theory of tangent categories, and it has not played any role in this paper. Nonetheless, it would be interesting to explore what aspects of the Goodwillie tangent structure allow for the cotangent complex and corresponding deformation theory to be developed.

There are other topics and questions regarding the Goodwillie tangent structure which we would like to have addressed in this paper. One such question is the extent to which the Goodwillie tangent structure on $\Catdiff$ is unique. We believe that indeed it is the unique (up to contractible choice) tangent $\infty$-category which extends Lurie's tangent bundle construction, but we do not try to give a proof of that conjecture here.

Another topic concerns what are known as `representable' tangent categories. It was observed by Rosick\'{y} that a model for synthetic differential geometry gives rise to a tangent category in which the tangent bundle functor is represented by an object with so-called `infinitesimal' structure; see~\cite[5.6]{cockett/cruttwell:2014}. It appears relatively easy to extend the definition of representable tangent structure to $\infty$-categories, but also to prove that the Goodwillie tangent structure is \emph{not} representable in this sense. In fact, Lurie's tangent bundle functor itself is not representable.

However, if we restrict the Goodwillie tangent structure to the subcategory of $\Catdiff$ consisting of $\infty$-toposes, e.g.\ see~\cite[6.1.0.4]{lurie:2009}, and the left exact colimit-preserving functors, then this restricted tangent structure is \emph{adjoint}, in the sense of~\cite[5.17]{cockett/cruttwell:2014}, to a representable tangent structure on $\Topinf$ (the $\infty$-category of $\infty$-toposes and geometric morphisms). The representing object for that tangent structure is the $\infty$-topos $T(\spaces) = \Exc(\finbased,\spaces)$, whose objects are parameterized spectra over arbitrary topological spaces, with infinitesimal structure arising directly from the Goodwillie tangent structure. Details of these claims are worked out by the third author in~\cite{ching:2021}. We wonder if such structure is related to work of Anel, Biedermann, Finster and Joyal~\cite{anel/biedermann/finster/joyal:2018},~\cite{anel/biedermann/finster/joyal:2023} on Goodwillie calculus for $\infty$-toposes.

Earlier in this introduction, we listed various topics from ordinary differential geometry that have been developed in the abstract setting of tangent categories, including vector bundles, connections, and curvature. For each of these topics, or any others that can be formulated in an abstract tangent category, we can speculate on what form they take in the Goodwillie tangent structure. Some ideas along those lines appear at the end of this paper in a section on \hyperref[sec:future]{Proposals for Future Work}.

There are two other directions for generalization that seem particularly worthy of exploration. One is to replace $\Catdiff$ in the Goodwillie tangent structure with a different $\infty$-cosmos in the sense of Riehl and Verity~\cite[1.2]{riehl/verity:2022}. The other is to look for tangent $\infty$-categories that encode other versions of functor calculus, such as the `manifold calculus' of Goodwillie and Weiss~\cite{weiss:1999,goodwillie/weiss:1999}, or the `orthogonal calculus' of Weiss~\cite{weiss:1995}. Other lines of inquiry arise from finding counterparts in tangent categories of concepts that are already established in Goodwillie calculus, such as Heuts's work~\cite{heuts:2021} on Goodwillie towers for (pointed compactly-generated) $\infty$-\emph{categories}, instead of functors, and work of Greg Arone and the third author~\cite{arone/ching:2011} on chain rules and the role of operads in functor calculus.

Finally, we might speculate on how the Goodwillie tangent structure of this paper fits into the much bigger programme of `higher differential geometry' developed by Schreiber~\cite[4.1]{schreiber:2013}, or into the framework of homotopy type theory~\cite{hott:2013}, though we don't have anything concrete to say about these possible connections.

\subsection*{Background on $\infty$-categories; notation and conventions}

This paper is written largely in the language of $\infty$-categories, as developed by Lurie in the books~\cite{lurie:2009} and~\cite{lurie:2017}, based on the quasi-categories of Boardman and Vogt~\cite{boardman/vogt:1973}. However, for much of the paper, specifics of that theory are not particularly important, and other models for $(\infty,1)$-categories could easily be used instead, especially if the reader is interested only in the main ideas of this work rather than the technical details. We expect that most of the mathematics in this paper could be written in the framework of `synthetic' higher categories being developed by Cisinski, Cnossen, Nguyen, and Walde~\cite{cisinski/cnossen/nguyen/walde:2026}, although that work does not address $\infty$-bicategories. We include in Appendices~\ref{sec:infty} and~\ref{sec:monoidal} some background information about $\infty$-categories and $\infty$-bicategories, as well as the notions of monoidal $\infty$-category and modules over such, which are central to our approach. 

All the basic concepts needed to define tangent $\infty$-categories, and to describe the underlying data of the Goodwillie tangent structure, will be familiar to readers versed in ordinary category theory. Those ideas include limits and colimits (in particular, pullbacks and pushouts), adjunctions, monoidal structures, and functor categories. To follow the details of our constructions, however, the reader will need a close acquaintance with simplicial sets and, especially, their relationship to categories via the nerve construction.

Throughout this paper, we take the perspective that categories, $\infty$-categories, and indeed $(\infty,2)$-categories, are all really the same sort of thing. We often do not distinguish notationally between any of these types of object unless we want to highlight the different roles they play. In particular, we often identify a category $\bcat{X}$ with its nerve, an $\infty$-category.

One exception to this convention is the category of simplicial sets itself, which we denote as $\sset$, and some related categories (such as `marked' or `scaled' simplicial sets) which we describe in the course of this paper. We will make some use of the various model structures on these categories, including the Quillen and Joyal model structures on $\sset$. There are $\infty$-categories associated to these model categories, but we will use separate notation to denote those $\infty$-categories when we need to invoke them.

One of the most confusing aspects of our theory is that $\infty$-categories (and $(\infty,2)$-categories) play several different roles in this paper at different `levels':
\begin{enumerate} \itemsep=5pt
	\item we can define \textbf{tangent objects in} any (potentially very large) $(\infty,2)$-category $\mathbf{C}$ (Definition~\ref{def:tangent-object});
	\item taking $\mathbf{C}$ to be the $(\infty,2)$-category of (large) $\infty$-categories $\mathbf{CAT}_\infty$, we get a notion of \textbf{tangent structure on} a specific $\infty$-category $\bcat{C} \in \mathbf{C}$ (Definition~\ref{def:tangent-structure-infty});
	\item taking $\bcat{C}$ to be an $\infty$-category of (smaller) $\infty$-categories, such as $\Catdiff$, we get the \textbf{tangent bundle} $T\cat{C}$ of a specific $\infty$-category $\cat{C} \in \bcat{C}$ (Definition~\ref{def:tan-cat}).
\end{enumerate}
We can also take $\cat{C}$ in (3) to be an $\infty$-category $\mathscr{C}\mathrm{at}_\infty$ of (even smaller) $\infty$-categories, in which case we would also have the \textbf{tangent space at} an $\infty$-category $Y$, i.e.\ $T_{Y}\cat{C}$. These tangent spaces are studied in~\cite{harpaz/nuiten/prasma:2018} but do not play any particular role for us here.

We distinguish between these three uses for $\infty$- and $(\infty,2)$-categories by applying the fonts $\mathbf{C},\bcat{C},\cat{C}$\index{c@$\cat{C},\bcat{C},\mathbf{C}$, notation for $\infty$-(bi)categories of various sizes} as indicated in the list above. In particular, we use these different fonts to signify the size restrictions that are implicit in our hierarchy. To be precise, we assume, where necessary, the existence of inaccessible cardinals that determine the different `sizes' of the $\infty$-categories described in (1), (2) and (3) above; see also~\cite[1.2.15]{lurie:2009} for a discussion of this foundational issue. Beyond the requirement to keep these three levels separate, size issues do not play a significant role in this paper.

One final (and important) comment on notation, especially for readers familiar with the papers of Cockett and Cruttwell: in this paper, we use `algebraic' order for writing composition, as opposed to the diagrammatic order employed in many papers in the tangent category literature. So, for morphisms $f:A \to B$ and $g:B \to C$, we write $gf$ for the composite morphism $A \to C$. Because of this choice, some of the expressions we use look different to those appearing in a corresponding place in~\cite{cockett/cruttwell:2014}.

\subsection*{Acknowledgements}

The initial research on this project was carried out during a visit by the third author to the University of Calgary, supported by the Pacific Institute for the Mathematical Sciences. Conversations with many people contributed to this paper, including Thomas Blom, James Cranch, Robin Cockett, Geoffrey Cruttwell, Jonathan Gallagher, Thomas Goodwillie, Marcello Lanfranchi, JS Lemay, Rory Lucyshyn-Wright, Jacob Lurie, Benjamin MacAdam, and Florian Schwarz. We are grateful to Thomas Nikolaus and Maxime Ramzi for pointing out a significant error in the first version of this work. We would also like to thank an anonymous referee for comments and suggestions which improved the paper. This project was largely inspired by work of the first author joint with Brenda Johnson, Christina Osborne, Emily Riehl and Amelia Tebbe~\cite{bauer/johnson/osborne/riehl/tebbe:2018} which was undertaken at the Banff International Research Station during the Women in Topology programme.

\part{Tangent Structures on $\infty$-Categories}

The goal of Part 1 is to extend the notion of tangent category of Cockett and Cruttwell to an $\infty$-categorical context, and to begin the study of these `tangent $\infty$-categories'.

We refer the reader to the paper of Cockett and Cruttwell~\cite{cockett/cruttwell:2014} for a detailed introduction to the theory of tangent categories. The definition that we use in this paper, however, is based on an alternative characterization of tangent structure due to Leung~\cite[14.1]{leung:2017}, and it extends to tangent functors and natural transformations via work of Garner~\cite[Thm. 9]{garner:2018}. That characterization is in terms of a certain monoidal category of `Weil-algebras', and in Chapter~\ref{sec:tancat} we recall that category and the corresponding notion of tangent structure on a category.

In Chapter~\ref{sec:weil} we turn to $\infty$-categories. Our first goal is to define an $\infty$-categorical version of the monoidal category of Weil-algebras, which we then use to give a definition of tangent $\infty$-category in Chapter~\ref{sec:tangent-infty-category}. We start to develop the theory of tangent $\infty$-categories in Chapter~\ref{sec:tangent-functor} with a definition of tangent functor, the appropriate notion of morphism between tangent $\infty$-categories. In Chapter~\ref{sec:differential} we consider `differential objects' in a tangent $\infty$-category.

In the last topic in this part, Chapter~\ref{sec:infty2}, we extend the definition of tangent structure to $(\infty,2)$-categories and other kinds of objects.

\chapter{Weil-Algebras and Tangent Categories} \label{sec:tancat}

We start by recalling the category of Weil-algebras used in Leung's characterization of tangent structures~\cite{leung:2017}.

\begin{definition} \label{def:crig}
Let $\CRig$\index{CRig@$\CRig$, category of commutative semirings} denote the category of commutative semirings (or rigs). The set $\N$\index{N@$\N$, the set of natural numbers including $0$} of natural numbers, including $0$, is an initial object in $\CRig$. There is a monoidal structure on $\CRig$ given by the tensor product $\otimes$ with unit object $\N$. \end{definition}

\begin{example} \label{ex:weil-prod}
For each finite set $J$, let $W^J$\index{WJ@$W^J$, the trivial commutative semiring generated by a finite set $J$} denote the commutative semiring whose elements are formal sums
\[ m + \sum_{j \in J} m_j x_j \]
where $m,m_j \in \N$ and $x_j$ is a formal generator corresponding to $j \in J$, with multiplication determined by the relation $x_j x_{j'}= 0$ for all $j, j' \in J$. When $J = \{1,\dots,n\}$ (or some other ordered set of size $n$), we can write\index{Wn@$W^n$, the trivial commutative semiring on $n$ generators}
\[ W^J \isom W^n := \N[x_1,\dots,x_n]/(x_ix_j)_{i,j = 1}^{n}. \]
We often denote $W^1 = \N[x]/(x^2)$ just as $W$\index{W@$W = W^1 = \N[x]/(x^2)$, the Weil-algebra on $1$ generator}, and when $J = \varnothing$ is the empty set, we have $W^{\varnothing} = \N$.
\end{example}

\begin{definition} \label{def:weil1}
Let $\Weil$\index{Weil@$\Weil$, the category of Weil-algebras} be the full subcategory of $\CRig$ whose objects are those of the form
\[ W^{J_1} \otimes \dots \otimes W^{J_r} \]
for some finite sequence of nonempty finite sets $J_1,\dots,J_r$. We include the case of the empty sequence, in which case this tensor product is equal to $\N$, the unit for the monoidal structure. We refer to these objects as \emph{Weil-algebras}\index{Weil-algebra}, and to the morphisms in $\Weil$ as \emph{morphisms of Weil-algebras}\index{Weil-algebra!morphism}.

The subcategory $\Weil$ is also a monoidal category under the tensor product, and it is (monoidally) equivalent to that denoted `$\N{-}\mathbf{Weil_1}$' in~\cite{leung:2017}.
\end{definition}

\begin{remark} \label{rem:gen}
Let $A = W^{J_1} \otimes \dots \otimes W^{J_r}$ be an object of $\Weil$. If we write $n_i = |J_i|$ and $n = n_1+\dots+n_r$, then we can denote the generators of $A$ (as a semiring) by
\[ x_1,\dots,x_n. \]
Note that this notation implicitly involves choosing an ordering on each finite set $J_i$. For small numbers of generators, we might use $x,y,z$, or, if we are discussing multiple Weil-algebras, we might use $y_1,\dots,y_m$ to make it clear which object the generators refer to. 

We can then also write
\[ A = \N[x_1,\dots,x_n]/(x_ix_j \; | \; i \sim j) \]
where $\sim$ is an equivalence relation on $\{1,\dots,n\}$ whose equivalence classes are in bijection with the sets $J_1,\dots,J_r$. The additive structure on the Weil-algebra $A$ has a basis consisting of the nonzero \emph{monomials}, i.e.\ products of generators, no two from the same set $J_i$.

We also note that any object $A$ of $\Weil$ is an augmented commutative semiring, with augmentation
\[ \epsilon: A \to \N \]
given by $\epsilon(j) = 0$ for all $j \in J_1,\dots,J_r$.
\end{remark}

There are certain pullback squares in $\Weil$ that play a crucial role in the definition of tangent structure.

\begin{lemma}[{\cite[3.14]{leung:2017}}] \label{lem:pb1}
For finite sets $J,J'$, let $J \sqcup J'$ denote their disjoint union. Then there is a pullback square in $\Weil$ of the form
\[ \begin{diagram}
		\node{W^{J \sqcup J'}} \arrow{e} \arrow{s} \node{W^J} \arrow{s} \\
		\node{W^{J'}} \arrow{e} \node{\N}
\end{diagram} \]
in which the horizontal and vertical maps are given by $x_{j'} \mapsto 0$ (for $j' \in J'$) and $x_j \mapsto 0$ (for $j \in J$) respectively. We refer to these diagrams as the \index{foundational pullback}\emph{foundational} pullbacks in $\Weil$.
\end{lemma}

\begin{lemma} \label{lem:pb2}
There is a pullback square in $\Weil$ of the form
\[ \begin{diagram}
		\node{W^2} \arrow{e,t}{\mu} \arrow{s,l}{\epsilon} \node{W \otimes W} \arrow{s,r}{1_W \otimes \epsilon} \\
		\node{\N} \arrow{e,t}{\eta} \node{W}
\end{diagram} \]
where $\mu: \N[x,y]/(x^2,xy,y^2) \to \N[a,b]/(a^2,b^2)$ is given by
\[ \mu(x) = ab, \quad \mu(y) = b. \]
We refer to this square as the \index{vertical lift pullback}\emph{vertical lift} pullback in $\Weil$. Collectively, the pullback diagrams in Lemmas~\ref{lem:pb1} and~\ref{lem:pb2} are the \index{tangent pullback}\emph{tangent pullbacks}.
\end{lemma}
\begin{proof}
A cone over the maps $1_W \otimes \epsilon$ and $\eta$ consists of a morphism
\[ \phi: \N[z_1,\dots,z_k]/(z_i z_j \; | \; i \sim j) \to \N[a,b]/(a^2,b^2) \]
such that each $\phi(z_i)$ is a sum of monomials $ab$ and $b$ (but not $a$). The corresponding lift
\[ \tilde{\phi}: \N[z_1,\dots,z_k]/(z_i z_j \; | \; i \sim j) \to \N[x,y]/(x^2,xy,y^2) \]
is given by replacing $ab$ with $x$ and $b$ with $y$ in the formula for each $\phi(z_i)$.
\end{proof}

\subsection*{Tangent categories}

Having established the monoidal category $\Weil$, we can now give a definition of tangent structure on a category $\bcat{X}$. Leung showed in~\cite[14.1]{leung:2017} that Cockett and Cruttwell's original definition in~\cite{cockett/cruttwell:2014} can be interpreted as an \emph{action} on $\bcat{X}$ by the monoidal category $\Weil$ in the following sense.

\begin{definition} \label{def:actegory}
Let $\bcat{X}$ be a category. A \emph{$\Weil$-module structure}\index{Weil-module@$\Weil$-module} on $\bcat{X}$ is a (strong) monoidal functor
\[ T^\bullet: \Weil \to \Fun(\bcat{X},\bcat{X}) \]
where the monoidal structure on $\Fun(\bcat{X},\bcat{X})$ is given by composition of functors. Corresponding to $T^\bullet$ is an \emph{action map} which we usually denote
\[ \odot: \Weil \times \bcat{X} \to \bcat{X}. \]
The structure encoded by $\odot$ is also termed a `$\Weil$-actegory'; see~\cite{capucci/gavranovic:2022}.
\end{definition}

\begin{definition}[Tangent categories] \label{def:tangent-category}
Let $\bcat{X}$ be a category. A \emph{tangent structure}\index{tangent structure!on a category} on $\bcat{X}$ is a $\Weil$-module structure
\[ T^\bullet: \Weil \to \Fun(\bcat{X},\bcat{X}) \]
such that for each $X \in \bcat{X}$, the corresponding action map
\[ - \odot X : \Weil \to \bcat{X} \]
sends the pullback squares in Lemmas~\ref{lem:pb1} and~\ref{lem:pb2} to pullback squares in $\bcat{X}$.

A \emph{tangent category}\index{tangent category} consists of a category $\bcat{X}$ together with a tangent structure on $\bcat{X}$.
\end{definition}

\begin{remark} \label{rem:pointwise-limit}
The condition that each action map $- \odot X$ preserves tangent pullbacks implies that the functor $T^\bullet$ as a whole preserves tangent pullbacks: if a diagram in a functor category such as $\Fun(\bcat{X},\bcat{X})$ has pointwise limits for each indexing object, then the diagram has a limit in the functor category which agrees with each of the pointwise limits. However, the converse does not hold: it is possible for a functor category to have a limit which is not given by pointwise limits (when some of those pointwise limits do not exist). Therefore, the pointwise condition is required in Definition~\ref{def:tangent-category} in order that it matches up with Cockett and Cruttwell's original notion.
\end{remark}

\begin{remark}
For each Weil-algebra $A$, we write\index{TA@$T^A$, the action of a Weil-algebra $A$ in a tangent structure}
\[ T^A := T^\bullet(A) = A \odot - : \bcat{X} \to \bcat{X} \]
for the corresponding endofunctor on the underlying category $\bcat{X}$. In particular, when $A = W$, the endofunctor $T^W$, which we often denote by $T$, is the \emph{tangent bundle functor}\index{tangent bundle functor} for the tangent structure. In this formulation, the monoidal structure on the functor $T^\bullet$ is given by isomorphisms
\[ \alpha: T^{A \otimes A'} \isom T^AT^{A'}, \quad \lambda: T^{\N} \isom 1_{\bcat{X}} \]
which satisfy appropriate coherence conditions.
\end{remark}

We now recall how Definition~\ref{def:tangent-category} reduces to Cockett and Cruttwell's original definition of tangent category~\cite[2.1]{cockett/cruttwell:2018}.

\begin{remark} \label{rem:T}
Let $(\bcat{X},T^\bullet)$ be a tangent category. The functors $T^A: \bcat{X} \to \bcat{X}$ for Weil-algebras $A$ are determined by the single functor $T = T^W: \bcat{X} \to \bcat{X}$ in the following way:
\begin{itemize}
	\item since $T^\bullet$ is monoidal, we have, for the unit object $\N$ of the monoidal structure on $\Weil$, an isomorphism
	\[ T^\N \isom 1_{\bcat{X}} \]
	with the identity functor on $\bcat{X}$;
	\item since the module action map $\odot$ is required to preserve the foundational pullbacks, we have, for any positive integer $n$ and any $X \in \bcat{X}$:
	\[ T_n(X) := T^{W^n}(X) \isom T(X) \times_{X} \dots \times_{X} T(X); \]
	the wide pullback of $n$ copies of the `projection' map $p_X: T(X) \to X$ corresponding to the unique map $\epsilon: W \to \N$; for an arbitrary finite set $J$, we also write $T_J := T^{W^J}$;
	\item since any Weil-algebra $A \in \Weil$ is of the form $W^{J_1} \otimes \dots \otimes W^{J_r}$, the monoidal condition then implies
	\[ T^A \isom T_{J_1} \cdots T_{J_r}. \]
\end{itemize}
We also write $T^n := T \cdots T$ for the $n$-fold composite of $T$, which corresponds to the action of the Weil-algebra $W^{\otimes n} = \N[x_1,\dots,x_n]/(x_1^2,\dots,x_n^2)$.

The main content of Leung's result is that the values of a tangent structure $T^\bullet$ on morphisms in $\Weil$ are determined by those values on five specific morphisms which correspond to the five natural transformations appearing in Cockett and Cruttwell's definition:
\begin{itemize}
	\item corresponding to the unique map $\epsilon: W \to \N$ is the \emph{projection}\index{tangent structure maps!projection},
	\[ p_T: T \to I; \]
	\item corresponding to the unit map $\eta: \N \to W$ is the \emph{zero section}\index{tangent structure maps!zero section}
	\[ 0_T: I \to T; \]
	\item corresponding to the map
	\[ \phi: \N[x,y]/(x^2,xy,y^2) \to \N[z]/(z^2); \quad x \mapsto z, \; y \mapsto z, \]
	is the \emph{addition}\index{tangent structure maps!addition}
	\[ +_T: T_2 \to T; \]
	\item corresponding to the symmetry map
	\[ \sigma: \N[x,y]/(x^2,y^2) \to \N[x,y]/(x^2,y^2); \quad x \mapsto y, \; y \mapsto x, \]
	is the \emph{flip}\index{tangent structure maps!flip}
	\[ c_T: T^2 \to T^2; \]
	\item corresponding to the map
	\[ \delta: \N[z]/(z^2) \to \N[x,y]/(x^2,y^2); \quad z \mapsto xy, \]
	is the \emph{vertical lift}\index{tangent structure maps!vertical lift}
	\[ \ell_T: T \to T^2. \]
\end{itemize}
Finally, the requirement that the tangent structure maps $- \odot X$ preserve the vertical lift pullback (\ref{lem:pb2}) corresponds to the condition that Cockett and Cruttwell refer to in~\cite[2.1]{cockett/cruttwell:2017} as the `Universality of the Vertical Lift'\index{universality of the vertical lift}, i.e.\ that for all $X \in \bcat{X}$ there is a pullback square (in $\bcat{X}$) of the form
\begin{equation} \label{eq:vert} \begin{diagram}
		\node{TX\times_X TX} \arrow{e} \arrow{s} \node{T(TX)} \arrow{s,r}{T(p)} \\
		\node{X} \arrow{e,t}{0} \node{TX.}
\end{diagram} \end{equation}
\end{remark}

\begin{examples} \label{ex:tangent-category}
Here are some of the standard examples of tangent categories. More examples, and more details, appear in the papers~\cite{cockett/cruttwell:2014,cockett/cruttwell:2018, garner:2018}.
\begin{enumerate} \itemsep=5pt
	\item Let \index{Mfld@$\Mfld$, category of finite-dimensional smooth manifolds}$\bcat{X} = \Mfld$, the category of finite-dimensional smooth manifolds and smooth maps, and let $T: \Mfld \to \Mfld$ be the ordinary tangent bundle functor. Then there is a tangent structure on $\Mfld$ with tangent bundle functor $T$ and projection map given by the usual bundle projections $TM \to M$. The zero and addition maps come from the vector bundle structure on $TM$, and the flip and vertical lift can be defined directly in terms of tangent vectors~\cite[2.2(i)]{cockett/cruttwell:2018}.
	\item The category $\Sch$\index{Sch@$\Sch$, category of schemes} of schemes has a tangent structure in which the tangent bundle functor $T: \Sch \to \Sch$ on a scheme $X$ is the vector bundle associated to the $\mathcal{O}_X$-algebra of K\"{a}hler differentials of $X$ (over $\mathrm{Spec}\;\mathbb{Z}$):
	\[ T(X) = \mathrm{Spec}\;\mathrm{Sym}\;\Omega_{X/\mathbb{Z}}. \]
	See~\cite[Ex. 2(iii)]{garner:2018} for further details.
	\item The category $\CRig$ of commutative semirings has a tangent structure given by $T^A(R) := A \otimes R$, i.e.\ with tangent bundle functor $T(R) = R[x]/(x^2)$. This structure restricts to a tangent structure on the full subcategory $\CRing$\index{CRing@$\CRing$, category of commutative rings} of commutative rings.
	\item The tangent structure on $\Sch$ of (2) restricts to a tangent structure on affine schemes, i.e.\ on $\CRing^{op}$, in which the tangent bundle functor $U: \CRing^{op} \to \CRing^{op}$ is such that $U^{op}$ is left adjoint to the tangent bundle functor $T: \CRing \to \CRing$ of (3).
	\item The category $\Weil$ has a tangent structure given by the canonical action of a monoidal category on itself.
	\item The category $\bcat{X}$ of `infinitesimally linear' objects in a model of synthetic differential geometry (SDG)\index{synthetic differential geometry (SDG)} has a tangent structure whose tangent bundle functor $T: \bcat{X} \to \bcat{X}$ is given by the exponential $T(C) = C^D$ where $D$ is an `object of infinitesimals' in $\bcat{X}$. See~\cite[5.1]{cockett/cruttwell:2014} for more details.
	\item Let $\bcat{X}$ be any category. Then there is a \emph{trivial tangent structure}\index{trivial tangent structure} on $\bcat{X}$ in which $T^A: \bcat{X} \to \bcat{X}$ is the identity functor for every Weil-algebra $A$.
	\item Let $\bcat{X}$ be a tangent category with tangent bundle functor $T$, and let $\bcat{J}$ be a category. Then the functor category $\Fun(\bcat{J},\bcat{X})$ has a tangent structure whose tangent bundle functor is composition with $T$.
	\item Any cartesian differential category (CDC)\index{cartesian differential category (CDC)} in the sense of~\cite{blute/cockett/seely:2009} has a tangent structure in which the tangent bundle functor is given by $T(A) = A \times A$; see~\cite[4.7]{cockett/cruttwell:2014}.
\end{enumerate}
\end{examples}

\subsection *{Tangent functors and natural transformations}

Cockett and Cruttwell introduce in~\cite[2.7]{cockett/cruttwell:2014} `strong morphisms of tangent structure' between tangent categories. Garner~\cite[Thm.\ 9]{garner:2018} identifies those morphisms with the `$\Weil$-functors' between $\Weil$-actegories, that is, those functors that commute (up to coherent natural isomorphism) with the actions of the monoidal category $\Weil$. We use Garner's result as the basis for the following definition, which unfolds the notion of functor between actegories; see also~\cite[3.3.2]{capucci/gavranovic:2022}.

\begin{definition} [Tangent functor] \label{def:tangent-1-functor}
Let $(\bcat{X},T^\bullet)$ and $(\bcat{Y},U^\bullet)$ be tangent categories in the sense of Definition~\ref{def:tangent-category}. A \emph{tangent functor}\index{tangent functor!between tangent categories} from $(\bcat{X},T^\bullet)$ to $(\bcat{Y},U^\bullet)$ consists of a functor $F: \bcat{X} \to \bcat{Y}$, together with natural isomorphisms
\[ \alpha_A: FT^A \arrow{e,t}{\isom} U^AF \]
for each Weil-algebra $A$, such that the following diagrams commute:
\begin{enumerate} \itemsep=5pt
\item for each morphism $\phi: A \to A'$ in $\Weil$,
\[ \begin{diagram}
	\node{FT^A} \arrow{e,t}{\alpha_A} \arrow{s,l}{FT^\phi} \node{U^AF} \arrow{s,r}{U^\phi F} \\
	\node{FT^{A'}} \arrow{e,t}{\alpha_{A'}} \node{U^{A'}F \, ;}
\end{diagram} \]
\item for each pair of Weil-algebras $A,A'$,
\[ \begin{diagram}
	\node{FT^AT^{A'}} \arrow{e,t}{\alpha_A T^{A'}} \arrow{s,l}{\isom} \node{U^A F T^{A'}} \arrow{e,t}{U^A \alpha_{A'}} \node{U^AU^{A'}F} \arrow{s,r}{\isom} \\
	\node{FT^{A \otimes T'}} \arrow[2]{e,t}{\alpha_{A \otimes A'}} \node[2]{U^{A \otimes A'}F}
\end{diagram} \]
where the vertical maps are based on the isomorphisms coming from the tangent structures $T^\bullet$ and $U^\bullet$ respectively;
\item for the unit Weil-algebra $\N$:
\[ \begin{diagram}
	\node{FT^{\N}} \arrow[2]{e,t}{\alpha_{\N}} \arrow{se,b}{\isom} \node[2]{U^{\N}F} \arrow{sw,b}{\isom} \\
	\node[2]{F}
\end{diagram} \]
where the diagonal maps are based on the isomorphisms $T^\N \isom 1_{\mathbb{X}}$ and $U^{\N} \isom 1_{\mathbb{Y}}$ coming from the tangent structures $T^\bullet$ and $U^\bullet$ respectively.
\end{enumerate}
We denote a tangent functor as $(F,\alpha)$, or often simply as $F$, leaving the natural isomorphisms $\alpha_A$ understood.
\end{definition}

\begin{definition}
For a tangent category $(\bcat{X},T^\bullet)$ the \emph{identity tangent functor on $\bcat{X}$}\index{tangent functor!identity} is $(1_{\bcat{X}},1)$ where $1_A: T^A \to T^A$ is the identity natural transformation. For tangent functors $(F,\alpha): (\bcat{X},T^\bullet) \to (\bcat{Y},U^\bullet)$ and $(G,\alpha'): (\bcat{Y},U^\bullet) \to (\bcat{Z},V^\bullet)$, their composite is the tangent functor $(GF,\alpha'')$ where $\alpha''_A$ is the composite natural transformation
\[\dgTEXTARROWLENGTH=3em GFT^A \arrow{e,t}{G\alpha_A} GU^AF \arrow{e,t}{\alpha'_AF} V^AGF. \]
\end{definition}

\begin{remark}
Cockett and Cruttwell also introduce a `lax' version of these tangent functors, in which the natural transformations $\alpha_A$ are no longer required to be isomorphisms, and there is a corresponding `oplax' notion using instead natural transformations
\[ \alpha'_A: U^AF \to FT^A \]
that satisfy corresponding conditions.
\end{remark}

We also have a notion of natural transformation between tangent functors.

\begin{definition} \label{def:tangent-transformation}
Let $(F,\alpha), (G,\alpha'): \bcat{X} \to \bcat{Y}$ be tangent functors as in Definition~\ref{def:tangent-1-functor}. A \emph{tangent transformation}\index{tangent transformation!between tangent functors on ordinary tangent categories} from $(F,\alpha)$ to $(G,\alpha')$ is a natural transformation
\[ \beta: F \to G \]
such that for each Weil-algebra $A$ the following diagam commutes:
\[ \begin{diagram}
	\node{FT^A} \arrow{e,t}{\alpha_A} \arrow{s,l}{\beta T^A} \node{U^AF} \arrow{s,r}{U^A \beta} \\
	\node{GT^A} \arrow{e,t}{\alpha'_A} \node{U^AG \, .}
\end{diagram} \]
Identity and composite tangent transformations are the ordinary identity and composite natural transformations.
\end{definition}

The definitions described above can be assembled into a $2$-category.

\begin{definition} \label{def:tan}
Let $\Tan$\index{Tan@$\Tan$, the $2$-category of ordinary tangent categories} denote the $2$-category whose objects are tangent categories, $1$-morphisms are tangent functors, and $2$-morphisms are tangent transformations, with identity and composition as described above. We refer to $\Tan$ as \emph{the $2$-category of tangent categories}.
\end{definition}

The main goal of the first part of this paper is to construct an $\infty$-categorical version of the $2$-category $\Tan$. In the next chapter, we begin that process by introducing an $\infty$-categorical version of the category of Weil-algebras.

\chapter{A Monoidal $\infty$-Category of Weil-Algebras} \label{sec:weil}

We now turn to $\infty$-categories and begin the process of setting up a theory of tangent structures in that context. The main goal of this chapter is to define an $\infty$-categorical analogue $\Weilinfty$ of the monoidal category of Weil-algebras of Definition~\ref{def:weil1}. The monoidal $\infty$-category $\Weilinfty$ plays the same role in the theory of tangent $\infty$-categories (to be defined in Chapter~\ref{sec:tangent-infty-category}) as $\Weil$ plays in that of ordinary tangent categories.

We will see from the precise definition below that the monoidal $\infty$-category $\Weilinfty$ has its homotopy category isomorphic to the ordinary category $\Weil$. That is, $\Weilinfty$ and $\Weil$ have the same objects (the ordinary Weil-algebras of Definition~\ref{def:weil1}), and the path components of the mapping space
\[ \Hom_{\Weilinfty}(A,A') \]
correspond to the ordinary Weil-algebra morphisms from $A$ to $A'$. This mapping space is the classifying space of a groupoid of what we call `labelled' Weil-algebra morphisms. In other words, $\Weilinfty$ is the nerve of a bicategory whose mapping categories are those groupoids, and the main goal of the first part of this chapter is to define this bicategory; see Proposition~\ref{prop:weil2}.

The construction of $\Weilinfty$ in this chapter highlights the close relationship with $\Weil$, but an alternative (and equivalent) construction can be given in terms of the $E_\infty$-semirings of May~\cite{may:1972}, therein called `$E_\infty$-ring spaces'. These objects can be viewed as homotopical analogues of the ordinary commutative semirings which appear in the definition of $\Weil$. We prove in Appendix~\ref{sec:e-infty} that our two approaches to the $\infty$-category $\Weilinfty$ are equivalent.

We now turn to the construction of the bicategory of Weil-algebras which underlies our definition of $\Weilinfty$.

\subsection*{A bicategory of Weil-algebras}

The goal of this section is to construct a bicategory $\Weilinfty$ whose objects are Weil-algebras, and whose mapping categories are groupoids of what we will call `labelled' Weil-algebra morphisms. Our construction of this object is inspired by work of Cranch~\cite{cranch:2010} on $E_\infty$-monoids, and we will return to that connection in Chapter~\ref{sec:differential}.

To motivate our definition, consider, for example, the Weil-algebra morphism
\begin{equation} \label{eq:example-label} \phi: \N[x,y]/(x^2,y^2) \to \N[u,v]/(u^2,uv,v^2) \end{equation}
given by
\[ \phi(1) = 1, \quad \phi(x) = 2u+v, \quad \phi(y) = v, \quad \phi(xy) = 2uv. \]
A $1$-morphism in $\Weilinfty$ consists of a Weil-algebra morphism such as $\phi$ together with extra information that allows us to distinguish the two copies of $u$ which appear in the expression $\phi(x)$ by `labelling' each of those two copies with a different element of some finite set of labels. In this case, there will then be a non-trivial $2$-morphism in the bicategory $\Weilinfty$, which permutes those two copies. 

Here is how we make this idea explicit. Let $\phi: A \to A'$ be a Weil-algebra morphism, and write $M_A$ for the set of nonzero monomials in the Weil-algebra $A$. Then a \emph{labelling} of $\phi$ will consist of a finite set $K$ (the set of \emph{labels}) together with functions
\[ s: K \to M_A, \quad t: K \to M_{A'} \]
such that for each $x_\alpha \in M_A$ we have
\[ \phi(x_\alpha) = \sum_{k \in s^{-1}(x_\alpha)} t(k). \]
The set $K$ of labels has additional structure corresponding to the multiplicative property of $\phi$: if $a$ is a label for a monomial in $\phi(x)$ and $b$ is a label for a monomial in $\phi(y)$, then it is natural to use $ab$ as a label for one of the corresponding monomials in the product $\phi(xy) = \phi(x)\phi(y)$. Note, however, that $ab$ only makes sense as a label if the product of the monomials $t(a)$ and $t(b)$ is nonzero in $A'$. In this way, we can identify each of the labels in $K$ with a product of some set of labels $k_1,\dots,k_r$ for which each $s(k_i)$ is one of the generators of the Weil-algebra $A$. 

For example, a labelling of the Weil-algebra morphism of (\ref{eq:example-label}) could be written in the form
\[ K = \{1, a_1, a_2, a_3, b, a_1b, a_2b\} \]
with
\[ s(a_1) = s(a_2) = s(a_3) = x, \quad s(b) = y \]
and
\[ t(a_1) = t(a_2) = u, \quad t(a_3) = v, \quad t(b) = v \]
where both $s$ and $t$ are extended multiplicatively to all of $K$. This observation gives $K$ the structure of a \emph{partial commutative monoid}, as in the following definition.

\begin{definition} \label{def:FPCM}
A \index{partial commutative monoid}\emph{partial commutative monoid} consists of a set $K$, a subset $K_2 \subseteq K \times K$, and a function
\[ K_2 \to K; \quad (x,y) \mapsto x \cdot y, \]
with the following properties:
\begin{enumerate}
	\item for all $x,y \in K$, $x \cdot y$ is defined (that is, $(x,y) \in K_2$) if and only if $y \cdot x$ is defined, in which case
	\[ x \cdot y = y \cdot x; \]
	\item there is an element $1_K \in K$ such that for all $x \in K$, $1_K \cdot x$ is defined, and
	\[ 1_K \cdot x = x; \]
	\item for all $x,y,z \in K$, $(x \cdot y) \cdot z$ is defined if and only if $x \cdot (y \cdot z)$ is defined, in which case
	\[ (x \cdot y) \cdot z = x \cdot (y \cdot z), \]
	and we write this element as $xyz$. (In particular, when an expression of the form $y_1 \cdots y_r$ is defined in $K$, the product of any subset of these elements is also defined.)
\end{enumerate}
For partial commutative monoids $K, K'$, a \index{partial commutative monoid!homomorphism}\emph{homomorphism} from $K$ to $K'$ is a function $f: K \to K'$ such that
\begin{enumerate} \itemsep=5pt
	\item $f(1_K) = 1_{K'}$, and
	\item for all $x,y \in K$, if $x \cdot y$ is defined, then $f(x) \cdot f(y)$ is defined, and
	\[ f(x \cdot y) = f(x) \cdot f(y). \]
\end{enumerate}
Let \index{FPCM@$\FPCM$, the category of finite partial commutative monoids}$\FPCM$ denote the category of \emph{finite} partial commutative monoids (i.e. those for which the underlying set is finite) and their homomorphisms, with identities and composition given by the ordinary identities and composition of functions.
\end{definition}

\begin{example}
For a Weil-algebra $A$, the set $M_A$\index{MA@$M_A$, the partial commutative monoid of nonzero monomials in the Weil-algebra $A$} of nonzero monomials in $A$ is a partial commutative monoid where the product $x_\alpha x_\beta$ of two monomials is defined, and in that case equals the ordinary product of those monomials, if and only if that product is nonzero in the Weil-algebra $A$.
\end{example}

We can now give a precise definition of the $1$-morphisms in our desired bicategory $\Weilinfty$.

\begin{definition} \label{def:weil-mor}
Let $A,A'$ be Weil-algebras. A \index{labelled Weil-algebra morphism}\emph{labelled Weil-algebra morphism} from $A$ to $A'$ is a diagram of finite partial commutative monoids
\begin{equation} \label{eq:span} \begin{diagram} \node[2]{K} \arrow{sw,t}{s} \arrow{se,t}{t} \\
	\node{M_A} \node[2]{M_{A'}}
\end{diagram} \end{equation}
such that
\begin{enumerate} \itemsep=5pt
	\item for $k \in K$ and a factorization
	\[ s(k) = m_1 \cdots m_r \]
	in $M_A$, there is a unique factorization
	\[ k = k_1 \cdots k_r \]
	in $K$ such that $s(k_i) = m_i$, and
	\item for $k,k' \in K$, the product $k k'$ is defined in $K$ if and only if $t(k) t(k')$ is defined in $M_{A'}$, i.e.\ is a nonzero monomial in $A'$.
\end{enumerate}
The \emph{underlying Weil-algebra morphism} $\phi:A \to A'$ associated to a span of the form (\ref{eq:span}) is given on the generator $x_\alpha$ of $A$ by
\begin{equation} \label{eq:label} \phi(x_\alpha) = \sum_{k \in s^{-1}(x_\alpha)} t(k). \end{equation}
\end{definition}

\begin{example} \label{ex:canonical}
For any morphism of Weil-algebras $\phi: A \to A'$ in which each generator of $A$ is mapped to either a single monomial of $A'$ or to zero, there is a canonical labelling given by a diagram of the form
\[ \begin{diagram}
	\node[2]{K} \arrow{sw,L} \arrow{se,t}{\phi |_{K}} \\
	\node{M_A} \node[2]{M_{A'}}
\end{diagram} \]
where $K \into M_A$ is the inclusion of the subset consisting of those monomials which do not map to $0$ in $A'$. For example, the identity morphism on the Weil-algebra $A$ has a canonical labelling given by the \emph{identity span}
\[ \begin{diagram}
	\node[2]{M_A} \arrow{sw,=} \arrow{se,=} \\
	\node{M_A} \node[2]{M_A \, .}
\end{diagram} \]
\end{example}

\begin{remark}
A labelling for a Weil-algebra morphism $\phi: A \to A'$ is determined, up to canonical isomorphism, by the labels $k \in K$ for which $s(k)$ is a generator of the Weil-algebra $A$. That is, a labelling is determined by labels for the monomials appearing in the expressions $\phi(x_i)$ where $x_1,\dots,x_n$ are the generators of $A$.
\end{remark}

\begin{remark} \label{rem:1}
Applying condition (1) in Definition~\ref{def:weil-mor} with $r = 0$, we deduce that the only element $k \in K$ for which $s(k) = 1$ is $k = 1$.

Similarly, a consequence of (2) is that the only element $k \in K$ for which $t(k) = 1$ is $k = 1$: suppose $t(k) = 1$; then $t(k) \cdot t(k) = 1 \cdot 1 = 1$ is nonzero in $A'$, and so by (2), $k \cdot k$ is defined in $K$. Then $s(k) s(k) = s(k \cdot k)$ is an element of $M_A$, but the only nonzero square monomial in the Weil-algebra $A$ is $1$, so $s(k) = 1$ and hence $k= 1$ as already noted.
\end{remark}

\begin{lemma} \label{lem:weil-mor}
Consider a labelled Weil-algebra morphism as in Definition~\ref{def:weil-mor}. The map $\phi: A \to A'$ defined in (\ref{eq:label}) is indeed a Weil-algebra morphism.
\end{lemma}
\begin{proof}
Since we have defined $\phi$ on an additive basis and extended linearly, $\phi$ evidently preserves the additive structure. It remains to show that $\phi$ also preserves the multiplicative structure. Firstly, we have (by Remark~\ref{rem:1})
\[ \phi(1) = \sum_{k \in s^{-1}(1)} t(k) = t(1) = 1. \]
To show that $\phi$ is multiplicative, it is sufficient to show that $\phi(x_\alpha x_\beta) = \phi(x_\alpha)\phi(x_\beta)$ for any two nonzero monomials $x_\alpha,x_\beta$ in $A$. We split the argument into two cases, depending on whether the product $x_\alpha x_\beta$ is zero or not in $A$.

Suppose $x_\alpha x_\beta$ is not zero in $A$, and hence is an element of $M_A$. We will establish a one-to-one correspondence between the set $s^{-1}(x_\alpha x_\beta)$ and the set of pairs $(k',k'') \in s^{-1}(x_\alpha) \times s^{-1}(x_\beta)$ such that $t(k')t(k'')$ is nonzero in $A'$. (Note that if $\phi(x_\alpha x_\beta) = 0$ in $A'$, then our argument will imply that both of these sets are empty.)

Take first an element $k \in K$ such that $s(k) = x_\alpha x_\beta$. By~\ref{def:weil-mor}(1), there is a unique factorization $k = k'k''$ such that $s(k') = x_\alpha$ and $s(k'') = x_\beta$. So we have a uniquely determined pair $(k',k'') \in s^{-1}(x_\alpha) \times s^{-1}(x_\beta)$ and $t(k')t(k'') = t(k)$ is a nonzero monomial in $A'$.

Conversely, suppose $(k',k'')$ is such a pair. Then the product $t(k')t(k'')$ is defined in $M_{A'}$, and so the product $k = k'k''$ is defined in $K$ by~\ref{def:weil-mor}(2), and satisfies $s(k) = s(k')s(k'') = x_\alpha x_\beta$. This construction establishes the one-to-one correspondence. We therefore have
\[ \phi(x_\alpha x_\beta) = \sum_{s(k) = x_\alpha x_\beta} t(k) = \sum_{s(k') = x_\alpha, s(k'') = x_\beta} t(k')t(k'') = \phi(x_\alpha)\phi(x_\beta). \]

For the second case, suppose that $x_\alpha x_\beta = 0$ in $A$. Our definition of $\phi$ automatically preserves the additive structure, so we then have $\phi(x_\alpha x_\beta) = 0$. Suppose $\phi(x_\alpha)\phi(x_\beta)$ is not zero. Then it includes a nonzero term $t(k')t(k'')$ for some $k' \in s^{-1}(x_\alpha)$ and $k'' \in s^{-1}(x_\beta)$. It follows that the product $k'k''$ is defined in $K$, and then $s(k'k'') = s(k')s(k'') = x_\alpha x_\beta$, which is a contradiction since the product $x_\alpha x_\beta$ is not defined in $M_A$. So $\phi(x_\alpha)\phi(x_\beta) = 0 = \phi(x_\alpha x_\beta)$.

This completes the check that $\phi$ preserves the multiplicative structure, and so $\phi: A \to A'$ is a Weil-algebra morphism.
\end{proof}

We now turn to the $2$-morphisms in the bicategory $\Weilinfty$, which are all isomorphisms.

\begin{definition} \label{def:relabel}
Let $A,A'$ be Weil-algebras, and let $(K,s,t)$ and $(K',s',t')$ be labelled Weil-algebra morphisms from $A$ to $A'$. A \emph{relabelling isomorphism} from $(K,s,t)$ to $(K',s',t')$ is an isomorphism of partial commutative monoids of the form $\alpha: K \isom K'$ such that the following diagram commutes
\[ \begin{diagram}
	\node[2]{K} \arrow{sw,t}{s} \arrow{se,t}{t} \arrow[2]{s,lr}{\isom}{\alpha} \\
	\node{M_A} \node[2]{M_{A'}} \\
	\node[2]{K'} \arrow{nw,b}{s'} \arrow{ne,b}{t'}
\end{diagram} \]
The identity function is a relabelling isomorphism from $(K,s,t)$ to itself, and a composite of relabelling isomorphisms is another relabelling isomorphism. Thus, the labelled Weil-algebra morphisms from $A$ to $A'$, together with the relabelling isomorphisms between them, form a groupoid which we denote
\[ \Hom_{\Weilinfty}(A,A'). \]
\end{definition}

The next step in defining the bicategory $\Weilinfty$ is to describe horizontal composition, which tells us how to label the composite $\phi'\phi$ of two Weil-algebra morphisms given a labelling of each of $\phi'$ and $\phi$.

\begin{definition} \label{def:span-comp}
For Weil-algebras $A$, $A'$, and $A''$, we define a functor
\[ \circ: \Hom_{\Weilinfty}(A',A'') \times \Hom_{\Weilinfty}(A,A') \to \Hom_{\Weilinfty}(A,A'') \]
by mapping the pair
\[ ((K',s',t'), (K,s,t)) \]
to the span in $\FPCM$ from $M_A$ to $M_{A''}$ given by the composite maps in the diagram
\begin{equation} \label{eq:span-comp} \begin{diagram}
	\node[3]{L} \arrow{sw} \arrow{se} \\
	\node[2]{K} \arrow{sw,t}{s} \arrow{se,t}{t} \node[2]{K'} \arrow{sw,t}{s'} \arrow{se,t}{t'} \\
	\node{M_A} \node[2]{M_{A'}} \node[2]{M_{A''}}
\end{diagram} \end{equation}
in which the middle square is a pullback in $\FPCM$. The specific choice of pullback is unimportant, but to be definite, we will use the standard fibre product
\[ L := K \times_{M_{A'}} K' = \{(k,k') \; | \; t(k) = s'(k') \} \]
and define
\[ s'': L \to M_A, \quad s''(k,k') := s(k) \]
and
\[ t'': L \to M_{A''}, \quad t''(k,k') := t'(k'). \]
We give $L$ the structure of a partial commutative monoid in which the product $(k,k') \cdot (l,l')$ is defined if and only if both $k \cdot l$ and $k' \cdot l'$ are defined in $K$ and $K'$ respectively. In that case, we set
\[ (k,k') \cdot (l,l') := (k \cdot l, k' \cdot l'). \]
We then define the desired functor $\circ$ on objects by setting
\[ (K',s',t') \circ (K,s,t) := (L,s'',t''). \]
This construction makes $L$ into a pullback in the category $\FPCM$. Relabelling isomorphisms for $K'$ and $K$ induce a relabelling of the pullback $L$ in a canonical way, so that $\circ$ becomes a functor as required.
\end{definition}

\begin{lemma} \label{lem:span-comp}
The underlying Weil-algebra morphism of $(K',s',t') \circ (K,s,t)$ is equal to the composite of the underlying Weil-algebra morphisms $\phi$ and $\phi'$ of $(K,s,t)$ and $(K',s',t')$ respectively.
\end{lemma}
\begin{proof}
For $x_\alpha \in M_A$, we have
\[ \begin{split} \phi'(\phi(x_\alpha)) &= \phi' \left( \sum_{k \in s^{-1}(x_\alpha)} t(k) \right) \\
	 &= \sum_{k \in s^{-1}(x_\alpha)} \phi'(t(k)) \\
	 & = \sum_{k \in s^{-1}(x_\alpha)} \sum_{k' \in s'^{-1}(t(k))} t'(k') \\
	 & = \sum_{(k,k') : s(k) = x_\alpha, t(k) = s'(k')} t'(k')\end{split} \]
which confirms that the span $(L,s'',t'')$ in Definition~\ref{def:span-comp} has $\phi'\phi$ as its underlying Weil-algebra morphism.
\end{proof}

The data of a bicategory also includes isomorphisms that witness the associativity and unitality of horizontal composition, so we define those now.

\begin{definition} \label{def:span-assoc}
Consider labelled Weil-algebra morphisms $(K,s,t): A \to A'$, $(K',s',t'): A' \to A''$, and $(K'',s'',t''): A'' \to A'''$. Then each of the composites
\[ ((K'',s'',t'') \circ (K',s',t')) \circ (K,s,t) \quad \text{and} \quad (K'',s'',t'') \circ ((K',s',t') \circ (K,s,t)) \]
corresponds to a diagram of the form
\begin{equation} \label{eq:span-assoc} \begin{diagram}
		\node[4]{M} \arrow{sw} \arrow{se} \\
		\node[3]{L} \arrow{sw} \arrow{se} \node[2]{L'} \arrow{sw} \arrow{se} \\
		\node[2]{K} \arrow{sw,t}{s} \arrow{se,t}{t} \node[2]{K'} \arrow{sw,t}{s'} \arrow{se,t}{t'} \node[2]{K''} \arrow{sw,t}{s''} \arrow{se,t}{t''} \\
		\node{M_A} \node[2]{M_{A'}} \node[2]{M_{A''}} \node[2]{M_{A'''}}
	\end{diagram} \end{equation}
in which all three squares are pullbacks. The universal property of the pullbacks gives a unique isomorphism between them which commutes with the structure maps, i.e.\ is a relabelling isomorphism between the two iterated composites. We take that map as the required associator isomorphism.

Composing $(K,s,t)$ with the identity spans on $A$ and $A'$ we get diagrams of the form
\[ \begin{diagram}
	\node[3]{L} \arrow{sw} \arrow{se} \\
	\node[2]{M_A} \arrow{sw,=} \arrow{se,=} \node[2]{K} \arrow{sw,t}{s} \arrow{se,t}{t} \\
	\node{M_A} \node[2]{M_A} \node[2]{M_{A'}}
\end{diagram} \]
and
\[ \begin{diagram}
	\node[3]{L} \arrow{sw} \arrow{se} \\
	\node[2]{K} \arrow{sw,t}{s} \arrow{se,t}{t} \node[2]{M_{A'}} \arrow{sw,=} \arrow{se,=} \\
	\node{M_A} \node[2]{M_{A'}} \node[2]{M_{A'}}
\end{diagram} \]
each of which admits a canonical isomorphism to $(K,s,t)$ itself, which we take to be the unit isomorphisms.
\end{definition}

\begin{proposition} \label{prop:weil2}
There is a \index{Weil-infinity@$\Weilinfty$!the bicategory of Weil-algebras and labelled Weil-algebra morphisms}bicategory $\Weilinfty$ whose objects are the Weil-algebras, mapping categories are the groupoids $\Hom_{\Weilinfty}(A,A')$, the identity $1$-morphisms are the identity spans of Example~\ref{ex:canonical}, horizontal composition is as in Definition~\ref{def:span-comp}, and associator and unit isomorphisms are as in Definition~\ref{def:span-assoc}.
\end{proposition}
\begin{proof}
We have to check the unit and pentagon axioms of~\cite[2.1.3]{johnson/yau:2021}. Each of those axioms asserts an equality between certain relabelling isomorphisms, each of which is determined by the universal property of a pullback. The uniqueness part of that universal property then confirms that the required equalities hold.
\end{proof}

\begin{lemma} \label{lem:hWeil}
The homotopy category of the bicategory $\Weilinfty$ is isomorphic to the category $\Weil$.
\end{lemma}
\begin{proof}
The objects of the homotopy category $h\Weilinfty$ are the objects of $\Weilinfty$, which are the Weil-algebras, so the same as the objects of $\Weil$. The morphisms in the homotopy category $h\Weilinfty$ are the connected components of the mapping groupoids in $\Weilinfty$, so we show that taking the underlying Weil-algebra morphism determines a bijection
\[ U: \pi_0\Hom_{\Weilinfty}(A,A') \isom \Hom_{\Weil}(A,A'). \]
To see that $U$ is well-defined, we observe that if $\alpha: (K,s,t) \arrow{e,t}{\isom} (K',s',t')$ is a relabelling isomorphism, and $x_\alpha \in M_A$ and $y_\beta \in M_{A'}$, then the number of elements in $K$ with $s(k) = x_\alpha$ and $t(k) = y_\beta$ is equal to the number of elements of $K'$ with $s'(k') = x_\alpha$ and $t'(k') = y_\beta$. Hence, the underlying Weil-algebra morphisms of $(K,s,t)$ and $(K',s',t')$ are equal.

To see that $U$ is injective, suppose that $(K,s,t)$ and $(K',s',t')$ are labelled Weil-algebra morphisms with the same underlying morphism $\phi: A \to A'$. For each generator $x_i$ of the Weil-algebra $A$, and each monomial $y_\beta \in M_{A'}$, we know that the subsets
\[ (s \times t)^{-1}(x_i,y_\beta) \quad \text{and} \quad (s' \times t')^{-1}(x_i,y_\beta) \]
have the same number of elements, equal to the coefficient of the monomial $y_\beta$ in the expression $\phi(x_i)$. Choose any bijections $\alpha_{i,\beta}$ between these sets, and define $\alpha: K \to K'$ by
\[ \alpha(k_1 \cdots k_r) := \alpha(k_1) \cdots \alpha(k_r) = \alpha_{i_1,\beta_1}(k_1) \cdots \alpha_{i_r,\beta_r}(k_r) \]
where $k = k_1 \cdots k_r$ is the unique decomposition of $k$ into a product such that $s(k_j) = x_{i_j}$ and $t(k_j) = y_{\beta_j}$. Then $\alpha$ is a relabelling isomorphism between $(K,s,t)$ and $(K',s',t')$.

To see that $U$ is surjective, take a Weil-algebra morphism $\phi: A \to A'$. Choose any set $K_1$ of distinct labels for the nonzero monomials in the expressions $\phi(x_i)$ for the generators $x_i$ of $A$ with $s: K_1 \to M_A$ and $t: K_1 \to M_{A'}$ the associated functions. Let $K$ be the set of subsets $\{ k_1,\dots,k_r\} \subseteq K_1$ such that $s(k_1) \cdots s(k_r)$ and  $t(k_1) \cdots t(k_r)$ are nonzero in $M_A$ and $M_{A'}$ respectively. Extend $s$ and $t$ multiplicatively to $K$. We make $K$ into a partial commutative monoid in which the product of two subsets is defined if and only if the subsets are disjoint, in which case the product is the union. The maps $s$ and $t$ are then partial monoid homomorphisms, and $(K,s,t)$ is a labelled Weil-algebra morphism whose underlying morphism is $\phi$.

Finally, we note that $U$ preserves composition by Lemma~\ref{lem:span-comp} and the underlying morphism of the identity span is the identity morphism. So $U$ determines an isomorphism of categories.
\end{proof}

\subsection*{An $\infty$-category of Weil-algebras}

We can now introduce our definition of the $\infty$-category $\Weilinfty$.

\begin{definition} \label{def:duskin-weil}
The $\infty$-category \index{Weil-infinity@$\Weilinfty$!the $\infty$-category of Weil-algebras and labelled Weil-algebra morphisms}$\Weilinfty$ is the Duskin nerve of the bicategory also denoted $\Weilinfty$ in Proposition~\ref{prop:weil2}. See Definition~\ref{def:duskin-nerve} for a description of the nerve of a bicategory.
\end{definition}

It will be convenient to have an explicit description of the $n$-simplexes in the $\infty$-category $\Weilinfty$. It turns out that those $n$-simplexes can be expressed as diagrams in $\FPCM$ indexed by the following categories.

\begin{definition} \label{def:Jn}
Let $\mathsf{J}_n$ denote the poset of \emph{intervals}\index{Jn@$\mathsf{J}_n$, the poset of intervals in $[n]$} in $[n] = \{0,1,\dots,n\}$, ordered by reverse inclusion. That is, the objects of $\mathsf{J}_n$ are the sets of the form
\[ [i,j] = \{k \in [n] \; | \; i \leq k \leq j \} \]
with a morphism $[i,j] \to [i',j']$ if and only if $i \leq i' \leq j' \leq j$. An order-preserving function $f: [n] \to [m]$ determines a functor $f_*: \mathsf{J}_n \to \mathsf{J}_m$ by
\[ f_*([i,j]) := [f(i),f(j)]. \]
\end{definition}

\begin{proposition} \label{prop:weil}
An $n$-simplex in $\Weilinfty$ can be expressed as a diagram $\alpha: \mathsf{J}_n \to \FPCM$ such that
\begin{enumerate} \setcounter{enumi}{-1}
	\item for each $j = 0,\dots,n$: $\alpha([j,j]) = M_{A_j}$ for some Weil-algebra $A_j$;
	\item for each $i < j$ in $[n]$, the diagram
	\[ \begin{diagram}
		\node[2]{\alpha([i,j])} \arrow{sw} \arrow{se} \\
		\node{\alpha([i,i])} \node[2]{\alpha([j,j])}
	\end{diagram} \]
	is a labelled Weil-algebra morphism from $A_i$ to $A_j$;
	\item for each $i < i' \leq j' < j$ in $[n]$, the square diagram
	\[ \begin{diagram}
		\node[2]{\alpha([i,j])} \arrow{sw} \arrow{se} \\
		\node{\alpha([i,j'])} \arrow{se} \node[2]{\alpha([i',j])} \arrow{sw} \\
		\node[2]{\alpha([i',j'])}
	\end{diagram} \]
	is a pullback in $\FPCM$.
\end{enumerate}
The simplicial structure maps in $\Weilinfty$ are given by precomposition with the functors labelled $f_*$ in Definition~\ref{def:Jn}.
\end{proposition}
\begin{proof}
A $0$-simplex of $\Weilinfty$ is a Weil-algebra $A$, to which we associate the diagram $\alpha: \mathsf{J}_0 \to \FPCM$ given by 
\[ \alpha([0,0]) = M_A. \]
A $1$-simplex is a labelled Weil-algebra morphism, which is precisely a diagram $\alpha: \mathsf{J}_1 \to \FPCM$ satisfying condition (1). For $n \geq 2$, the proposition follows by the same argument as in~\cite[4.5]{cranch:2010}.
\end{proof}

\begin{remark}
Our construction of the $\infty$-category $\Weilinfty$ fits into the general framework of `span $\infty$-categories', see~\cite[Sec.\ 2]{haugseng/hebestreit/linskens/nuiten:2023a}, introduced by Barwick as `effective Burnside $\infty$-categories', see~\cite[5.10]{barwick:2017}. In particular, let
\[ \FPCM_{\Weil} \subseteq \FPCM \]
denote the full subcategory of $\FPCM$ consisting of those finite partial commutative monoids of the form $M_A$ for a Weil-algebra $A$. Write
\[ \FPCM_{\Weil}^s, \FPCM_{\Weil}^t \]
for the subcategories of $\FPCM_{\Weil}$ consisting of, respectively, those morphisms $s$ satisfying condition (1) of Definition~\ref{def:weil-mor}, and those morphisms $t$ satisfying condition (2) of Definition~\ref{def:weil-mor}. Then
\[ (\FPCM_{\Weil}, \FPCM_{\Weil}^s,\FPCM_{\Weil}^t) \]
is an adequate triple in the sense of~\cite[2.1]{haugseng/hebestreit/linskens/nuiten:2023a}, and $\Weilinfty$ is the associated span $\infty$-category
\[ \Weilinfty \homeq \mathrm{Span}(\FPCM_{\Weil}, \FPCM_{\Weil}^s,\FPCM_{\Weil}^t). \]
\end{remark}

\subsection*{Some limits and colimits in $\Weilinfty$}

We now establish some basic facts about the $\infty$-category $\Weilinfty$. The mapping spaces in $\Weilinfty$ are the classifying spaces of the groupoids which make up the bicategory $\Weilinfty$. We can identify those classifying spaces as follows.

\begin{lemma} \label{lem:hom-weil}
	Let $A,A'$ be Weil-algebras. For each Weil-algebra morphism $\phi: A \to A'$, we can write
	\[ \phi(x_i) = \sum_{y \in M_{A'}} n_{\phi,i,y} \; y \]
	where $x_i$ is one of the generators of $A$, $y$ ranges over the distinct nonzero monomials in $A'$, that is, the elements of $M_{A'}$, and $n_{\phi,i,y}$ is a non-negative integer. Then the mapping space in the $\infty$-category $\Weilinfty$ is a disjoint union of components indexed by the Weil-algebra morphisms, each of which is a product of classifying spaces of symmetric groups:
	\[ \Hom_{\Weilinfty}(A,A') \homeq \coprod_{\phi: A \to A'} \; \prod_{i,y} B\Sigma_{n_{\phi,i,y}}. \]
\end{lemma}
\begin{proof}
	It follows from Lemma~\ref{lem:hWeil} that the components of the mapping space $\Hom_{\Weilinfty}(A,A')$ are in one-to-one correspondence with the Weil-algebra morphisms from $A$ to $A'$. The component corresponding to a given $\phi: A \to A'$ is the classifying space of the group of automorphisms of any specific labelling $(K,s,t)$ of $\phi$ in the groupoid $\Hom_{\Weilinfty}(A,A')$. Such an automorphism $\alpha$ is determined by its values on those elements $k \in K$ such that $s(k)$ is a generator of $A$, where it can freely permute the $n_{\phi,i,y}$ elements in the preimage $(s \times t)^{-1}(x_i,y)$. Therefore, the overall automorphism group is the product of the symmetric groups on each of those integers.	
\end{proof}

We can use Lemma~\ref{lem:hom-weil} to show that the Weil-algebra $\N$ is both initial and terminal in the $\infty$-category $\Weilinfty$.

\begin{lemma}
	The $\infty$-category $\Weilinfty$ is pointed with null object $\N$.
\end{lemma}
\begin{proof}
	Since $\N$ is a null object in the $1$-category $\Weil$, the mapping spaces $\Hom_{\Weilinfty}(A,\N)$ and $\Hom_{\Weilinfty}(\N,A)$ have a unique component. The unique morphism $\phi: A \to \N$ is given by $\phi(x_i) = 0$ for all generators $x_i$ of $A$, and so by Lemma~\ref{lem:hom-weil}, $\Hom_{\Weilinfty}(A,\N) \homeq *$. Since $\N$ itself has zero generators, it also follows from~\ref{lem:hom-weil} that $\Hom_{\Weilinfty}(\N,A)$ is contractible. Therefore, $\N$ is both terminal and initial in $\Weilinfty$.
\end{proof}

Recall that a crucial role in the theory of ordinary tangent categories is played by certain pullback diagrams in the category $\Weil$ (Lemmas~\ref{lem:pb1},~\ref{lem:pb2}). We now show that the corresponding diagrams in the $\infty$-category $\Weilinfty$ are also pullbacks.

\begin{proposition} \label{prop:tangent-pb}
Each of the tangent pullback squares in Lemmas~\ref{lem:pb1} and~\ref{lem:pb2} is the underlying diagram of a pullback in the $\infty$-category $\Weilinfty$ in which the morphisms are given the canonical labellings of Example~\ref{ex:canonical}. We refer to these pullbacks as the \emph{tangent pullbacks} in $\Weilinfty$.
\end{proposition}
\begin{proof}
We use the following criterion for a diagram in an $\infty$-category $\Weilinfty$ to be a pullback, which is a consequence of~\cite[4.2.4.1]{lurie:2009}: that for any object $A$, applying the mapping space construction $\Hom_{\Weilinfty}(A,-)$ yields a homotopy pullback square of spaces.

Now consider the foundational pullback diagram in Lemma~\ref{lem:pb1}. Each Weil-algebra morphism in that diagram has a canonical labelling described in Example~\ref{ex:canonical}. The fact that the diagram commutes in $\Weil$ implies that those canonical labellings can be uniquely extended to a commutative diagram in the $\infty$-category $\Weilinfty$ of the same form.

Now take the corresponding diagram of spaces given by applying the mapping space construction:
\begin{equation} \label{eq:hom-pb1} \begin{diagram}
		\node{\Hom_{\Weilinfty}(A, W^{J \sqcup J'})} \arrow{s} \arrow{e} \node{\Hom_{\Weilinfty}(A, W^J)} \arrow{s} \\
		\node{\Hom_{\Weilinfty}(A, W^{J'})} \arrow{e} \node{\Hom_{\Weilinfty}(A,\N)}
\end{diagram} \end{equation}
Since the original diagram in ~\ref{lem:pb1} is a pullback in $\Weil$, we know that the components of $\Hom_{\Weilinfty}(A, W^{J \sqcup J'})$ correspond to the components of the pullback of the remainder of the diagram. So fix such a component, i.e.\ a Weil-algebra morphism $\phi: A \to W^{J \sqcup J'}$.

Each monomial appearing as one of the summands in one of the expressions $\phi(x_i) \in W^{J \sqcup J'}$ is either an element of $J$ or an element of $J'$. (The constant monomial $1$ cannot appear because $\phi$ is a semiring homomorphism, and there are no non-trivial products in $M_{W^{J \sqcup J'}}$.) The corresponding expressions for the projected maps $A \to W^J$ and $A \to W^{J'}$ are given simply by deleting the terms in $J'$ or $J$, respectively. For each $i$, we therefore have a (homotopy) pullback square
\[ \begin{diagram} \dgARROWLENGTH=1.5em
	\node{\prod_{j \in J \sqcup J'} B\Sigma_{n_{i,j}}} \arrow{e} \arrow{s} \node{\prod_{j \in J} B\Sigma_{n_{i,j}}} \arrow{s} \\
	\node{\prod_{j \in J'} B\Sigma_{n_{i,j}}} \arrow{e} \node{*}
\end{diagram} \]
Taking the product over all generators $x_i$ of $A$, and then the disjoint union over all Weil-algebra morphisms $\phi: A \to W^{J \sqcup J'}$, we obtain by Lemma~\ref{lem:hom-weil} the desired homotopy pullback diagram (\ref{eq:hom-pb1}).
			
Similarly, for the vertical lift diagram in Lemma~\ref{lem:pb2}, we want to show that the corresponding diagram of spaces
\[ \begin{diagram}
	\node{\Hom_{\Weilinfty}(A,W^2)} \arrow{e} \arrow{s} \node{\Hom_{\Weilinfty}(A,W \otimes W)} \arrow{s} \\
	\node{\Hom_{\Weilinfty}(A,\N)} \arrow{e} \node{\Hom_{\Weilinfty}(A,W)}
\end{diagram} \]
is a homotopy pullback for any Weil-algebra $A$.

With $x,y,a,b$ as in Lemma~\ref{lem:pb2} we have, for each generator $x_i$ of $A$, a homotopy pullback of the form
\[ \begin{diagram}
	\node{B\Sigma_{n_{i,x}} \times B\Sigma_{n_{i,y}}} \arrow{e} \arrow{s} \node{B\Sigma_{n_{i,a}} \times B\Sigma_{n_{i,b}} \times B\Sigma_{n_{i,ab}}} \arrow{s} \\
	\node{*} \arrow{e} \node{B\Sigma_{n_{i,a}}}
\end{diagram} \]
where $n_{i,x} = n_{i,ab}$ and $n_{i,y} = n_{i,b}$, the top horizontal map matches up the corresponding terms, and the right-hand vertical map is the projection. Taking the product over all $i$ and the disjoint union over all Weil-algebra morphisms $A \to W^2$, we get the desired homotopy pullback of mapping spaces. So the vertical lift diagram is a pullback in $\Weilinfty$.
\end{proof}
		
Similarly, we verify that the tensor product of Weil-algebras realizes the coproduct in the $\infty$-category $\Weilinfty$.

\begin{proposition} \label{prop:coprod-weil}
	Let $A$ and $A'$ be Weil-algebras. Then $A \otimes A'$ is the coproduct of $A$ and $A'$ in the $\infty$-category $\Weilinfty$.
\end{proposition}
\begin{proof}
	The inclusion maps $A \to A \otimes A'$ and $A' \to A \otimes A'$ admit canonical labellings in the sense of Example~\ref{ex:canonical}. We check that for any Weil-algebra $A''$, the induced map
	\[ \Hom_{\Weilinfty}(A \otimes A', A'') \to \Hom_{\Weilinfty}(A,A'') \times \Hom_{\Weilinfty}(A',A'') \]
	is a homotopy equivalence. The tensor product is the coproduct in the homotopy category $\Weil$, so this map is a bijection on components. For a given Weil-algebra morphism $\phi: A \otimes A' \to A''$, the sequence of positive integers appearing in the expressions $\phi(x_i)$ for generators $x_i$ of $A \otimes A'$ is clearly the disjoint union of the sequences for the corresponding maps $A \to A''$ and $A' \to A''$. The claim follows again by Lemma~\ref{lem:hom-weil}.
\end{proof}

\subsection*{The Strict Monoidal Structure on $\Weilinfty$}

We saw in Chapter~\ref{sec:tancat} that the monoidal structure on the category of Weil-algebras given by tensor product plays an essential role in the definition of tangent category, and we now turn to the corresponding structure for the $\infty$-category $\Weilinfty$, which on objects is also the tensor product of Weil-algebras.

As with ordinary monoidal categories, a monoidal structure on an $\infty$-category is only associative or unital up to equivalence. However, in the $\infty$-categorical context, coherence conditions for those equivalences, such as MacLane's pentagon condition, hold only up to even higher equivalences, and so on. The combinatorics of those equivalences can be encoded in a number of ways, and we describe a few different, though equivalent, approaches in Appendix~\ref{sec:mon}.

However, it is convenient for us to have a \emph{strict} monoidal structure on $\Weilinfty$; that is, one in which the monoidal product on the underlying $\infty$-category is strictly associative and unital. We will use the following definition of strict monoidal $\infty$-category, which can be made without any of the more sophisticated approaches alluded to above.

\begin{definition} \label{def:mon-quasi}
A \index{strict monoidal $\infty$-category}\index{monoidal $\infty$-category!strict}\emph{strict monoidal $\infty$-category} is a simplicial monoid for which the underlying simplicial set is an $\infty$-category. 
\end{definition}

We will make the $\infty$-category $\Weilinfty$ into a simplicial monoid for which the monoidal structure on objects is given by the tensor product of Weil-algebras. Notice that each nonzero monomial in the tensor product $A \otimes A'$ is uniquely a tensor of a nonzero monomial in $A$ with a nonzero monomial in $A'$. In other words, there is an isomorphism of partial commutative monoids of the form
\[ M_{A \otimes A'} \isom M_A \times M_{A'}. \]
Thus, we can model the tensor product of Weil-algebras via the cartesian product in $\FPCM$, though there is one technical wrinkle in carrying out that plan: the usual model for the cartesian product of sets (and hence of partial commutative monoids) is not strictly associative or unital. For example, any one-element set is a unit up to isomorphism, but there is no object that is a \emph{strict} unit. In order to obtain a strict monoidal structure on $\Weilinfty$, we therefore need to introduce a different model for the product of sets which is strictly monoidal. The details of that construction are unimportant for us, but for definiteness we refer to work of Schauenburg, who proved the following proposition.

\begin{proposition}[{\cite[4.3]{schauenburg:2001}}] \label{prop:set-prod}
There is a strict monoidal structure $\times$ on the category $\Fin$ of finite sets, with strict unit given by a fixed one-element set which we denote $\{1\}$, such that $B \times C$ is a categorical product of $B$ and $C$ for any sets $B,C$.
\end{proposition}

\begin{remark}
We note that Schauenburg's proof in~\cite[4.3]{schauenburg:2001} does not work as stated there for the category of all finite sets. The functor $R$ constructed in the course of the proof is not injective on objects: if $(J_1,\dots,J_r)$ is any string of finite sets for which some $J_i = \varnothing$, then $R(J_1,\dots,J_r) = \varnothing$. However, if we restrict to the category of nonempty finite sets, then the proof is valid and determines a strict monoidal structure $\times$ which models the cartesian product. We extend to all finite sets by taking $J \times \varnothing = \varnothing = \varnothing \times J$ for all $J$.
\end{remark}

\begin{remark} \label{rem:prod-not}
The set $B \times C$, where $\times$ is as in Proposition~\ref{prop:set-prod}, is not necessarily equal to the set of ordered pairs $(b,c)$ with $b \in B$ and $c \in C$. However, the universal property of the product determines a specific bijection between those sets. We can therefore refer to elements of $B \times C$ via the notation $(b,c)$, with the understanding that we mean the element of $B \times C$ which corresponds to that ordered pair under this specific bijection. Similarly, we write elements of the product of three sets in the form $(b,c,d)$, and so on.
\end{remark}

\begin{remark} \label{rem:prod-fpcm}
If $B$ and $C$ are finite partial commutative monoids, then we can form the product set $B \times C$ using the strict monoidal structure of Proposition~\ref{prop:set-prod}. That set has a canonical partial monoid structure in which the product $(b,c) \cdot (b',c')$ is defined, and is equal to $(b \cdot b', c \cdot c')$, if and only if both of those products are defined in $B$ and $C$ respectively. Notice here that we are using the notation for elements of the set $B \times C$ introduced in Remark~\ref{rem:prod-not}. The resulting construction defines a strict monoidal structure on $\FPCM$ with unit object given by the set $\{1\}$, the chosen unit for the operation $\times$, with its unique partial monoid structure.
\end{remark}

\begin{example} \label{ex:MA}
In this paper, the most important objects in $\FPCM$ are the finite partial commutative monoids $M_A$ given by the sets of nonzero monomials in a Weil-algebra $A$. It follows from our previous observations that for a Weil-algebra $A = W^{J_1} \otimes \dots \otimes W^{J_r}$, there is an isomorphism
\[ M_A \isom M_{W^{J_1}} \times \dots \times M_{W^{J_r}}. \]
From now on, we take this isomorphism as our \emph{definition} of the finite partial commutative monoid $M_A$.

To be precise, for a nonempty finite set $J$, we define the finite partial commutative monoid:
\[ M_{W^J} := J \sqcup \{1\} \]
with only those products involving the identity element, such as $j \cdot 1 = j$, defined in $M_{W^J}$. We then set
\[ M_A := M_{W^{J_1}} \times \dots \times M_{W^{J_r}} \]
where $\times$ denotes the product on $\FPCM$ described in Remark~\ref{rem:prod-fpcm}.

This choice means that in $\FPCM$ we now have an \emph{equality}
\[ M_{A \otimes A'} = M_A \times M_{A'} \]
for Weil-algebras $A$ and $A'$, and $M_{\N}$ is \emph{equal} to the unit object $\{1\}$ for the strict monoidal structure on $\FPCM$.
\end{example}

We have now set things up so that the tensor product of Weil-algebras is given exactly by the chosen product $\times$ on the category $\FPCM$. We therefore have the following strict monoidal structure on the $\infty$-category $\Weilinfty$. 

\begin{definition} \label{def:weil-tensor}
Given $n$-simplexes in $\Weilinfty$, i.e.\ functors $\alpha,\beta : \mathsf{J}_n \to \FPCM$ satisfying the conditions of Proposition~\ref{prop:weil}, we define $\alpha \otimes \beta: \mathsf{J}_n \to \FPCM$ by
\[ (\alpha \otimes \beta)([i,j]) := \alpha([i,j]) \times \beta([i,j]) \]
where $\times$ is the strict monoidal product on $\FPCM$ of Remark~\ref{rem:prod-fpcm}. More precisely, $\alpha \otimes \beta$ is the composite functor
\[ \dgTEXTARROWLENGTH=3em \mathsf{J}_n \arrow{e,t}{\Delta} \mathsf{J}_n \times \mathsf{J}_n \arrow{e,t}{\alpha \times \beta} \FPCM \times \FPCM \arrow{e,t}{\times} \FPCM. \]
In particular, on objects of $\Weilinfty$ the operation $\otimes$ is given by the ordinary tensor product of Weil-algebras.
\end{definition}

\begin{proposition} \label{prop:weil-tensor}
The operation $\otimes$ of Definition~\ref{def:weil-tensor} makes $\Weilinfty$ into a strict monoidal $\infty$-category, which we denote $\Weilinfty$\index{Weil-infinity@$\Weilinfty$!monoidal structure}.
\end{proposition}
\begin{proof}
The main item to check is that the functor $\alpha \otimes \beta: \mathsf{J}_n \to \FPCM$ is indeed an $n$-simplex in $\Weilinfty$ when $\alpha$ and $\beta$ are. So we have to verify conditions (0)--(2) of Proposition~\ref{prop:weil}. Condition (0) is satisfied because, as described in Example~\ref{ex:MA}, we have
\[ (\alpha \otimes \beta)([i,i]) = \alpha([i,i]) \times \beta([i,i]) = M_{A_i} \times M_{A'_i} = M_{A_i \otimes A'_i}. \]
Condition (2) follows from the fact that a product of pullback diagrams is also a pullback.

For condition (1), we have to show that given labelled Weil-algebra morphisms
\[ \begin{diagram} 
	\node[2]{K} \arrow{sw,t}{s} \arrow{se,t}{t} \node[3]{K'} \arrow{sw,t}{s'} \arrow{se,t}{t'} \\
	\node{M_{A_0}} \node[2]{M_{A_1}} \node{M_{A'_0}} \node[2]{M_{A'_1}}
\end{diagram} \]
the product diagram
\[ \begin{diagram}
	\node[2]{K \times K'} \arrow{sw,t}{(s,s')} \arrow{se,t}{(t,t')} \\
	\node{M_{A_0 \otimes A'_0}} \node[2]{M_{A_1 \otimes A'_1}}
\end{diagram} \]
is also a labelled Weil-algebra morphism, i.e.\ satisfies conditions (1)--(2) of Definition~\ref{def:weil-mor}.

Take a pair $(k,k') \in K \times K'$ and suppose that
\[ s(k) = m_1 \cdots m_r, \quad s'(k') = m'_1 \cdots m'_r \]
so that $(s(k),s'(k'))$ has a factorization in $M_{A_0 \otimes A'_0}$. We know there are corresponding unique factorizations
\[ k = k_1 \cdots k_r, \quad k' = k'_1 \cdots k'_r \]
such that $s(k_i) = m_i$ and $s(k'_i) = m'_i$. In other words, we have a factorization of $(k,k')$ as the product of the pairs $(k_i,k'_i)$ for $i = 1,\dots,r$. For uniqueness, note that any other factorization of $(k,k')$ comprises corresponding factorizations of $k$ and $k'$, which must be the unique factorizations above. This verifies condition (1) of~\ref{def:weil-mor}.

Now suppose that $(t,t')(k_1,k'_1) \cdot (t,t')(k_2,k'_2)$ is defined in $M_{A_1 \otimes A'_1}$. That object is equal to
\[ (t(k_1)t(k_2), t'(k'_1)t'(k'_2)), \]
and so $t(k_1)t(k_2)$ is defined in $M_{A_1}$, and $t'(k'_1)t'(k'_2)$ is defined in $M_{A'_1}$. It follows that $k_1k_2$ is defined in $K$, and $k'_1 k'_2$ is defined in $K'$, and so
\[ (k_1,k'_1) \cdot (k_2,k'_2) = (k_1k_2,k'_1k'_2) \]
is defined in $K \times K'$. That verifies condition (2), completing the proof that $\alpha \otimes \beta$ is an $n$-simplex in $\Weilinfty$.

We now note that $\otimes$ commutes with the simplicial structure maps in $\Weilinfty$ and so defines a map of simplicial sets
\[ \otimes: \Weilinfty \times \Weilinfty \to \Weilinfty. \]
It remains to check that $\otimes$ is strictly associative and unital, with $M_{\N} = \{1\}$ as the unit object. These claims follow immediately from the corresponding facts about the product operation $\times$ on $\FPCM$.
\end{proof}

\begin{lemma} \label{lem:hWeil-monoidal}
The strict monoidal $\infty$-category $\Weilinfty$ induces a strict monoidal structure on the homotopy category $h\Weilinfty \isom \Weil$, which agrees with the ordinary tensor product of Weil-algebras.
 \end{lemma}
\begin{proof}
We have to show that for labelled Weil-algebra morphisms $(K,s,t): A_0 \to A_1$ and $(K',s',t'): A'_0 \to A'_1$, with underlying Weil-algebra morphisms $\phi$ and $\phi'$ respectively, the tensor product $(K,s,t) \otimes (K',s',t'): A_0 \otimes A'_0 \to A_1 \otimes A'_1$ has underlying Weil-algebra morphism $\phi \otimes \phi'$. That tensor product is given by the cartesian product of spans
\[ \begin{diagram}
	\node[2]{K \times K'} \arrow{sw,t}{s \times s'} \arrow{se,t}{t \times t'} \\
	\node{M_{A_0} \times M_{A'_0}} \node[2]{M_{A_1} \times M_{A'_1}}
\end{diagram} \]
so its underlying morphism is given by
\[ x_{\alpha} \otimes x_{\alpha'} \mapsto \sum_{(k,k') \in (s \times s')^{-1}(x_\alpha \otimes x_{\alpha'})} (t \times t')(k,k'). \]
According to Example~\ref{ex:MA} we are identifying the tensor $x_\alpha \otimes x_{\alpha'}$ with the pair $(x_\alpha,x_{\alpha'}) \in M_{A_0} \times M_{A'_0}$, so this sum can be written
\[ \sum_{k \in s^{-1}(x_\alpha)} \sum_{k' \in s'^{-1}(x_{\alpha'})} t(k) \otimes t'(k') \]
which is equal to $(\phi \otimes \phi')(x_\alpha \otimes x_{\alpha'})$.
\end{proof}

This completes our construction of the monoidal $\infty$-category $\Weilinfty$, which is central to our definition of tangent $\infty$-category in Chapter~\ref{sec:tangent-infty-category}. In the remainder of this chapter, we describe an alternative approach to the definition of $\Weilinfty$ based on the notion of $E_\infty$-semiring.

\subsection*{Weil-algebras and $E_\infty$-semirings}

Recall that the $1$-category $\Weil$ of Definition~\ref{def:weil1} is a full subcategory of $\CRig$, the category of commutative semirings. It turns out that, similarly, the $\infty$-category $\Weilinfty$ can be identified with a full subcategory of a corresponding $\infty$-category $\ERig$ of `$E_\infty$-semirings'. An $E_\infty$-semiring is a homotopical generalization of the notion of commutative semiring: roughly speaking, it consists of a topological space with two binary operations (addition and multiplication) together with homotopies which express the usual semiring axioms. Those homotopies are subject to coherence conditions, which hold only up to higher-order homotopies and so on.

A precise definition of $E_\infty$-semirings was given first by May in~\cite{may:1972}, where they were termed `$E_\infty$-ring-spaces'. Gepner, Groth, and Nikolaus give in~\cite{gepner/groth/nikolaus:2015} a conceptual characterization of the $\infty$-category $\ERig$, which we recall now.

Let $\EMon$\index{EMon@$\EMon$, the $\infty$-category of $E_\infty$-monoids} denote the $\infty$-category of $E_\infty$-monoids; see, for example~\cite[Sec.\ 1]{gepner/groth/nikolaus:2015}. An object in $\EMon$ can be viewed as a topological space together with a single binary operation that is unital, associative, and commutative, all up to homotopy (and higher coherent homotopies). Gepner, Groth, and Nikolaus prove in~\cite[5-1]{gepner/groth/nikolaus:2015} that there is a canonical symmetric monoidal structure $\otimes$ on $\EMon$, which plays the role of the tensor product of ordinary commutative monoids. That monoidal structure is characterized by the condition that the free $E_\infty$-monoid functor $\spaces \to \EMon$ is symmetric monoidal with respect to the cartesian product on $\spaces$. In particular, the monoidal unit in $\EMon$ is the free $E_\infty$-monoid on a one-point space, which is\index{N-infinity@$\N_\infty$, the free $E_\infty$-monoid on one generator}
\[ \N_\infty := B\Fin^{\isom} \isom \bigsqcup_{n \geq 0} B\Sigma_n, \]
the classifying space of the groupoid of finite sets and bijections, with binary operation induced by disjoint union.

In any symmetric monoidal $\infty$-category, there is a notion of `$E_\infty$-algebra' consisting of an object $X$ together with a `multiplication' map $\mu: X \otimes X \to X$ which is unital, associative, and commutative (again, up to homotopy and higher coherent homotopies). Gepner, Groth, and Nikolaus then define $E_\infty$-semirings as follows.

\begin{definition}[$E_\infty$-semiring] \label{def:E-semiring}
An \emph{$E_\infty$-semiring}\index{E-infinity-semiring@$E_\infty$-semiring} is an $E_\infty$-algebra in the symmetric monoidal $\infty$-category $\EMon$. The $\infty$-category $\ERig$\index{ERig@$\ERig$, the $\infty$-category of $E_\infty$-semirings}, of $E_\infty$-semirings, inherits the symmetric monoidal tensor product from $\EMon$, with unit $\N_\infty$. It follows from~\cite[3.2.4.7]{lurie:2017} that $\otimes$ is also the coproduct in $\ERig$.
\end{definition}

\begin{example} \label{ex:e-infty}
Let $M$ be a finite partial commutative monoid. In Definition~\ref{def:sigma(M)} we will give the precise construction of an $E_\infty$-semiring $\N_\infty^M$ whose underlying $E_\infty$-monoid is free on the finite set $M$, with multiplicative structure given by the map of $E_\infty$-monoids
\begin{equation} \label{eq:sigma(M)-mult} \N_\infty^M \otimes \N_\infty^M \homeq \N_\infty^{M \times M} \arrow{e,t}{\mu} \N_\infty^M, \end{equation}
where the map $\mu$ is determined by the partial monoid structure on $M$ in the following way. On the term of $\N_\infty^{M \times M}$ corresponding to a pair $(m_1,m_2)$, the map $\mu$ is either trivial (if $m_1 \cdot m_2$ is undefined in $M$) or is projection onto the term of $\N_\infty^M$ given by the product $m_1 \cdot m_2$ (if this product is defined). 
\end{example}

\begin{definition}[$E_\infty$-Weil-algebras] \label{def:infty-Weil}
An \index{E-infinity-Weil-algebra@$E_\infty$-Weil-algebra}\emph{$E_\infty$-Weil-algebra} is an $E_\infty$-semiring of the form $\N_\infty^{M_A}$, where $M_A$ is the finite partial commutative monoid of nonzero monomials in an ordinary Weil-algebra $A$.
\end{definition}

\begin{remark}
Applying the path components functor $\pi_0$ to an $E_\infty$-semiring, we get an ordinary commutative semiring. For an $E_\infty$-Weil-algebra $\N_\infty^{M_A}$, we have
\[ \pi_0(\N_\infty^{M_A}) \isom A. \]
\end{remark}

There is a one-to-one correspondence between the $E_\infty$-Weil-algebras of Definition~\ref{def:infty-Weil} and the objects of $\Weilinfty$, that is, the ordinary Weil-algebras of Definition~\ref{def:weil1}. We will prove in Appendix~\ref{sec:e-infty} that the mapping spaces in $\Weilinfty$ are equivalent to the spaces of morphisms between those $E_\infty$-Weil-algebras. More precisely, we have the following result identifying $\Weilinfty$ with a full monoidal subcategory of $\ERig$.

\begin{theorem} \label{thm:comparison}
There is a fully faithful symmetric monoidal functor
\[ \sigma: \Weilinfty \to \ERig \]
which on objects sends the Weil-algebra $A$ to the $E_\infty$-Weil-algebra $\N_\infty^{M_A}$. 
\end{theorem}

\chapter{Tangent $\infty$-Categories} \label{sec:tangent-infty-category}

With the monoidal $\infty$-category $\Weilinfty$ in hand, we now turn to tangent $\infty$-categories. The goal of this chapter is to give our definition of tangent $\infty$-category and to explore basic examples of such objects. Our approach is a direct generalization to $\infty$-categories of the notion of ordinary tangent category in Definition~\ref{def:tangent-category}. We start with the definition of a \emph{module} over a monoidal $\infty$-category.

\begin{definition}[$\Weilinfty$-module] \label{def:module}
Let $\bcat{X}$ be an $\infty$-category. A \emph{$\Weilinfty$-module structure}\index{Weil-infinity@$\Weilinfty$!module over} on $\bcat{X}$ is a (strong) monoidal functor
\[ T^\bullet: \Weilinfty \to \Fun(\bcat{X},\bcat{X}) \]
where $\Weilinfty$ has the monoidal structure of Proposition~\ref{prop:weil-tensor}, and $\Fun(\bcat{X},\bcat{X})$ has monoidal structure given by composition of functors. We think of $T^\bullet$ as making $\bcat{X}$ into a module over the monoid $\Weilinfty$ via a corresponding action map
\[ \odot: \Weilinfty \times \bcat{X} \to \bcat{X}, \]
and we refer to the pair $(\bcat{X},T^\bullet)$ as a \emph{$\Weilinfty$-module $\infty$-category}.
\end{definition}

\begin{remark}
We noted in the previous section that the theory of monoidal $\infty$-categories can be approached in a number of different, though equivalent, ways. Each of those ways yields a different, but equivalent, notion of monoidal functor between monoidal $\infty$-categories, and hence a different way to define $\Weilinfty$-module $\infty$-categories via Definition~\ref{def:module}. Several of those approaches are described in Appendix~\ref{sec:mon}. One can also give various definitions of the notion of module over a monoidal $\infty$-category, and some of these are described in Appendix~\ref{sec:module}. Those other approaches are, in many cases, more amenable to describing functors between module $\infty$-categories, which will be important in the next chapter for defining tangent functors between tangent $\infty$-categories.
\end{remark}

Our notion of tangent $\infty$-category is now a straightforward translation of Definition~\ref{def:tangent-category} into the $\infty$-category context.

\begin{definition}[Tangent $\infty$-category] \label{def:tangent-structure-infty}
Let $\bcat{X}$ be an $\infty$-category. A \index{tangent structure!on an $\infty$-category}\emph{tangent structure} on $\bcat{X}$ is a $\Weilinfty$-module structure 
\[ T^\bullet: \Weilinfty \to \Fun(\bcat{X},\bcat{X}) \]
such that for each object $X \in \bcat{X}$, the corresponding functor
\[ - \odot X: \Weilinfty \to \bcat{X} \]
preserves the tangent pullbacks of Proposition~\ref{prop:tangent-pb}. We say that the pair $(\bcat{X},T^\bullet)$ is a \emph{tangent $\infty$-category}\index{tangent infinity-category@tangent $\infty$-category}.
\end{definition}

\begin{remark}
Roughly speaking, a tangent structure on the $\infty$-category $\bcat{X}$ includes the following information: for each Weil-algebra $A$, an endofunctor
\[ T^A: \bcat{X} \to \bcat{X}, \]
for each labelled Weil-algebra morphism $\phi: A \to A'$ a natural transformation
\[ T^\phi: T^A \to T^{A'}, \]
and natural equivalences
\begin{equation} \label{eq:tangent-monoidal} T^{\N} \homeq 1_{\bcat{X}}, \quad T^{A \otimes A'} \homeq T^AT^{A'} \end{equation}
for Weil-algebras $A,A'$. These data are subject to higher coherences, and the exact way those coherences are organized depends on the choice of explicit model for the notion of $\Weilinfty$-module $\infty$-category.  However, we will show in Corollary~\ref{cor:strict} that every tangent $\infty$-category is equivalent (in a suitable sense) to one that is \emph{strict}, that is, where the equivalences in (\ref{eq:tangent-monoidal}) are identities. 
\end{remark}

\begin{remark}
Much of Remark~\ref{rem:T} extends to the $\infty$-categorical case. The endofunctor $T^A: \bcat{X} \to \bcat{X}$ for any Weil-algebra $A$ is determined, up to equivalence, by the tangent bundle functor $T = T^W$. For a positive integer $n$ we write $T_n := T^{W^n}$ and for a finite set $J$ we similarly write $T_J := T^{W^J}$. It follows from the condition on preservation of tangent pullbacks that
\[ T_n(X) \homeq T(X) \times_X \dots \times_X T(X) \]
for any $n \geq 1$, just as in an ordinary tangent category. It then follows from the equivalences in (\ref{eq:tangent-monoidal}) that if $A = W^{n_1} \otimes \dots \otimes W^{n_r}$, then
\[ T^A(X) \homeq T_{n_1} \cdots T_{n_r}(X). \]
A tangent structure on an $\infty$-category $\bcat{X}$ also entails the five natural transformations $p,0,+,c,\ell$ described in~\ref{rem:T}. Each of those natural transformations is based on a morphism of Weil-algebras which admits a canonical labelling and hence determines a morphism in $\Weilinfty$ via the process of Example~\ref{ex:canonical}. The tangent structure associates to that morphism a corresponding natural transformation between endofunctors on $\bcat{X}$.

However, in place of the strictly commutative diagrams in Cockett and Cruttwell's definition of tangent category~\cite[2.1]{cockett/cruttwell:2014}, a tangent $\infty$-category includes higher-level coherence data that establishes the commutativity of those diagrams up to homotopy.

For example, a tangent structure on a $1$-category $\bcat{X}$ gives each projection map $p_C: TC \to C$ the structure of a commutative monoid in the slice category $\bcat{X}_{/C}$. (Such a structure makes $p$ into an \emph{additive bundle} in the terminology of~\cite{cockett/cruttwell:2014}.) On the other hand, a tangent structure on an $\infty$-category $\bcat{X}$ gives each projection map the structure of an $E_\infty$-monoid, that is, an operation which is associative and commutative only up to higher coherent equivalences.
\end{remark}

\subsection*{Examples of tangent $\infty$-categories}

With our main definition set, we can now describe some of the basic examples of tangent $\infty$-categories. We start by showing that ordinary tangent categories are also tangent $\infty$-categories.

\begin{definition}[Ordinary tangent categories] \label{def:nerve-tangent}
Let $(\bcat{X},U^\bullet)$ be an ordinary tangent category in the sense of Definition~\ref{def:tangent-category}, i.e.\ $U^\bullet$ is a (strong) monoidal functor
\[ U^\bullet: \Weil \to \Fun(\bcat{X},\bcat{X}). \]
We define a $\Weilinfty$-module structure on $\bcat{X}$ by the composite
\[ \dgTEXTARROWLENGTH=3em n^*(U^\bullet): \Weilinfty \arrow{e,t}{n} \Weil \arrow{e,t}{U^\bullet} \Fun(\bcat{X},\bcat{X}) \]
where $n: \Weilinfty \to \Weil$ is the map of (strict) simplicial monoids which sends a simplex in $\Weilinfty$ to the underlying sequence of Weil-algebra morphisms. That map is transpose to the monoidal isomorphism of Lemma~\ref{lem:hWeil-monoidal} under the adjunction between the homotopy category and nerve constructions between strict monoidal categories and strict monoidal $\infty$-categories (and strict monoidal functors in each case).
\end{definition}

\begin{lemma} \label{lem:nerve-tangent}
Let $\bcat{X}$ be a tangent category. Then $\bcat{X}$, equipped with the $\Weilinfty$-module structure of Definition~\ref{def:nerve-tangent}, is a tangent $\infty$-category.
\end{lemma}
\begin{proof}
Let $X$ be an object of $\bcat{X}$. We have to show that the composite
\[ \dgTEXTARROWLENGTH=3em  \Weilinfty \arrow{e,t}{n} \Weil \arrow{e,t}{U^\bullet(X)} \bcat{X}. \]
preserves tangent pullbacks. The functor $n$ preserves tangent pullbacks by Proposition~\ref{prop:tangent-pb}, and the second map preserves tangent pullbacks since $\bcat{X}$ is an ordinary tangent category.
\end{proof}

Conversely, any tangent structure (in the sense of Definition~\ref{def:tangent-structure-infty}) on the nerve of an ordinary category arises in the manner of Lemma~\ref{lem:nerve-tangent}. 

\begin{lemma} \label{lem:category-tangent}
Let $(\bcat{X},T^\bullet)$ be a tangent $\infty$-category whose underlying $\infty$-category is (the nerve of) an ordinary category. Then $T^\bullet$ is monoidally equivalent to $n^*(U^\bullet)$ for some ordinary tangent structure $U^\bullet: \Weil \to \Fun(\bcat{X},\bcat{X})$ in the sense of Definition~\ref{def:tangent-category}.
\end{lemma}
\begin{proof}
Since $\bcat{X}$ is an ordinary category, the (strict) monoidal $\infty$-category $\Fun(\bcat{X},\bcat{X})$ is the nerve of the ordinary (strict) monoidal category $\Fun(\bcat{X},\bcat{X})$. Under the adjunction of Proposition~\ref{prop:mon-nerve}, (the monoidal equivalence class of) the monoidal functor $U^\bullet$ corresponds to (the monoidal equivalence class of) a monoidal functor
\[ U^\bullet: \Weil \isom h\Weilinfty \to \Fun(\bcat{X},\bcat{X}). \]
In other words, $T^\bullet$ is monoidally equivalent to $n^*(U^\bullet)$, where $U^\bullet$ is the composite monoidal functor
\[ \dgTEXTARROWLENGTH=2.5em \Weil \isom h\Weilinfty \arrow{e,t}{hT^\bullet} h\Fun(\bcat{X},\bcat{X}) \isom \Fun(\bcat{X},\bcat{X}). \]
Since pullbacks in $\bcat{X}$ viewed as an $\infty$-category $\bcat{X}$ are the same as pullbacks in the ordinary category $\bcat{X}$, the $\Weil$-module structure map $U^\bullet$ is a tangent structure on $\bcat{X}$.
\end{proof}

\begin{remark}
Lemma~\ref{lem:category-tangent} gives us one sense in which tangent $\infty$-categories reduce to ordinary tangent categories. A complete description of their relationship is contained in Proposition~\ref{prop:tangent-functor}, where tangent functors (and tangent natural transformations) are also taken into account.
\end{remark}

\begin{warning} \label{warn:tangent}
For any tangent $\infty$-category $\bcat{X}$, we can apply the homotopy category construction to the $\Weilinfty$-action map to get a map
\[ h\Weilinfty \times h\bcat{X} \to h\bcat{X}, \]
which makes $h\bcat{X}$ into a $\Weil$-module via the monoidal isomorphism $\Weil \isom h\Weilinfty$ of Lemma~\ref{lem:hWeil-monoidal}. However, this module structure is not necessarily a tangent structure on the category $h\bcat{X}$ because the tangent pullbacks in $\bcat{X}$ do not necessarily determine pullbacks in $h\bcat{X}$.
\end{warning}

\begin{example}[Trivial tangent structures]
Let $\bcat{X}$ be an arbitrary $\infty$-category. There is a trivial $\Weilinfty$-module structure on $\bcat{X}$ for which the structure map $I^\bullet: \Weilinfty \to \Fun(\bcat{X},\bcat{X})$ is the constant monoidal functor with value the identity functor on $\bcat{X}$. This action map preserves the tangent pullbacks since a square of identity morphisms is a pullback. So this $\Weilinfty$-module is a tangent structure, which we refer to as the \index{trivial tangent structure}\emph{trivial tangent structure} on $\bcat{X}$.
\end{example}

\begin{example}[Pointwise tangent structures on functor $\infty$-categories] \label{ex:pointwise}
Let $(\bcat{X},T^\bullet)$ be a tangent $\infty$-category, and let $S$ be any simplicial set. There is a monoidal functor
\[ *: \Fun(\bcat{X},\bcat{X}) \to \Fun(\Fun(S,\bcat{X}),\Fun(S,\bcat{X})) \]
which, on objects, sends an endofunctor $F: \bcat{X} \to \bcat{X}$ to the map
\[ F_*: \Fun(S,\bcat{X}) \to \Fun(S,\bcat{X}); \quad G \mapsto FG. \]
The functor $*$ is monoidal (with respect to composition) because
\[ (FF')_*(G) = FF'G = F_*(F'G) = F_*F'_*(G). \]
The composite map
\[ \dgTEXTARROWLENGTH=3em \Weilinfty \arrow{e,t}{T^\bullet} \Fun(\bcat{X},\bcat{X}) \arrow{e,t}{*} \Fun(\Fun(S,\bcat{X}),\Fun(S,\bcat{X})) \]
therefore makes $\Fun(S,\bcat{X})$ into a $\Weilinfty$-module.

We wish to show that for each diagram $D: S \to \bcat{X}$, the corresponding $\Weilinfty$-action map $- \odot D: \Weilinfty \to \Fun(S,\bcat{X})$ preserves the tangent pullbacks. That action map is given by
\[ (A \odot D)(s) = A \odot (D(s)) \]
where $\odot$ on the right-hand side denotes the action map corresponding to the $\Weilinfty$-module structure on $\bcat{X}$, and we know that for each $s \in S$, the functor $- \odot D(s): \Weilinfty \to \bcat{X}$ preserves tangent pullbacks. A square diagram in $\Fun(S,\bcat{X})$ which determines a pullback in $\bcat{X}$ for each $s \in S$ is a pullback square in $\Fun(S,\bcat{X})$ by~\cite[5.1.2.2]{lurie:2009}. It follows that $- \odot D$ preserves pullbacks too. Therefore $\Fun(S,\bcat{X})$, with the above $\Weilinfty$-action, is a tangent $\infty$-category. We refer to this as the \index{pointwise tangent structure}\emph{pointwise tangent structure} on $\Fun(S,\bcat{X})$.
\end{example}

\begin{lemma} \label{lem:transfer}
Let $i: \bcat{X} \weq \bcat{Y}$ be an equivalence of $\infty$-categories, and let $T: \bcat{X} \to \bcat{X}$ be the tangent bundle functor for a tangent structure on $\bcat{X}$. Then there is a tangent structure on $\bcat{Y}$ whose underlying tangent bundle functor $\bcat{Y} \to \bcat{Y}$ is equivalent to $iTi^{-1}$.
\end{lemma}
\begin{proof}
We show that $i$ induces an equivalence of monoidal $\infty$-categories
\[ \Fun(\bcat{X},\bcat{X}) \homeq \Fun(\bcat{Y},\bcat{Y}) \]
whose underlying functor is $i(-)i^{-1}$. This claim follows from the methods described in~\cite[4.7.1]{lurie:2017} in which a universal property is established for a monoidal $\infty$-category of endomorphisms such as $\End(\bcat{X}) := \Fun(\bcat{X},\bcat{X})$. That universal property says that $\End(\bcat{X})$ is the terminal object in the $\infty$-category $\bCatinf[\bcat{X}]$ in which an object is a pair $(\bcat{F},\eta)$ consisting of an $\infty$-category $\bcat{F}$ and a functor $\eta: \bcat{F} \times \bcat{X} \to \bcat{X}$. Any such terminal object has a unique associative algebra structure in $\bCatinf[\bcat{X}]$ which induces the structure of an associative algebra in $\bCatinf$, i.e.\ a monoidal $\infty$-category. (To be precise, this claim is an application of~\cite[4.7.1.40]{lurie:2017} with $\cat{C} = \cat{M} = \bCatinf$, the $\infty$-category of $\infty$-categories, and $M = \bcat{X}$.)

We therefore have an equivalence of $\infty$-categories
\[ \bCatinf[\bcat{X}] \weq \bCatinf[\bcat{Y}]; \quad (\bcat{F},\eta) \mapsto (\bcat{F},i(-)i^{-1} \circ \eta) \]
which preserves the terminal objects and their algebra structures. In particular, $(\End(\bcat{X}),i(-)i^{-1})$ is a terminal object in $\bCatinf[\bcat{Y}]$, so is equivalent, as an associative algebra, to $(\End(\bcat{Y}),\mathrm{id})$, yielding the desired monoidal equivalence.

Composing the given tangent structure map $T^\bullet: \Weilinfty \to \End(\bcat{X})$ with the equivalence constructed above determines the necessary tangent structure on the $\infty$-category $\bcat{Y}$.
\end{proof}

\begin{example}[Strict tangent structures] \label{ex:tangent-strict}
Recall that we have constructed $\Weilinfty$ as a \emph{strict} monoidal $\infty$-category, i.e.\ a simplicial monoid whose underlying simplicial set is an $\infty$-category. One way to get an action of $\Weilinfty$ on another $\infty$-category $\bcat{X}$ is to have a \emph{strict} action of that simplicial monoid. Therefore we will say that a \emph{strict} $\Weilinfty$-module consists of an $\infty$-category $\bcat{X}$ together with an action map
\[ \Weilinfty \times \bcat{X} \to \bcat{X} \]
which is \emph{strictly} associative and unital. In other words, $\bcat{X}$ is a module, in the usual strict sense, over the simplicial monoid $\Weilinfty$. If the corresponding monoidal functor $\Weilinfty \to \Fun(\bcat{X},\bcat{X})$ is a tangent structure, then we say that $\bcat{X}$ is a \index{strict tangent $\infty$-category}\emph{strict tangent $\infty$-category}. We will prove in the next section that any tangent $\infty$-category is equivalent (in a suitable sense) to a strict one.
\end{example}

\begin{example}[Tangent structure on Weil-algebras] \label{ex:weil-tangent}
Any monoidal $\infty$-category is a module over itself in a canonical way. In particular, for the strict monoidal $\infty$-category $\Weilinfty$, that module structure is given by the (strictly associative and unital) multiplication map
\[ \otimes: \Weilinfty \times \Weilinfty \to \Weilinfty. \]
For any Weil-algebra $A$, the functor $- \otimes A: \Weilinfty \to \Weilinfty$ preserves tangent pullbacks. That claim follows from the corresponding fact in the $1$-category $\Weil$, by Lemma~\ref{lem:hom-weil} and the argument used in the proof of Proposition~\ref{prop:tangent-pb}. It follows that $\Weilinfty$ is a strict tangent $\infty$-category in the sense of Example~\ref{ex:tangent-strict}.
\end{example}

We can generalize Example~\ref{ex:weil-tangent} in the following way.

\begin{proposition} \label{prop:monoidal-tangent}
Let $U: \Weilinfty \to \bcat{X}$ be a monoidal functor from $\Weilinfty$ to another monoidal $\infty$-category $\bcat{X}$. Suppose also that the underlying functor $U: \Weilinfty \to \bcat{X}$ preserves the tangent pullbacks, and the monoidal structure on $\bcat{X}$ preserves pullbacks in each variable. Then $\bcat{X}$ has a tangent structure given by the composite
\[ \dgTEXTARROWLENGTH=3em \Weilinfty \arrow{e,t}{U} \bcat{X}^{\odot} \arrow{e,t}{C} \Fun(\bcat{X},\bcat{X}) \]
where $C$ is given by the canonical action of $\bcat{X}$ on itself.
\end{proposition}
\begin{proof}
The given conditions ensure that the underlying functor of this composite preserves the tangent pullbacks, so it yields a tangent structure.
\end{proof}

\subsection*{Tangent structures on $E_\infty$-ring spectra}

We now use Proposition~\ref{prop:monoidal-tangent} to show that the ordinary tangent structures on the category of commutative rings (and its opposite) generalize to tangent structures on the $\infty$-category of $E_\infty$-ring-spectra (and its opposite).

\begin{example}[Tangent structure on $E_\infty$-ring-spectra] \label{ex:ring-spec}
Let $\ERig$ be the $\infty$-category of $E_\infty$-semirings of Gepner, Groth, and Nikolaus~\cite{gepner/groth/nikolaus:2015}, and \index{E-infinity-ring-spectra@$\E\spectra$!the $\infty$-category of $E_\infty$-ring spectra}$\E\spectra$ the $\infty$-category of $E_\infty$-ring-spectra with monoidal structure given by smash product. In Appendix~\ref{sec:e-infty} we prove that $\Weilinfty$ is equivalent to a monoidal full subcategory of $\ERig$. It also follows from~\cite[5.1, 7.5]{gepner/groth/nikolaus:2015} that there is a monoidal functor
\[ \ERig \to \E\spectra \]
given by composing the `ring completion' and `associated ring spectrum' functors. Restricting to $\Weilinfty$ we obtain a monoidal functor
\[ \hat{(-)}: \Weilinfty \to \E\spectra. \]

We claim that this monoidal functor satisfies the hypotheses of Proposition~\ref{prop:monoidal-tangent}. Pullbacks in $\E\spectra$ are calculated in the underlying category of spectra, so it is sufficient to show that for each tangent pullback in $\Weilinfty$, the corresponding diagram of spectra is a pullback. Each spectrum in one of those diagrams is a (bi)product of finitely many copies of the sphere spectrum, and each map is either a projection or an inclusion. It is simple to verify that each such diagram is a pullback.

Finally, we observe that pullback and pushout squares of spectra coincide, and the smash product commutes with pushouts in each variable, so it also preserves pullbacks. This completes the verification of the conditions of Proposition~\ref{prop:monoidal-tangent}, and so there is a tangent structure on $\E\spectra$\index{E-infinity-ring-spectra@$\E\spectra$!tangent structure} given by
\[ T^A(R) := \hat{A} \smsh R. \]
In particular, the underlying tangent bundle functor is given on an $E_\infty$-ring-spectrum $R$ by the square-zero extension of the $R$-module $R$:
\[ T(R) = \hat{W} \smsh R \isom R \ltimes R. \]
\end{example}

\begin{example}[Another tangent structure on $E_\infty$-ring spectra] \label{ex:ring-spec-op}
As noted above, for each Weil-algebra $A$, the underlying spectrum of the $E_\infty$-ring spectrum $\hat{A}$ is a (bi)product of finitely many copies of the sphere spectrum. The smash product map
\[ \hat{A} \smsh -: \E\spectra \to \E\spectra \]
is of the form $E \mapsto E^{\oplus n}$, which preserves all limits (since they are calculated in the underlying $\infty$-category of spectra). It follows from the Adjoint Functor Theorem~\cite[5.5.2.9]{lurie:2009} that $\hat{A} \smsh -$ admits a \emph{left} adjoint which we will denote
\[ S(\hat{A},-): \E\spectra \to \E\spectra. \]
The universal properties of the adjunctions make this construction functorial in $A$, i.e.\ we have a functor
\[ S(\hat{(-)},-): \Weilinfty^{op} \times \E\spectra \to \E\spectra.\]
There are natural equivalences
\[ \begin{split} \Hom_{\E\spectra}(S(\widehat{A \otimes A'},R),R') & \homeq \Hom_{\E\spectra}(R,\widehat{A \otimes A'} \smsh R') \\
	& \homeq \Hom_{\E\spectra}(R,\hat{A} \smsh \hat{A'} \smsh R') \\
	& \homeq \Hom_{\E\spectra}(S(\hat{A},R),\hat{A'} \smsh R') \\
	& \homeq \Hom_{\E\spectra}(S(\hat{A'},S(\hat{A},R)),R')
\end{split} \]
and so
\[ S(\widehat{A \otimes A'},R) \homeq S(\hat{A'},S(\hat{A},R)). \]
In other words, the monoidal structure on $\hat{(-)}$ determines a monoidal functor
\[ \Weilinfty^{op} \to \Fun(\E\spectra,\E\spectra); \quad A \mapsto S(\hat{A},-) \]
and so, taking opposites, a monoidal functor
\[ U: \Weilinfty \to \Fun(\E\spectra^{op},\E\spectra^{op}). \]
We claim that $U$ preserves tangent pullbacks, i.e. that for any $E_\infty$-ring spectrum $R$, the functor
\[ U(-)(R) = S(\hat{(-)},R) \]
takes the tangent pullbacks in $\Weilinfty$ to pushouts in $\E\spectra$. We can check that condition by showing for any $E_\infty$-ring spectra $E'$, the functor
\[ \Hom_{\E\spectra}(S(\hat{(-)},E),E') \homeq \Hom_{\E\spectra}(E,\hat{(-)} \smsh E') \]
takes the tangent pullbacks in $\Weilinfty$ to pullbacks of mapping spaces. Since $\hat{(-)}$ takes tangent pullbacks to pullbacks in $\E\spectra$, and $- \smsh E'$ preserves pullbacks, this condition holds. It follows that $U$ determines a tangent structure on the opposite $\infty$-category $\E\spectra^{op}$\index{E-infinity-ring-spectra@$\E\spectra$!dual tangent structure}, which is analogous to the tangent structure on $\CRing^{op}$ of Example~\ref{ex:tangent-category}. The underlying tangent bundle on an $E_\infty$-ring-spectrum $R$ is given by the free symmetric $R$-algebra
\[ U(R) = \mathrm{Sym}_R(L_R) \]
where $L_R$ is the \emph{absolute cotangent complex} on $R$ of~\cite[7.3.2.14]{lurie:2017}.
\end{example}

\subsection*{Tangent structure on derived manifolds}

The prototypical example of an ordinary tangent category is the standard structure on $\Mfld$, the category of smooth manifolds and smooth maps, with the ordinary smooth tangent bundle functor. Spivak introduced in~\cite{spivak:2010} an $\infty$-category $\dMfld$ of \emph{derived manifolds}, which contains $\Mfld$ as a discrete subcategory, and which allows for pullbacks along non-transverse pairs of smooth maps. We show here that the standard tangent structure on $\Mfld$ extends to a tangent structure on the $\infty$-category $\dMfld$.

For this construction, we use a version of $\dMfld$ that is characterized by the following universal property of Carchedi and Steffens~\cite{carchedi/steffens:2019}; see there for more details.

\begin{definition} \label{def:dMfld}
Let \index{dMfld@$\dMfld$, the $\infty$-category of derived manifolds}$\dMfld$ denote an idempotent-complete $\infty$-category with finite limits that admits a functor $i: \Mfld \to \dMfld$ with the following universal property: for any other idempotent-complete $\infty$-category with finite limits $\bcat{C}$, the functor $i$ induces an equivalence
\[ i^*: \Fun^{\mathrm{lex}}(\dMfld,\bcat{C}) \weq \Fun^{\pitchfork}(\Mfld,\bcat{C}) \]
between the $\infty$-categories of finite-limit-preserving functors (on the left-hand side) and functors that preserve the \emph{transverse pullbacks} and terminal object of $\Mfld$ (on the right-hand side). In particular, $i$ preserves those transverse pullbacks and the terminal object. An explicit model for $\dMfld$ is given by the opposite of the $\infty$-category of (homotopically finitely presented) simplicial $C^\infty$-rings~\cite[5.4]{carchedi/steffens:2019}.
\end{definition}

\begin{lemma} \label{lem:dMfld}
There is a monoidal functor
\[ \tilde{\bullet}: \Fun^{\pitchfork}(\Mfld,\Mfld) \to \Fun^{\mathrm{lex}}(\dMfld,\dMfld) \]
such that the following diagram commutes
\[ \begin{diagram}
	\node{\Fun^{\pitchfork}(\Mfld,\Mfld)} \arrow{e,t}{\tilde{\bullet}} \arrow{se,b}{i_*} \node{\Fun^{\mathrm{lex}}(\dMfld,\dMfld)} \arrow{s,lr}{i^*}{\homeq} \\
	\node[2]{\Fun^{\pitchfork}(\Mfld,\dMfld)}
\end{diagram} \]
where $i^*$ is the equivalence of Definition~\ref{def:dMfld}, with $\bcat{C} = \dMfld$, and $i_*$ is given by postcomposition with $i$.
\end{lemma}
\begin{proof}
A formal argument can be made using ideas from~\cite[4.7.1]{lurie:2017} in a similar manner to the proof of Lemma~\ref{lem:transfer}. Here we sketch an informal approach. Given $F,G: \Mfld \to \Mfld$ that preserve the transverse pullbacks and terminal object, the universal property on $\dMfld$ implies that $iF$ and $iG$ factor uniquely (up to contractible choice) as
\[ iF \homeq \tilde{F}i, \quad iG \homeq \tilde{G}i \]
for $\tilde{F},\tilde{G}: \dMfld \to \dMfld$. We then have
\[ iFG \homeq \tilde{F}iG \homeq \tilde{F}\tilde{G}i \]
which implies that there is a canonical equivalence $\widetilde{FG} \homeq \tilde{F}\tilde{G}$.
\end{proof}

\begin{proposition}[Tangent structure on derived manifolds] \label{prop:dMfld}
The standard tangent structure on $\Mfld$ induces, by composition with $\tilde{\bullet}$, a tangent structure on the $\infty$-category $\dMfld$.
\end{proposition}
\begin{proof}
First, we argue that the structure map for the ordinary tangent structure on $\Mfld$, that is
\[ T^\bullet: \Weil \to \Fun(\Mfld,\Mfld), \]
factors via $\Fun^{\pitchfork}(\Mfld,\Mfld)$. It is sufficient to note that the tangent bundle functor $T: \Mfld \to \Mfld$ preserves the transverse pullbacks and terminal object, and this claim follows directly from the definition of transversality.

Composing $T^\bullet$ with the monoidal functor $\tilde{\bullet}$ described in Lemma~\ref{lem:dMfld}, and with the truncation functor from $\Weilinfty$ to its homotopy category $\Weil$, we obtain a monoidal functor
\[ \tilde{T}^\bullet: \Weilinfty \to \Fun^{\mathrm{lex}}(\dMfld,\dMfld). \]
Since the target of $\tilde{T}^\bullet$ is a full subcategory of $\Fun(\dMfld,\dMfld)$, it is now sufficient to show that $\tilde{T}^\bullet$ preserves the tangent pullbacks. Equivalently, we must show that the composite
\[ \dgTEXTARROWLENGTH=2em \Weilinfty \arrow{e,t}{T^\bullet} \Fun^{\pitchfork}(\Mfld,\Mfld) \arrow{e,t}{i_*} \Fun^{\pitchfork}(\Mfld,\dMfld) \]
preserves those pullbacks. The first map does (as $T^\bullet$ is a tangent structure on $\Mfld$), and the second map does too (as each tangent pullback is transverse by~\cite[Ex. 4.4(ii)]{cockett/cruttwell:2018}, and $i$ preserves transverse pullbacks).
\end{proof}

\chapter{Tangent Functors} \label{sec:tangent-functor}

We now turn to functors between tangent $\infty$-categories. For ordinary tangent categories, Cockett and Cruttwell introduce in~\cite[2.7]{cockett/cruttwell:2014} a notion of `strong morphism of tangent structure' which consists of a functor on the underlying categories that commutes, up to natural isomorphism, with those tangent structures. Garner shows in~\cite[Thm.\ 9]{garner:2018} that Cockett and Cruttwell's definition is equivalent to the standard notion of functor between $\Weil$-actegories; see~\cite[3.3.2]{capucci/gavranovic:2022}. In particular, Garner shows that for two tangent categories $\bcat{X},\bcat{Y}$, the category of tangent functors $\bcat{X} \to \bcat{Y}$ (whose morphisms are tangent transformations) is equivalent to the category of $\Weil$-functors and $\Weil$-transformations~\cite[3.3.9]{capucci/gavranovic:2022}. We use Garner's result as the basis for our definition of tangent functor (and tangent natural transformation) between tangent $\infty$-categories.

Recall that a tangent structure on an $\infty$-category $\bcat{X}$ consists of a suitable action of the monoidal $\infty$-category $\Weilinfty$ on $\bcat{X}$. A tangent functor from $\bcat{X} \to \bcat{Y}$ is then a functor $F: \bcat{X} \to \bcat{Y}$ that commutes with those actions up to equivalence (and higher coherent equivalences). A tangent transformation between tangent functors $F,G: \bcat{X} \to \bcat{Y}$ is a natural transformation $F \to G$ that is itself coherent with those equivalences.

We make the above notions precise in Appendix~\ref{sec:module} where, for a monoidal $\infty$-category $\bcat{W}$, we describe an $\infty$-bicategory $\Mod_{\bcat{W}}(\mathbf{CAT}_\infty)$ whose objects are the $\bcat{W}$-module $\infty$-categories; see Definition~\ref{def:mod4}. We therefore make the following definition, which we view as the central object of the first part of this paper.

\begin{definition} \label{def:TANCAT}
The \emph{$\infty$-bicategory of tangent $\infty$-categories}\index{Tan-infinity@$\Tan_\infty$, the $\infty$-bicategory of tangent $\infty$-categories} is the full sub-bicategory
\[ \Tan_\infty \subseteq \Mod_{\Weilinfty}(\mathbf{CAT}_\infty) \]
whose objects are those $\Weilinfty$-modules that are tangent $\infty$-categories in the sense of Definition~\ref{def:tangent-structure-infty}, i.e.\ for which the $\Weilinfty$-action maps preserve the tangent pullbacks.
\end{definition}

In particular, we obtain the following notions.

\begin{definition} \label{def:tangent-functor}
Let $\bcat{X},\bcat{Y}$ be tangent $\infty$-categories. The \emph{$\infty$-category of tangent functors from $\bcat{X}$ to $\bcat{Y}$}\index{tangent functor!$\infty$-category of} is the mapping $\infty$-category\index{Tan-infinity-XY@$\Tan_\infty(\bcat{X},\bcat{Y})$, the $\infty$-category of tangent functors between tangent $\infty$-categories $\bcat{X},\bcat{Y}$}
\[ \Tan_\infty(\bcat{X},\bcat{Y}) := \Hom_{\Mod_{\Weilinfty}(\mathbf{CAT}_\infty)}(\bcat{X},\bcat{Y}) \]
associated to the $\infty$-bicategory $\Mod_{\Weilinfty}(\mathbf{CAT}_\infty)$. A \emph{tangent functor}\index{tangent functor!between tangent $\infty$-categories} from $\bcat{X}$ to $\bcat{Y}$ is an object, and a \emph{tangent transformation}\index{tangent transformation!between tangent functors on tangent $\infty$-categories} is a morphism, of this $\infty$-category.
\end{definition}

\begin{remark}
Roughly speaking, a tangent functor between two tangent $\infty$-categories $(\bcat{X},T^\bullet)$ and $(\bcat{Y},U^\bullet)$ consists of an underlying functor $F: \bcat{X} \to \bcat{Y}$, together with natural equivalences
\[ \alpha_A: FT^A \homeq U^AF \]
for each Weil-algebra $A$. Those equivalences must be compatible with morphisms of Weil-algebras, and with the monoidal structure on $\Weilinfty$ and action of $\Weilinfty$ on $\bcat{X}$ and $\bcat{Y}$. That is, for each morphism $\phi: A \to A'$ of Weil-algebras, the following diagram commutes up to equivalence:
\[ \begin{diagram}
\node{FT^{A}} \arrow{e,t}{\alpha_A} \arrow{s,l}{FT^\phi} \node{U^AF} \arrow{s,r}{U^\phi F} \\
\node{FT^{A'}} \arrow{e,t}{\alpha_{A'}} \node{U^{A'}F}
\end{diagram} \]
for each pair of Weil-algebras $A,A'$, the following diagram commutes up to equivalence
\[ \begin{diagram} 
\node{FT^AT^{A'}} \arrow{s,l}{\sim} \arrow{e,t}{\alpha_A T^{A'}} \node{U^AFT^{A'}} \arrow{e,t}{U^A \alpha_{A'}} \node{U^AU^{A'}F} \arrow{s,r}{\sim} \\
\node{FT^{A \otimes A'}} \arrow[2]{e,t}{\alpha_{A \otimes A'}} \node[2]{U^{A \otimes A'}F}
\end{diagram} \]
and the following diagram also commutes up to equivalence
\[ \begin{diagram}
\node{FT^{\N}} \arrow{se,b}{\isom} \arrow[2]{e,t}{\alpha_{\N}} \node[2]{U^{\N}F} \arrow{sw,b}{\isom} \\
\node[2]{F}
\end{diagram} \]
Needless to say, those equivalences are coherent up to still-higher equivalences and so on. The precise form of all this data depends on a specific choice of model for the $\infty$-bicategory $\Tan_\infty$ and its mapping $\infty$-categories. We describe one concrete approach in Corollary~\ref{cor:tangent-functor} below.
\end{remark}

\begin{remark}
A tangent transformation between tangent functors $(F,\alpha),(G,\beta): (\bcat{X},T^\bullet) \to (\bcat{Y},U^\bullet)$ consists of a natural transformation $\gamma: F \to G$ such that for each Weil-algebra $A$, the following squares commute up to equivalence
\[ \begin{diagram}
\node{FT^A} \arrow{e,t}{\alpha_A} \arrow{s,l}{\gamma T^A} \node{U^AF} \arrow{s,r}{U^A \gamma} \\
\node{GT^A} \arrow{e,t}{\beta_A} \node{U^A G}
\end{diagram} \]
and those equivalences are compatible with morphisms of Weil-algebras and the monoidal structure up to higher equivalence, and so on.
\end{remark}

\subsection*{Tangent functors between ordinary tangent categories}

Our first goal is to show that the notion of tangent functor in Definition~\ref{def:tangent-functor} recovers Cockett and Cruttwell's `strong morphisms of tangent structure' between ordinary tangent categories~\cite[2.7]{cockett/cruttwell:2014}. In fact, we will show that for any two ordinary tangent categories $\bcat{X},\bcat{Y}$, the $\infty$-category $\Tan_\infty(\bcat{X},\bcat{Y})$ is equivalent to the (nerve of the) category of (strong) tangent functors from $\bcat{X}$ to $\bcat{Y}$ and tangent transformations, in the sense of Definition~\ref{def:tangent-1-functor}.

In fact, we can state our result in terms of a fully faithful embedding of $\infty$-bicategories.

\begin{proposition} \label{prop:tangent-functor}
Let $\Tan$ denote the $2$-category of tangent categories, tangent functors, and tangent transformations, viewed as an $\infty$-bicategory via the simplicial nerve. Then the construction of Definition~\ref{def:nerve-tangent} underlies a fully-faithful embedding of $\infty$-bicategories
\[ \mathbf{N}:  \Tan \to \Tan_\infty \]
whose image is spanned by those tangent $\infty$-categories for which the underlying $\infty$-category is (the nerve of) an ordinary category.
\end{proposition}
\begin{proof}
Garner~\cite[Thm.\ 9]{garner:2018} identifies $\Tan$ with a full sub-$2$-category of the $2$-category of $\Weil$-modules, $\Weil$-functors, and $\Weil$-natural transformations, which (by Proposition~\ref{prop:nerve-mod-bicat}) is equivalent to a full subcategory of $\Mod_{\Weil}(\mathbf{CAT}_\infty)$. Recall that we have a monoidal functor $n: \Weilinfty \to \Weil$, which displays $\Weil$ as the homotopy category of $\Weilinfty$, and which (by Definition~\ref{def:n-module}) induces a map of $\infty$-bicategories
\[ n^*: \Mod_{\Weil}(\mathbf{CAT}_\infty) \to \Mod_{\Weilinfty}(\mathbf{CAT}_\infty). \]
It follows from Lemma~\ref{lem:nerve-tangent} that the composite
\[ \mathbf{N}: \Tan \into \Mod_{\Weil}(\mathbf{CAT}_\infty) \arrow{e,t}{n^*} \Mod_{\Weilinfty}(\mathbf{CAT}_\infty) \]
takes values in the full subcategory $\Tan_\infty$. It follows from Proposition~\ref{prop:homotopy-monoidal-module} that $n^*$, and hence $\mathbf{N}$, is fully faithful.
\end{proof}

\begin{corollary} \label{cor:nerve}
Let $\bcat{X},\bcat{Y}$ be ordinary tangent categories. Then the $\infty$-category
\[ \Tan_\infty(\bcat{X},\bcat{Y}) \]
is equivalent to the nerve of the ordinary category of (strong) tangent morphisms from $\bcat{X}$ to $\bcat{Y}$ and tangent transformations between them.
\end{corollary}

We view Proposition~\ref{prop:tangent-functor} as justification that our theory of tangent $\infty$-categories subsumes the theory of ordinary tangent categories developed by Cockett and Cruttwell.

In fact, we have the following stronger version of Corollary~\ref{cor:nerve}.

\begin{proposition} \label{prop:Y-functor}
Let $\bcat{X}$ be a tangent $\infty$-category, and $\bcat{Y}$ an ordinary tangent category. Recall that the homotopy category $h\bcat{X}$ becomes a $\Weil$-module in the sense of~\ref{def:actegory}, although this $\Weil$-module is not necessarily a tangent category. Nonetheless, the $\infty$-category
\[ \Tan_\infty(\bcat{X},\bcat{Y}) \]
is equivalent to the nerve of the ordinary category of $\Weil$-functors $h\bcat{X} \to \bcat{Y}$ and $\Weil$-natural transformations between them. If $h\bcat{X}$ is a tangent category, then this is equivalent to the category
\[ \Tan(h\bcat{X},\bcat{Y}) \]
of tangent functors $h\bcat{X} \to \bcat{Y}$ and tangent transformations of Definition~\ref{def:tangent-1-functor}.
\end{proposition}
\begin{proof}
This result follows immediately from Proposition~\ref{prop:Y-module-functor}.
\end{proof}

\subsection*{Tangent equivalences}

The $\infty$-bicategory $\Tan_\infty$ determines a notion of equivalence between tangent $\infty$-categories.

\begin{definition} \label{def:tangent-equiv}
A tangent functor $F: \bcat{X} \to \bcat{Y}$ is a \emph{tangent equivalence}\index{tangent equivalence} if it is an equivalence in the $\infty$-bicategory $\Tan_\infty$, that is, if there is a tangent functor $G: \bcat{Y} \to \bcat{X}$ and tangent natural equivalences (i.e.\ invertible tangent transformations) $GF \homeq I_{\bcat{X}}$ and $FG \homeq I_{\bcat{Y}}$.
\end{definition}

Since the notion of tangent equivalence depends only on the invertible tangent transformations, it can be defined in terms of the $(\infty,1)$-core of the $\infty$-bicategory $\Tan_\infty$. That core is a full subcategory of the core $\Mod_{\Weilinfty}(\bCatinf)$ of the $\infty$-bicategory of $\Weilinfty$-modules. Proposition~\ref{prop:mon-mod-bicat} therefore has the following consequence.

\begin{proposition} \label{prop:mod-weilinfty}
The $\infty$-category $\Tan_\infty^{\sim}$ of tangent $\infty$-categories and tangent functors is equivalent to a full subcategory of the $\infty$-category associated to the projective model structure on the category of marked $\Weilinfty$-modules $\mathsf{Mod}_{\Weilinfty}^+$.
\end{proposition}
\begin{proof}
In Proposition~\ref{prop:mon-mod-bicat} we construct an equivalence of $\infty$-categories
\[ \mathsf{Mod}^{+,\circ}_{\Weilinfty} \weq \Mod_{\Weilinfty}(\bCatinf). \]
Passing to the full subcategories of tangent $\infty$-categories, i.e.\ those $\Weilinfty$-modules which satisfy the tangent pullback condition, we get the desired equivalence.
\end{proof}

\begin{corollary} \label{cor:strict}
Every tangent $\infty$-category is tangent equivalent to a strict tangent $\infty$-category.
\end{corollary}

\begin{corollary} \label{cor:equiv}
A tangent functor is a tangent equivalence if and only if its underlying functor is an equivalence of $\infty$-categories.
\end{corollary}
\begin{proof}
Proposition~\ref{prop:mod-weilinfty} tells us that equivalences in $\Tan_\infty^{\sim}$ (and hence in $\Tan_\infty$) come from weak equivalences between fibrant objects in the projective model structure on $\mathsf{Mod}^+_{\Weilinfty}$. Those projective weak equivalences are precisely the morphisms which are weak equivalences in the underlying category of marked simplicial sets, i.e.\ equivalences of $\infty$-categories.
\end{proof}

\subsection*{Tangent functors between strict tangent $\infty$-categories}

In Corollary~\ref{cor:strict} we showed that any tangent $\infty$-category is tangent equivalent to a strict tangent $\infty$-category $\bcat{X}$, i.e.\ one given by a strictly associative and unital action of the simplicial monoid $\Weilinfty$ on the simplicial set $\bcat{X}$. In this section, we describe how to construct tangent functors between two strict tangent $\infty$-categories. Since tangent functors are merely $\Weilinfty$-functors between the underlying $\Weilinfty$-module $\infty$-categories, we can use Remark~\ref{rem:bar} to give an explicit description of the data required.

\begin{corollary} \label{cor:tangent-functor}
Let $(\bcat{X},T^\bullet)$ and $(\bcat{Y},U^\bullet)$ be strict tangent $\infty$-categories. Then, up to natural (tangent) isomorphism, a tangent functor from $\bcat{X}$ to $\bcat{Y}$\index{tangent functor!between strict tangent $\infty$-categories} is specified by the following data. For each $n \geq 0$, we have a map of simplicial sets
\[ F_n: \Weilinfty^n \times \bcat{X} \to \Fun(\Delta^n,\bcat{Y}). \]
These maps satisfy the following conditions:
\[ d_i(F_n(A_n,\dots,A_1,X)) = \begin{cases} F_{n-1}(A_n,\dots,A_2,T^{A_1}(X)) & \text{if $i = 0$}; \\
	F_{n-1}(A_n,\dots,A_{i+1} \otimes A_i,\dots A_1,X) & \text{if $0 < i < n$}; \\
	U^{A_n}(F_{n-1}(A_{n-1},\dots,A_1,X)) & \text{if $i = n$}; \end{cases} \]
\[ s_j(F_n(A_n,\dots,A_1,X)) = F_{n+1}(A_n,\dots,A_{j+1},\N,A_j,\dots,A_1,X) \quad \text{for $0 \leq j \leq n$}; \]
and in addition the morphism $F_1(A,X)$ is an equivalence in $\bcat{Y}$ for all $A \in \Weilinfty$ and $X \in \bcat{X}$.
\end{corollary}

Each of the equations appearing in Corollary~\ref{cor:tangent-functor} is shorthand for a certain commutative diagram of simplicial sets. The maps $d_i$ and $s_j$ refer to the maps
\[ d_i: \Fun(\Delta^n,\bcat{Y}) \to \Fun(\Delta^{n-1},\bcat{Y}), \quad s_j: \Fun(\Delta^n,\bcat{Y}) \to \Fun(\Delta^{n+1},\bcat{Y}) \]
induced by the standard coface and codegeneracy maps $d^i: \Delta^{n-1} \to \Delta^n$ and $s^j: \Delta^{n+1} \to \Delta^n$. 
 
\begin{remark} \label{rem:tangent-functor}
Further unwrapping the description in Corollary~\ref{cor:tangent-functor}, we can describe a tangent functor $\bcat{X} \to \bcat{Y}$ by the following data:
\begin{itemize} \itemsep=5pt
\item a functor $F: \bcat{X} \to \bcat{Y}$
\item a natural transformation, between functors $\Weilinfty \times \bcat{X} \to \bcat{Y}$, of the form
\[ \alpha_\bullet: FT^\bullet \to U^\bullet F \]
such that for each Weil-algebra $A$, the component
\[ \alpha_A: FT^A \weq U^AF \]
is a natural equivalence, and the component
\[ \alpha_{\N}: F = FT^{\N} \weq U^{\N}F = F \]
is the identity natural transformation on $F$;
\item a $2$-simplex in the $\infty$-category of functors $\Weilinfty \times \Weilinfty \to \Fun(\bcat{X},\bcat{Y})$ of the form
\[ \begin{diagram}
	\node{FT^\bullet T^{\bullet'}} \arrow[2]{e,t}{\alpha_{\bullet \otimes \bullet'}} \arrow{se,b}{\alpha_\bullet T^{\bullet'}} \node[2]{U^\bullet U^{\bullet'} F} \\
	\node[2]{U^\bullet F T^{\bullet'}} \arrow{ne,b}{U^\bullet \alpha_{\bullet'}}
\end{diagram} \]
which yields a degenerate $2$-simplex in $\Fun(\Weilinfty \times \bcat{X},\bcat{Y})$ when evaluated with either $\bullet$ or $\bullet'$ equal to $\N$;
\item for each $n \geq 3$, an $n$-simplex in $\Fun(\Weilinfty^n \times \bcat{X},\bcat{Y})$ whose boundary is determined by the previous data following the same pattern as above, and which determines degenerate $n$-simplexes when evaluated at $\N$ in any factor of $\Weilinfty$.
\end{itemize}
\end{remark}

\begin{remark}
The description in Remark~\ref{rem:tangent-functor} is closely related to Cockett and Cruttwell's original definition of (strong) morphism between tangent categories~\cite[2.7]{cockett/cruttwell:2014} and reduces to it when $\bcat{X}$ and $\bcat{Y}$ are ordinary (strict) tangent categories. Dropping the condition that each $\alpha_A$ is a natural equivalence, we also obtain a definition corresponding to Cockett and Cruttwell's more general (non-strong) morphisms of tangent structure. We could refer to these objects as `lax' tangent functors, though we do not develop that notion further in this paper.
\end{remark}

\chapter{Differential Objects in Cartesian Tangent $\infty$-Categories} \label{sec:differential}
Having established the basic definitions of tangent $\infty$-category and tangent functor, we now start to generalize the broader theory developed by Cockett and Cruttwell to the $\infty$-categorical setting. In this chapter, we consider the notion of `differential object' introduced in~\cite[Def. 4.8]{cockett/cruttwell:2014} (and strengthened slightly in~\cite[3.1]{cockett/cruttwell:2018}) to describe the connection between tangent categories and cartesian differential categories.

Roughly speaking, the differential objects in a tangent category $\bcat{X}$ are the tangent \emph{spaces}, that is, the fibres $T_xM$ of the tangent bundle projections over (generalized) points $x: * \to M$ of objects in $\bcat{X}$. In the tangent category of smooth manifolds, these differential objects are the Euclidean vector spaces $\R^n$; in the tangent category of schemes, they are the affine spaces $\mathbb{A}^n = \operatorname{Spec} \mathbb{Z}[x_1,\dots,x_n]$.

Our first goal in this chapter is to generalize Cockett and Cruttwell's definition of differential objects to tangent $\infty$-categories. We do that by identifying a tangent $\infty$-category $\E$ which `represents' differential objects in a tangent $\infty$-category $\bcat{X}$ in the sense that those differential objects correspond to (finite-product-preserving) tangent functors from $\E$ to $\bcat{X}$.

We then prove Proposition~\ref{prop:tangent-space-differential}, our analogue of~\cite[4.15]{cockett/cruttwell:2014}, which provides a canonical structure of a differential object on each tangent space in a tangent $\infty$-category. In Corollary~\ref{cor:tangent-space-differential}, we note that conversely any object that admits a differential structure is equivalent to a tangent space.

Finally in this chapter we turn to the connection between differential objects and the cartesian differential categories of Blute, Cockett and Seely~\cite{blute/cockett/seely:2009}. In~\ref{thm:cartesian-differential-structure}, we prove a generalization of~\cite[4.11]{cockett/cruttwell:2014} by showing that from every (cartesian) tangent $\infty$-category $\bcat{X}$ we can construct a cartesian differential category whose objects are the differential objects of $\bcat{X}$, and whose morphisms come from those in the homotopy category of $\bcat{X}$.

Cockett and Cruttwell describe differential objects in the context of tangent categories that are \emph{cartesian} in the following sense.

\begin{definition} \label{def:cartesian}
	A tangent $\infty$-category $(\bcat{X},T)$ is \emph{cartesian}\index{cartesian tangent $\infty$-category} if $\bcat{X}$ admits finite products (including a terminal object which we denote $*$), and the tangent bundle functor $T: \bcat{X} \to \bcat{X}$ preserves those finite products.	
\end{definition}

We now recall from~\cite[3.1]{cockett/cruttwell:2018} the definition of differential object.

\begin{definition} \label{def:diff-obj}
	A \emph{differential object}\index{differential object!in a cartesian tangent category} in a cartesian tangent category $(\bcat{X},T)$ consists of
	\begin{itemize}
		\item an object $D$ in $\bcat{X}$ (the \emph{underlying object}),
		\item morphisms $\sigma: D \times D \to D$ and $\zeta: * \to D$ that provide $D$ with the structure of a commutative monoid in $\bcat{X}$, and
		\item a morphism $\hat{p}: T(D) \to D$,
	\end{itemize}
	such that
	\begin{itemize}
		\item the map $\langle p,\hat{p} \rangle: T(D) \to D \times D$ is an isomorphism, and
		\item the five diagrams listed in~\cite[3.1]{cockett/cruttwell:2018} commute.
	\end{itemize}
	A \emph{linear morphism of differential objects}\index{linear morphism of differential objects} is a morphism of underlying objects that commutes with the differential structure maps $\sigma$, $\zeta$ and $\hat{p}$. We denote by $\Diff(\bcat{X})$ the category of differential objects in $\bcat{X}$ and their morphisms.
\end{definition}

\begin{warning}
	Cockett and Cruttwell in~\cite[4.11]{cockett/cruttwell:2014} use the notation $\Diff(\bcat{X})$ for a category whose morphisms include \emph{all} morphisms between underlying objects, not only those that are linear. We will write $\widehat{\Diff}(\bcat{X})$ when we need to refer to that larger category.
\end{warning}

In order to translate Definition~\ref{def:diff-obj} into the $\infty$-categorical context, we first observe that the commutative monoid structure on a differential object should be replaced, in the context of $\infty$-categories, with the structure of an $E_\infty$-monoid. We recall that there is an $\infty$-category $\E$, which plays the role of a Lawvere theory for $E_\infty$-monoids in the sense that $E_\infty$-monoids in an $\infty$-category $\bcat{X}$ correspond to finite-product-preserving functors $\E \to \bcat{X}$.

\subsection*{The Lawvere theory for $E_\infty$-monoids}

Cranch showed in~\cite{cranch:2010} that the $\infty$-category $\E$ can be constructed as the nerve of a bicategory of spans of finite sets.

\begin{definition} \label{def:E}
Let $\E$\index{E-infinity@$\E$!the bicategory of spans of finite sets} denote the bicategory in which
\begin{itemize} \itemsep=10pt
	\item an object is a finite set;
	\item a $1$-morphism from $J$ to $J'$ is a span of finite sets of the form
	\begin{equation} \label{eq:spanfin} \begin{diagram}
	\node[2]{K} \arrow{sw,t}{s} \arrow{se,t}{t} \\
	\node{J} \node[2]{J'}
	\end{diagram} \end{equation}
	\item a $2$-morphism between two such spans is a commutative diagram of the form
	\[ \begin{diagram}
		\node[2]{K} \arrow{sw,t}{s} \arrow[2]{s,lr}{\alpha}{\isom} \arrow{se,t}{t} \\
		\node{J} \node[2]{J'} \\
		\node[2]{K'} \arrow{nw,b}{s'} \arrow{ne,b}{t'}
	\end{diagram} \]
	where $\alpha: K \to K'$ is a bijection;
	\item composition of $2$-morphisms is composition of such bijections, and the identity $2$-morphism on a span is the identity map $1_K$;
	\item the identity $1$-morphism on a finite set $J$ is the identity span $J \arrow{e,=} J \arrow{e,=} J$;
	\item horizontal composition is given by the composition of spans in the manner described for $\Weilinfty$ in Definition~\ref{def:span-comp}, forming pullbacks in the category of finite sets;
	\item associator and unit isomorphisms are determined by the universal property of pullbacks, in the manner described for $\Weilinfty$ in Definition~\ref{def:span-assoc};
\end{itemize}
Since all $2$-morphisms are invertible, the Duskin nerve (Definition~\ref{def:duskin-nerve}) of $\E$ is an $\infty$-category, which we also denote $\E$\index{E-infinity@$\E$!the $\infty$-category of spans of finite sets}. Notice that the construction of $\E$ mirrors that of $\Weilinfty$, but with spans of finite partial commutative monoids (between Weil-algebras) replaced by spans of finite sets. As with $\Weilinfty$, we can give an explicit description of the $n$-simplexes in the $\infty$-category $\E$ associated to the above bicategory.
\end{definition}

\begin{proposition} \label{prop:E}
An $n$-simplex in the $\infty$-category $\E$ consists of a diagram $\beta: \mathsf{J}_n \to \Fin$, where $\mathsf{J}_n$ is the poset of intervals in $[n] = \{0,\dots,n\}$ introduced in Definition~\ref{def:Jn}, such that
\begin{itemize} 
	\item for each $i < i' \leq j' < j$ in $[n]$, the square diagram
	\[ \begin{diagram}
	\node[2]{\beta([i,j])} \arrow{sw} \arrow{se} \\
	\node{\beta([i,j'])} \arrow{se} \node[2]{\beta([i',j])} \arrow{sw} \\
	\node[2]{\beta([i',j'])}
\end{diagram} \]
is a pullback in $\Fin$.
\end{itemize}
The simplicial structure maps in $\E$ are given by precomposition with the functors labelled $f_*$ in Definition~\ref{def:Jn}.
\end{proposition}

We also make $\E$ into a strict monoidal $\infty$-category via the objectwise cartesian product of finite sets, in the same manner as $\Weilinfty$.

\begin{definition} \label{def:E-mon}
Let $\otimes: \E \times \E \to \E$ be the operation given on $n$-simplexes by
\[ \beta \otimes \beta'([i,j]) := \beta([i,j]) \times \beta'([i,j]) \]
where $\times$ denotes the strict monoidal cartesian product on the category of sets chosen in Proposition~\ref{prop:set-prod}.
\end{definition}

\begin{proposition} \label{prop:E-mon}
The operation $\otimes$ of Definition~\ref{def:E-mon} makes $\E$ into a strict monoidal $\infty$-category.
\end{proposition}
\begin{proof}
Since a product of pullbacks is a pullback, the operation $\otimes$ determines a binary operation on $\E$ as claimed. Since we have chosen a strictly associative and unital model for the cartesian product of finite sets, that operation makes $\E$ into a simplicial monoid.
\end{proof}

Just as the homotopy category of $\Weilinfty$ is Leung's ordinary (monoidal) category of Weil-algebras, so the homotopy category of $\E$ is a familiar object.

\begin{proposition}
There is a monoidal equivalence between the homotopy category of the strict monoidal $\infty$-category $\E$ and the category of free finitely-generated commutative monoids, and monoid homomorphisms, under tensor product, which we denote $\N^\bullet$.
\end{proposition}
\begin{proof}
The desired equivalence maps the finite set $J$ (an object of $\E$) to the free commutative monoid generated by $J$, which we denote $\N^J$, and the monoidal property is reflected in the natural isomorphism
\[ \N^{J \times J'} \isom \N^J \otimes \N^{J'}. \]
Morphisms in the homotopy category of $\E$ from $J$ to $J'$ are isomorphism classes of spans of the form (\ref{eq:spanfin}). Such a span determines a monoid homomorphism $\phi: \N^J \to \N^{J'}$ by the formula
\[ \phi(j) := \sum_{k \in s^{-1}(j)} t(k) \]
for each $j \in J$, extended linearly to all of $\N^J$. Here we view an element of $\N^J$ as a finite formal sum of elements of $J$. Conversely, given such a homomorphism, choosing distinct labels for the elements of each of the sums of the form $\phi(j)$ determines a span from $J$ to $J'$.
\end{proof}

\begin{remark}
It is clear from Definition~\ref{def:E} that $\E$ is self-opposite via the functor which reverses a span. The corresponding equivalence between $\N^\bullet$ and its opposite sends a homomorphism $\N^J \to \N^{J'}$ to its transpose, i.e.\ the homomorphism $\N^{J'} \to \N^J$ given by the transpose matrix.
\end{remark}

Given the similarity in the construction of the $\infty$-categories $\E$ (Definition~\ref{def:E}) and of $\Weilinfty$ (Proposition~\ref{prop:weil}), it is not surprising that there is a more direct connection between them.

\begin{definition} \label{def:UWE}
Let $U_{\FPCM}: \FPCM\to \Fin$ be the forgetful functor from finite partial commutative monoids to finite sets. Then there is a functor
\[ U: \Weilinfty \to \E; \quad \alpha \mapsto U_{\FPCM} \circ \alpha \]
which sends an $n$-simplex $\alpha$, that is, a diagram of finite partial commutative monoids, to the same diagram considered in the category of finite sets. In particular, for a Weil-algebra $A$ viewed as an object of $\Weilinfty$, we have
\[ U(A) = M_A. \]
It follows from our definition of the strict monoidal structures on $\Weilinfty$ and $\E$, both given by our chosen strict monoidal product of finite sets, that $U$ is a strict monoidal functor, i.e.\ a map of simplicial monoids. On objects we have
\[ U(A \otimes A') = M_{A \otimes A'} = M_A \times M_{A'} = U(A) \times U(A'). \]
\end{definition}

\begin{proposition}
The strict monoidal functor $U: \Weilinfty \to \E$ makes $\E$ into a (strict) tangent $\infty$-category, which we denote $(\E,U^\bullet)$\index{E-infinity@$\E$!tangent structure}, with
\[ U^A(M) = U(A) \otimes M. \]
\end{proposition}
\begin{proof}
We have to show that for each $J$ in $\E$, the functor
\[ U(-) \otimes J : \Weilinfty \to \E \]
preserves the tangent pullbacks. 

We use the criterion that a diagram in an $\infty$-category such as $\E$ is a pullback if and only if for every object $J'$, applying the functor $\Hom_{\E}(J',-)$ yields a homotopy pullback of mapping spaces. 

Thus, for the foundational pullbacks (\ref{lem:pb1}), we have to show that for any finite sets $J,J'$, and any $m,n \geq 1$, the following diagram is a homotopy pullback of spaces
\[ \begin{diagram}
	\node{\Hom_{\E}(J',M_{W^{m+n}} \times J)} \arrow{e} \arrow{s} \node{\Hom_{\E}(J',M_{W^m} \times J)} \arrow{s} \\
	\node{\Hom_{\E}(J',M_{W^n} \times J)} \arrow{e} \node{\Hom_{\E}(J',J)}
\end{diagram} \]
This claim follows from the same type of argument as in the proof of Proposition~\ref{prop:tangent-pb}. We can calculate the mapping spaces in $\E$ following the same argument as in Lemma~\ref{lem:hom-weil}: they are given by
\[ \Hom_{\E}(J',J) \homeq \prod_{\phi: \N^{J'} \to \N^J} \prod_{j' \in J',j \in J} B\Sigma_{\phi(j'),j} \]
where $\Sigma_{\phi(j'),j}$ denotes the symmetric group on the (non-negative) integer coefficient of $j$ in the expression $\phi(j') \in \N^J$, for a monoid homomorphism $\phi: \N^{J'} \to \N^J$.

For the vertical lift pullback (\ref{lem:pb2}), we similarly have to show that for any finite sets $J,J'$, the following is a homotopy pullback
\[ \begin{diagram}
	\node{\Hom_{\E}(J', M_{W^2} \times J)} \arrow{e} \arrow{s} \node{\Hom_{\E}(J', M_{W \otimes W} \times J)} \arrow{s} \\
	\node{\Hom_{\E}(J',J)} \arrow{e} \node{\Hom_{\E}(J', M_W \times J)}
\end{diagram} \]
and again this follows in the same way.
\end{proof}

\subsection*{Differential objects in a tangent $\infty$-category}

Now we return to the notion of differential objects in a tangent $\infty$-category $\bcat{X}$. An $E_\infty$-monoid structure on an object in $\bcat{X}$ corresponds to a finite-product-preserving functor $\mathscr{D}: \E \to \bcat{X}$. The additional structure needed to make that object into a differential object now corresponds to the structure of a tangent functor on $\mathscr{D}$. Hence, we have the following definition.

\begin{definition} \label{def:differential-structure}
Let $(\bcat{X},T^\bullet)$ be a cartesian tangent $\infty$-category. A \emph{differential object}\index{differential object!in a cartesian tangent $\infty$-category} in $\bcat{X}$ is a tangent functor
\[ \mathscr{D}: (\E,U^\bullet) \to (\bcat{X},T^\bullet) \]
for which the underlying functor $\E \to \bcat{X}$ preserves finite products. We also say that $\mathscr{D}$ is a \emph{differential structure} on the object $\mathscr{D}(*)$, where $*$ denotes a one-element set.

We define the \emph{$\infty$-category of differential objects} in $\bcat{X}$\index{Diff(X)@$\Diff(\bcat{X})$, $\infty$-category of differential objects in the cartersian tangent $\infty$-category $\bcat{X}$}:
\[ \Diff(\bcat{X}) :=  \Tan^{\times}_\infty(\E,\bcat{X}) \subseteq \Tan_\infty(\E,\bcat{X}) \]
to be the full subcategory of the $\infty$-category of tangent functors $\E \to \bcat{X}$ (see Definition~\ref{def:tangent-functor}), whose objects are those tangent functors for which the underlying functor preserves finite products.
\end{definition}

\subsection*{Differential objects in an ordinary tangent category}

The first main result of this chapter is to show that Definition~\ref{def:differential-structure} reduces to Definition~\ref{def:diff-obj} in the case that $\bcat{X}$ is an ordinary tangent category.

\begin{proposition} \label{prop:diff}
Let $(\bcat{X}, T)$ be a cartesian tangent category. Then the $\infty$-category $\Diff(\bcat{X})$ is an ordinary category, which is equivalent to the category of differential objects in $\bcat{X}$ and their linear morphisms, as described in Definition~\ref{def:diff-obj}.
\end{proposition}
\begin{proof}
The homotopy category of the $\infty$-category $\E$ is equivalent to the category $\N^\bullet$ of free finitely-generated commutative monoids. That homotopy category inherits a $\Weil$-module structure from the $\Weilinfty$-module structure on $\E$ and in this case, the tangent pullbacks in $\E$ \emph{do} pass to pullbacks in $\N^\bullet$. So $\N^\bullet$ is an ordinary tangent category, with structure given by the ordinary tensor product of commutative monoids.
\[ A \odot \N^J := A \otimes \N^J. \]
When $\bcat{X}$ is an ordinary tangent category, it then follows from Proposition~\ref{prop:Y-functor} that
\[ \Tan_\infty(\E,\bcat{X}) \isom \Tan(\N^\bullet,\bcat{X}). \]
Hence the full subcategory $\Diff(\bcat{X}) \subseteq \Tan_\infty(\E,\bcat{X})$ is equivalent to the category of finite-product-preserving (ordinary) tangent functors $\N^\bullet \to \bcat{X}$ and their tangent transformations.

Thus, our goal reduces to a claim purely within the ordinary tangent category theory of Cockett and Cruttwell, which we prove in the next proposition.
\end{proof}

\begin{proposition}
Let $\bcat{X}$ be an arbitrary tangent category. Then there is an equivalence between the category of differential objects in $\bcat{X}$ and linear morphisms and the category of finite-product-preserving strong tangent functors $\N^\bullet \to \bcat{X}$ and tangent transformations.
\end{proposition}
\begin{proof}
To build the required equivalence, notice that the object $\N$ in the category $\N^\bullet$ has a differential structure given by its ordinary monoid structure, and with $\hat{p}: T_{\otimes}(\N) = W \otimes \N = W \to \N$ given by the monoid (but not algebra) map $a+bx \mapsto b$.

Any product-preserving tangent functor preserves differential objects, and the components of a tangent transformation are morphisms that commute with the differential structure maps $\sigma$, $\zeta$ and $\hat{p}$. We therefore have a functor $e$ from the category of product-preserving tangent functors to the category of differential objects.

To see that $e$ is fully faithful, take two product-preserving tangent functors $\mathscr{D},\mathscr{D'}: \N^\bullet \to \bcat{X}$, i.e. strong tangent morphisms (in the sense of~\cite[2.7]{cockett/cruttwell:2014}) $(\mathscr{D},\alpha)$ and $(\mathscr{D}',\alpha')$. A tangent transformation $\beta: (\mathscr{D},\alpha) \to\nolinebreak (\mathscr{D}',\alpha')$ is a natural transformation $\beta:\nolinebreak \mathscr{D} \to\nolinebreak \mathscr{D}'$ such that the following diagram commutes
\begin{equation} \label{eq:beta} \begin{diagram}
		\node{\mathscr{D}T_{\otimes}} \arrow{s,l}{\beta T_{\N^\bullet}} \arrow{e,t}{\alpha} \node{T\mathscr{D}} \arrow{s,r}{T\beta} \\
		\node{\mathscr{D}'T_{\otimes}} \arrow{e,t}{\alpha'} \node{T\mathscr{D}'}
\end{diagram} \end{equation}
where $T$ denotes the tangent bundle functor on $\bcat{X}$.

Since all the functors in (\ref{eq:beta}) preserve finite products, that diagram commutes if and only if it does so on its $\N$-component, and the natural transformation $\beta$ is uniquely determined by the component $\beta_{\N}: \mathscr{D}(\N) \to \mathscr{D}'(\N)$, which must be a map of commutative monoids in $\bcat{X}$.

Therefore, the natural transformation $\beta$ is a tangent transformation (i.e.\ the diagram above commutes) if and only if $\beta_{\N}$ commutes with the differential structure map $\hat{p}$. Thus $e$ determines a bijection between the tangent transformations $(\mathscr{D},\alpha) \to (\mathscr{D}',\alpha')$ and morphisms of differential objects $\mathscr{D}(\N) \to \mathscr{D}'(\N)$.

To see that $e$ is essentially surjective, consider a differential object $(D,\sigma,\zeta,\hat{p})$ in $\bcat{X}$. The commutative monoid $(D,\sigma,\zeta)$ determines (uniquely up to isomorphism) a product-preserving functor $\mathscr{D}: \N^\bullet \to \bcat{X}$ with $\mathscr{D}(\N^J) \isom D^J$.

To make $\mathscr{D}$ into a tangent functor, we have to define a natural isomorphism
\[ \alpha: \mathscr{D}T_{\otimes} \to T\mathscr{D} \]
that commutes with the tangent structure maps. We define the component
\[ \alpha_{\N}: \mathscr{D}T_{\otimes}(\N) = \mathscr{D}(W) \isom D^{\{1,x\}} \weq T(D) = T\mathscr{D}(\N) \]
of $\alpha$ to be the inverse of the isomorphism $\langle p, \hat{p} \rangle$ of Definition~\ref{def:diff-obj}. Other components are determined by the naturality requirement and the fact that all functors involved preserve finite products. Thus, we obtain the natural isomorphism $\alpha$.

To see that $(\mathscr{D},\alpha)$ is a tangent functor, we check the conditions of~\cite[2.7]{cockett/cruttwell:2014}. We will write out the proof for the commutative diagram involving the vertical lift $\ell$; the other conditions are much easier to verify. We must show that the following diagram commutes:
\[ \begin{diagram}
	\node{\mathscr{D}T_{\otimes}} \arrow[2]{e,tb}{\alpha}{\isom} \arrow{s,l}{\mathscr{D}\ell_{\otimes}} \node[2]{T\mathscr{D}} \arrow{s,r}{\ell \mathscr{D}} \\
	\node{\mathscr{D}T_{\otimes}^2} \arrow{e,tb}{\alpha T_{\otimes}}{\isom} \node{T\mathscr{D}T_{\otimes}} \arrow{e,tb}{T\alpha}{\isom} \node{T^2\mathscr{D}}
\end{diagram} \]
Since all functors in this diagram preserve finite products, it is sufficient to look at the diagram of components at the object $\N \in \N^\bullet$. This diagram in $\bcat{X}$ takes the form
\[ \begin{diagram} \dgARROWLENGTH=5em
	\node{D^2} \arrow{s,l}{\langle \pi_1,\zeta !,\zeta !,\pi_2 \rangle} \node[2]{T(D)} \arrow[2]{w,tb}{\langle p, \hat{p} \rangle}{\isom} \arrow{s,r}{\ell} \\
	\node{D^4} \node{T(D)^2} \arrow{w,tb}{\langle p\pi_1,\hat{p}\pi_1,p \pi_2,\hat{p}\pi_2 \rangle}{\isom} \node{T^2(D)} \arrow{w,tb}{\langle T(p),T(\hat{p}) \rangle}{\isom}
\end{diagram} \]
where we have identified the terms $D^2$, $T(D)^2$ and $D^4$ by choosing an order on the basis elements of $W$ (we use $\{1,x\}$) and $W \otimes W$ (we use $\{1,x_1,x_2,x_1x_2\}$), and we are writing $!$ for the unique map to the terminal object. We check that this diagram commutes by looking at each of the four components in turn:
\[ pT(p)\ell = p0p = p, \quad \hat{p}T(p)\ell = \hat{p}0p = \zeta !, \quad pT(\hat{p})\ell = \hat{p}p_T\ell = \zeta !, \quad \hat{p}T(\hat{p})\ell = \hat{p}. \]
These equations follow from the definition of differential object in~\cite[3.1]{cockett/cruttwell:2018}. In particular, the last equation is precisely the `extra' axiom that was added to the definition of differential object in~\cite{cockett/cruttwell:2018}.

The remaining four diagrams in~\cite[2.7]{cockett/cruttwell:2014} are verified in a similar manner. To see that $(\mathscr{D},\alpha)$ is a strong morphism of tangent categories, we also need to show that the functor $\mathscr{D}: \N^J \mapsto D^J$ preserves the tangent pullbacks in $\N^\bullet$. Each of those pullbacks can be explicitly described using basis elements, and its preservation follows from the form of $\mathscr{D}$.

Finally we check that the differential structure on $e(\mathscr{D}) = \mathscr{D}(\N) = D$ induced by that on $\N$ agrees with the original structure $(D,\sigma,\zeta,\hat{p})$. For the commutative monoid structure, this is a standard fact about models for the Lawvere theory $\N^\bullet$. For the projection map $\hat{p}$, we have to check that the map
\[ \dgTEXTARROWLENGTH=3em T(D) = T\mathscr{D}(\N) \arrow{e,t}{\alpha_\N^{-1}} \mathscr{D}T_{\otimes}(\N) = \mathscr{D}(W) \arrow{e,t}{\mathscr{D}(\hat{p}_{\N})} \mathscr{D}(\N) = D \]
agrees with $\hat{p}$, which it does. Thus $e$ is essentially surjective and is an equivalence of categories as claimed.
\end{proof}

\subsection*{Differential structure on tangent spaces}

The next main goal of this chapter is to identify differential objects in a tangent $\infty$-category with the `tangent spaces'. This part of the paper parallels~\cite[Sec.\ 4.4]{cockett/cruttwell:2014}.

\begin{definition} \label{def:tangent-space}
	Let $(\bcat{X},T^\bullet)$ be a cartesian tangent $\infty$-category. For a pointed object in $\bcat{X}$, that is, a morphism $x: * \to C$ in $\bcat{X}$ where $*$ denotes a terminal object, the \emph{tangent space to $C$ at $x$}\index{TAE@$T_xC$, tangent space at $x$ to an object $C$ in a cartesian tangent $\infty$-category}\index{tangent space} is the pullback
	\[ \begin{diagram}
		\node{T_xC} \arrow{e} \arrow{s} \node{T(C)} \arrow{s,r}{p_{C}} \\
		\node{*} \arrow{e,t}{x} \node{C}
	\end{diagram} \]
	if that pullback exists and is preserved by the functor $T^A: \bcat{X} \to \bcat{X}$ for each Weil-algebra $A$. We write $\bcat{X}^T_*$ for the subcategory of $\Fun(\Delta^1,\bcat{X})$ whose objects are morphisms in $\bcat{X}$ of the form $x: * \to C$ with the property that the tangent space to $C$ at $x$ exists.
\end{definition}

We now wish to construct a differential object in $\bcat{X}$ whose underlying object is the tangent space $T_xC$, i.e.\ a product-preserving tangent functor $\cat{T}_xC: \E \to \bcat{X}$ with $* \mapsto T_xC$. We will use the presentation of tangent functors from Corollary~\ref{cor:tangent-functor} to describe $\cat{T}_xC$, assuming that we have replaced $\bcat{X}$ with an equivalent strict tangent $\infty$-category.

Our first goal is to describe the underlying functor $\cat{T}_xC$. Since it must preserve products, that functor is necessarily given on objects by
\[ J \mapsto (T_xC)^J \]
for each finite set $J$. The following lemma shows how to relate these finite products to the tangent structure on $\bcat{X}$.

\begin{lemma} \label{lem:W}
For a finite set $J$, there is an equivalence
\[ (T_xC)^J \homeq T^{W^J}_xC \]
where $T^A_xC$, for a Weil-algebra $A$, denotes the pullback
\begin{equation} \label{eq:TxC} \begin{diagram}
\node{T^A_xC} \arrow{s} \arrow{e} \node{T^AC} \arrow{s} \\
\node{*} \arrow{e,t}{x} \node{C}
\end{diagram} \end{equation}
in which the right-hand vertical map is given by the augmentation $A \to \N$ of the Weil-algebra $A$. Moreover, the pullback (\ref{eq:TxC}) with $A = W^J$ is preserved by $T^{A'}$ for any Weil-algebra $A'$.
\end{lemma}
\begin{proof}
If $J$ is empty, then both $(T_xC)^J$ and $T^{W^J}_xC = T^{\N}_xC$ are equal to the terminal object $*$. If $|J| = 1$, then both sides are equal to the tangent space $T_xC$. For arbitrary $J$, we have $T^{W^J}C = TC \times_C \dots \times _C TC$ (with $|J|$ factors of $TC$) and the diagram (\ref{eq:TxC}) with $A = W^J$ is an iterated pullback of $|J|$ copies of the diagram (\ref{eq:TxC}) with $A = W = W^1$ over the diagram (\ref{eq:TxC}) with $A = \N = W^0$. It follows that $T^{W^J}_xC$ is a product of $|J|$ copies of $T_xC$. 

The second claim then follows from the fact that $T^{A'}$ preserves finite products (as the tangent structure on $\bcat{X}$ is cartesian) and also preserves the pullback square defining $T_xC$ (Definition~\ref{def:tangent-space}).
\end{proof}

\begin{definition} \label{def:TxC}
We define $\cat{T}_xC: \E \to \bcat{X}$\index{TAF@$\cat{T}_xC$, the differential object whose underlying object is the tangent space $T_xC$} by
\[ \cat{T}_xC(J) := T^{W^J}_xC \]
where $W^J$ denotes the Weil-algebra $\N \oplus \N^J$, the square-zero extension of $\N$ by $\N^J$. In other words, $\cat{T}_xC$ is the composite of functors
\[ \E \arrow{e,t}{W} \Weilinfty \arrow{e,t}{T^\bullet_xC} \bcat{X}. \]
Since $W(\N^J) = W^J$, this functors is finite-product-preserving. (In fact, it is the composite of two finite-product-preserving functors.)
\end{definition}

It remains to show that the functor $\cat{T}_xC: \E \to \bcat{X}$ commutes, up to coherent equivalence, with the tangent structures on $\E$ and $\bcat{X}$. We continue to follow the outline described in Corollary~\ref{cor:tangent-functor} by constructing the required map $\Delta^1 \times \Weilinfty \times \E \to \bcat{X}$, that is, a natural equivalence
\[ \alpha: (\cat{T}_xC)U^\bullet \weq T^\bullet(\cat{T}_xC). \]
The component of $\alpha$ for $A \in \Weilinfty$ and $M \in \E$ is a morphism
\[ \alpha_{A,M}: T^{W(A \otimes M)}_xC \to T^A(T^{W(M)}_xC) \]
which we construct as follows. Central to the definition is the following natural map of Weil-algebras.

\begin{definition} \label{def:W}
For a Weil-algebra $A$, with corresponding partial monoid $M_A$, and a finite set $J$, let $\phi_{A,J}$ denote the span in $\FPCM$ of the form
\[ \begin{diagram}
\node[2]{(M_A \times J)_+} \arrow{sw,t}{=} \arrow{se,t}{t} \\
\node{(M_A \times J)_+} \node[2]{_AM \times J_+}
\end{diagram} \]
where $t: (M_A \times J)_+ \to M_A \times J_+$ is the map of finite partial commutative monoids defined by
\[ t(m,j) = (m,j), \quad t(1) = (1,1). \]
\end{definition}

\begin{lemma} \label{lem:W-1}
The span $\phi_{M_A,J}$ of Definition~\ref{def:W} is a labelled Weil-algebra morphism
\[ \phi_{A,J}: W(A \otimes \N^J) \to A \otimes W(\N^J). \]
Moreover, this map is natural in both $A$ and $J$.
\end{lemma}
\begin{proof}
We have to check the conditions of Definition~\ref{def:weil-mor}. Condition (1) is satisfied by the identity on $(M_A \times J)_+$. For condition (2), we note that only products with the identity element are defined in $(M_A \times J)_+$. For any two non-identity elements $(m,j)$ and $(m'.,j')$, the product of $t(m,j)$ and $t(m',j')$ is defined if and only if both $m \cdot m'$ and $j \cdot j'$ are defined in $M_A$ and $J_+$ respectively. However, there are no non-identity products in $J_+$, so $t(m,j) \cdot t(m',j')$ is undefined too.

For the naturality claims, we have to define a suitable naturality square. So consider a labelled Weil-algebra morphism $\alpha: A \to A'$ given by the following span in $\FPCM$:
\[ \begin{diagram}
\node[2]{K} \arrow{sw,t}{s} \arrow{se,t}{t} \\
\node{M_A} \node[2]{M_{A'}}
\end{diagram} \]
and a morphism $\beta: J \to J'$ in $\E$ given by a span of finite sets:
\[ \begin{diagram}
\node[2]{L} \arrow{sw,t}{f} \arrow{se,t}{g} \\
\node{J} \node[2]{J'}
\end{diagram} \]
Then desired naturality square in $\Weilinfty$, of the form
\[ \begin{diagram}
\node{W(A \otimes \N^J)} \arrow{e,t}{\phi_{A,J}} \arrow{s,l}{W(\alpha \otimes \beta)} \node{A \otimes W(\N^J)} \arrow{s,r}{\alpha \otimes W(\beta)} \\
\node{W(A' \otimes \N^{J'})} \arrow{e,t}{\phi_{A',J'}} \node{A' \otimes W(\N^{J'})}
\end{diagram} \]
is given by an isomorphism of finite partial commutative monoids
\[ (M_A \times J)_+ \times_{(M_A \times J_+)} (K \times L_+) \isom (K \times L)_+ \times_{(M_{A'} \times J')_+} (M_{A'} \times J')_+. \]
Each side is canonically isomorphic to $(K \times L)_+$, so we have the desired naturality.
\end{proof}

\begin{definition} \label{def:W-mor}
Consider the following morphism of cospans in $\bcat{X}$:
\begin{equation} \label{eq:cospan} \begin{diagram}
\node{*} \arrow{e,t}{x} \arrow{s,l}{\sim} \node{C} \arrow{s} \node{T^{W(A \otimes \N^J)}C} \arrow{w} \arrow{s,r}{T^{\phi_{A,J}}C} \\
\node{T^A(*)} \arrow{e,t}{T^A(x)} \node{T^A(C)} \node{T^AT^{W(\N^J)}C} \arrow{w}
\end{diagram} \end{equation}
The left-hand square is induced by acting on the morphism $x: * \to C$ in $\bcat{X}$ by the morphism $\N \to A$ in $\Weilinfty$, and the left-hand vertical map is an equivalence since the tangent structure on $\bcat{X}$ is cartesian. The right-hand square is induced by acting on the object $C$ with the naturality square for $\phi$ of the form
\begin{equation} \label{eq:W-nat} \begin{diagram}
\node{\N = W(A \otimes *)} \arrow{s,l}{\phi_{A,*}} \node{W(A \otimes \N^J)} \arrow{w,t}{W(A \otimes f)} \arrow{s,r}{\phi_{A,J}} \\
\node{A = A \otimes W(*)} \node{A \otimes W(\N^J)} \arrow{w,t}{A \otimes W(f)}
\end{diagram} \end{equation}
where $* = \N_{\infty}^0$ denotes the trivial $E_\infty$-monoid and $f: \N^J \to *$ is the unique map. 

Let
\[ \alpha_{A,J}: T^{W(A \otimes \N^J)}_xC \to T^A(T^{W(\N^J)}_xC) \]
by the morphism induced by taking limits of the cospans in the diagram (\ref{eq:cospan}), using the fact that $T^A$ preserves the pullback square (\ref{eq:TxC}) in Lemma~\ref{lem:W}. All these constructions are natural in $A$ and $J$, so that we obtain a natural transformation $\alpha: (\cat{T}_xC)U^\bullet \to T^\bullet(\cat{T}_xC)$ as required.
\end{definition}

\begin{lemma} \label{lem:W-mor}
The morphism $\alpha_{A,\N^J}$ of Definition~\ref{def:W-mor} is an equivalence.
\end{lemma}
\begin{proof}
Since the left-hand vertical map is an equivalence, it is sufficient to show that the right-hand square in the diagram (\ref{eq:cospan}) is a pullback. Let us denote that right-hand square by $S_C(A,\N^J)$. In other words, $S_C(A,\N^J)$ is given by acting on $C$, via the tangent structure on $\bcat{X}$, by the square diagram in $\Weilinfty$ of the form
\begin{equation} \label{eq:W-mor} \begin{diagram}
\node{W(A \otimes \N^J)} \arrow{e,t}{\phi_{A,J}} \arrow{s} \node{A \otimes W(\N^J)} \arrow{s} \\
\node{W(A \otimes *)} \arrow{e,t}{\phi_{A,*}} \node{A \otimes W(*)}
\end{diagram} \end{equation}
It is in this proof that, crucially, we use the fact that the vertical lift pullbacks in $\Weilinfty$ are preserved by the tangent structure. Indeed, the square $S_C(W,\N)$ is precisely the vertical lift pullback given by acting on $C$ by the vertical lift pullback in $\Weilinfty$ from Proposition~\ref{prop:tangent-pb}. We will prove that $S_C(A,\N^J)$ is a pullback in $\bcat{X}$ for all $A \in \Weilinfty$ and all $J \in \E$.

When $A = \N$, the vertical maps in (\ref{eq:W-mor}) are equivalences, so the square is a pullback. When $J = \varnothing$, the horizontal maps are equivalences, so it is a pullback.

When $A = W$, $|J|=1$, the square is a vertical lift pullback as already noted.

We now claim that the square $S_C(W^n,\N^J)$ is a pullback of $n$ copies of the square $S_C(W,\N^J)$ over $S_C(\N,\N^J)$. We have to show that
\[ T^{W(W^n \otimes \N^J)}C \homeq T^{W(W \otimes \N^J)} \times _{T^{W(\N \otimes \N^J)}C} \dots \times_{T^{W(\N \otimes \N^J)}C} T^{W(W \otimes \N^J)}C \]
and
\[ T^{W^n \otimes W(\N^J)}C \homeq T^{W \otimes W(\N^J)} \times _{T^{\N \otimes W(\N^J)}C} \dots \times_{T^{\N \otimes W(\N^J)}C} T^{W \otimes W(\N^J)}C, \]
which both follow from the fact that the tangent structure preserves foundational pullbacks. We can therefore conclude that $S_C(W^n,\N)$ is a pullback in $\bcat{X}$.

Next we claim that the square $S_C(A,\N^k)$ is a pullback of $k$ copies of the square $S_C(A,\N)$ over $S_C(A,*)$. We have to show that
\[ T^{W(A \otimes \N^k)}C \homeq T^{W(A \otimes \N)}C \times_{T^{W(A \otimes *)}C} \dots \times_{T^{W(A \otimes *)}C} T^{W(A \otimes \N)}C \]
and
\[ T^{A \otimes W(\N^k)}C \homeq T^{A \otimes W(\N)}C \times_{T^{A \otimes W(*)}C} \dots \times_{T^{A \otimes W(*)}C} T^{A \otimes W(\N)}C, \]
which also follow from the fact that the tangent structure preserves foundational pullbacks (for $T^AC$ in the latter case).

We can now deduce that $S_C(W^n,\N^J)$ is a pullback for all $n \geq 2$ and all $J$ (and all $C$ in $\bcat{X}$). So now consider the square $S_C(A \otimes A',\N^J)$ for Weil-algebras $A,A'$. We can write this square as follows
\[ \begin{diagram}
\node{T^{W(A \otimes A' \otimes \N^J)}C} \arrow{e} \arrow{s} \node{T^{A \otimes W(A' \otimes \N^J)}C} \arrow{e} \arrow{s} \node{T^{A \otimes A' \otimes W(\N^J)}C} \arrow{s} \\
\node{T^{W(A \otimes A' \otimes \N^J)}C} \arrow{e}  \node{T^{A \otimes W(A' \otimes *)}C} \arrow{e} \node{T^{A \otimes A' \otimes W(*)}C}
\end{diagram} \]
The right-hand square is $S_{T^AC}(A',\N^J)$, and the left-hand square is $S_C(A,A' \otimes N)$. By induction, we can assume both of these are pullbacks, and so we deduce that the composed square is a pullback too. Thus, we can show that $S_C(A,\N^J)$ is a pullback for all $A,J,C$.
\end{proof}

We now have the basic data of a tangent functor from $\E$ to $\bcat{X}$, that is, a functor $\cat{T}_xC: \E \to \bcat{X}$ together with a natural equivalence
\[ \alpha: (\cat{T}_xC)U^\bullet \to T^\bullet(\cat{T}_xC). \]
Our goal is now to extend that data to a tangent functor from $\E$ to $\bcat{X}$ in the sense of the previous chapter. We will assume (without loss of generality) that $\bcat{X}$ is a strict tangent $\infty$-category, and hence a (fibrant) marked $\Weilinfty$-module. By Corollary~\ref{cor:tangent-functor}, our goal is to construct a map of marked $\Weilinfty$-modules of the form
\[ \cat{T}_xC: B\E\to \bcat{X}, \]
which extends the map of Definition~\ref{def:TxC}, and where $B\E$ denotes the simplicial bar construction on the marked $\Weilinfty$-module $\E$. 

We start by extending Definition~\ref{def:W} to simplexes of higher dimension.

\begin{definition} \label{def:W-n}
Let $A_n,\dots,A_1$ be Weil-algebras, and take $J \in \E$. We construct an $n$-simplex $\sigma_{A_n,\dots,A_1,J}$ in $\Weilinfty$ of the form
\[ W(A_n \otimes \dots \otimes A_1 \otimes \N^J) \to A_n \otimes W(A_{n-1} \otimes \dots \otimes \N^J) \to \dots \to A_n \otimes \dots \otimes A_1 \otimes W(\N^J). \]
in which each non-degenerate edge is the morphism
\[ A_n \otimes \dots \otimes A_{j+1} \otimes \phi_{A_j \otimes \dots \otimes A_{i+1}, A_i \otimes \dots \otimes A_1 \otimes \N^J} \]
for some $0 \leq i < j \leq n$. To be precise, recall from Proposition~\ref{prop:weil} that an $n$-simplex in $\Weilinfty$ can be identified with a diagram $\mathsf{J}_n \to \FPCM$, where $\mathsf{J_n}$ is the poset of intervals in $[0,n] = \{0,1,\dots,n\}$. We take
\[ \sigma_{A_n,\dots,A_1,J}([i,j]) := M_{A_n} \times \dots \times M_{A_{j+1}} \times (M_{A_j} \times \dots \times M_{A_1} \times J)_+ \]
with structure maps all of the form $(M \times K)_+ \to M \times K_+$ as in~\ref{def:W}. This definition reduces to Definition~\ref{def:W} in the case $n = 1$.  Moreover, the construction above is natural in $A_n,\dots,A_1$, and $J$, and extends to a functor
\[ \sigma_n : \Weilinfty^n \times \E \to \Fun(\Delta^n,\Weilinfty). \]
\end{definition}

Our next task is to build a higher-dimensional version of diagram (\ref{eq:cospan}). We start with the naturality squares (\ref{eq:W-nat}).

\begin{definition} \label{def:tilde-sigma}
Define a functor
\[ \tilde{\sigma}_n: \Weilinfty^n \times \E\to \Fun(\Delta^n,\Fun(\Delta^1,\bcat{X})) \]
by
\[ \tilde{\sigma}_n(A_n,\dots,A_1,J) := (\sigma_{A_n,\dots,A_1,J} \to \sigma_{A_n,\dots,A_1,\varnothing}) \]
where the right-hand side is the map of $n$-simplexes induced by the canonical morphism $J \to \varnothing$ in $\E$, given by the cospan $J \leftarrow \varnothing = \varnothing$ (i.e.\ the morphism from $J$ to the terminal object of $\E$). 
\end{definition}

\begin{lemma} \label{lem:W-BE}
The maps $\tilde{\sigma}_n: \Weilinfty^n \times \E \to \Fun(\Delta^n,\Fun(\Delta^1,\Weilinfty))$ of Definition~\ref{def:tilde-sigma} satisfy the face and degeneracy conditions described in Corollary~\ref{cor:tangent-functor} and hence determine a map of $\Weilinfty$-modules
\[ \tilde{\sigma}: B\E \to \Fun(\Delta^1,\Weilinfty), \]
which on vertices is given by
\[ (A,J) \mapsto (A \otimes W^J \to A) \]
where the right-hand side here is given by the canonical map $W^J \to \N$. Here $B\E$ is the simplicial bar construction on the $\Weilinfty$-module $\E$.
\end{lemma}
\begin{proof}
From Definition~\ref{def:W-n} we have the following equations:
\[ d^i(\sigma_{A_n,\dots,A_1,J}) =
	\begin{cases} \sigma_{A_n,\dots,A_2,M_{A_1} \times J} & \text{if $i = 0$}; \\ 
		\sigma_{A_n,\dots,A_{i+1} \otimes A_i,\dots,A_1,J} & \text{if $0 < i < n$}; \\
		A_n \otimes \sigma_{A_{n-1},\dots,A_1,J} & \text{if $i = n$}; \end{cases}\]
and
\[ s^i(\sigma_{A_n,\dots,A_1,J}) = \sigma_{A_n,\dots,A_{i+1},\N,A_i,\dots,A_1,J}. \]
These facts tell us that the maps $\sigma_n$ of Definition~\ref{def:W-n} satisfy the required conditions, and it follows that the maps $\tilde{\sigma}_n$ do too. 
\end{proof}

\begin{remark}
The map $\tilde{\sigma}$ does not send marked edges in $B\E$ to equivalences in $\Fun(\Delta^1,\Weilinfty)$, and so does not determine a tangent functor in the sense of Definition~\ref{def:tangent-functor}. (In fact, we can view $\sigma$ as a `lax' tangent functor analogous to Cockett and Cruttwell's `morphisms of tangent structure' of~\cite[2.7]{cockett/cruttwell:2014}, but we will not pursue that analogy here.)
\end{remark}

\begin{definition}
Given a pointed object $x: * \to C$ in $\bcat{X}$ for which the tangent space exists, and a morphism $\phi: A' \to A$ in $\Weilinfty$, consider the cospan in $\bcat{X}$ of the form
\[ \begin{diagram}
\node[2]{T^{A'}C} \arrow{s,r}{T^\phi(C)} \\
\node{T^A(*)} \arrow{e,t}{T^A(x)} \node{T^A(C)}
\end{diagram} \]
This construction determines a functor
\[ T^\bullet_*: \bcat{X}^T_* \times \Fun(\Delta^1,\Weilinfty) \to \Fun(\lrcorner,\bcat{X}) \]
where $\lrcorner$ is the indexing category for a cospan, i.e.\ two morphisms with common target. Lemma~\ref{lem:W} tells us that $T^\bullet_*$ takes values in the full subcategory
\[ \Fun^T(\lrcorner,\bcat{X}) \subseteq \Fun(\lrcorner,\bcat{X}) \]
consisting of those cospans which admit a pullback in $\bcat{X}$ and such that this pullback is preserved by $T^A$ for all $A \in \Weilinfty$. Notice that the functor $T^\bullet_*$ is a map of (strict) $\Weilinfty$-modules, where the $\Weilinfty$-action on the left-hand side comes from the pointwise action on $\Weilinfty$, and on the right-hand side from the pointwise action on $\bcat{X}$.

Combining $T^\bullet_*$ with the map $\tilde{\sigma}$ of Lemma~\ref{lem:W-BE}, we get a map of $\Weilinfty$-modules
\[ \cat{T}^{\lrcorner}: \bcat{X}^T_* \times B\E \to \Fun^T(\lrcorner,\bcat{X}). \]
On vertices, the map is given by
\[ (x:* \to C, (A,J)) \mapsto (T^A(*) \to T^A(C) \leftarrow T^A(T^{W^J}(C)). \]
Notice that Lemma~\ref{lem:W} tells us the pullback of this cospan is $T^A((T_xC)^J)$. Also notice that applying this map of simplicial sets to a given $x: * \to C$ and the (marked) edge in $B\E$ of the form
\[ (A,J) \to (\N,M_A \times J) \]
yields precisely diagram (\ref{eq:cospan}) of Definition~\ref{def:W-mor}.
\end{definition}

We would now like to compose $\cat{T}^{\lrcorner}$ with a suitable functor $\lim: \Fun(\lrcorner,\bcat{X}) \to \bcat{X}$ to get the desired maps $\cat{T}_xC: B\E \to \bcat{X}$. However, there is a technical wrinkle here, which is that we need those maps to commute \emph{strictly} with the $\Weilinfty$-actions, so we have to be careful to construct the pullbacks in a way that makes this happen.

\begin{definition} \label{def:TxC-diff}
Recall that we are writing $\Fun^T(\lrcorner,\bcat{X})$ for the full subcategory of $\Fun(\lrcorner,\bcat{X})$ consisting of those cospans which admit a pullback in $\bcat{X}$ which is preserved by $T^A$ for all $A \in \Weilinfty$. We give $\Fun^T(\lrcorner,\bcat{X})$ a marking where a morphism is marked if it induces an equivalence on pullbacks, and we give $\Fun^T(\lrcorner,\bcat{X})$ the pointwise $\Weilinfty$-module structure coming from the $\Weilinfty$-action on $\bcat{X}$. These choices make $\Fun^T(\lrcorner,\bcat{X})$ into a marked $\Weilinfty$-module.

Similarly, let $\Fun^T(\Delta^1 \times \Delta^1,\bcat{X})$ denote the full subcategory of square diagrams in $\bcat{X}$ which are pullbacks that are preserved by $T$, with the pointwise $\Weilinfty$-module structure, and where a morphism is marked if it is an equivalence on the initial vertex of $\Delta^1 \times \Delta^1$. Then $\Fun^T(\Delta^1 \times \Delta^1,\bcat{X})$ is also a marked $\Weilinfty$-module.

Now consider the map
\[ r: \Fun^T(\Delta^1 \times \Delta^1,\bcat{X}) \to \Fun^T(\lrcorner,\bcat{X}) \]
given by restriction along $\lrcorner \subseteq \Delta^1 \times \Delta^1$. The map $r$ is a trivial fibration in the Joyal model structure on simplicial sets, by~\cite[4.3.2.15]{lurie:2009}, and it is a map of marked $\Weilinfty$-modules. Thus $r$ is a trivial fibration in the projective model structure on $\ModWeil$.

Note that we also have a map of marked $\Weilinfty$-modules
\[ e: \Fun^T(\Delta^1 \times \Delta^1,\bcat{X}) \to \bcat{X} \]
given by evaluating at the initial vertex.

Finally, consider the following diagram of marked $\Weilinfty$-modules
\[ \begin{diagram}
\node[2]{\Fun^T(\Delta^1 \times \Delta^1,\bcat{X})} \arrow{s,lr,A}{\sim}{r} \arrow{e,t}{e} \node{\bcat{X}} \\
\node{\bcat{X}^T_* \times B\E} \arrow{e,b}{\cat{T}^{\lrcorner}} \arrow{ne,t,..}{\bar{\cat{T}}} \node{\Fun^T(\lrcorner,\bcat{X})}
\end{diagram} \]
Since $B\E$ is cofibrant, so is $\cat{X}^T_* \times B\E$, and hence a lift $\bar{\cat{T}}$ exists (and is unique up to equivalence). We can then define
\[ \cat{T}: \bcat{X}^T_* \times B\E \to \bcat{X} \]
to be the composite $e \circ \bar{\cat{T}}$.
\end{definition}

\begin{proposition} \label{prop:tangent-space-differential}
Let $\bcat{X}$ be a cartesian tangent $\infty$-category. Then there is a functor
\[ \cat{T}_*: \bcat{X}_*^{T} \to \Diff(\bcat{X}) \]
which maps a pointed object $x: * \to C$ that admits a tangent space $T_xC$ to a differential object in $\bcat{X}$ whose underlying object is that tangent space.
\end{proposition}
\begin{proof}
It follows from Lemma~\ref{lem:W-mor} that transpose to the map $\cat{T}$ of Definition~\ref{def:TxC-diff} is a functor
\[ \cat{T}_*: \bcat{X}^T_* \to \Hom_{\ModWeil}(B\E,\bcat{X}) \homeq \mathbf{Tan_\infty}(\E,\bcat{X}). \]
Lemma~\ref{lem:W} implies that $\cat{T}$ takes values in $\Diff(\bcat{X})$.
\end{proof}

\begin{corollary} \label{cor:tangent-space-differential}
Let $\bcat{X}$ be a cartesian tangent $\infty$-category. Then an object $D$ in $\bcat{X}$ admits a differential structure if and only if $D$ is equivalent to some tangent space $T_xC$.
\end{corollary}
\begin{proof}
The if direction follows from Proposition~\ref{prop:tangent-space-differential}. To see the only if direction, it is sufficient to note that if $D$ admits a differential structure, then we have a diagram
\[ \begin{diagram}
	\node{D} \arrow{e,t}{\langle \zeta,1 \rangle} \arrow{s,l}{\sim} \node{D \times D} \arrow{s,lr}{\sim}{\alpha_{\N}} \\
	\node{T_{\zeta}D} \arrow{e} \arrow{s} \node{T(D)} \arrow{s,r}{p} \\
	\node{*} \arrow{e,t}{\zeta} \node{D}
\end{diagram} \]
and the top-left map is an equivalence because the bottom square and composite squares are both pullbacks in $\bcat{X}$.
\end{proof}

\begin{remark}
Proposition~\ref{prop:tangent-space-differential} implies that a morphism $f: C \to D$ in $\bcat{X}$ induces a map
\[ T_xf: T_xC \to T_{f(x)}D \]
that preserves the differential structures on these tangent spaces. This map $T_xf$ is the analogue of the ordinary derivative of a map $f$ at a point $x$, and preservation of the differential structure is the analogue of this derivative being a linear map in a setting where there is no precise version of the vector space structure on an ordinary tangent space.
\end{remark}

\subsection*{Tangent $\infty$-categories and cartesian differential categories}

One of the motivations for Cockett and Cruttwell to study differential objects was to make a connection between cartesian tangent categories and the cartesian differential categories of Blute, Cockett and Seely~\cite{blute/cockett/seely:2009}. Roughly speaking, they show that for a cartesian tangent category $\bcat{X}$ in which every object has a canonical differential structure, there is a corresponding cartesian differential structure on $\bcat{X}$, and that every cartesian differential category arises in this way. We provide a generalization of that result to tangent $\infty$-categories.

\begin{definition}
	Let $\bcat{X}$ be a cartesian tangent $\infty$-category. We define a category $\widehat{h\Diff}(\bcat{X})$\index{hDiffX@$\widehat{h\Diff}(\bcat{X})$, the cartesian differential category of differential objects in a cartesian tangent $\infty$-category $\bcat{X}$} with objects the differential objects of $\bcat{X}$, and with morphisms from $\mathscr{D}$ to $\mathscr{D}'$ given by morphisms in the homotopy category of $\bcat{X}$ between the underlying objects $\mathscr{D}(*)$ and $\mathscr{D}'(*)$. This construction is not the homotopy category of $\Diff(\bcat{X})$ because we are now, like Cockett and Cruttwell, including \emph{all} morphisms between the underlying objects, not only those that commute with the differential structures.
\end{definition}

A full definition of \emph{cartesian differential structure} on a category $\bcat{X}$ can be found in~\cite[2.1.1]{blute/cockett/seely:2009}. Such a structure consists of an assignment, to each morphism $f: A \to B$ in $\bcat{X}$, a morphism $D(f): A \times A \to B$, referred to as the \emph{derivative} of $f$, satisfying a list of seven axioms.

The canonical example for the category whose objects are the Euclidean spaces $\mathbb{R}^n$, and whose morphisms are the smooth maps $f: \mathbb{R}^m \to \mathbb{R}^n$, has derivative $\nabla(f): \mathbb{R}^m \times \mathbb{R}^m \to \mathbb{R}^n$ given by the map
\[ \nabla(f)(a,v) = D_af(v) \]
where $D_af: \R^m \to \R^n$ denotes the derivative of $f$ at the point $a \in \R^m$.  In this example, the seven axioms describe standard properties of multivariable calculus such as the linearity of the total derivative and the equality of mixed partial derivatives.

The following result generalizes~\cite[4.11]{cockett/cruttwell:2014}. 

\begin{theorem} \label{thm:cartesian-differential-structure}
	Let $\bcat{X}$ be a cartesian tangent $\infty$-category. There is a cartesian differential structure on $\widehat{h\Diff}(\bcat{X})$ in which the monoid structure on an object is that inherited from its differential structure, and the derivative of a morphism $f: A \to B$ is given by the composite
	\[ \dgTEXTARROWLENGTH=3.5em \nabla(f): A \times A \arrow{e,tb}{\langle p_A,\hat{p}_A \rangle^{-1}}{\sim} TA \arrow{e,t}{T(f)} TB \arrow{e,t}{\hat{p}_B} B. \]
	where $\hat{p}_A$ and $\hat{p}_B$ are determined by the differential structures on $A$ and $B$ respectively.
\end{theorem}
\begin{proof}
	Let $\bcat{X}_d \subseteq \bcat{X}$ be the full subcategory of $\bcat{X}$ whose objects are those that admit a differential structure. We claim that $\bcat{X}_d$ is a tangent subcategory, i.e.\ is closed under the action of $\Weilinfty$. To see this claim, we note that, if $D$ admits a differential structure, $T^A(D) \homeq D^{M_A}$. A finite product of differential objects has a canonical differential structure, so $\bcat{X}_d$ is closed under finite products, and it follows that $\bcat{X}_d$ is a cartesian tangent subcategory of $\bcat{X}$.
	
	We now claim that, unusually, the homotopy category $h\bcat{X}_d$ inherits a tangent structure from that on $\bcat{X}_d$. The $\Weilinfty$-action on $\bcat{X}_d$ passes to an action of $\Weil$ on $h\bcat{X}_d$, so the key is to show that each of the foundational and vertical lift pullbacks in $\bcat{X}_d$ is also a pullback in $h\bcat{X}_d$.
	
	First, consider an object $D \in \bcat{X}_d$ and finite sets $J,J'$. We have to show that
	\[ \begin{diagram}
		\node{T^{W^{J \sqcup J'}}(D)} \arrow{e} \arrow{s} \node{T^{W^J}(D)} \arrow{s} \\
		\node{T^{W^{J'}}(D)} \arrow{e} \node{D}
	\end{diagram} \]
	is a pullback in $h\bcat{X}_d$. A differential structure on $D$ determines equivalences in the $\infty$-category $\bcat{X}_d$, and hence isomorphisms in $h\bcat{X}_d$ of the form
	\[ T^A(D) = T^A \cat{D}(\ast) \homeq \cat{D}(A \otimes \ast) = \cat{D}(M_A) \homeq \cat{D}(\ast)^{M_A} = D^{M_A} \]
	where the penultimate equivalence comes from the fact that $\cat{D}$ preserves finite products, which are given in $\E$ by the disjoint union of sets. Note also that products in an $\infty$-category pass to the homotopy category.  The above equivalences are natural in $A$, which implies that we can replace the square above with a corresponding diagram
	\[ \begin{diagram}
			\node{D^{J \sqcup J' \sqcup \{1\}}} \arrow{e} \arrow{s} \node{D^{J \sqcup \{1\}}} \arrow{s} \\
			\node{D^{J' \sqcup\{1\}}} \arrow{e} \node{D}
	\end{diagram} \]
	where each map is the evident projection. A square of this form is a pullback in any category with finite products.
	
	For the vertical lift axiom, we have to consider the following diagram in $h\bcat{X}_d$
	\[ \begin{diagram} \dgARROWLENGTH=4em
		\node{D^3} \arrow{e,t}{\langle \pi_1,\pi_2,\zeta !,\pi_3 \rangle} \arrow{s,l}{\pi_3} \node{D^4} \arrow{s,r}{\langle \pi_3,\pi_4 \rangle} \\
		\node{D} \arrow{e,t}{\langle \zeta !,1 \rangle} \node{D^2}
	\end{diagram} \]
	where $\zeta: * \to D$ is the zero map for the differential structure on $D$. A diagram of this form is a pullback in any category (regardless of the map $\zeta$).
	
	We have completed the check that $h\bcat{X}_d$ has a cartesian tangent structure, and we now apply~\cite[4.11]{cockett/cruttwell:2014} to this structure. This provides a cartesian differential structure on the category $\widehat{\Diff}(h\bcat{X}_d)$ whose objects are the differential objects in $h\bcat{X}_d$, and whose morphisms are morphisms in $h\bcat{X}_d$ between the underlying objects.
	
	We complete the proof of our theorem by constructing an embedding of the category $\widehat{h\Diff}(\bcat{X})$ into $\widehat{\Diff}(h\bcat{X}_d)$. Given a differential object in $\bcat{X}$ with underlying object $D$, we obtain a differential structure on $D$ in $h\bcat{X}_d$ in the same manner as in the proof of Proposition~\ref{prop:diff}. Since morphisms in each case are homotopy classes of maps in $\bcat{X}$, we see that this construction determines a fully faithful functor
	\[ V: \widehat{h\Diff}(\bcat{X}) \to \widehat{\Diff}(h\bcat{X}_d). \]
	Since $\widehat{h\Diff}(\bcat{X})$ is closed under finite products, we obtain the desired cartesian differential structure by restriction along $V$. The formula for the derivative $\nabla(f)$ follows from that in $\widehat{\Diff}(h\bcat{X}_d)$.
\end{proof}

\chapter{Tangent Structures in and on an $(\infty,2)$-Category} \label{sec:infty2}

The goal of this chapter is to generalize our notion of `tangent $\infty$-category' to $(\infty,2)$-categories in two different ways. Firstly, we observe that one can easily extend Definition~\ref{def:tangent-structure-infty} to a notion of \emph{tangent object} \textbf{in} any $(\infty,2)$-category, and that this context provides the natural setting for such a notion. Secondly, we apply that general notion to give a definition of \emph{tangent $(\infty,2)$-category}, i.e.\ a tangent structure \textbf{on} an $(\infty,2)$-category, or tangent object in the $(\infty,2)$-category of $(\infty,2)$-categories. In that last case, we can also specialize to those $(\infty,2)$-categories which are the Duskin nerve of an ordinary bicategory, and thus we obtain a notion of \emph{tangent bicategory}.

Lanfranchi~\cite{lanfranchi:2023a} has independently made an equivalent definition of tangent object in the context of ordinary $2$-categories. Our notion extends that definition to $(\infty,2)$-categories.

We refer to Appendix~\ref{sec:infty-bicat} for background about the theory of $(\infty,2)$-categories, for which we use Lurie's notion of `$\infty$-bicategory' from~\cite{lurie:kerodon} as our principal model.

\subsection*{Tangent objects in an $\infty$-bicategory}

Our definition of tangent object in an arbitrary $\infty$-bicategory is based on the observation that a monoidal $\infty$-category can be identified with an $\infty$-bicategory that has a single object. Thus, associated to $\Weilinfty$ is an $\infty$-bicategory $\uWI$. Moreover, a $\Weilinfty$-module $\infty$-category then corresponds to a map of $\infty$-bicategories
\[ \uWI \to \mathbf{CAT}_\infty \]
into the $\infty$-bicategory of $\infty$-categories. This bicategorical approach to modules over a monoidal $\infty$-category is described in more detail in Appendix~\ref{sec:module}.

There is one subtlety to the definition, which Lanfranchi attributes to Lucyshyn-Wright; see~\cite[4.2]{lanfranchi:2023a}. We need a bicategorical interpretation of the condition that a tangent structure map
\[ T^\bullet: \Weilinfty \to \Fun(\bcat{X},\bcat{X}) \]
sends the tangent pullbacks in $\Weilinfty$ to `pointwise' pullbacks in $\Fun(\bcat{X},\bcat{X})$, i.e.\ those which are given pointwise by pullbacks in $\bcat{X}$; see Remark~\ref{rem:pointwise-limit}. To state that condition, we use the following definition from~\cite[4.1]{lanfranchi:2023a}.

\begin{definition} \label{def:pointwise-limit}
Let $\mathbf{C}$ be an $\infty$-bicategory, and let $D: \mathbb{I} \to \mathbf{C}(x,y)$ be a limit diagram in one of the mapping $\infty$-categories of $\mathbf{C}$. Then $D$ is a \emph{pointwise}\index{pointwise limit in an $\infty$-bicategory} limit if it is preserved by the induced functor
\[  \mathbf{C}(f,y): \mathbf{C}(x,y) \to \mathbf{C}(x',y) \]
for every $1$-morphism $f: x' \to x$ in $\mathbf{C}$.
\end{definition}

\begin{example} \label{ex:pointwise-limit}
When $\mathbf{C} = \mathbf{CAT}_\infty$, a limit diagram in
\[ \mathbf{CAT}_\infty(\bcat{X},\bcat{Y}) = \Fun(\bcat{X},\bcat{Y}) \]
is pointwise if and only if it is given by a limit in $\bcat{Y}$ for each object in $\bcat{X}$. If the $\infty$-category $\bcat{Y}$ admits all limits of a certain shape, then limits of that shape in $\mathbf{CAT}_\infty(\bcat{X},\bcat{Y})$ will automatically be pointwise, but in general that need not be the case.
\end{example}

We can now state our definition of tangent object in an $\infty$-bicategory using the notion from Definition~\ref{def:mod-infty-bicat} of a module over a monoidal $\infty$-category in an arbitrary $\infty$-bicategory.

\begin{definition} \label{def:tangent-object}
Let $\mathbf{C}$ be an $\infty$-bicategory. A \emph{tangent object}\index{tangent object in an $\infty$-bicategory} in $\mathbf{C}$ is a $\Weilinfty$-module $\bcat{X}$ in $\mathbf{C}$, for which the corresponding monoidal functor
\[ T^\bullet: \Weilinfty \to \mathbf{C}(\bcat{X},\bcat{X}) \]
sends the tangent pullbacks in $\Weilinfty$ to pointwise pullbacks in $\mathbf{C}(\bcat{X},\bcat{X})$.

The \emph{$\infty$-bicategory of tangent objects in $\mathbf{C}$}\index{TanC@$\Tan(\mathbf{C})$, the $\infty$-bicategory of tangent objects in the $\infty$-bicategory $\mathbf{C}$} is the full sub-bicategory
\[ \Tan(\mathbf{C}) \subseteq \Mod_{\Weilinfty}(\mathbf{C}) \]
whose objects are the tangent objects in the sense above.
\end{definition}

\begin{remark}
A $\Weilinfty$-module in $\mathbf{C}$ corresponds to a functor of $\infty$-bicategories
\[ \bcat{X}(-): \uWI \to \mathbf{C}, \]
where $\uWI$ is the one-object $\infty$-bicategory with mapping $\infty$-category
\[ \uWI(\bullet,\bullet) := \Weilinfty. \]
The functor $\bcat{X}(-)$ then corresponds to a monoidal functor
\[ T^\bullet: \Weilinfty \to \mathbf{C}(\bcat{X}(\bullet),\bcat{X}(\bullet)), \]
which provides the action of $\Weilinfty$ on $\bcat{X} = \bcat{X}(\bullet)$.
\end{remark}

\begin{example}
Taking $\mathbf{C} = \mathbf{CAT}_{\infty}$, the $\infty$-bicategory of $\infty$-categories, we have
\[ \Tan(\mathbf{CAT}_\infty) = \Tan_\infty \]
the $\infty$-bicategory of tangent $\infty$-categories previously introduced.
\end{example}

\begin{example}
Taking $\mathbf{C} = \mathbf{CAT}$, the simplicial nerve of the $2$-category of ordinary categories, functors, and natural transformations, we have
\[ \Tan(\mathbf{CAT}) \homeq \Tan \]
the (simplicial nerve of the) $2$-category of ordinary tangent categories, tangent functors, and tangent natural transformations.
\end{example}

\begin{example} \label{ex:tangent-infty-bicategories}
Taking $\mathbf{C} = \mathbf{CAT}_{(\infty,2)}$, the $\infty$-bicategory of $\infty$-bicategories, we obtain an $\infty$-bicategory
\[ \Tan_{(\infty,2)} := \Tan(\mathbf{CAT}_{(\infty,2)}) \]
of \emph{tangent $\infty$-bicategories}.\index{TanCAT2@$\Tan_{(\infty,2)}$, the $\infty$-bicategory of tangent $\infty$-bicategories}\index{tangent $\infty$-bicategory}
\end{example}

\begin{example}
Let $\mathbf{C}$ be the full sub-bicategory $\mathbf{CAT}_2 \subseteq \mathbf{CAT}_{(\infty,2)}$ whose objects are the (Duskin nerves of) ordinary bicategories. Then we obtain an $\infty$-bicategory
\[ \Tan_{2} := \Tan(\mathbf{CAT}_2) \]
of \emph{tangent bicategories}\index{Tan2@$\Tan_2$, the $\infty$-bicategory of tangent bicategories}\index{tangent bicategory}. For a bicategory $\bcat{B}$, the $\infty$-category $\Hom_{\mathbf{CAT}_2}(\bcat{B},\bcat{B})$ is a (strict) monoidal $(2,1)$-category, and a tangent structure on $\bcat{B}$ in this sense can be identified with a functor of monoidal $(2,1)$-categories of the form
\[ \Weilinfty \to \Hom_{\mathbf{CAT}_2}(\bcat{B},\bcat{B}) \]
which maps the tangent pullbacks to pointwise pullbacks.
\end{example}

\subsection*{Tangent $\infty$-bicategories}

Our goal in the rest of this chapter is to elaborate on the notion of a tangent $\infty$-bicategory associated to Example~\ref{ex:tangent-infty-bicategories}. According to Definition~\ref{def:tangent-object}, a tangent structure on the $\infty$-bicategory $\bcat{B}$ consists of a map of $\infty$-bicategories
\[ \uWI \to \mathbf{CAT}_{(\infty,2)} \]
which sends $\bullet$ to $\bcat{B}$, and hence induces a monoidal functor
\[ T^\bullet: \Weilinfty \to \Fun(\bcat{B},\bcat{B})^{\sim}. \]
That functor is required to map the tangent pullbacks to pointwise pullbacks. In order to understand this condition, we need to get a handle on the pointwise pullbacks, i.e.\ on pullbacks in $\infty$-categories of the form $\Fun(\bcat{X},\bcat{B})^{\sim}$. 

We do not have a complete characterization of these pointwise pullbacks, but we can give a sufficient condition for a diagram in $\Fun(\bcat{B},\bcat{B})$ to be a pointwise pullback, which will be enough for the application we have in mind. That condition is based on the following notion of pullback in an $\infty$-bicategory $\bcat{Y}$. 

\begin{definition} \label{def:2-pullback}
	Let $\bcat{Y}$ be an $\infty$-bicategory, and consider a commutative diagram:
	\[ \begin{diagram}
		\node{C} \arrow{e} \arrow{s} \arrow{se} \node{C_1} \arrow{s} \\
		\node{C_2} \arrow{e} \node{C_0}
	\end{diagram} \]
	in $\bcat{Y}$. We say that this diagram is a \emph{homotopy $2$-pullback}\index{homotopy 2-pullback in an $\infty$-bicategory} if, for every $D \in \bcat{Y}$, the induced diagram of mapping $\infty$-categories
	\[ \begin{diagram}
		\node{\bcat{Y}(D,C)} \arrow{e} \arrow{s} \node{\bcat{Y}(D,C_1)} \arrow{s} \\
		\node{\bcat{Y}(D,C_2)} \arrow{e} \node{\bcat{Y}(D,C_0)}
	\end{diagram} \]
	is a pullback of $\infty$-categories (that is, a homotopy pullback in the marked model structure on marked simplicial sets).
\end{definition}

Our condition is given by the following lemma.

\begin{lemma} \label{lem:htpy-pb}
Let $\bcat{X},\bcat{Y}$ be $\infty$-bicategories, and let $D$ denote a commutative diagram
\[ \begin{diagram}
	\node{F} \arrow{e} \arrow{s} \node{F_0} \arrow{s} \\
	\node{F_1} \arrow{e} \node{F_2}
\end{diagram} \]
in $\Fun(\bcat{X},\bcat{Y})^{\sim}$. Suppose that for each object $c \in \bcat{X}$, the diagram $D(c)$ given by
\[ \begin{diagram}
	\node{F(c)} \arrow{e} \arrow{s} \node{F_0(c)} \arrow{s} \\
	\node{F_1(c)} \arrow{e} \node{F_2(c)}
\end{diagram} \]
is a homotopy $2$-pullback in $\bcat{Y}$. Then $D$ is a pointwise pullback in the $\infty$-bicategory $\mathbf{CAT}_\infty$.
\end{lemma}
\begin{proof}
Let $G: \bcat{X}' \to \bcat{X}$ be an arbitrary map of $\infty$-bicategories. We have to show that the diagram $DG$ given by
\[ \begin{diagram}
	\node{FG} \arrow{e} \arrow{s} \node{F_0G} \arrow{s} \\
	\node{F_1G} \arrow{e} \node{F_2G}
\end{diagram} \]
is a pullback in the $\infty$-category $\Fun(\bcat{X}',\bcat{Y})^{\sim}$.

To get a handle on this $\infty$-category, we can identify it as a full subcategory of the $\infty$-category associated with a certain model category. First, recall that any $\infty$-bicategory is equivalent to the nerve of a category (strictly) enriched in $\infty$-categories (that is, in simplicial sets satisfying the inner horn conditions). Let us write
\[ \bcat{X}' \homeq N(\mathsf{X}'), \quad \bcat{Y} \homeq N(\mathsf{Y}) \]
where $\mathsf{X}',\mathsf{Y}$ are simplicial categories of that form. We will view $\mathsf{X}',\mathsf{Y}$ as enriched in \emph{marked} simplicial sets, where the marked edges in an $\infty$-category are (as usual) given by the equivalences.

Now consider the (ordinary) category $\Fun_{\msset}(\mathsf{X}' \times \mathsf{Y}^{op}, \msset)$ of $\msset$-enriched functors of the form
\[ \mathsf{X}' \times \mathsf{Y}^{op} \to \msset. \]
That category has a projective model structure in which a morphism (that is, a $\msset$-enriched natural transformation) is a weak equivalence or fibration if and only if each of its components is a weak equivalence or fibration (respectively) in the marked model structure on $\msset$. Associated with that projective model structure is an $\infty$-category
\[ \Fun_{\msset}(\mathsf{X}' \times \mathsf{Y}^{op}, \msset)^{\circ} \]
whose objects are the fibrant-cofibrant objects.

It follows from~\cite[A.3.4.14]{lurie:2009} that we can identify $\Fun(\bcat{X}',\bcat{Y})^{\sim}$ with the full subcategory of $\Fun_{\msset}(\mathsf{X}' \times \mathsf{Y}^{op}, \msset)^{\circ}$ whose objects are those (fibrant-cofibrant) $\msset$-enriched functors
\[ f: \mathsf{X}' \times \mathsf{Y}^{op} \to \msset \]
such that for each $x \in \mathsf{X}'$, the functor $f(x,-): \mathsf{Y}^{op} \to \msset$ is in the essential image of the Yoneda embedding for the $\msset$-enriched category $\mathsf{Y}$. In particular, a map of $\infty$-bicategories $H: \bcat{X}' \to \bcat{Y}$ corresponds to a functor of the form
\[ H': \mathsf{X}' \times \mathsf{Y}^{op} \to \msset; \quad H'(x,y) \homeq \Hom_{\bcat{Y}}(y,H(x)) \]

The diagram $DG$ now determines a diagram $(DG)'$ in the model category $\Fun_{\msset}(\mathsf{X}' \times \mathsf{Y}^{op}, \msset)$. For each $x \in \mathsf{X}'$ and $y \in \mathsf{Y}$, the diagram $(DG)'(x,y) \homeq D'(G(x),y)$ is, by hypothesis, a pullback of $\infty$-categories, i.e.\ a homotopy pullback in the marked model structure on $\msset$. It follows that $(DG)'$ is a homotopy pullback in the projective model structure on $\Fun_{\msset}(\mathsf{X}' \times \mathsf{Y}^{op}, \msset)$, and hence a pullback in the corresponding $\infty$-category
\[ \Fun_{\msset}(\mathsf{X}' \times \mathsf{Y}^{op}, \msset)^{\circ}. \]
Therefore $DG$ is the pullback in the full subcategory $\Fun(\bcat{X}',\bcat{Y})^{\sim}$, as required.
\end{proof}

Lemma~\ref{lem:htpy-pb} then has the following immediate consequence.

\begin{proposition} \label{prop:2-pullback}
Let $\bcat{B}$ be an $\infty$-bicategory, and let
\[ T^\bullet: \Weilinfty \to \Fun(\bcat{B},\bcat{B})^{\sim} \]
be a monoidal functor. Suppose that for each $X \in \bcat{B}$, the functor $T^\bullet(-)(X): \Weilinfty \to \bcat{B}$ maps the tangent pullbacks to homotopy $2$-pullbacks in $\bcat{B}$. Then $T^\bullet$ determines a tangent structure on $\bcat{B}$.
\end{proposition}

\begin{corollary} \label{cor:tangent-bicat-quasi}
Let $(\bcat{B},T^\bullet)$ be a tangent $\infty$-bicategory that satisfies the condition of Proposition~\ref{prop:2-pullback}. Then $T^\bullet$ determines a tangent structure on the core $\infty$-category $\bcat{B}^{\sim}$.
\end{corollary}
\begin{proof}
The core construction $(-)^{\sim}: \mathbf{CAT}_{(\infty,2)} \to \mathbf{CAT}_\infty$ determines a monoidal functor
\[ (-)^{\sim}: \Fun(\bcat{B},\bcat{B})^{\sim} \to \Fun(\bcat{B}^{\sim},\bcat{B}^{\sim}). \]
Composing with the structure map $T^\bullet$ we get the required action $(T^\bullet)^{\sim}$ of $\Weilinfty$ on the $\infty$-category $\bcat{B}^{\sim}$.

To see that this action is a tangent structure on $\bcat{B}^{\sim}$, let $X$ be an arbitrary object of $\bcat{B}^{\sim}$, i.e.\ an object of $\bcat{B}$. By assumption $T^\bullet(-)(X)$ maps tangent pullbacks in $\Weilinfty$ to homotopy $2$-pullbacks in $\bcat{B}$. The core groupoid functor $(-)^{\simeq}$ maps homotopy pullbacks in the marked model structure to homotopy pullbacks in the Quillen model structure. It follows that $(-)^{\sim}$ maps a homotopy $2$-pullback in $\bcat{B}$ to a pullback in the $\infty$-category $\bcat{B}^{\sim}$. It follows that $(T^\bullet)^{\sim}(-)(X)$ preserves the tangent pullbacks.
\end{proof}

\part{A Tangent $\infty$-Category of $\infty$-Categories}

\chapter{Goodwillie Calculus and the Tangent Bundle Functor} \label{sec:goodwillie}

Up to this point, we have been developing the general theory of tangent $\infty$-categories, but now we start to look at the specific tangent structure which motivated this work. The construction and study of this structure, which we term the \emph{Goodwillie tangent structure}, occupies the remainder of this paper. We start with a brief review of Goodwillie's calculus of functors, ideas from which will permeate the definitions and proofs to come. Most of this chapter is based on Lurie's development~\cite[Ch. 6]{lurie:2017} of Goodwillie's ideas in the general context of $\infty$-categories.

The central notion in Goodwillie calculus is the \emph{Taylor tower} of a functor; that is, a sequence of approximations to the functor that plays the role of the Taylor series in ordinary calculus. To describe the Taylor tower, we recall Goodwillie's analogues of polynomials in the context of functors. We need the following preliminary definition.

\begin{definition} \label{def:cube}
	An \emph{$n$-cube}\index{cube in an $\infty$-category} in an $\infty$-category $\cat{C}$ is a diagram $\mathcal{X}$ in $\cat{C}$ indexed by the poset $\mathcal{P}[n]$ of subsets of $[n] = \{1,\dots,n\}$. Such a cube is \emph{strongly cocartesian}\index{cube in an $\infty$-category!strongly cocartesian} if each $2$-dimensional face is a pushout, and is \emph{cartesian}\index{cube in an $\infty$-category!cartesian} if the cube as a whole is a limit diagram, i.e.\ the map
	\[ \mathcal{X}(\varnothing) \to \holim_{\varnothing \neq S \subseteq [n]} \mathcal{X}(S) \]
	is an equivalence in $\cat{C}$.
\end{definition}

\begin{definition}[{Goodwillie~\cite[3.1]{goodwillie:1991}}] \label{def:excisive}
	Let $F: \cat{C} \to \cat{D}$ be a functor between $\infty$-categories. We say that $F$ is \emph{$n$-excisive}\index{n-excisive@$n$-excisive functor} if it maps each strongly cocartesian $(n+1)$-cube in $\cat{C}$ to a cartesian cube in $\cat{D}$. In particular, $F$ is $1$-excisive (or, simply, \emph{excisive})\index{excisive functor} if it maps pushout squares in $\cat{C}$ to pullback squares in $\cat{D}$, and $F$ is $0$-excisive if and only if it is constant (up to equivalence). We write\index{Excn@$\Exc^n$, $\infty$-category of $n$-excisive functors}
	\[ \Exc^n(\cat{C},\cat{D}) \subseteq \Fun(\cat{C},\cat{D}) \]
	for the full subcategory of the functor $\infty$-category whose objects are the $n$-excisive functors.
\end{definition}

\begin{remark}
	Motivation for this definition comes from algebraic topology; the property of being excisive is closely related to the excision property in the Eilenberg-Steenrod axioms for homology, and excisive functors on the $\infty$-category of pointed topological spaces are closely related to homology theories.
\end{remark}

One of Goodwillie's key constructions is of a universal $n$-excisive approximation for functors between suitable $\infty$-categories. In order to state a condition for this approximation to exist, we introduce the following definition from Lurie~\cite[6.1.1.6]{lurie:2017}.

\begin{definition} \label{def:diff-cat}
	An $\infty$-category $\cat{C}$ is \emph{differentiable}\index{differentiable $\infty$-category} if it admits finite limits and sequential colimits (i.e.\ colimits along countable sequences of composable morphisms), and those limits and colimits commute.
\end{definition}

\begin{remark}
	We caution the reader to distinguish carefully between the words `differentiable', as applied to an $\infty$-category in the above definition, and `differential', which we introduced in Chapter~\ref{sec:differential} to refer to the structure inherent on a tangent space in any tangent $\infty$-category. In Chapter~\ref{sec:catdiff-diff} we make that distinction particularly challenging by considering differential objects in the $\infty$-category of differentiable $\infty$-categories, or `differential differentiable $\infty$-categories' if you prefer.
	
	We did consider introducing a new term for what Lurie calls `differentiable' $\infty$-categories in order to avoid the clash of terminology described above. For example, we might have called them `Goodwillie' $\infty$-categories in order to emphasise their role in the theory of Goodwillie calculus. They are also closely related to the `precontinuous categories' of~\cite[1.2]{adamek/lawvere/rosicky:2003}, and we might have introduced the term `finitely precontinuous $\infty$-category' for the notion relevant here.
	
	Ultimately, we decided to adopt Lurie's terminology despite the potential for confusion. While that terminology does not yet appear to be well-established across the literature, our work involves extensive references to~\cite{lurie:2017}, where the term `differentiable' was introduced and which itself is a widely-read work. We found it preferable to retain the same language and to trust that the reader will carefully distinguish between the terms `differentiable' and `differential'.
\end{remark}

The following result is due to Lurie~\cite[6.1.1.10]{lurie:2017} in this generality, but its proof is based on Goodwillie's original construction from~\cite[Sec.\ 1]{goodwillie:2003}.

\begin{proposition} \label{prop:Pn}
	Let $\cat{C},\cat{D}$ be $\infty$-categories such that $\cat{C}$ has finite colimits and a terminal object, and $\cat{D}$ is differentiable. Then the inclusion
	\[ \Exc^n(\cat{C},\cat{D}) \into \Fun(\cat{C},\cat{D}) \]
	admits a left adjoint $P_n$ which preserves finite limits.\index{Pn@$P_n$, $n$-excisive approximation}
\end{proposition}

Proposition~\ref{prop:Pn} implies that a functor $F: \cat{C} \to \cat{D}$ admits a universal $n$-excisive approximation map $p_n: F \to P_nF$, with the property that any map $F \to G$, where $G$ is $n$-excisive, factors uniquely (i.e.\ up to contractible choice) as $F \to P_nF \to G$.

\begin{definition} \label{def:Taylor}
	Let $F:\cat{C} \to \cat{D}$ be a functor between $\infty$-categories that satisfies the conditions in Proposition~\ref{prop:Pn}. The \emph{Taylor tower}\index{Taylor tower of a functor} of $F$ is the sequence of natural transformations
	\[ F \to \dots \to P_nF \to P_{n-1}F \to \dots \to P_1F \to P_0F = F(*) \]
	determined by the universal property of each $P_nF$ and the observation that an $n$-excisive functor is also $(n+1)$-excisive, for each $n$.
\end{definition}

We need the following generalization of~\cite[3.1]{arone/ching:2011}.

\begin{lemma} \label{lem:PnGF}
	Let $F: \cat{C} \to \cat{D}$ and $G: \cat{D} \to \cat{E}$ be functors between $\infty$-categories that satisfy the conditions of Proposition~\ref{prop:Pn}, and suppose that $F$ preserves the terminal object in $\cat{C}$. Then the map
	\[ P_n(GF) \to P_n((P_nG)F) \]
	induced by $p_n: G \to P_nG$, is an equivalence.
\end{lemma}
\begin{proof}
	Following Goodwillie~\cite{goodwillie:2003}, Lurie defines~\cite[6.1.1.27]{lurie:2017}
	\[ P_nF := \colim_k T_n^kF \]
	where
	\[ T_nF(X) \homeq \lim_{\varnothing \neq S \subseteq [n]} F(C_S(X)) \]
	and $C_S(X) \homeq \hocofib(\bigvee_S X \to X)$ is a model for the `join' of the object $X \in \cat{C}$ with a finite set $S$. Since $F$ preserves the terminal object, we get canonical maps
	\[ C_S(F(X)) \to F(C_S(X)) \]
	and hence
	\[ T_n(G)F \to T_n(GF) \]
	and therefore also a map (commuting with the canonical maps from $GF$) of the form
	\[ P_n(G)F \to P_n(GF). \]
	Since $P_n(GF)$ is $n$-excisive, this must factor via a map
	\[ P_n((P_nG)F) \to P_n(GF) \]
	which we claim to be an inverse to the map in the statement of the Lemma. It follows from the universal property of $P_n$ that the composite
	\[ P_n(GF) \to P_n((P_nG)F) \to P_n(GF) \]
	is an equivalence. In the following diagram
	\[ \begin{diagram}
		\node{P_n(GF)} \arrow{e} \arrow{s} \node{P_n((P_nG)F)} \arrow{e} \arrow{s} \node{P_n(GF)} \arrow{s} \\
		\node{P_n((P_nG)F)} \arrow{e} \node{P_nP_n((P_nG)F)} \arrow{e} \node{P_n((P_nG)F)}
	\end{diagram} \]
	the horizontal composites are equivalences, and the middle vertical map is an equivalence. It follows that the end map is a retract of an equivalence, hence is itself an equivalence.
\end{proof}

A significant role in our later constructions is played by a multivariable version of Goodwillie's calculus, so we describe that briefly, too. More details are in~\cite[6.1.3]{lurie:2017}.

\begin{definition} \label{def:multi-excisive}
	Let $F: \cat{C}_1 \times \cat{C}_2 \to \cat{D}$ be a functor between $\infty$-categories such that each of $\cat{C}_1,\cat{C}_2$ admits finite colimits and a terminal object, and $\cat{D}$ is differentiable. We say that $F$ is \emph{$n_1$-excisive in its first variable} (and similarly for its second variable) if, for all $X_2 \in \cat{C}_2$, the functor $F(-,X_2): \cat{C}_1 \to \cat{D}$ is $n_1$-excisive. We say that $F$ is \emph{$(n_1,n_2)$-excisive}\index{n1n2-excisive@$(n_1,n_2)$-excisive functor} if it is $n_i$-excisive in its $i$-th variable, for each $i$. We write\index{Excn1n2@$\Exc^{n_1,n_2}(\cat{C}_1 \times \cat{C}_2,\cat{D})$, the $\infty$-category of $(n_1,n_2)$-excisive functors}
	\[ \Exc^{n_1,n_2}(\cat{C}_1 \times \cat{C}_2,\cat{D}) \subseteq \Fun(\cat{C}_1 \times \cat{C}_2,\cat{D}) \]
	for the full subcategory on the $(n_1,n_2)$-excisive functors. The inclusion of this subcategory has a left adjoint $P_{n_1,n_2}$ given by applying the functor $P_{n_i}$ of Proposition~\ref{prop:Pn} to each variable in turn, keeping the other variable constant.
\end{definition}

We now begin our description of the Goodwillie tangent structure. First, we introduce the underlying $\infty$-category for this structure.

\begin{definition} \label{def:Catdiff-informal}
	Let $\Catinf$\index{Cat01@$\Catinf$, the $\infty$-category of $\infty$-categories} be Lurie's model~\cite[3.0.0.1]{lurie:2009} for the $\infty$-category of $\infty$-categories\footnote{Recall that we are not being explicit about the various size restrictions on our $\infty$-categories, but for this definition to make sense, we of course require the objects in the $\infty$-category $\Catinf$ to be restricted to be smaller than $\Catinf$ itself.}, and let $\Catdiff \subseteq \Catinf$\index{Cat0diff@$\Catdiff$, the $\infty$-category of differentiable $\infty$-categories} be the subcategory whose objects are the differentiable $\infty$-categories and whose morphisms are the functors that preserve sequential colimits.
\end{definition}

The tangent bundle functor for the Goodwillie tangent structure is an endofunctor $T: \Catdiff \to \Catdiff$ defined by Lurie in~\cite[7.3.1.10]{lurie:2017}. That definition, and much of the rest of this paper, depends heavily on the $\infty$-category $\finbased$ of `finite pointed spaces', so let us be entirely explicit about that object.

\begin{definition} \label{def:based}
	We say that a Kan complex is \emph{finite}\index{finite Kan complex} if it is weak homotopy equivalent to a simplicial set with finitely many nondegenerate simplexes. Let $\finsset$ be the simplicial category in which
	\begin{itemize}
	\item objects are the finite pointed Kan complexes\index{finite pointed Kan complex}, i.e.\ pairs $(X,x_0)$ where $X$ is a finite Kan complex and $x_0$ is a vertex;
	\item mapping simplicial sets are given by the simplicial subsets
	\[ \Hom_{\finsset}((X,x_0),(Y,y_0)) \subseteq \Hom_{\sset}(X,Y) \]
	of the ordinary simplicial mapping spaces, in which the vertices are those maps $f:X \to Y$ such that $f(x_0) = y_0$.
	\end{itemize}
	Let $\finbased$\index{Sa@$\finbased$, the $\infty$-category of finite pointed spaces} be the simplicial nerve of the simplicial category $\finsset$. Since the mapping spaces in $\finsset$ are Kan complexes, the simplicial nerve $\finbased$ is an $\infty$-category by~\cite[1.1.5.10]{lurie:2009}.
\end{definition}

\begin{definition}[{\cite[7.3.1.10]{lurie:2017}}] \label{def:tan-cat}
	Let $\cat{C}$ be a differentiable $\infty$-category. The \emph{tangent bundle}\index{TAB@$T(\cat{C})$, the tangent bundle on a differentiable $\infty$-category $\cat{C}$} on $\cat{C}$ is the $\infty$-category
	\[ T(\cat{C}) := \Exc(\finbased,\cat{C}) \]
	of excisive functors $\finbased \to \cat{C}$ (i.e.\ those that map pushout squares in $\finbased$ to pullback squares in $\cat{C}$).
	
	The tangent bundle on $\cat{C}$ is equipped with a projection map
	\[ p_{\cat{C}}: T(\cat{C}) \to \cat{C}; \quad L \mapsto L(*) \]
	given by evaluating an excisive functor $L$ at the one-point space $*$.
\end{definition}

\begin{remark}
	Lurie actually makes Definition~\ref{def:tan-cat} for a slightly different class of $\infty$-categories: those that are \emph{presentable} in the sense of~\cite[5.5.0.1]{lurie:2009}. There is a significant overlap between the presentable and differentiable $\infty$-categories, including any compactly-generated $\infty$-category and any $\infty$-topos; see~\cite[6.1.1]{lurie:2017}.
	
	It seems very likely that most of the rest of this paper could be made with the $\infty$-category $\Catdiff$ replaced by an $\infty$-category $\Catpr$ of presentable $\infty$-categories in which the morphisms are those functors that preserve all \emph{filtered} colimits, not merely the sequential colimits. For example, Proposition~\ref{prop:Pn} holds when $\cat{D}$ is a presentable $\infty$-category, and \cite[7.3.1.14]{lurie:2017} implies that $T$ can be extended to a functor $T: \Catpr \to \Catpr$. However, much of our argument is based on results from~\cite[Sec. 6]{lurie:2017} which is written in the context of differentiable $\infty$-categories, so we follow that lead.
\end{remark}

\begin{definition} \label{def:TF}
	Let $F: \cat{C} \to \cat{D}$ be a functor between differentiable $\infty$-categories. We define a functor
	\[ T(F): T(\cat{C}) \to T(\cat{D}) \]
	by the formula
	\[ L \mapsto P_1(FL). \]
\end{definition}

It follows from Lemma~\ref{lem:T(F)} below that the constructions of Definitions~\ref{def:tan-cat} and~\ref{def:TF} can be made into a functor
\[ T: \Catdiff \to \Catdiff \]
so that the projection maps $p_\cat{C}$ form a natural transformation $p$ from $T$ to the identity functor $I$.

The following result is the main theorem of this paper.

\begin{theorem} \label{thm:main}
	The tangent bundle functor $T: \Catdiff \to \Catdiff$ and the projection map $p: T \to I$ extend to a tangent structure (in the sense of Definition~\ref{def:tangent-structure-infty}) on the $\infty$-category $\Catdiff$. We refer to this structure as the \emph{Goodwillie tangent structure}\index{Goodwillie tangent structure}.
\end{theorem}

We believe that the tangent structure described in Theorem~\ref{thm:main} is \emph{unique} (i.e.\ that the space of tangent structures that extend $T$ and $p$ is contractible) though we will not prove that claim here.

The proof of Theorem~\ref{thm:main} occupies the next two chapters, concluding with the proof of Theorem~\ref{thm:T}. In Chapter~\ref{sec:underlying} we introduce the basic constructions and definitions which underlie the tangent structure; that is, we describe how to construct, from a Weil-algebra $A$ and a differentiable $\infty$-category $\cat{C}$, a new $\infty$-category $T^A(\cat{C})$, and how these constructions interact with labelled morphisms of Weil-algebras and (sequential-colimit-preserving) functors of $\infty$-categories.

In Chapter~\ref{sec:tangent-structure-catdiff} we turn those constructions into an actual tangent structure on the $\infty$-category $\Catdiff$. It turns out to be more convenient to define that tangent structure on an $\infty$-category $\RelCatdiff$ that is equivalent to $\Catdiff$, and whose definition is based on the notion of \emph{relative} $\infty$-category which we describe in Definition~\ref{def:rel}. Having established a tangent structure on $\RelCatdiff$, we transfer that structure to $\Catdiff$ using Lemma~\ref{lem:transfer}.

In Chapter~\ref{sec:catdiff-diff} we begin the study of the Goodwillie tangent structure by identifying its differential objects in the sense of Definition~\ref{def:differential-structure}. Those differential objects turn out to be precisely the stable $\infty$-categories. That fact is of no surprise given that the tangent bundle construction~\ref{def:tan-cat} is formed precisely so that its tangent spaces are the stable $\infty$-categories $\spectra(\cat{C}_{/X})$. Nonetheless, this observation is a check that our tangent structure is acting as intended.

The characterization of differential objects as stable $\infty$-categories also confirms the intuition, promoted by Goodwillie, that in the analogy between functor calculus and the ordinary calculus of manifolds, one should view the category of spectra as playing the role of Euclidean space.

Further developing that analogy, in Chapter~\ref{sec:jet} we characterize the \emph{$n$-excisive} functors for $n > 1$, as corresponding to a notion of `$n$-jet' of a smooth map between manifolds. The construction of a Taylor tower itself does not precisely fit into the framework of tangent $\infty$-categories because it involves non-invertible natural transformations. We therefore show in Chapter~\ref{sec:infty2-catdiff} that the Goodwillie tangent structure on $\Catdiff$ extends to a tangent structure, in the sense of Example~\ref{ex:tangent-infty-bicategories}, on an $\infty$-bicategory $\CATdiff$ whose underlying $\infty$-category is $\Catdiff$. That tangent $\infty$-bicategory is the natural setting for Goodwillie calculus.

\chapter{The Goodwillie Tangent Structure: Underlying Data} \label{sec:underlying}

Our goal in this chapter is to introduce the basic data of the Goodwillie tangent structure on the $\infty$-category $\Catdiff$ without paying attention to the higher coherence information needed to obtain an actual tangent $\infty$-category. Thus, we define the tangent structure on objects and morphisms in $\Weilinfty$ and $\Catdiff$, and prove basic lemmas concerning functoriality and the preservation of pullbacks, including the crucial vertical lift axiom (Proposition~\ref{prop:vertical-lift}).

\subsection*{Tangent structure on objects}

We start by defining $T^A(\cat{C})$ for a Weil-algebra $A$ and a differentiable $\infty$-category $\cat{C}$. When $A = W$, this definition reduces precisely to Definition~\ref{def:tan-cat}.

\begin{definition} \label{def:A-excisive}
	Let $A = W^{n_1} \otimes \dots \otimes W^{n_r}$ be an object in the $\infty$-category $\Weilinfty$ with $n = n_1+\dots+n_r$ generators. We write
	\[ \finbased^n := \finbased \times \dots \times \finbased \]
	for the product of $n$ copies of the $\infty$-category $\finbased$ (with $\finbased^0 = *$). Let
	\begin{equation} \label{eq:TA} T^A(\cat{C}) = \Exc^A(\finbased^n,\cat{C}) := \Exc^{1,\dots,1}(\finbased^{n_1} \times \dots \times \finbased^{n_r},\cat{C}) \end{equation}
	be the full subcategory\index{TAC@$T^A(\cat{C})$, the action of a Weil-algebra $A$ on a differentiable $\infty$-category $\cat{C}$} of the functor $\infty$-category $\Fun(\finbased^n,\cat{C})$ that consists of those functors
	\[ L: \finbased^{n_1} \times \dots \times \finbased^{n_r} \to \cat{C} \]
	that are excisive (i.e.\ take pushouts to pullbacks) in each of their $r$ variables individually; see Definition~\ref{def:multi-excisive}. We say that a functor $L$ with this property is \emph{$A$-excisive}\index{A-excisive functor@$A$-excisive functor, for a Weil-algebra $A$}
\end{definition}

It is crucial in Definition~\ref{def:A-excisive} that we view $L$ as a functor of $r$ variables, each of the form $\finbased^{n_j}$, according to the description of the Weil-algebra $A$ as a tensor product of terms of the form $W^{n_j}$.

\begin{examples} \label{ex:TA}
	We can identify some particular examples of $T^A(\cat{C})$ to show that we are on the right track to define a tangent structure.
	\begin{enumerate} \itemsep=10pt
		\item For $A = \N$, the unit object for the monoidal structure on $\Weilinfty$, we get the identity functor:
		\[ T^{\N}(\cat{C}) = \Fun(\finbased^0,\cat{C}) \isom \cat{C}. \]
		\item For $A = W = \N[x]/(x^2)$, we get the tangent bundle functor from Definition~\ref{def:tan-cat}:
		\[ T^W(\cat{C}) = T(\cat{C}) = \Exc(\finbased,\cat{C}). \]
		\item For $A = W \otimes W = \N[x,y]/(x^2,y^2)$, we get the $\infty$-category of functors $\finbased^2 \to \cat{C}$ that are excisive \emph{in each variable individually}. In the notation of~\cite[6.1.3.1]{lurie:2017} we can write
		\[ T^{W \otimes W}(\cat{C}) = \Exc^{1,1}(\finbased^2,\cat{C}), \]
		and we have $T^{W \otimes W}(\cat{C}) \homeq T^W(T^W(\cat{C}))$ as required in a tangent structure; see~\cite[6.1.3.3]{lurie:2017}.
		\item \label{item:Wn} For $A = W^2 = \N[x,y]/(x^2,xy,y^2)$, we get the $\infty$-category of functors $\finbased^2 \to \cat{C}$ that are excisive when viewed as a functor of \emph{one} variable:
		\[ T_2(\cat{C}) := T^{W^2}(\cat{C}) = \Exc(\finbased^2,\cat{C}). \]
		For this definition to satisfy the pullback conditions in a tangent structure, we need to have an equivalence of $\infty$-categories
		\[ \Exc(\finbased^2,\cat{C}) \homeq \Exc(\finbased,\cat{C}) \times_{\cat{C}} \Exc(\finbased,\cat{C}), \]
		and indeed there is such an equivalence under which an excisive functor $L: \finbased^2 \to \cat{C}$ corresponds to the pair of excisive functors
		\[ (L(*,-),L(-,*)), \]
		and the pair $(L_1,L_2)$ with $L_1(*) = L_2(*)$ corresponds to the excisive functor given by the fibre product
		\[ (X,Y) \mapsto L_1(X) \times_{L_1(*) = L_2(*)} L_2(Y). \]
		That claim is proved in the next lemma.
	\end{enumerate}
\end{examples}

\begin{lemma} \label{lem:A-excisive}
	Let $\cat{S}_1,\cat{S}_2$ be $\infty$-categories each with finite colimits and a terminal object $*$, and let $\cat{C}$ be a differentiable $\infty$-category. Then a functor
	\[ L: \cat{S}_1 \times \cat{S}_2 \to \cat{C} \]
	is excisive (as a functor of one variable) if and only if
	\begin{enumerate}
		\item $L$ is excisive in each variable individually; and
		\item The morphisms $X_1 \to *$ and $X_2 \to *$ determine equivalences
		\[ L(X_1,X_2) \homeq L(X_1,*) \times_{L(*,*)} L(*,X_2) \]
		for all $X_1 \in \cat{S}_1$, $X_2 \in \cat{S}_2$.
	\end{enumerate}
\end{lemma}
\begin{proof}
	Key to this result is the fact that a pushout square in the $\infty$-category $\cat{S}_1 \times \cat{S}_2$ is a square for which each component is a pushout in its respective $\cat{S}_i$; see~\cite[5.1.2.3]{lurie:2009}. 
	
	Suppose that $L$ is excisive. Applying that condition to a pushout square in $\cat{S}_1 \times \cat{S}_2$ of the form
	\[ \begin{diagram}
		\node{(X_0,Y)} \arrow{e} \arrow{s} \node{(X_1,Y)} \arrow{s} \\
		\node{(X_2,Y)} \arrow{e} \node{(X_{12},Y)}
	\end{diagram} \]
	consisting of an arbitrary pushout in $\cat{S}_1$ and a fixed object in $\cat{S}_2$, we deduce that $L$ is excisive in its $\cat{S}_1$ variable. Similarly for its $\cat{S}_2$ variable.
	
	Applying the condition that $L$ is excisive to a pushout square in $\cat{S}_1 \times \cat{S}_2$ of the form
	\[ \begin{diagram}
		\node{(X_1,X_2)} \arrow{e} \arrow{s} \node{(X_1,*)} \arrow{s} \\
		\node{(*,X_2)} \arrow{e} \node{(*,*)}
	\end{diagram} \]
	we deduce condition (2).
	
	Conversely, suppose that $L$ satisfies conditions (1) and (2), and consider a diagram
	\[ \begin{diagram}
		\node{(X_0,Y_0)} \arrow{e} \arrow{s} \node{(X_1,Y_1)} \arrow{s} \\
		\node{(X_2,Y_2)} \arrow{e} \node{(X_{12},Y_{12})}
	\end{diagram} \]
	that is a pushout in each component.
	
	Consider the following diagram in $\cat{C}$:
	\[ \begin{diagram}
		\node{L(X_0,Y_0)} \arrow{e} \arrow{s} \node{L(X_0,Y_1)} \arrow{e} \arrow{s} \node{L(X_1,Y_1)} \arrow{s} \\
		\node{L(X_2,Y_0)} \arrow{e} \arrow{s} \node{L(X_2,Y_1)} \arrow{e} \arrow{s} \node{L(X_{12},Y_1)} \arrow{s} \\
		\node{L(X_2,Y_2)} \arrow{e} \node{L(X_2,Y_{12})} \arrow{e} \node{L(X_{12},Y_{12})}
	\end{diagram} \]
	The top-right and bottom-left squares are pullbacks because $L$ is excisive in each variable individually, so it is sufficient to show that the top-left and bottom-right squares are also pullbacks, since then the whole square is a pullback by standard pasting properties of pullbacks (see~\cite[4.4.2.1]{lurie:2009}).
	
	For the top-left square (the bottom-right is similar), consider the diagram
	\[ \begin{diagram}
		\node{L(X_0,Y_0)} \arrow{e} \arrow{s} \node{L(X_0,Y_1)} \arrow{e} \arrow{s} \node{L(X_0,*)} \arrow{s} \\
		\node{L(X_2,Y_0)} \arrow{e} \arrow{s} \node{L(X_2,Y_1)} \arrow{e} \arrow{s} \node{L(X_2,*)} \arrow{s} \\
		\node{L(*,Y_0)} \arrow{e} \node{L(*,Y_1)} \arrow{e} \node{L(*,*)}
	\end{diagram} \]
	Condition (2) implies that the bottom-right square, the bottom half, right-hand half, and entire square are all pullbacks. From that, it follows by a succession of applications of~\cite[4.4.2.1]{lurie:2009} that each individual square is a pullback too, including the top-left as required.
\end{proof}

A crucial role in the following sections will be played by a universal $A$-excisive approximation. The next proposition is a generalization of~\ref{prop:Pn}, which corresponds to the case $A = W$.

\begin{proposition} \label{prop:PA}
	Let $A$ be a Weil-algebra with $n$ generators, and let $\cat{C}$ be a differentiable $\infty$-category. Then there is a functor\index{PA@$P_A$, $A$-excisive approximation}
	\[ P_A: \Fun(\finbased^n,\cat{C}) \to \Exc^A(\finbased^n,\cat{C}) \]
	that is left adjoint to the inclusion and preserves finite limits. Moreover, the $\infty$-category $T^A(\cat{C}) = \Exc^A(\finbased^n,\cat{C})$ is differentiable with finite limits and sequential colimits all calculated objectwise in $\cat{C}$. We write
	\[ p_A: F \to P_A(F) \]
	for the corresponding universal $P_A$-approximation map.
\end{proposition}
\begin{proof}
	The first part is an example of the multivariable excisive approximation construction described in Definition~\ref{def:multi-excisive}. It follows from the definition of excisive, and the fact that $\cat{C}$ is differentiable, that the subcategory $\Exc^A(\finbased^n,\cat{C})$ of $\Fun(\finbased^n,\cat{C})$ is closed under finite limits and sequential colimits which are computed objectwise in $\cat{C}$, and hence commute. Therefore $\Exc^A(\finbased^n,\cat{C})$ is differentiable.
\end{proof}

The following property of $P_A$ is used in the proof of Lemma~\ref{lem:T1}.

\begin{lemma} \label{lem:PA}
	Let $F: \cat{C} \to \cat{D}$ be a functor between differentiable $\infty$-categories that preserves both finite limits and sequential colimits. Then the following diagram commutes (up to natural equivalence)
	\[ \begin{diagram}
		\node{\Fun(\finbased^n,\cat{C})} \arrow{s,l}{F_*} \arrow{e,t}{P_A} \node{\Exc^A(\finbased^n,\cat{C})} \arrow{s,r}{F_*} \\
		\node{\Fun(\finbased^n,\cat{D})} \arrow{e,t}{P_A} \node{\Exc^A(\finbased^n,\cat{D})}
	\end{diagram} \]
	where each functor $F_*$ denotes post-composition with $F$.
\end{lemma}
\begin{proof}
	This claim is a multivariable version of \cite[6.1.1.32]{lurie:2017} and follows from the explicit construction of the excisive approximation, as in the proof of Lemma~\ref{lem:PnGF}.
\end{proof}

\subsection*{Tangent structure on morphisms in  $\Catdiff$}

We now turn to the action of the tangent structure functors $T^A$ on morphisms in $\Catdiff$, i.e.\ functors $F: \cat{C} \to \cat{D}$ between differentiable $\infty$-categories that preserve sequential colimits. This action is the obvious extension of that described in~\ref{def:TF} for the tangent bundle functor $T$.

\begin{definition} \label{def:TAF}
	Let $A$ be a Weil-algebra, and let $F: \cat{C} \to \cat{D}$ be a functor between differentiable $\infty$-categories that preserves sequential colimits. We define
	\[ T^A(F): T^A(\cat{C}) \to T^A(\cat{D}) \]
	to be the composite
	\[ \Exc^A(\finbased^n,\cat{C}) \subseteq \Fun(\finbased^n,\cat{C}) \arrow{e,t}{F_*} \Fun(\finbased^n,\cat{D}) \arrow{e,t}{P_A} \Exc^A(\finbased^n,\cat{D}). \]
	That is, we have
	\begin{equation} \label{eq:T^AF} T^A(F)(L) := P_A(FL) \end{equation}
	for an $A$-excisive functor $L: \finbased^n \to \cat{C}$.
\end{definition}

\begin{lemma}
	Let $A$ be a Weil-algebra, and let $F: \cat{C} \to \cat{D}$ be a sequential-colimit-preserving functor between differentiable $\infty$-categories. Then $T^A(F)$ also preserves sequential colimits.
\end{lemma}
\begin{proof}
	Each of the functors $F_*,P_A$ in Definition~\ref{def:TAF} preserves sequential colimits, which, in all cases, are computed pointwise in the $\infty$-categories $\cat{C}$ and $\cat{D}$.
\end{proof}

We now check that Definition~\ref{def:TAF} makes $T^A$ into a functor $\Catdiff \to \Catdiff$, at least up to higher equivalence.

\begin{lemma} \label{lem:T(F)}
	Let $A$ be a Weil-algebra. Then:
	\begin{enumerate} \itemsep=10pt
		\item for the identity functor $I_\cat{C}$ on a differentiable $\infty$-category $\cat{C}$, there is a natural equivalence
		\[ I_{T^A(\cat{C})} \weq T^A(I_{\cat{C}}) \]
		given by the maps $p_A: L \weq P_A(L)$ for $A$-excisive $L: \finbased^n \to \cat{C}$;
		\item for sequential-colimit-preserving functors $F: \cat{C} \to \cat{D}$ and $G: \cat{D} \to \cat{E}$ between differentiable $\infty$-categories, there is a natural equivalence
		\[ T^A(GF) \weq T^A(G)T^A(F) \]
		which comprises natural equivalences
		\[ P_A(GFL) \weq P_A(GP_A(FL)) \]
		induced by the $A$-excisive approximation map $p_A: FL \to P_A(FL)$ for $L \in \Exc^A(\finbased^n,\cat{C})$.
	\end{enumerate}
\end{lemma}
\begin{proof}
	Part (1) is a standard property of a localization functor such as $P_A$. Part (2) is more substantial. When $A = W$, so that $P_A = P_1$, this result is proved by Lurie in~\cite[7.3.1.14]{lurie:2017}. (That result is in the context of presentable $\infty$-categories and functors that preserve all filtered colimits, but the proof works equally well for differentiable $\infty$-categories and functors that only preserve sequential colimits. Fundamentally this result relies on the Klein-Rognes~\cite{klein/rognes:2002} Chain Rule as extended to $\infty$-categories by Lurie in~\cite[6.2.1.24]{lurie:2017}.)
	
	In particular, Lurie's argument shows that when the functor $G: \cat{D} \to \cat{E}$ between differentiable $\infty$-categories preserves sequential colimits, we have an equivalence
	\[ P_1(GH) \to P_1(GP_1(H)) \]
	for any $H: \finbased \to \cat{D}$. That argument extends to the case of $H: \finbased^n \to \cat{D}$, simply by replacing $\finbased$ with $\finbased^n$ and the null object $*$ with $(*,\dots,*)$. This observation provides the desired result when $A = W^n$.
	
	Now consider an arbitrary Weil-algebra $A = W^{n_1} \otimes \dots \otimes W^{n_r}$. Then we have
	\[ P_A = P_1^{(r)} \cdots P_1^{(1)} \]
	where $P_1^{(j)}$ is the ordinary excisive approximation applied to a functor
	\[ \finbased^{n_1} \times \dots \times \finbased^{n_r} \to \cat{D} \]
	in its $\finbased^{n_j}$ variable. For any $H: \finbased^n \to \cat{D}$, the map $p_A: H \to P_A(H)$ factors as a composite
	\[ \dgTEXTARROWLENGTH=2.5em  H \arrow{e,t}{p_1^{(1)}} P_1^{(1)}H \arrow{e,t}{p_1^{(2)}} P_1^{(2)}P_1^{(1)}H \arrow{e,t}{p_1^{(3)}} \cdots \arrow{e,t}{p_1^{(r)}} P_AH \]
	of excisive approximations in each variable separately. Applying $G$ to the $j$-th map yields a $P_1^{(j)}$-equivalence by the extended version of Lurie's argument, and hence a $P_A$-equivalence as required. Thus, we have an equivalence
	\begin{equation} \label{eq:PAPA} P_A(GH) \weq P_A(GP_A(H)). \end{equation}
	Taking $H$ to be $FL$ yields the desired result.
\end{proof}

\subsection*{Tangent structure on morphisms in $\Weilinfty$}

We next address the functoriality of our putative tangent structure in the $\Weilinfty$ variable. Let $\phi:A \to A'$ be a morphism in $\Weilinfty$, i.e.\ a labelled morphism of Weil-algebras, and let $\cat{C}$ be a differentiable $\infty$-category. We construct a (sequential-colimit-preserving) functor
\[ T^\phi(\cat{C}): T^A(\cat{C}) \to T^{A'}(\cat{C}) \]
as follows. Roughly speaking, $T^\phi$ is the map
\[ \Exc^A(\finbased^n,\cat{C}) \to \Exc^{A'}(\finbased^{n'},\cat{C}) \]
given by precomposition with a suitable functor
\[ \tilde{\phi}: \finbased^{n'} \to \finbased^n \]
whose definition mirrors the Weil-algebra morphism $\phi$.

\begin{definition} \label{def:phi}
Let $\phi: A \to A'$ be a labelled morphism between Weil-algebras with $n$ and $n'$ generators respectively. Let $M_A$ and $M_{A'}$ be the partial monoids of nonzero monomials in $A$ and $A'$, respectively. Recall from Definition~\ref{def:weil-mor} that $\phi$ is a span of finite partial commutative monoids of the form
\[ \begin{diagram}
	\node[2]{K} \arrow{sw,t}{s} \arrow{se,t}{t} \\
	\node{M_A} \node[2]{M_{A'}}
\end{diagram} \]
and that the underlying ordinary morphism of Weil-algebras is given by
\[ \phi(x_i) = \sum_{k \in s^{-1}(x_i)} t(k) \]
for each generator $x_i$ of the Weil-algebra $A$.

We define a simplicially-enriched functor
\[ \tilde{\phi}: \finsset^{n'} \to \finsset^n \]
by setting
\[ \tilde{\phi}(Y_1,\dots,Y_{n'}) := (X_1,\dots,X_n) \]
where
\[ X_i := \bigvee_{k \in s^{-1}(x_i)} Y_{t(k)} \]
where a monomial $t(k) = y_{j_1} \cdots y_{j_r}$ in $A'$ determines a corresponding smash product
\[ Y_{t(k)} := Y_{j_1} \smsh \dots \smsh Y_{j_r}. \]
Note that if $s^{-1}(x_i)$ is empty, i.e.\ if $\phi(x_i) = 0$, then we set $X_i := *$. 

The above construction extends to mapping spaces by way of the usual simplicial enrichment of the wedge sum and smash product of pointed simplicial sets. Thus, we obtain the desired simplicially-enriched functor $\tilde{\phi}$.

We now take the simplicial nerve of $\tilde{\phi}$ to get a functor of $\infty$-categories, which we denote also as
\[ \tilde{\phi}: \finbased^{n'} \to \finbased^{n}. \]
On objects, this functor is given by the same formulas above.
\end{definition}

\begin{definition} \label{def:Tphi}
	Let $\phi: A \to A'$ be a labelled Weil-algebra morphism, and let $\cat{C}$ be a differentiable $\infty$-category. Then we define\index{TAD@$T^\phi$, action of a labelled Weil-algebra morphism $\phi$ in the Goodwillie tangent structure}
	\[ T^\phi(\cat{C}): T^A(\cat{C}) \to T^{A'}(\cat{C}); \quad L \mapsto P_{A'}(L\tilde{\phi}), \]
	that is, $T^{\phi}(\cat{C})$ is the composite
	\[ \dgTEXTARROWLENGTH=2.5em \Exc^A(\finbased^n,\cat{C}) \subseteq \Fun(\finbased^n,\cat{C}) \arrow{e,t}{\tilde{\phi}^*} \Fun(\finbased^{n'},\cat{C}) \arrow{e,t}{P_{A'}} \Exc^{A'}(\finbased^{n'},\cat{C}) \]
	where $\tilde{\phi}^*$ denotes precomposition with $\tilde{\phi}$.
\end{definition}

\begin{examples} \label{ex:Tphi}
	Using Definition~\ref{def:Tphi}, we can now work out how the fundamental natural transformations from Cockett and Cruttwell's definition of tangent structure~\cite[2.3]{cockett/cruttwell:2014} manifest in our case.
	\begin{enumerate} \itemsep=10pt
		\item Let $\epsilon: W \to \mathbb{N}$ be the augmentation. Then $\tilde{\epsilon}: * \to \finbased$ is the functor that picks out the null object $*$, and so the \emph{projection}
		\[ p := T^\epsilon: T(\cat{C}) \to \cat{C} \]
		can be identified with the evaluation map
		\[ \Exc(\finbased,\cat{C}) \to \cat{C}; \quad L \mapsto L(*), \]
		as in Definition~\ref{def:tan-cat}.
		\item Let $\eta: \mathbb{N} \to W$ be the unit map. Then $\tilde{\eta}: \finbased \to *$ is, of course, the constant map, and so the \emph{zero section}
		\[ 0 := T^\eta: \cat{C} \to T(\cat{C}) \]
		can be identified with the map
		\[ \cat{C} \to \Exc(\finbased,\cat{C}); \quad C \mapsto \mathrm{const}_C \]
		that picks out the constant functors.
		\item Let $\phi: W^2 \to W$ be the addition map given by $x \mapsto x$ and $y \mapsto x$. Then $\tilde{\phi}: \finbased \to \finbased^2$ is the diagonal $X \mapsto (X,X)$, and so the \emph{addition}
		\[ + := T^\phi: T^{W^2}(\cat{C}) \to T(\cat{C}) \]
		is given by the map
		\[ \Exc(\finbased^2,\cat{C}) \to \Exc(\finbased,\cat{C}); \quad L \mapsto [X \mapsto L(X,X)]. \]
		Under the equivalence $T^{W^2}(\cat{C}) \homeq T(\cat{C}) \times_{\cat{C}} T(\cat{C})$ of Example~\ref{ex:TA}(\ref{item:Wn}), we can identify $+$ with the fibrewise product map
		\[ \begin{split} \Exc(\finbased,\cat{C}) \times_{\cat{C}} \Exc(\finbased,\cat{C}) \to \Exc(\finbased,\cat{C}); \\ (L_1,L_2) \mapsto L_1(-) \times_{L_1(*) = L_2(*)} L_2(-). \end{split} \]
		\item Let $\sigma: W \otimes W \to W \otimes W$ be the symmetry map. Then $\tilde{\sigma}: \finbased^2 \to \finbased^2$ is given by $(X,Y) \mapsto (Y,X)$, and the \emph{flip}
		\[ c := T^\sigma: T^2(\cat{C}) \to T^2(\cat{C}) \]
		is the symmetry map
		\[ \Exc^{1,1}(\finbased^2,\cat{C}) \to \Exc^{1,1}(\finbased^2,\cat{C}); \quad L \mapsto [(X,Y) \mapsto L(Y,X)]. \]
		\item Let $\delta: W \to W \otimes W$ be the map $x \mapsto xy$. Then $\tilde{\delta}: \finbased^2 \to \finbased$ is the smash product $(X,Y) \mapsto X \smsh Y$, and the \emph{vertical lift}
		\[ \ell := T^\delta: T(\cat{C}) \to T^2(\cat{C}) \]
		is the map
		\[ \Exc(\finbased,\cat{C}) \to \Exc^{1,1}(\finbased^2,\cat{C}); \quad L \mapsto [(X,Y) \mapsto L(X \smsh Y)]. \]
	\end{enumerate}
\end{examples}

\begin{remark}
	In each of the cases in Examples~\ref{ex:Tphi}, we did not need to apply $P_{A'}$ at the end because precomposition with $\tilde{\phi}$ preserved the excisiveness property. However, this is not always the case; for example, if $\phi: W \otimes W \to W$ is the fold map given by $x_1,x_2 \mapsto x$, then $\tilde{\phi}(X) = (X,X)$, but the functor $X \mapsto L(X,X)$ is not typically excisive when $L$ is excisive in each variable individually.
\end{remark}

We should check that $T^\phi$ is a morphism in $\Catdiff$, and it turns out to have a stronger property too.

\begin{lemma}
	The functor $T^\phi: T^A(\cat{C}) \to T^{A'}(\cat{C})$ preserves finite limits and sequential colimits.
\end{lemma}
\begin{proof}
	This result follows immediately from Definition~\ref{def:Tphi}, since $P_{A'}$ preserves these (co)limits, which in each case are calculated objectwise in $\cat{C}$.
\end{proof}

We now check functoriality in the $\Weilinfty$ factor, at least up to higher equivalence.

\begin{lemma} \label{lem:Tphi}
	Let $\cat{C}$ be a differentiable $\infty$-category. Then:
	\begin{enumerate} \itemsep=10pt
		\item for the identity morphism $1_A$ on a Weil-algebra $A$ we have $\tilde{1}_A = I_{\finbased^n}$, the identity functor on $\finbased^n$, and so there is a natural equivalence
		\[ I_{T^A(\cat{C})} \weq T^{1_A}(\cat{C}) \]
		given by the maps $p_A: L \weq P_A(L)$ for $A$-excisive $L: \finbased^n \to \cat{C}$;
		\item for a $2$-simplex in $\Weilinfty$, consisting of labelled Weil-algebra morphisms $\phi_1: A \to A'$ and $\phi_2: A' \to A''$, and a labelling on the composite $\phi_2\phi_1$, together with extra data relating these labellings, which we recall below, there is a natural equivalence
		\[ T^{\phi_2\phi_1}(\cat{C}) \weq T^{\phi_2}(\cat{C})T^{\phi_1}(\cat{C}) \]
		given by a composite of two maps of the form
		\[ \dgTEXTARROWLENGTH=2em P_{A''}((-)\widetilde{\phi_2\phi_1}) \arrow{e,t}{(i)} P_{A''}((-)\tilde{\phi_1}\tilde{\phi_2}) \arrow{e,t}{(ii)} P_{A''}(P_{A'}((-)\tilde{\phi_1})\tilde{\phi_2}) \]
		induced (i) by a natural transformation $\alpha: \widetilde{\phi_2\phi_1} \to \tilde{\phi_1}\tilde{\phi_2}$, described in Definition~\ref{def:alpha} below, and (ii) by the universal $A'$-excisive approximation map $p_{A'}$.
	\end{enumerate}
\end{lemma}

\begin{definition} \label{def:alpha}
	Let $\sigma$ be a $2$-simplex in $\Weilinfty$, i.e.\ a diagram of finite partial commutative monoids of the form
	\[ \begin{diagram}
		\node[3]{K''} \arrow{sw,t}{s''} \arrow{se,t}{t''} \\
		\node[2]{K} \arrow{sw,t}{s} \arrow{se,t}{t} \node[2]{K'} \arrow{sw,t}{s'} \arrow{se,t}{t'} \\
		\node{M_A} \node[2]{M_{A'}} \node[2]{M_{A''}}
	\end{diagram} \]
	satisfying the conditions of Proposition~\ref{prop:weil}, and suppose the Weil-algebras $A,A',A''$ have $n,n',n''$ generators respectively.
	
	The simplicially-enriched functor $\widetilde{\phi_2\phi_1}: \finsset^{n''} \to \finsset^n$ is given by
	\[ \widetilde{\phi_2\phi_1}(Z_1,\dots,Z_{n''}) = (X_1,\dots,X_n) \]
	where
	\[ X_i = \Wdge_{k'' \; : \; ss''(k'') = x_i} Z_{t't''(k'')}. \]
	On the other hand, the composite functor $\tilde{\phi_1}\tilde{\phi_2}$ is given by
	\[ \tilde{\phi_1}\tilde{\phi_2}(Z_1,\dots,Z_{n''}) \isom (X'_1,\dots,X'_n) \]
	where
	\[ X'_i = \Wdge_{k \; : \; s(k) = x_i} \left(\Wdge_{k'_1,\dots,k'_r \in K'_1 \; : \; t(k) = s'(k'_1) \cdots s'(k'_r)} Z_{t'(k_1)} \smsh \dots \smsh Z_{t'(k_r)} \right). \]
	We define a simplicially-enriched natural transformation
	\[ \alpha: \widetilde{\phi_2\phi_1} \to \tilde{\phi_1}\tilde{\phi_2} \]
	by the maps $X_i \to X'_i$ which are given, on the factor corresponding to $k'' \in K''$ such that $s's''(k'') = x_i$, by the canonical isomorphism
	\[ \alpha_{k''}: Z_{t't''(k'')} \to Z_{t'(k'_1)} \smsh \dots \smsh Z_{t'(k'_r)} \]
	into the factor given by $k = s''(k'')$ and $k'_1 \cdots k'_r = t''(k'')$, noting that
	\[ t't''(k'') = t'(k'_1 \cdots k'_r) = t'(k'_1) \cdots t'(k'_r). \]
	Note that no monomial in a Weil-algebra has a repeated factor, so there is indeed only one canonical isomorphism $\alpha_{k''}$ which matches up the factors in the two smash products.
	
	We now take simplicial nerves of the simplicial natural transformation $\alpha$ to get a natural transformation between functors $\finbased^{n''} \to \finbased^{n}$ (in the $\infty$-categorical sense), for which we use the same notation,
	\[ \alpha: \widetilde{\phi_2\phi_1} \to \tilde{\phi_1}\tilde{\phi_2}. \]
\end{definition}

\begin{remark}
	The pullback condition in Proposition~\ref{prop:weil} for the $2$-simple $\sigma$ identifies the set $K''_i$ with pairs $(k,k'_1 \cdots k'_r)$ such that $s(k) = x_i$, and $t(k) = s'(k'_1) \cdots s'(k'_r)$, \emph{and such that} the product $k'_1 \cdots k'_r$ is defined in $K'$. However, there are other terms appearing in the wedge sum that forms $X'_i$: those where $k'_1 \cdots k'_r$ is not an element of $K'$ because $t'(k'_1) \cdots t'(k'_r)$ is zero in $A''$.
	
	It follows that $\alpha$ is the inclusion of the wedge summand in a decomposition
	\begin{equation} \label{eq:decomp-phi} \tilde{\phi_1}\tilde{\phi_2} \homeq \widetilde{\phi_2\phi_1} \wdge \zeta(\phi_1,\phi_2) \end{equation}
	where $\zeta(\phi_1,\phi_2)(Z_1,\dots,Z_{n''})$ is a finite wedge sum of terms, each of the form
	\[ Z_{k_1} \smsh \dots \smsh Z_{k_q} \]
	where the monomial $z_{k_1} \cdots z_{k_q}$ is zero in the Weil-algebra $A''$.
\end{remark}

\begin{proof}[{Proof of Lemma~\ref{lem:Tphi}}]
	The calculation of $\tilde{1}_A$ in part (1) follows immediately from Definition~\ref{def:phi}. Part (2) is much more substantial and will occupy the next several pages. We prove that each of the two maps $(i)$ and $(ii)$ is an equivalence, and both of those facts are important later in the paper.
	
	For $(i)$, let $L: \finbased^{n} \to \cat{C}$ be an $A$-excisive functor. We have to show that the map
	\begin{equation} \label{eq:alpha} P_{A''}(L\alpha): P_{A''}(L\widetilde{\phi_2\phi_1}) \to P_{A''}(L\tilde{\phi_1}\tilde{\phi_2}) \end{equation}
	is an equivalence. The idea here is that the difference between $\widetilde{\phi_2\phi_1}$ and $\tilde{\phi_1}\tilde{\phi_2}$, i.e.\ the factor $\zeta(\phi_1,\phi_2)$ in (\ref{eq:decomp-phi}), consists of terms that make no contribution after taking the $A''$-excisive approximation.
	
	By induction on the number of wedge summands in $\zeta(\phi_1,\phi_2)$, we reduce to showing that for any functor $\beta: \finbased^{n''} \to \finbased^{n}$, the inclusion $\beta \to \beta \wdge \xi$ induces an equivalence
	\[ P_{A''}(L\beta) \to P_{A''}(L(\beta \wdge \xi)) \]
	where
	\[ \xi(Z_1,\dots,Z_{n''}) \homeq (*,\dots,*,Z_{k_1} \smsh \dots \smsh Z_{k_q},*,\dots,*) \]
	for some sequence of indices $k_1,\dots,k_r$ such that $z_{k_i}z_{k_{i'}} = 0$ in $A''$ for some $i \neq i'$. Equivalently, we show that the collapse map $\beta \wdge \xi \to \beta$ induces an equivalence in the other direction.
	
	Writing $L: \finbased^{n_1} \times \dots \times \finbased^{n_r} \to \cat{C}$ according to the decomposition of the Weil-algebra $A = W^{n_1} \otimes \dots \otimes W^{n_r}$, we can assume without loss of generality that $r = 1$ (i.e.\ $A = W^n$) since $\xi$ is only non-trivial in one variable. In that case, $L: \finbased^n \to \cat{C}$ is excisive. Applying $L$ to the pushout diagram in $\Fun(\finbased^{n''},\finbased^n)$ of the form
	\[ \begin{diagram}
		\node{\beta \vee \xi} \arrow{e} \arrow{s} \node{\beta} \arrow{s} \\
		\node{\xi} \arrow{e} \node{*}
	\end{diagram} \]
	and taking $P_{A''}$ (which preserves pullbacks), we get a pullback square
	\[ \begin{diagram}
		\node{P_{A''}(L(\beta \wdge \xi))} \arrow{e} \arrow{s} \node{P_{A''}(L\beta)} \arrow{s} \\
		\node{P_{A''}(L\xi)} \arrow{e} \node{P_{A''}(L(*))}
	\end{diagram} \]
	so it is sufficient to show that the bottom map is an equivalence. Replacing $\cat{C}$ with the slice $\infty$-category $\cat{C}_{/L(*)}$, we can assume that $L$ is reduced. Our goal is then to show that the functor
	\begin{equation} \label{eq:LZ} (Z_1,\dots,Z_{n''}) \mapsto L(*, \dots, *, Z_{k_1} \smsh \dots \smsh Z_{k_q}, *, \dots, *) \end{equation}
	has trivial $A''$-excisive part whenever $z_{k_i}z_{k_{i'}} = 0$ in $A''$ for some $i \neq i'$.
	
	If we write
	\[ \finbased^{n''} = \finbased^{n''_1} \times \dots \times \finbased^{n''_s} \]
	according to the decomposition $A'' = W^{n''_1} \otimes \dots \otimes W^{n''_s}$, then the condition that $z_{k_i}z_{k_{i'}} = 0$ implies that the variables $Z_{k_i}$ and $Z_{k_{i'}}$ are in the same factor $\finbased^{n''_j}$, and it is sufficient to show that the functor (\ref{eq:LZ}) has trivial excisive approximation with respect to that factor. Without loss of generality, we can take $A'' = W^{n''}$ so that we simply have to show that a functor $F: \finbased^{n''} \to \cat{C}$ of the form (\ref{eq:LZ}) has trivial excisive approximation whenever $q \geq 2$.
	
	Since $F$ is reduced in each of the variables $Z_{k_i}$, it is, by~\cite[6.1.3.10]{lurie:2017}, $q$-reduced when viewed as a functor of all of those variables. Since $q \geq 2$, it follows that $F$ has trivial excisive approximation with respect to those variables. It follows from Lemma~\ref{lem:A-excisive} that the excisive approximation of a functor $\finbased^{n''} \to \cat{C}$ factors via its excisive approximation with respect to any subset of its variables, so $F$ also has trivial excisive approximation as a functor $\finbased^{n''} \to \cat{C}$, as desired. This completes the proof that the map $(i)$ is an equivalence.
	
	For $(ii)$, we prove the following: for any $G: \finbased^{n'_1} \times \dots \times \finbased^{n'_r} \to \cat{C}$ and any Weil-algebra morphism $\phi: A' \to A''$, there is an equivalence
	\begin{equation} \label{eq:pp} P_{A''}(G\tilde{\phi}) \weq P_{A''}(P_{A'}(G)\tilde{\phi}) \end{equation}
	induced by the $A'$-excisive approximation map $p_{A'}: G \to P_{A'}G$.
	
	It is sufficient to show that excisive approximation with respect to each of the $r$ variables induces an equivalence. Thus we can reduce to the case that $A' = W^{n'}$, in which case $P_{A'} = P_1$. Replacing $\cat{C}$ with the slice $\infty$-category $\cat{C}_{G(*)}$ of objects over and under $G(*)$, we can also assume that $\cat{C}$ is a pointed $\infty$-category and that $G$ is reduced, i.e.\ $G(*) \homeq *$. So we have to show that for any reduced $G: \finbased^{n'} \to \cat{C}$ the map
	\begin{equation} \label{eq:PA1G} P_{A''}(G\tilde{\phi}) \to P_{A''}((P_1G)\tilde{\phi}) \end{equation}
	is an equivalence.
	
	We break this proof into two parts. First, we use downward induction on the Taylor tower of $G$ to reduce to the case $G$ is $m$-homogeneous for some $m \geq 2$. From there, we use the specific form of $\tilde{\phi}$ and the fact that $\phi$ is an algebra homomorphism to show (\ref{eq:PA1G}) is an equivalence by direct calculation.
	
	To start the induction we apply Lemma~\ref{lem:PnGF} which tells us that, since $\tilde{\phi}$ is reduced, the universal $n$-excisive approximation $p_n: G \to P_nG$ induces an equivalence
	\begin{equation} \label{eq:PnPnGphi} P_n(G\tilde{\phi}) \weq P_n((P_nG)\tilde{\phi}) \end{equation}
	for any $n$.
	
	If $n \geq s$, then an $A''$-excisive functor $\finbased^{n''_1} \times \dots \times \finbased^{n''_s} \to \cat{C}$ is $n$-excisive by~\cite[6.1.3.4]{lurie:2017}. It follows from (\ref{eq:PnPnGphi}) that $p_n$ induces an equivalence
	\begin{equation} \label{eq:PAPnG} P_{A''}(G\tilde{\phi}) \weq P_{A''}((P_nG)\tilde{\phi}). \end{equation}
	Next, consider the fibre sequences
	\[ P_mG \to P_{m-1}G \to R_mG \]
	provided by~\cite[6.1.2.4]{lurie:2017} (Goodwillie's delooping theorem for homogeneous functors) in which $R_mG$ is $ m$-homogeneous. We show below that
	\begin{equation} \label{eq:PAH} P_{A''}(H\tilde{\phi}) \homeq * \end{equation}
	for any functor $H: \finbased^{n'} \to \cat{C}$ that is $m$-homogeneous for some $m \geq 2$. Since $P_{A''}$ preserves fibre sequences, we then deduce that for $m \geq 2$
	\[ P_{A''}((P_mG)\tilde{\phi}) \weq P_{A''}((P_{m-1}G)\tilde{\phi}), \]
	which combined with (\ref{eq:PAPnG}) implies the desired result (\ref{eq:PA1G}).
	
	So our goal now is (\ref{eq:PAH}). By~\cite[6.1.2.9]{lurie:2017}, we can assume that $\cat{C}$ is a stable $\infty$-category, and then Lemma~\ref{lem:hom-class} below gives a classification of $m$-homogeneous functors $\finbased^{n'} \to \cat{C}$. That result tells us that $H$ is a finite product of terms of the form
	\begin{equation} \label{eq:LY} (Y_1,\dots,Y_{n'}) \mapsto L(Y_1^{\smsh m_1} \smsh \dots \smsh Y_{n'}^{\smsh m_{n'}})_{h(\Sigma_{m_1} \times \dots \times \Sigma_{m_{n'}})} \end{equation}
	for linear $L: \finbased \to \cat{C}$ and $m = m_1+\dots+m_{n'}$. Since $P_{A''}$ preserves finite products, we can reduce our goal (\ref{eq:PAH}) to the case that $H$ is of the form in (\ref{eq:LY}).
	
	The next part of the proof depends essentially on the fact that $\phi: A' \to A''$ is an algebra homomorphism, and we start by describing what that condition implies about the functor $\tilde{\phi}$.
	
	Recall that we have reduced to the case 
	\[ A' = W^{n'} = \N[y_1,\dots,y_{n'}]/(y_iy_{i'})_{i,i' = 1,\dots,n'} \]
	and let us denote the generators of
	\[ A'' = W^{n''_1} \otimes \dots \otimes W^{n''_s} \]
	as $z_{j,1}, \dots, z_{j,n''_j}$ for $j = 1,\dots,s$. Similarly, we denote inputs to the corresponding functor
	\[ \tilde{\phi}: \finbased^{n''_1} \times \dots \times \finbased^{n''_s} \to \finbased^{n'} \]
	with the notation $Z_{j,k}$. We write $\un{Z}$ for the object of $\finbased^{n''_1} \times \dots \times \finbased^{n''_s}$ with these entries.
	
	Since $\phi$ is an algebra homomorphism, we have, for any $i,i'$:
	\[ \phi(y_i)\phi(y_{i'}) = \phi(y_iy_{i'}) = \phi(0) = 0. \]
	There are no terms in $A''$ with negative coefficients, so this means that the product of any monomial in $\phi(y_i)$ and any monomial in $\phi(y_{i'})$ must contain a factor of the form $z_{j,k}z_{j,k'}$ for some $j,k,k'$.
	
	Translating this observation into a statement about the map $\tilde{\phi}$, and writing $\tilde{\phi}_1, \dots, \tilde{\phi}_{n'}$ for the components of that map, we see that for all $i,i' = 1,\dots,n'$:
	\begin{equation} \label{eq:phiZ} \tilde{\phi}_{i}(\un{Z}) \smsh \tilde{\phi}_{i'}(\un{Z}) \homeq \Wdge Z_{j_1,k_1} \smsh \dots \smsh Z_{j_t,k_t}, \end{equation}
	a finite wedge sum of terms each of which satisfies $j_l = j_{l'}$ for some $l \neq l'$.
	
	Now consider the functor $H\tilde{\phi}: \finbased^{n''} \to \cat{C}$ where $H$ is of the form (\ref{eq:LY}) with $m = m_1+\dots+m_{n'} \geq 2$. Applying $H$ to $\tilde{\phi}$ that satisfies (\ref{eq:phiZ}), and since the linear functor $L$ takes finite wedge sums to products, we obtain a decomposition of the form
	\begin{equation} \label{eq:HphiZ} H\tilde{\phi}(\un{Z}) \homeq \prod L(Z_{j_1,k_1} \smsh \dots \smsh Z_{j_t,k_t}) _{hG}\end{equation}
	a product of terms where $L: \finbased \to \cat{C}$ is linear, $j_l = j_{l'}$ for some $l \neq l'$, and $G$ is a subgroup of $\Sigma_m$ that acts by permuting some of the factors $Z_{j,k}$.
	
	By Lemma~\ref{lem:hom-class} below, each of the terms in (\ref{eq:HphiZ}) is $r_j$-homogeneous in the variable $\finbased^{n''_j}$ where $r_j$ is the number of times that the index $j$ appears in the list $j_1,\dots,j_t$. We know that there is some $j$ such that $r_j \geq 2$, which implies that each of these terms has trivial $A''$-excisive approximation. We therefore obtain the desired (\ref{eq:PAH}):
	\[ P_{A''}(H\tilde{\phi}) \homeq *. \]
	Putting this calculation together with our earlier induction, we get the equivalence (\ref{eq:pp}), which completes the proof that map $(ii)$ is an equivalence, once we have proved the following lemma.
\end{proof}

\begin{lemma} \label{lem:hom-class}
	Let $\cat{C}$ be a stable $\infty$-category. A functor $H: \finbased^k \to \cat{C}$ is $m$-homogeneous if and only if it can be written in the form
	\[ H(Y_1,\dots,Y_k) \homeq \prod_{m_1+\dots+m_k = m} L_{(m_1,\dots,m_k)}(Y_1^{\smsh m_1} \smsh \dots\smsh Y_k^{m_k})_{h(\Sigma_{m_1} \times \dots \times \Sigma_{m_k})} \]
	for a collection of linear (i.e.\ reduced and excisive) functors
	\[ L_{(m_1,\dots,m_k)}: \finbased \to \cat{C} \]
	indexed by the ordered partitions of $m$ into $k$ non-negative integers. The right-hand side of this equivalence involves the (homotopy) coinvariants for the action of the group $\Sigma_{m_1} \times \dots \times \Sigma_{m_k}$ that permutes each of the smash powers $Y_i^{\smsh m_i}$.
\end{lemma}
\begin{proof}
	Suppose $H$ is $m$-homogeneous, and let $\un{Y} = (Y_1,\dots,Y_k)$ be an object of $\finbased^k$. By~\cite[6.1.4.14]{lurie:2017}, we can write
	\[ H(\un{Y}) \homeq M(\un{Y},\dots,\un{Y})_{h\Sigma_m} \]
	where $M: (\finbased^k)^m \to \cat{C}$ is symmetric multilinear. By induction on Lemma~\ref{lem:A-excisive}, a linear functor $\finbased^k \to \cat{C}$ is of the form
	\[ L_1(Y_1) \times \dots \times L_k(Y_k) \]
	and it follows that the symmetric multilinear $M$ can be written as
	\[ M(\un{Y}_1,\dots,\un{Y}_m) \homeq \prod_{1 \leq j_1,\dots,j_m \leq k} L_{j_1,\dots,j_m}(Y_{1,j_1},\dots,Y_{m,j_m}) \]
	where the symmetric group $\Sigma_m$ acts by permuting the indexes $j_1,\dots,j_m$ as well as the inputs of each multilinear functor $L_{j_1,\dots,j_m}: \finbased^m \to \cat{C}$.
	
	We therefore get
	\[ H(\un{Y}) \homeq \prod_{m_1+\dots+m_k = m} L_{1,\dots,1,\dots,k,\dots,k}(Y_1, \dots, Y_1, \dots, Y_k, \dots, Y_k)_{h(\Sigma_{m_1} \times \dots \times \Sigma_{m_k})} \]
	where the index $i$ is repeated $m_i$ times in each term. Finally, by induction using the result of~\cite[1.4.2.22]{lurie:2017}, multilinear functors $\finbased^m \to \cat{C}$ factor via the smash product $\smsh: \finbased^m \to \finbased$ yielding the desired expression.
	
	Conversely, suppose $H$ is of the given form. It is sufficient to show that each functor of the form
	\[ F(Y_1,\dots,Y_k) \mapsto L(Y_1^{\smsh m_1} \smsh \dots \smsh Y_k^{\smsh m_k}) \]
	with $L$ linear, is $m$-homogeneous, since the finite product and coinvariants constructions preserve homogeneity of functors with values in a stable $\infty$-category. It can be shown directly from the definition that $F$ is $m_i$-excisive in its $i$-th variable, so that $F$ is $m$-excisive by~\cite[6.1.3.4]{lurie:2017}. To see that $F$ is also $m$-reduced, we apply~\cite[6.1.3.24]{lurie:2017} by directly calculating the relevant cross-effects of $F$.
\end{proof}

Finally, for this section, we note the following compatibility between natural transformations of the form $\alpha$, which will be crucial for establishing the higher homotopy coherence of our constructions in Chapter~\ref{sec:tangent-structure-catdiff}.

\begin{lemma} \label{lem:alpha}
	Let $\sigma$ be a $3$-simplex in $\Weilinfty$ with edges given by the labelled Weil-algebra morphisms $\phi_1,\phi_2,\phi_3,\phi_{12},\phi_{23},\phi_{123}$. Then the natural transformations determined by applying Definition~\ref{def:alpha} to the $2$-faces of $\sigma$ form a strictly commutative diagram, in the category of simplicial functors $\finsset^{n_3} \to \finsset^{n_0}$ and simplicial natural transformations, of the form:
	\[ \begin{diagram}
		\node{\widetilde{\phi_3\phi_2\phi_1}} \arrow{e,t}{\alpha_{12,3}} \arrow{s,l}{\alpha_{1,23}} \node{\widetilde{\phi_2\phi_1}\tilde{\phi}_3} \arrow{s,r}{\alpha_{1,2}\tilde{\phi}_3} \\
		\node{\tilde{\phi}_1\widetilde{\phi_3\phi_2}} \arrow{e,t}{\tilde{\phi}_1\alpha_{2,3}} \node{\tilde{\phi}_1\tilde{\phi}_2\tilde{\phi}_3}
	\end{diagram} \]
	Similarly, any $m$-simplex in $\Weilinfty$, for $m \geq 3$, determines a strictly commutative $(m-1)$-cube of functors $\finsset^{n_m} \to \finsset^{n_0}$.
\end{lemma}
\begin{proof}
	Recall from Proposition~\ref{prop:weil} that the $3$-simplex $\sigma$ can be represented by a diagram of finite partial commutative monoids indexed by the poset of intervals in $\{0,1,2,3\}$. The components of $\widetilde{\phi_3\phi_2\phi_1}$ are wedge summands indexed by elements $k \in \sigma([0,3])$. The pullback conditions in Proposition~\ref{prop:weil} imply that each such $k$ is uniquely determined by elements $k_1 \in \sigma([0,1]), k_2 \in \sigma([1,2]), k_3 \in \sigma([2,3])$ which satisfy $t(k_1) = s'(k_2)$ and $t'(k_2) = s''(k_3)$. Following through Definition~\ref{def:alpha}, we see that each composite map in the square above maps the wedge summand corresponding to $k$ by the canonical isomorphism to the wedge summand of $\tilde{\phi_1}\tilde{\phi_2}\tilde{\phi_3}$ corresponding to the triple $(k_1,k_2,k_3)$. Therefore, the diagram strictly commutes.
	
	For an $m$-simplex, $\sigma$, the same argument implies that each $2$-dimensional face of the corresponding $(m-1)$-cube commutes, hence so does the whole cube.
\end{proof}

\subsection*{Basic tangent structure properties}

We will show in Chapter~\ref{sec:tangent-structure-catdiff} that the constructions of Definitions~\ref{def:A-excisive},~\ref{def:TAF} and~\ref{def:Tphi} extend to a functor
\[ T: \Weilinfty \times \Catdiff \to \Catdiff. \]
Firstly, however, we will check that our definitions so far satisfy the various conditions needed to determine a tangent structure on $\Catdiff$, starting with the observation that $T$ corresponds to an action of the strict monoidal $\infty$-category $\Weilinfty$ on $\Catdiff$.

\begin{lemma} \label{lem:monoidal}
	Let $\cat{C}$ be a differentiable $\infty$-category. Then:
	\begin{enumerate} \itemsep=10pt
		\item there is a natural \emph{isomorphism} $T^{\mathbb{N}}(\cat{C}) \isom \cat{C}$ given by evaluation at the unique point in $\finbased^0 = *$;
		\item for Weil-algebras $A,A'$, there is a natural \emph{isomorphism}
		\[ T^{A'}(T^{A}(\cat{C})) \isom T^{A \otimes A'}(\cat{C}) \]
		given by restricting the isomorphism
		\[ \Fun(\finbased^{n'},\Fun(\finbased^{n},\cat{C})) \isom \Fun(\finbased^{n+n'},\cat{C}) \]
		to the subcategories of suitably excisive functors.
	\end{enumerate}
\end{lemma}
\begin{proof}
	These results follow immediately from the definition of $T^A(\cat{C})$ in (\ref{eq:TA}).
\end{proof}

We now verify that the tangent pullbacks in $\Weilinfty$ (see Proposition~\ref{prop:tangent-pb}) are preserved by the action map $T$.

\begin{lemma} \label{lem:foundational}
	Let $J,J'$ be disjoint finite sets, and $\cat{C}$ a differentiable $\infty$-category. Then there is a pullback of $\infty$-categories of the form
	\[ \begin{diagram}
		\node{T_{J \sqcup J'}(\cat{C})} \arrow{e} \arrow{s} \node{T_{J}(\cat{C})} \arrow{s} \\
		\node{T_{J'}(\cat{C})} \arrow{e} \node{\cat{C}}
	\end{diagram} \]
	where the horizontal maps are induced by the projection $T_{J'} \to 1$, and the vertical by $T_J \to 1$.
\end{lemma}
\begin{proof}
	To prove this lemma, we should first specify precisely what we mean by a `pullback of $\infty$-categories'. For the purposes of this proof, we take that condition to mean that the given diagram is a homotopy pullback in the Joyal model structure on simplicial sets, and hence a pullback in the $\infty$-category $\Catinf$. (In the proof of Theorem~\ref{thm:T} in the next chapter, we show that this diagram is also a pullback in the subcategory $\Catdiff$.)
	
	We can write the desired diagram in the following form
	\[ \begin{diagram}
		\node{\Exc(\finbased^{J \sqcup J'},\cat{C})} \arrow{e,t}{p'_1} \arrow{s,l}{p_1} \node{\Exc(\finbased^J,\cat{C})} \arrow{s,r}{p} \\
		\node{\Exc(\finbased^{J'},\cat{C})} \arrow{e,t}{p'} \node{\cat{C}}
	\end{diagram} \]
	where each map is given by evaluating either the $J$-indexed variables or the $J'$-indexed variables at $*$.
	
	First, we argue that the map $p$ is a fibration in the Joyal model structure. This is true for the corresponding projection
	\[ \Fun(\finbased^J,\cat{C}) \to \cat{C} \]
	because $* \to \finbased^J$ is a cofibration between cofibrant objects, and the Joyal model structure on simplicial sets is closed monoidal. Since $\Exc(\finbased^J,\cat{C})$ is a full subcategory of $\Fun(\finbased^J,\cat{C})$ that is closed under equivalences, $p$ is also a fibration by~\cite[2.4.6.5]{lurie:2009}.
	
	It is now sufficient to show that the induced map
	\[ \pi: \Exc(\finbased^J \times \finbased^{J'},\cat{C}) \to \Exc(\finbased^J,\cat{C}) \times_{\cat{C}} \Exc(\finbased^{J'},\cat{C}) \]
	is an equivalence of $\infty$-categories.
	
	To see this, we define an inverse map
	\[ \iota: \Exc(\finbased^J,\cat{C}) \times_{\cat{C}} \Exc(\finbased^{J'},\cat{C}) \to \Exc(\finbased^J \times \finbased^{J'},\cat{C}) \]
	that sends a pair $(F_1,F_2)$ of excisive functors with the property that $F_1(*) = F_2(*)$, to the functor
	\[ (X_1,X_2) \mapsto F_1(X_1) \times_{F_2(*)} F_2(X_2). \]
	It is simple to check that $\pi\iota \homeq \mathrm{id}$, and it follows from Lemma~\ref{lem:A-excisive} that $\iota\pi \homeq \mathrm{id}$.
\end{proof}

\subsection*{The vertical lift axiom}

The last substantial condition we need in order that the constructions earlier in this chapter underlie a tangent structure $T$ on the $\infty$-category $\Catdiff$ is that $T$ preserves the vertical lift pullback (\ref{lem:pb2}) in the category $\Weilinfty$. Applying Definition~\ref{def:Tphi} to the Weil-algebra morphisms that appear in that pullback square, and using Lemma~\ref{lem:A-excisive}, we reduce to the following result.

\begin{proposition} \label{prop:vertical-lift}
	Let $\cat{C}$ be a differentiable $\infty$-category. Then the following diagram is a pullback of $\infty$-categories:
	\[ \begin{diagram}
		\node{T(\cat{C}) \times_{\cat{C}} T(\cat{C})} \arrow{s,l}{p \times p} \arrow{e,t}{v} \node{T(T(\cat{C}))} \arrow{s,r}{T(p)} \\
		\node{\cat{C}} \arrow{e,t}{0} \node{T(\cat{C})}
	\end{diagram} \]
	where $v$ sends the pair $(L_1,L_2)$ of excisive functors to the $(1,1)$-excisive functor $\finbased^2 \to \cat{C}$ given by
	\[ (X,Y) \mapsto L_1(X \smsh Y) \times_{L_1(*) = L_2(*)} L_2(Y). \]
\end{proposition}
\begin{proof}
	As in the proof of Lemma~\ref{lem:foundational}, our goal is to show that this diagram is a homotopy pullback in the Joyal model structure on simplicial sets, and the argument for that lemma also shows that $T(p)$ is a fibration, so it is sufficient to show that the induced map
	\[ f: \Exc(\finbased,\cat{C}) \times_{\cat{C}} \Exc(\finbased,\cat{C}) \to \cat{C} \times_{T\cat{C}} T^2(\cat{C}) \]
	is an equivalence of $\infty$-categories.
	
	Using Lemma~\ref{lem:A-excisive} again, we can write $f$ as a map
	\[ f: \Exc(\finbased^2,\cat{C}) \to \Exc^{1,1_*}(\finbased^2,\cat{C}); \quad L \mapsto [(X,Y) \mapsto L(X \smsh Y,Y)] \]
	from the $\infty$-category of excisive functors $\finbased^2 \to \cat{C}$ to the $\infty$-category of functors $M: \finbased^2 \to \cat{C}$ that are excisive in each variable individually and reduced in the second variable (in the sense that the map $M(X,*) \weq M(*,*)$ is an equivalence for each $X \in \finbased$).
	
	To show that $f$ is an equivalence of $\infty$-categories, we describe an explicit homotopy inverse
	\[ g: \Exc^{1,1_*}(\finbased^2,\cat{C}) \to \Exc(\finbased^2,\cat{C}). \]
	For a functor $M: \finbased^2 \to \cat{C}$ we set
	\begin{equation} \label{eq:vert-inv} g(M)(X,Y) := M(X,S^0) \times_{M(*,S^0)} M(*,Y) \end{equation}
	where the map $M(*,Y) \to M(*,S^0)$ used to construct this pullback is induced by the null map $Y \to S^0$ (that maps every point in $Y$ to the basepoint in $S^0$).
	
	We should show that if $M$ is excisive in each variable and reduced in $Y$, then $g(M)$ is excisive as a functor $\finbased^2 \to \cat{C}$. Consider a pushout square in $\finbased^2$
	\[ \begin{diagram}
		\node{(X,Y)} \arrow{e} \arrow{s} \node{(X_1,Y_1)} \arrow{s} \\
		\node{(X_2,Y_2)} \arrow{e} \node{(X_0,Y_0)}
	\end{diagram} \]
	consisting of individual pushout squares in each variable. Applying $g(M)$ to this pushout square we get the square in $\cat{C}$ given by
	\[ \begin{diagram}
		\node{M(X,S^0) \times_{M(*,S^0)} M(*,Y)} \arrow{e} \arrow{s}
		\node{M(X_1,S^0) \times_{M(*,S^0)} M(*,Y_1)} \arrow{s} \\
		\node{M(X_2,S^0) \times_{M(*,S^0)} M(*,Y_2)} \arrow{e}
		\node{M(X_0,S^0) \times_{M(*,S^0)} M(*,Y_0)}
	\end{diagram} \]
	Since $M$ is excisive in each variable, this is a pullback of pullback squares, and hence is itself a pullback. Therefore, $g(M)$ is excisive.
	
	We then have, for $L \in \Exc(\finbased^2,\cat{C})$:
	\[ \begin{split} g(f(L))(X,Y) &= f(L)(X,S^0) \times_{f(L)(*,S^0)} f(L)(*,Y) \\
		&= L(X \smsh S^0,S^0) \times_{L(* \smsh S^0,S^0)} L(* \smsh Y,Y) \\
		&\homeq L(X,S^0) \times_{L(*,S^0)} L(*,Y).
	\end{split} \]
	Consider the diagram
	\[ \begin{diagram}
		\node{L(X,Y)} \arrow{e} \arrow{s} \node{L(X,*)} \arrow{e} \arrow{s} \node{L(X,S^0)} \arrow{e} \arrow{s} \node{L(X,*)} \arrow{s} \\
		\node{L(*,Y)} \arrow{e} \node{L(*,*)} \arrow{e} \node{L(*,S^0)} \arrow{e} \node{L(*,*)}
	\end{diagram} \]
	The first and third squares are pullbacks because $L$ is excisive. The composite of the second and third squares is also a pullback. Hence, the second square is a pullback by~\cite[4.4.2.1]{lurie:2009}, and so the composite of the first and second squares is also a pullback, again by~\cite[4.4.2.1]{lurie:2009}. This calculation implies that the natural map
	\[ L \to g(f(L)) \]
	is an equivalence.
	
	On the other hand we have, for $M \in \Exc^{1,1_*}(\finbased^2,\cat{C})$:
	\[ \begin{split}
		f(g(M))(X,Y) &= g(M)(X \smsh Y,Y) \\
		&= M(X \smsh Y,S^0) \times_{M(*,S^0)} M(*,Y).
	\end{split} \]
	We claim that $f(g(M))$ is equivalent to $M$. First note that $M$ factors via the pointed $\infty$-category $\cat{C}_{M(*,*)}$, of objects over/under $M(*,*)$. We can therefore assume without loss of generality that $\cat{C}$ is pointed and that $M(*,*)$ is a null object. Since the map $M(*,Y) \to M(*,S^0)$ is then null, we can write $f(g(M))$ as the product
	\begin{equation} \label{eq:fgM} f(g(M))(X,Y) \homeq \hofib[M(X \smsh Y,S^0) \to M(*,S^0)] \times M(*,Y). \end{equation}
	Now observe that $M \in \Exc^{1,1_*}(\finbased^2,\cat{C})$ can be viewed as a functor
	\[ \finbased \to \Exc_*(\finbased,\cat{C}); \quad X \mapsto M(X,-) \]
	from the \emph{pointed} $\infty$-category $\finbased$ to the \emph{stable} $\infty$-category $\Exc_*(\finbased,\cat{C}) \homeq \spectra(\cat{C})$ of linear functors $\finbased \to \cat{C}$. Any such functor splits off the image of the null object, so we have an equivalence
	\begin{equation} \label{eq:M} M(X,Y) \homeq \hofib[M(X,Y) \to M(*,Y)] \times M(*,Y). \end{equation}
	Comparing (\ref{eq:fgM}) and (\ref{eq:M}) we see that to get an equivalence $f(g(M)) \homeq M$, it is enough to produce an equivalence
	\[ D(X,Y) \homeq D(X \smsh Y,S^0) \]
	where $D: \finbased^2 \to \cat{C}$ is given by
	\[ D(X,Y) := \hofib(M(X,Y) \to M(*,Y)). \]
	This $D$ is given by reducing $M$ in its first variable, and it follows that $D$ is reduced and excisive in both variables, i.e.\ is multilinear. The equivalence we need is a consequence of the classification of multilinear functors $\finbased^2 \to \cat{C}$ for any $\infty$-category $\cat{C}$ with finite limits.
	
	To be explicit, it follows from \cite[1.4.2.22]{lurie:2017} that the evaluation map
	\[ \Exc^{1,1}_{*,*}(\finbased^2,\cat{C}) \to \Exc_*(\finbased,\cat{C}); \quad D \mapsto D(-,S^0) \]
	is an equivalence of $\infty$-categories, and it is clear that a one-sided inverse is given by
	\[ \Exc_*(\finbased,\cat{C}) \to \Exc^{1,1}_{*,*}(\finbased^2,\cat{C}); \quad F \mapsto F(X \smsh Y). \]
	Thus, this functor is a two-sided inverse, and we therefore get the desired natural equivalence
	\[ D(X,Y) \homeq D(X \smsh Y,S^0). \]
	This completes the proof that $M \homeq f(g(M))$, and hence the proof that $f$ is an equivalence of $\infty$-categories.
\end{proof}

\subsection*{Goodwillie tangent structure is cartesian}

The final thing we check in this section is that $\Catdiff$ has finite products given by the ordinary product of $\infty$-categories, which are preserved by the tangent bundle functor $T = T^W: \Catdiff \to \Catdiff$. This fact implies that the tangent structure on $\Catdiff$ is cartesian in the sense of Definition~\ref{def:cartesian}.

\begin{lemma} \label{lem:cartesian}
	For differentiable $\infty$-categories $\cat{C},\cat{C}'$, the product $\cat{C} \times \cat{C}'$ is differentiable, and is the product of $\cat{C}$ and $\cat{C}'$ in $\Catdiff$. Moreover, the projections determine an equivalence of $\infty$-categories
	\[ T(\cat{C} \times \cat{C}') \homeq T(\cat{C}) \times T(\cat{C}'). \]
\end{lemma}
\begin{proof}
	Since limits and colimits in the product $\infty$-category are detected in each term, the product $\cat{C} \times \cat{C}'$ is differentiable. A functor $\cat{D} \to \cat{C} \times \cat{C}'$, with $\cat{D}$ differentiable, preserves sequential colimits if and only if each component preserves sequential colimits. It follows that $\cat{C} \times \cat{C}'$ is the product in $\Catdiff$. Finally, the isomorphism
	\[ \Fun(\finbased,\cat{C}_1 \times \cat{C}_2) \isom \Fun(\finbased,\cat{C}_1) \times \Fun(\finbased,\cat{C}_2) \]
	restricts to the $\infty$-categories of excisive functors since a square in $\cat{C}_1 \times \cat{C}_2$ is a pullback if and only if it is a pullback in each factor.
\end{proof}

\chapter{The Goodwillie Tangent Structure: Formal Construction} \label{sec:tangent-structure-catdiff}
The constructions of Chapter~\ref{sec:underlying} contain the basic data and lemmas on which our tangent structure rests, but to have a tangent $\infty$-category we need to extend those definitions to an action of the monoidal $\infty$-category $\Weilinfty$ on a suitable $\infty$-category $\Catdiff$ of differentiable $\infty$-categories and sequential-colimit-preserving functors. We start by giving an explicit description of the $\infty$-category $\Catdiff$.

\subsection*{The $\infty$-category of differentiable $\infty$-categories}

We recall Lurie's model for the $\infty$-category of $\infty$-categories from~\cite[3.0.0.1]{lurie:2009}.

\begin{definition} \label{def:Catinf}
	Let $\Catinf$\index{Cat01@$\Catinf$, the $\infty$-category of $\infty$-categories} be the simplicial nerve~\cite[1.1.5.5]{lurie:2009} of the simplicial category whose objects are the $\infty$-categories, and for which the simplicial mapping spaces are the maximal Kan complexes inside the usual functor $\infty$-categories:
	\[ \Hom_{\Catinf}(\cat{C},\cat{D}) := \Fun(\cat{C},\cat{D})^{\homeq}, \]
	i.e.\ the subcategory of $\Fun(\cat{C},\cat{D})$ whose morphisms are the natural equivalences. An $n$-simplex in $\Catinf$ therefore consists of the following data:
	\begin{itemize}
		\item a sequence of $\infty$-categories $\cat{C}_0,\dots,\cat{C}_n$;
		\item for each $0 \leq i < j \leq n$, a map of simplicial sets
		\[ \lambda_{i,j}: \mathsf{P}_{i,j} \to \Fun(\cat{C}_i,\cat{C}_j) \]
		where $\mathsf{P}_{i,j}$ denotes (the nerve of) the poset of those subsets of
		\[ \{i,i+1,\dots,j-1,j\} \]
		that include both $i$ and $j$, ordered by inclusion;
	\end{itemize}
	such that
	\begin{enumerate}
		\item for each edge $e$ in $\mathsf{P}_{i,j}$, the natural transformation $\lambda_{i,j}(e)$ is an equivalence;
		\item for each $i < j < k$, the following diagram commutes:
		\[ \begin{diagram}
			\node{\mathsf{P}_{i,j} \times \mathsf{P}_{j,k}} \arrow{s,l}{\cup} \arrow{e,t}{\lambda_{i,j} \times \lambda_{j,k}} \node{\Fun(\cat{C}_i,\cat{C}_j) \times \Fun(\cat{C}_j,\cat{C}_k)} \arrow{s,r}{\circ} \\
			\node{\mathsf{P}_{i,k}} \arrow{e,t}{\lambda_{i,k}} \node{\Fun(\cat{C}_i,\cat{C}_k).}
		\end{diagram} \]
	\end{enumerate}
	In particular, a $2$-simplex in $\Catinf$ comprises three functors
	\[ \begin{diagram}
		\node{\cat{C}_0} \arrow{se,b}{F} \arrow[2]{e,t}{H} \node[2]{\cat{C}_2} \\
		\node[2]{\cat{C}_1} \arrow{ne,b}{G}
	\end{diagram} \]
	together with a natural equivalence $H \weq GF$.
\end{definition}

\begin{definition} \label{def:Catdiff}
	Let $\Catdiff \subseteq \Catinf$ be the maximal simplicial subset whose objects are the differentiable $\infty$-categories and whose morphisms are the functors that preserve sequential colimits. We refer to $\Catdiff$ as the \emph{$\infty$-category of differentiable $\infty$-categories}\index{Cat0diff@$\Catdiff$, the $\infty$-category of differentiable $\infty$-categories}.
	
	Note that $\Catdiff$ is the simplicial nerve of a simplicial category whose objects are the differentiable $\infty$-categories with simplicial mapping objects
	\[ \Hom_{\Catdiff}(\cat{C},\cat{D}) = \Fun_{\N}(\cat{C},\cat{D})^{\homeq}, \]
	the maximal Kan complex of the full subcategory $\Fun_{\N}(\cat{C},\cat{D}) \subseteq \Fun(\cat{C},\cat{D})$ of sequential-colimit-preserving functors. As is our convention, we do not distinguish notationally between the $\infty$-category $\Catdiff$ and this underlying simplicial category.
\end{definition}

\subsection*{Differentiable relative $\infty$-categories}

It turns out that $\Catdiff$ is not the most convenient $\infty$-category on which to build the Goodwillie tangent structure. In this section, we describe another $\infty$-category, denoted $\RelCatdiff$, which is equivalent to $\Catdiff$ and which is more amenable.

The objects of $\RelCatdiff$ are \emph{relative} $\infty$-categories. A \emph{relative category} is simply a category $\cat{C}$ together with a subcategory $\cat{W}$ of `weak equivalences' which contains all isomorphisms in $\cat{C}$. Associated to the pair $(\cat{C},\cat{W})$ is an $\infty$-category $\cat{C}[\cat{W}^{-1}]$ given by formally inverting the morphisms in $\cat{W}$. Barwick and Kan showed in~\cite{barwick/kan:2012} that any $\infty$-category can be obtained this way, so that relative categories are yet another model for $\infty$-categories.

Mazel-Gee~\cite{mazel-gee:2019} uses a nerve construction of Rezk to extend the localization construction to `relative $\infty$-categories' in which $\cat{C}$ and $\cat{W}$ are themselves allowed to be $\infty$-categories. In other words, the Rezk nerve describes a very general `calculus of fractions' for $\infty$-categories.

The reason that relative $\infty$-categories are convenient for us is that we can replace the $\infty$-category $\Exc^A(\finbased^{n},\cat{C})$ of $A$-excisive functors with the \emph{relative} $\infty$-category
\[ (\Fun(\finbased^{n},\cat{C}), P_A\cat{E}) \]
consisting of the full $\infty$-category of functors together with the subcategory consisting of those natural transformations that become equivalences on applying the $A$-excisive approximation $P_A$.

The benefit of this approach is that we can work directly with the functor $\infty$-categories $\Fun(\finbased^{n},\cat{C})$ without involving the explicit $P_A$-approximation functor. In particular, it makes the monoidal nature of our construction almost immediate.

\begin{definition} \label{def:rel}
	A \emph{relative $\infty$-category}\index{relative $\infty$-category} is a pair $(\cat{C},\cat{W})$ consisting of an $\infty$-category $\cat{C}$ and a subcategory $\cat{W} \subseteq \cat{C}$ that includes all equivalences in $\cat{C}$. (In particular $\cat{W}$ contains all the objects of $\cat{C}$.)
	
	A \emph{relative functor}\index{relative functor!between relative $\infty$-categories} $G: (\cat{C}_0,\cat{W}_0) \to (\cat{C}_1,\cat{W}_1)$ between relative $\infty$-categories is a functor $G: \cat{C}_0 \to \cat{C}_1$ such that $G(\cat{W}_0) \subseteq \cat{W}_1$.
	
	A \emph{natural transformation}\index{natural transformation!between relative functors} between relative functors $G,G': (\cat{C}_0,\cat{W}_0) \to (\cat{C}_1,\cat{W}_1)$ is a natural transformation between the functors $G,G': \cat{C}_0 \to \cat{C}_1$, that is, a functor $\alpha: \Delta^1 \times \cat{C}_0 \to \cat{C}_1$ that restricts to $G$ on $\{0\} \times \cat{C}_0$ and to $G'$ on $\{1\} \times \cat{C}_0$.
	
	We say that a natural transformation $\alpha$ is a \emph{relative equivalence}\index{relative equivalence} if for each $X \in \cat{C}_0$, the morphism $\alpha_X: G(X) \to G'(X)$ is in the subcategory $\cat{W}_1 \subseteq \cat{C}_1$. In this case, we also say that the natural transformation $\alpha$ \emph{takes values in $\cat{W}_1$}.
\end{definition}

\begin{example}
	Associated to any $\infty$-category is a \emph{minimal} relative $\infty$-category $(\cat{C},\cat{E}_{\cat{C}})$ where $\cat{E}_{\cat{C}}$ is the subcategory of equivalences in $\cat{C}$. A relative functor between minimal relative $\infty$-categories is just a functor between the underlying $\infty$-categories, and a relative equivalence between such functors is just a natural equivalence in the usual sense.
\end{example}

We now wish to restrict to those relative $\infty$-categories whose localization is a \emph{differentiable} $\infty$-category in the sense of Definition~\ref{def:diff-cat}. In fact, it will be convenient to consider not all such relative $\infty$-categories, but only those for which the localization is an exact reflective subcategory.

\begin{definition} \label{def:diffrel}
	We say that a relative $\infty$-category $(\cat{C},\cat{W})$ is \emph{differentiable}\index{differentiable relative $\infty$-category} if $\cat{C}$ is a differentiable $\infty$-category, and $\cat{W}$ is the subcategory of local equivalences for an exact localization functor $\cat{C} \to \cat{C}$ in the sense of~\cite[5.2.7]{lurie:2009}. In other words, there exists an adjunction of differentiable $\infty$-categories
	\[ f: \cat{C} \rightleftarrows \cat{D} : g \]
	such that
	\begin{itemize}
		\item $f$ preserves finite limits;
		\item $g$ is fully faithful and preserves sequential colimits;
		\item $\cat{W}$ is the subcategory of $f$-equivalences in $\cat{C}$, i.e.\ those morphisms that are mapped by $f$ to an equivalence in $\cat{D}$.
	\end{itemize}
	
	We will say that a relative functor $G: (\cat{C}_0,\cat{W}_0) \to (\cat{C}_1,\cat{W}_1)$ between differentiable relative $\infty$-categories is \emph{differentiable}\index{relative functor!differentiable} if its underlying functor $G: \cat{C}_0 \to \cat{C}_1$ preserves sequential colimits.
\end{definition}

Our task is now to build an $\infty$-category whose objects are the differentiable relative $\infty$-categories, and which is equivalent to $\Catdiff$. We start by giving a simplicial enrichment to the category of differentiable relative $\infty$-categories and differentiable relative functors.

\begin{definition} \label{def:relzcatdiff}
	Let $\RelzCatdiff$\index{RelzCatdiff@$\RelzCatdiff$, the simplicial category of differentiable relative $\infty$-categories} be the simplicial category in which:
	\begin{itemize}
		\item objects are the differentiable relative $\infty$-categories $(\cat{C},\cat{W})$;
		\item the simplicial mapping object $\Hom_{\RelzCatdiff}((\cat{C}_0,\cat{W}_0),(\cat{C}_1,\cat{W}_1))$ is given by the subcategory
		\[ \Fun^{\homeq}_{\N}((\cat{C}_0,\cat{W}_0),(\cat{C}_1,\cat{W}_1)) \subseteq \Fun(\cat{C}_0,\cat{C}_1) \]
		whose objects are the differentiable relative functors, and whose morphisms are the relative equivalences.
	\end{itemize}
\end{definition}

\begin{proposition} \label{prop:DK-relcatdiff}
	There is a Dwyer-Kan equivalence of simplicial categories (i.e.\ an equivalence in Bergner's model structure~\cite{bergner:2007})
	\[ M_0: \Catdiff \to \RelzCatdiff \]
	that sends each differentiable $\infty$-category $\cat{C}$ to the differentiable relative $\infty$-category $(\cat{C},\cat{E}_{\cat{C}})$ where $\cat{E}_{\cat{C}}$ is the subcategory of equivalences in $\cat{C}$, and on simplicial mapping objects is given by the equivalences (in fact, equalities):
	\[ \Fun^{\homeq}_{\N}(\cat{C}_0,\cat{C}_1) \weq \Fun^{\homeq}_{\N}((\cat{C}_0,\cat{E}_{\cat{C}_0}),(\cat{C}_1,\cat{E}_{\cat{C}_1})). \]
\end{proposition}
\begin{proof}
	First note that when $\cat{C}$ is a differentiable $\infty$-category, $(\cat{C},\cat{E}_{\cat{C}})$ is a differentiable relative $\infty$-category; the identity adjunction on $\cat{C}$ satisfies the conditions of Definition~\ref{def:diffrel}.
	
	Since $M_0$ is clearly fully faithful, it remains to show that it is essentially surjective on objects. So let $(\cat{C},\cat{W})$ be an arbitrary differentiable relative $\infty$-category. Then we know that $\cat{W}$ is the subcategory of $f$-equivalences for some adjunction of differentiable $\infty$-categories
	\[ f: \cat{C} \rightleftarrows \cat{D} :g \]
	that satisfies the conditions of Definition~\ref{def:diffrel}. The functors $f$ and $g$ both preserve sequential colimits, so they determine differentiable relative functors
	\[ f: (\cat{C},\cat{W}) \rightleftarrows (\cat{D},\cat{E}_{\cat{D}}) :g \]
	which we claim are isomorphisms in the homotopy category of the simplicial category $\RelzCatdiff$.
	
	The counit of the localizing adjunction $(f,g)$ is an equivalence $\epsilon: fg \weq 1_{\cat{D}}$ and therefore also a relative equivalence $fg \weq 1_{(\cat{D},\cat{E}_{\cat{D}})}$. Hence $fg = 1_{(\cat{D},\cat{E}_{\cat{D}})}$ in the homotopy category of $\RelzCatdiff$.
	
	It remains to produce a relative equivalence between $gf$ and $1_{(\cat{C},\cat{W})}$. The unit of the localizing adjunction $(f,g)$ is a natural transformation $1_{\cat{C}} \to gf$ between relative functors $(\cat{C},\cat{W}) \to (\cat{C},\cat{W})$. To test if this natural transformation is a relative equivalence, we must show that it becomes an equivalence in $\cat{D}$ after applying $f$. But $f\eta: f \to fgf$ is inverse to the equivalence $\epsilon f$, so this is indeed the case. Hence $gf = 1_{(\cat{C},\cat{W})}$ in the homotopy category too, and so $f$ and $g$ are isomorphisms as claimed. Thus $M_0$ is essentially surjective.
\end{proof}

Note that $\RelzCatdiff$ is enriched in $\infty$-categories but not in Kan complexes because a relative equivalence is not necessarily invertible. The simplicial nerve of $\RelzCatdiff$ is therefore not an $\infty$-category. In order to rectify this problem, we add inverses for the relative equivalences by taking a fibrant replacement for the simplicial mapping objects in $\RelzCatdiff$.

For this purpose, we use an explicit fibrant replacement for the Quillen model structure given by Kan's $\Exinf$ functor~\cite{kan:1957}. We do not give a complete definition of the functor $\Exinf: \sset \to \sset$, but here are the key properties from our point of view. For any simplicial set $Y$, the simplicial set $\Exinf(Y)$ is a fibrant replacement for $Y$ in the Quillen model structure. So $\Exinf(Y)$ is a Kan complex, and there is a natural map $r_Y: Y \weq \Exinf(Y)$ which is a weak equivalence and cofibration in that model structure.

The vertices of $\Exinf(Y)$ can be identified with the vertices of $Y$, and the edges of $\Exinf(Y)$ can be identified with zigzags
\[ y_0 \to y_1 \leftarrow y_2 \to \dots \leftarrow y_{2k} \]
of edges in $Y$, where we identify zigzags of different lengths by including additional identity morphisms on the right. The map $r$ sends an edge $y_0 \to y_1$ of $Y$ to the zigzag
\[ y_0 \to y_1 \arrow{e,=} y_1. \]
The functor $\Exinf$ preserves finite products and has simplicial enrichment coming from the composite
\[ \dgTEXTARROWLENGTH=2.5em X \times \Exinf(Y) \arrow{e,t}{r_X} \Exinf(X) \times \Exinf(Y) \isom \Exinf(X \times Y). \]
When $Y$ is an $\infty$-category, we can think of the Kan complex $\Exinf(Y)$ as a model for the `$\infty$-groupoidification' of $Y$ given by freely inverting the $1$-simplexes.

\begin{definition} \label{def:relcatdiff}
	Let $\RelCatdiff$ denote the simplicial category whose objects are the differentiable relative $\infty$-categories, and whose simplicial mapping spaces are the Kan complexes
	\[ \Hom_{\RelCatdiff}((\cat{C}_0,\cat{W}_0),(\cat{C}_1,\cat{W}_1)) := \Exinf \Fun^{\homeq}_{\N}((\cat{C}_0,\cat{W}_0),(\cat{C}_1,\cat{W}_1)) \]
	given by applying $\Exinf$ to the simplicial mapping objects in $\RelzCatdiff$. Composition in $\RelCatdiff$ is induced by that in $\RelzCatdiff$ using the fact that $\Exinf$ preserves finite products.
	
	There is a canonical functor
	\[ r: \RelzCatdiff \to \RelCatdiff \]
	given by the identity on objects and by inclusions of the form $r_Y: Y \to \Exinf(Y)$ on mapping spaces. Since $r_Y$ is a weak equivalence in the Quillen model structure, $r$ is a Dwyer-Kan equivalence.
	
	By construction, $\RelCatdiff$ is enriched in Kan complexes, and hence its simplicial nerve is an $\infty$-category which, following our usual convention, we also denote $\RelCatdiff$.\index{RelCatdiff@$\RelCatdiff$, the $\infty$-category of relative differentiable $\infty$-categories}
\end{definition}

\begin{corollary} \label{cor:DK-relcatdiff}
	There is an equivalence of $\infty$-categories
	\[ M: \Catdiff \weq \RelCatdiff \]
	given by composing the map $M_0$ from Proposition~\ref{prop:DK-relcatdiff} with the functor \[r: \RelzCatdiff \weq \RelCatdiff\] of Definition~\ref{def:relcatdiff}. 
\end{corollary}
\begin{proof}
	Each of $M_0$ and $r$ is a Dwyer-Kan equivalence, so their composite is a Dwyer-Kan equivalence between fibrant objects in the Bergner model structure. Taking simplicial nerves, we get an equivalence of $\infty$-categories.
\end{proof}

The equivalence $M$ tells us that we can use $\RelCatdiff$ as a model for the $\infty$-category of differentiable $\infty$-categories and sequential-colimit-preserving functors. We show in the next section that $\RelCatdiff$ admits the required tangent structure, which can then be transferred along $M$ to $\Catdiff$ using Lemma~\ref{lem:transfer}.

We conclude this section by giving an explicit description of the simplicial set $\RelCatdiff$, which will be useful in constructing our tangent structure.

\begin{remark} \label{rem:relcatdiff}
	An $n$-simplex $\lambda$ in the $\infty$-category $\RelCatdiff$ consists of the following data:
	\begin{itemize}
		\item a sequence of differentiable relative $\infty$-categories
		\[ (\cat{C}_0,\cat{W}_0), \dots, (\cat{C}_n,\cat{W}_n); \]
		\item for each $0 \leq i < j \leq n$, a functor
		\[ \lambda_{i,j}: \mathsf{P}_{i,j} \to \Exinf \Fun(\cat{C}_i,\cat{C}_j) \]
	\end{itemize}
	subject to the following conditions:
	\begin{enumerate}
		\item for each object $I \in \mathsf{P}_{i,j}$, $\lambda_{i,j}(I): \cat{C}_i \to \cat{C}_j$ is a differentiable relative functor;
		\item for each morphism $\iota: I \subseteq I'$ in $\mathsf{P}_{i,j}$, the edge
		\[ \lambda_{i,j}(\iota) \in \Exinf\Fun(\cat{C}_i,\cat{C}_j)_1 \]
		is a zigzag of relative equivalences in $\Fun(\cat{C}_i,\cat{C}_j)$;
		\item for each $i < j < k$, the following diagram commutes
		\[ \begin{diagram}
			\node{\mathsf{P}_{i,j} \times \mathsf{P}_{j,k}} \arrow{s,l}{\cup} \arrow{e,t}{\lambda_{i,j} \times \lambda_{j,k}} \node{\Exinf\Fun(\cat{C}_i,\cat{C}_j) \times \Exinf\Fun(\cat{C}_j,\cat{C}_k)} \arrow{s,r}{\Exinf(\circ)} \\
			\node{\mathsf{P}_{i,k}} \arrow{e,t}{\lambda_{i,k}} \node{\Exinf\Fun(\cat{C}_i,\cat{C}_k).}
		\end{diagram} \]
	\end{enumerate}
	Now let $\theta: [m] \to [n]$ be an order-preserving function. Then $\theta^*(\lambda)$ is the $m$-simplex given by the following data:
	\begin{itemize}
	\item the sequence of differentiable relative $\infty$-categories
	\[ (\cat{C}_{\theta(0)},\cat{W}_{\theta(0)}), \dots, (\cat{C}_{\theta(m)},\cat{W}_{\theta(m)}); \]
	\item for each $0 \leq i < j \leq m$, the functor
	\[ \lambda_{\theta(i),\theta(j)} \circ \theta(-): \mathsf{P}_{i,j} \to \Exinf \Fun(\cat{C}_{\theta(i)},\cat{C}_{\theta(j)}) \]
	where $\theta(-): \mathsf{P}_{i,j} \to \mathsf{P}_{\theta(i),\theta(j)}$ maps a subset $I \subseteq \{i,\dots,j\}$ to its image $\theta(I) \subseteq \{\theta(i),\dots,\theta(j)\}$.
	\end{itemize}	
\end{remark}

\subsection*{Tangent structure on differentiable relative $\infty$-categories}

We can now build a tangent structure on the $\infty$-category $\RelCatdiff$ by describing explicitly the corresponding action map
\[ T: \Weilinfty \times \RelCatdiff \to \RelCatdiff. \]

In order to make this a strict action of the monoidal $\infty$-category $\Weilinfty$, we need to be careful about one point. When we write a functor $\infty$-category of the form
\[ \Fun(\finbased^{n},\cat{C}) \]
we will actually mean the isomorphic simplicial set
\[ \Fun(\finbased,\Fun(\finbased,\dots,\Fun(\finbased,\cat{C})\dots)) \]
with $n$ iterations. It follows that the simplicial sets $\Fun(\finbased^{m},\Fun(\finbased^{n},\cat{C}))$ and $\Fun(\finbased^{m+n},\cat{C})$ are actually \emph{equal} not merely isomorphic.

With that warning in mind, we start by defining our desired functor $T$ on objects.

\begin{definition} \label{def:T0}
	Let $A$ be a Weil-algebra with $n$ generators, and $(\cat{C},\cat{W})$ a differentiable relative $\infty$-category. We define a relative $\infty$-category
	\[ T^A(\cat{C},\cat{W}) := (\Fun(\finbased^{n},\cat{C}),P_A\cat{W}) \]
	where $P_A\cat{W}$ is the subcategory of the functor $\infty$-category $\Fun(\finbased^{n},\cat{C})$ consisting of those morphisms, i.e.\ natural transformations $\beta: \Delta^1 \times \finbased^{n} \to \cat{C}$, that, after applying the $A$-excisive approximation functor $P_A$ from Proposition~\ref{prop:PA}, take values in $\cat{W}$. That is, for each $X \in \finbased^{n}$, the component $(P_A\beta)_X$ is in $\cat{W}$.
\end{definition}

\begin{lemma} \label{lem:T0}
	The relative $\infty$-category $T^A(\cat{C},\cat{W})$ is differentiable.
\end{lemma}
\begin{proof}
	Since $\cat{C}$ is differentiable, and limits and colimits in a functor $\infty$-category are calculated objectwise, it follows that $\Fun(\finbased^{n},\cat{C})$ is also differentiable. It therefore remains to show that $P_A\cat{W}$ is determined by a suitable localizing adjunction.
	
	Since $(\cat{C},\cat{W})$ is a differentiable relative $\infty$-category, there is an adjunction of differentiable $\infty$-categories
	\[ f: \cat{C} \rightleftarrows \cat{D} :g \]
	satisfying the conditions of Definition~\ref{def:diffrel}. Now consider the pair of adjunctions
	\[ \Fun(\finbased^{n},\cat{C}) \stackbin[\iota]{P_A}{\rightleftarrows} \Exc^A(\finbased^n,\cat{C}) \stackbin[g_*]{f_*}{\rightleftarrows} \Exc^A(\finbased^n,\cat{D}) \]
	where the maps $f_*$ and $g_*$ are given by composing with $f$ and $g$ respectively, noting that since these functors preserve finite limits, $f_*$ and $g_*$ preserve $A$-excisive functors.
	
	We verify that the composed adjunction $(f_*P_A,\iota g_*)$ satisfies the conditions of Definition~\ref{def:diffrel}:
	\begin{itemize}
		\item $P_A$ preserves finite limits by Proposition~\ref{prop:PA}, and $f_*$ preserves finite limits because $f$ does and those limits are calculated objectwise in $\Exc^A(\finbased^n,\cat{C})$ and $\Exc^A(\finbased^n,\cat{D})$ (also by~\ref{prop:PA});
		\item $\iota$ is fully faithful and preserves sequential colimits by~\ref{prop:PA}; $g_*$ is fully faithful because $g$ is, and because $\Fun(\finbased^{n},-)$ preserves fully faithful inclusion; $g_*$ preserves sequential colimits again because of~\ref{prop:PA};
		\item a morphism $\beta$ in $\Fun(\finbased^{n},\cat{C})$ is an $f_*P_A$-equivalence if and only if $P_A(\beta)$ is an $f_*$-equivalence if and only if (since equivalences in the $\infty$-category $\Exc^A(\finbased^n,\cat{D})$ are detected objectwise) $P_A(\beta)_X$ is an $f$-equivalence, i.e.\ in $\cat{W}$, for all $X \in \finbased^{n}$; thus the subcategory $P_A\cat{W}$ consists precisely of the $f_*P_A$-equivalences.
	\end{itemize}
	Thus $(\Fun(\finbased^n,\cat{C}),P_A\cat{W})$ is a differentiable relative $\infty$-category as claimed.
\end{proof}

Next, we define $T$ on morphisms.

\begin{definition} \label{def:T1}
	Let $\phi: A_0 \to A_1$ be a labelled morphism of Weil-algebras, and let $G: (\cat{C}_0,\cat{W}_0) \to (\cat{C}_1,\cat{W}_1)$ be a differentiable relative functor between differentiable relative $\infty$-categories.
	
	We let $T^{\phi}(G): (\Fun(\finbased^{n_0},\cat{C}_0),P_{A_0}\cat{W}_0) \to (\Fun(\finbased^{n_1},\cat{C}_1),P_{A_1}\cat{W}_1)$ be the relative functor
	\[ T^{\phi}(G): \Fun(\finbased^{n_0},\cat{C}_0) \to \Fun(\finbased^{n_1},\cat{C}_1); \quad L \mapsto GL\tilde{\phi} \]
	given by composition with the maps of simplicial sets $\tilde{\phi}: \finbased^{n_1} \to \finbased^{n_0}$ (of Definition~\ref{def:phi}) and $G: \cat{C}_0 \to \cat{C}_1$.
\end{definition}

\begin{lemma} \label{lem:T1}
	The functor $T^{\phi}(G)$ is a differentiable relative functor
	\[ T^{A_0}(\cat{C}_0,\cat{W}_0) \to T^{A_1}(\cat{C}_1,\cat{W}_1). \]
\end{lemma}
\begin{proof}
	We must show that $T^{\phi}(G)(P_{A_0}\cat{W}_0) \subseteq P_{A_1}\cat{W}_1$, so consider a morphism of $\Fun(\finbased^{n_0},\cat{C}_0)$, i.e.\ a natural transformation $L \to L'$, such that the induced map $(P_{A_0}L)(X) \to (P_{A_0}L')(X)$ is in $\cat{W}_0$ for all $X \in \finbased^{n_0}$.
	
	First, since $G$ is a relative functor, it follows that
	\[ G(P_{A_0}L)\tilde{\phi}(Y) \to G(P_{A_0}L')\tilde{\phi}(Y) \]
	is in $\cat{W}_1$ for any $Y \in \finbased^{n_1}$.
	
	Now recall that since $(\cat{C}_1,\cat{W}_1)$ is differentiable, there is an adjunction
	\[ f_1: \cat{C}_1 \rightleftarrows \cat{D}_1: g_1 \]
	satisfying the conditions of Definition~\ref{def:diffrel}. Thus $\cat{W}_1$ is the subcategory of $f_1$-equivalences, and so the map
	\[ f_1G(P_{A_0}L)\tilde{\phi} \to f_1G(P_{A_0}L')\tilde{\phi} \]
	is a natural equivalence between functors $\finbased^{n_1} \to \cat{D}_1$. Since $\cat{D}_1$ is differentiable, we can apply $P_{A_1}$ to this map to obtain a natural equivalence
	\[ P_{A_1}(f_1G(P_{A_0}L)\tilde{\phi}) \weq P_{A_1}(f_1G(P_{A_0}L')\tilde{\phi}). \]
	Since $f_1$ preserves both sequential colimits and finite limits, it commutes with the construction $P_{A_1}$ by Lemma~\ref{lem:PA}. Thus, we also have a natural equivalence
	\[ f_1P_{A_1}(G(P_{A_0}L)\tilde{\phi}) \weq f_1P_{A_1}(G(P_{A_0}L')\tilde{\phi}). \]
	In other words, the natural map
	\[ P_{A_1}(G(P_{A_0}L)\tilde{\phi}) \to P_{A_1}(G(P_{A_0}L')\tilde{\phi}) \]
	takes values in $\cat{W}_1$. But, applying the equivalences of (\ref{eq:PAPA}) and (\ref{eq:pp}), it follows that
	\[ P_{A_1}(GL\tilde{\phi}) \to P_{A_1}(GL'\tilde{\phi}) \]
	takes values in $\cat{W}_1$, so that
	\[ GL\tilde{\phi} \to GL'\tilde{\phi} \]
	is in $P_{A_1}\cat{W}_1$. So $T^{\phi}(G)$ is a relative functor as required.
	
	Finally, $T^{\phi}(G)$ preserves sequential colimits because $G$ does, and these colimits are calculated objectwise in the functor $\infty$-categories. So $T^{\phi}(G)$ is differentiable.
\end{proof}

Before moving on to the general case, it is worthwhile also to define $T$ explicitly on $2$-simplexes.

\begin{definition} \label{def:T2}
	Let $\phi$ be a $2$-simplex in $\Weilinfty$ consisting of, in part, a pair of labelled Weil-algebra morphisms
	\[ A_0 \arrow{e,t}{\phi_1} A_1 \arrow{e,t}{\phi_2} A_2 \]
	and a labelling on the composite $\phi_2\phi_1$.
	
	Let $\lambda$ be a $2$-simplex in $\RelCatdiff$. According to Remark~\ref{rem:relcatdiff}, $\lambda$ consists of a diagram of relative functors
	\[ \begin{diagram}
		\node{(\cat{C}_0,\cat{W}_0)} \arrow[2]{e,t}{H} \arrow{se,b}{F} \node[2]{(\cat{C}_2,\cat{W}_2)} \\
		\node[2]{(\cat{C}_1,\cat{W}_1)} \arrow{ne,b}{G}
	\end{diagram} \]
	together with an edge in $\Exinf \Fun(\cat{C}_0,\cat{C}_2)$, that is, a zigzag
	\[ \lambda_{0,1,2}: H \to E_1 \leftarrow \dots \to E_{2k-1} \leftarrow GF, \]
	in which each map is a relative equivalence.
	
	We define $T^{\phi}(\lambda)$ to be the $2$-simplex in $\RelCatdiff$ consisting of the corresponding diagram of relative functors
	\[ \begin{diagram}
		\node{T^{A_0}(\cat{C}_0,\cat{W}_0)} \arrow[2]{e,t}{T^{\phi_2\phi_1}(H)} \arrow{se,b}{T^{\phi_1}(F)} \node[2]{T^{A_2}(\cat{C}_2,\cat{W}_2)} \\
		\node[2]{T^{A_1}(\cat{C}_1,\cat{W}_1)} \arrow{ne,b}{T^{\phi_2}(G)}
	\end{diagram} \]
	together with the following zigzag of relative equivalences
	\begin{equation} \label{eq:zigzag} (H(-)\widetilde{\phi_2\phi_1}) \to (E_1(-)\tilde{\phi}_1\tilde{\phi}_2) \leftarrow \dots \to (E_{2k-1}(-)\tilde{\phi}_1\tilde{\phi}_2) \leftarrow (GF(-)\tilde{\phi}_1\tilde{\phi}_2) \end{equation}
	between $T^{\phi_2\phi_1}(H)$ and $T^{\phi_2}(G)T^{\phi_1}(F)$, in which each map is induced by the corresponding map in $\lambda_{0,1,2}$, and the first map, in addition, involves the natural transformation
	\[ \alpha: \widetilde{\phi_2\phi_1} \to \tilde{\phi}_1\tilde{\phi}_2 \]
	associated to the $2$-simplex $\phi$ by Definition~\ref{def:alpha}.
\end{definition}

\begin{lemma} \label{lem:T2}
	The construction of $T^{\phi}(\lambda)$ in Definition~\ref{def:T2} produces a $2$-simplex in $\RelCatdiff$.
\end{lemma}
\begin{proof}
	The only thing left to check is that each of the natural transformations in the zigzag (\ref{eq:zigzag}) is a relative equivalence, i.e.\ that each of these natural transformations takes values in $P_{A_2}\cat{W}_2$, that is, takes values in $\cat{W}_2$ after applying $P_{A_2}$.
	
	Let $\gamma: E \to E'$ be a relative equivalence between relative functors
	\[ E,E': (\cat{C}_0,\cat{W}_0) \to (\cat{C}_2,\cat{W}_2), \]
	i.e.\ for each $X \in \cat{C}_0$, the map $\gamma_X: E(X) \to E'(X)$ is in $\cat{W}_2$.
	
	It follows that for every functor $L: \finbased^{n_0} \to \cat{C}_0$, the map induced by $\gamma$
	\[ EL\tilde{\phi}_1\tilde{\phi}_2 \to E'L\tilde{\phi}_1\tilde{\phi}_2 \]
	takes values in $\cat{W}_2$. A similar argument to that in the proof of Lemma~\ref{lem:T1} implies that
	\[ P_{A_2}(EL\tilde{\phi}_1\tilde{\phi}_2) \to P_{A_2}(E'L\tilde{\phi}_1\tilde{\phi}_2) \]
	takes values in $\cat{W}_2$. Therefore, the map induced by $\gamma$,
	\[ E(-)\tilde{\phi}_1\tilde{\phi}_2 \to E'(-)\tilde{\phi}_1\tilde{\phi}_2 \]
	takes values in $P_{A_2}\cat{W}_2$ as required.
	
	This argument shows that each map in the zigzag (\ref{eq:zigzag}) after the first is a relative equivalence. To show that the first map is also a relative equivalence, we note that
	\[ P_{A_2}(E_1(-)\widetilde{\phi_2\phi_1}) \to P_{A_2}(E_1(-)\tilde{\phi}_1\tilde{\phi}_2) \]
	is an equivalence of the type described in (\ref{eq:alpha}), hence is in $\cat{W}_2$. Combined with the previous argument, we deduce that
	\[ H(-)\widetilde{\phi_2\phi_1} \to E_1(-)\tilde{\phi}_1\tilde{\phi}_2 \]
	is also a relative equivalence.
\end{proof}

We now extend our constructions above to simplexes of arbitrary dimension.

\begin{definition} \label{def:Tn}
	Let $\phi$ be an $n$-simplex in $\Weilinfty$. We note that $\phi$ consists of a sequence of labelled Weil-algebra morphisms
	\[ \dgTEXTARROWLENGTH=2em A_0 \arrow{e,t}{\phi_1} A_1 \arrow{e,t}{\phi_2} \cdots \arrow{e,t}{\phi_n} A_n \]
	together with compatible labellings on all composites.
	
	Also let $\lambda$ be an $n$-simplex in $\RelCatdiff$ as described in Remark~\ref{rem:relcatdiff}.
	
	We define $T^{\phi}(\lambda)$ to be the $n$-simplex in $\RelCatdiff$ consisting of the differentiable relative $\infty$-categories
	\[ T^{A_0}(\cat{C}_0,\cat{W}_0), \dots, T^{A_n}(\cat{C}_n,\cat{W}_n) \]
	and the functors
	\[ T^{\phi}(\lambda)_{i,j}: \mathsf{P}_{i,j} \to \Exinf \Fun(\Fun(\finbased^{n_i},\cat{C}_i),\Fun(\finbased^{n_j},\cat{C}_j)) \]
	constructed as follows.
	
	Firstly, from the $n$-simplex $\phi$, we construct for each $0 \leq i < j \leq n$ a functor (of ordinary categories)
	\[ \tilde{\phi}_{i,j}: \mathsf{P}_{i,j} \to \Fun_{\sset}(\finsset^{n_j},\finsset^{n_i}) \]
	where $\mathsf{P}_{i,j}$ is the poset of subsets of $\{i,i+1,\dots,j-1,j\}$ containing both $i$ and $j$, and $\Fun_{\sset}(\finsset^{n_j},\finsset^{n_i})$ is the category of simplicial functors and simplicial natural transformations.
	
	For the subset
	\[ I = \{i = i_0,i_1,\dots,i_k = j\} \subseteq \{i,i+1,\dots,j-1,j\}, \]
	we let $\tilde{\phi}_{i,j}(I)$ be the composite functor
	\begin{equation} \label{eq:tilde-phi} \widetilde{(\phi_{i_k} \cdots \phi_{i_{k-1}+1})} \cdots \widetilde{(\phi_{i_1} \cdots \phi_{i_0+1})}: \finsset^{n_j} \to \finsset^{n_i}. \end{equation}
	Here, where we write a `composite' of labelled Weil-algebra morphisms such as $\phi_{i_k} \cdots \phi_{i_{k-1}+1}$, we mean the corresponding edge of the $n$-simplex $\phi$ in $\Weilinfty$ with the labelling associated to it in $\phi$.
	
	For a morphism $e$ in $\mathsf{P}_{i,j}$, that is an inclusion $I \subseteq I'$ of subsets of $\{i,i+1,\dots,j-1,j\}$ that include $i$ and $j$, we have to produce a natural transformation
	\[ \tilde{\phi}_{i,j}(e): \tilde{\phi}_{i,j}(I) \to \tilde{\phi}_{i,j}(I'): \finbased^{n_j} \to \finbased^{n_i}. \]
	For example, when $e$ is the inclusion $\{0,2\} \subseteq \{0,1,2\}$, then $\tilde{\phi}_{0,2}(e)$ is required to be a natural transformation
	\[ \widetilde{\phi_2\phi_1} \to \tilde{\phi}_1 \tilde{\phi}_2 \]
	which we choose to be the map $\alpha$ of Definition~\ref{def:alpha}. For a general morphism $e$, the desired map $\tilde{\phi}_{i,j}(e)$ is obtained by (a composite of) generalizations of $\alpha$ to more than two factors. It follows from Lemma~\ref{lem:alpha} that, for inclusions $I \arrow{e,t}{e} I' \arrow{e,t}{e'} I''$, we have
	\[ \tilde{\phi}_{i,j}(e'e) = \tilde{\phi}_{i,j}(e')\tilde{\phi}_{i,j}(e) \]
	and so we obtain a functor $\tilde{\phi}_{i,j}: \mathsf{P}_{i,j} \to \Fun_{\sset}(\finsset^{n_j},\finsset^{n_i})$ as desired.
	
	We remark that the poset $\mathsf{P}_{i,j}$ is the shape of a $(j-i-1)$-cube, and the diagram $\tilde{\phi}_{i,j}$ is then precisely the cubical diagram referred to in Lemma~\ref{lem:alpha} associated to the $(j-i)$-dimensional face of the $n$-simplex $\sigma$ whose vertices are $A_i,\dots,A_j$.
	
	The functor $\tilde{\phi}_{i,j}$ corresponds to a simplicial functor
	\[ \mathsf{P}_{i,j} \times \finsset^{n_j} \to \finsset^{n_i}, \]
	and hence (taking simplicial nerves) a functor of $\infty$-categories
	\[ \mathsf{P}_{i,j} \times \finbased^{n_j} \to \finbased^{n_i}. \]
	We will also write $\tilde{\phi}_{i,j}$ for the corresponding functor
	\[ \tilde{\phi}_{i,j}: \mathsf{P}_{i,j} \to \Fun(\finbased^{n_j},\finbased^{n_i}). \]
	
	From the $n$-simplex $\lambda$, we also have a functor
	\[ \lambda_{i,j}: \mathsf{P}_{i,j} \to \Exinf \Fun(\cat{C}_i,\cat{C}_j) \]
	as described in Remark~\ref{rem:relcatdiff}. We then complete the definition of the $n$-simplex $T^{\phi}(\lambda)$ by defining $T^{\phi}(\lambda)_{i,j}$ to be the composite
	\[ \dgTEXTARROWLENGTH=3em \begin{split} \mathsf{P}_{i,j} &\arrow{e,t}{\langle \tilde{\phi}_{i,j}, \lambda_{i,j} \rangle} \Fun(\finbased^{n_j},\finbased^{n_i}) \times \Exinf \Fun(\cat{C}_i,\cat{C}_j) \\
		&\arrow{e,t}{\langle r, \mathrm{Id} \rangle} \Exinf \Fun(\finbased^{n_j},\finbased^{n_i}) \times \Exinf \Fun(\cat{C}_i,\cat{C}_j) \\
		&\arrow{e,t}{\isom} \Exinf (\Fun(\finbased^{n_j},\finbased^{n_i}) \times \Fun(\cat{C}_i,\cat{C}_j)) \\
		&\arrow{e,t}{\Exinf(c)} \Exinf \Fun(\Fun(\finbased^{n_i},\cat{C}_i),\Fun(\finbased^{n_j},\cat{C}_j)) \end{split} \]
	where
	\[ c: \Fun(\finbased^{n_j},\finbased^{n_i}) \times \Fun(\cat{C}_i,\cat{C}_j) \to \Fun(\Fun(\finbased^{n_i},\cat{C}_i),\Fun(\finbased^{n_j},\cat{C}_j)) \]
	is transpose to the composition map.
\end{definition}

\begin{lemma} \label{lem:Tn}
	The construction of Definition~\ref{def:Tn} defines an $n$-simplex $T^{\phi}(\lambda)$ in $\RelCatdiff$.
\end{lemma}
\begin{proof}
	We will verify each of the conditions in Remark~\ref{rem:relcatdiff}. First note that each $T^{A_i}(\cat{C}_i,\cat{W}_i)$ is a differentiable relative $\infty$-category by Lemma~\ref{lem:T0}.
	
	Now consider an object $I = \{i = i_0, i_1, \dots, i_k = j\} \in \mathsf{P}_{i,j}$. Then $T^{\phi}(\lambda)_{i,j}(I)$ is the functor
	\[ \Fun(\finbased^{n_i},\cat{C}_i) \to \Fun(\finbased^{n_j},\cat{C}_j) \]
	given by pre-composition with the functor
	\[ \widetilde{(\phi_{i_k} \cdots \phi_{i_{k-1}+1})} \cdots \widetilde{(\phi_{i_1} \cdots \phi_{i_0+1})}: \finbased^{n_j} \to \finbased^{n_i} \]
	and post-composition with the differentiable relative functor $\lambda_{i,j}(I): \cat{C}_i \to \cat{C}_j$. The argument of Lemma~\ref{lem:T1}, with $\tilde{\phi}$ replaced by the composite functor \[\widetilde{(\phi_{i_k} \cdots \phi_{i_{k-1}+1})} \cdots \widetilde{(\phi_{i_1} \cdots \phi_{i_0+1})},\] and $G$ replaced by $\lambda_{i,j}(I)$, implies that $T^{\phi}(\lambda)_{i,j}(I)$ is a differentiable relative functor. This verifies condition (1) of Remark~\ref{rem:relcatdiff} for our proposed $n$-simplex $T^{\phi}(\lambda)$.
	
	Next consider an edge $I \subseteq I'$ in $\mathsf{P}_{i,j}$. Then $T^{\phi}(\lambda)_{i,j}$ applied to that edge is a zigzag of natural transformations
	\[ T^{\phi}(\lambda)_{i,j}(I) \to \dots \leftarrow T^{\phi}(\lambda)_{i,j}(I'): \Fun(\finbased^{n_i},\cat{C}_i) \to \Fun(\finbased^{n_j},\cat{C}_j) \]
	induced by natural transformations of the type $\alpha: \widetilde{\phi_2\phi_1} \to \tilde{\phi}_1\tilde{\phi}_2$ and a zigzag $\lambda_{i,j}(I) \to \dots \leftarrow \lambda_{i,j}(I')$. The argument of Lemma~\ref{lem:T2}, slightly generalized, implies that each entry in this zigzag is a relative equivalence. This verifies condition (2) of Remark~\ref{rem:relcatdiff}.
	
	It remains to check condition (3), which is a large but easy diagram-chase in the category of simplicial sets, and which follows from the corresponding condition for the $n$-simplex $\lambda$, the naturality of $\Exinf$, and the fact that the following diagrams involving the (ordinary) functors $\tilde{\phi}_{i,j}$ commute:
	\[ \begin{diagram}
		\node{\mathsf{P}_{i,j} \times \mathsf{P}_{j,k}} \arrow{s,l}{\cup} \arrow{e,t}{\tilde{\phi}_{i,j} \times \tilde{\phi}_{j,k}} \node{\Fun_{\sset}(\finsset^{n_j},\finsset^{n_i}) \times \Fun_{\sset}(\finsset^{n_k},\finsset^{n_j})} \arrow{s,r}{\circ} \\
		\node{\mathsf{P}_{i,k}} \arrow{e,t}{\tilde{\phi}_{i,k}} \node{\Fun_{\sset}(\finsset^{n_k},\finsset^{n_i}).}
	\end{diagram} \]
	That last claim follows directly from the construction of $\tilde{\phi}_{i,j}$ in (\ref{eq:tilde-phi}).
\end{proof}

\begin{proposition} \label{prop:T}
	The construction of $T$ on simplexes in Definition~\ref{def:Tn} gives a strict action of the simplicial monoid $\Weilinfty$ on the simplicial set $\RelCatdiff$.
\end{proposition}
\begin{proof}
	We start by showing that the assignment $(\phi,\lambda) \mapsto T^{\phi}(\lambda)$ commutes with the simplicial structures. So take an order-preserving function $\theta: [m] \to [n]$. We want to show that
	\[ \theta^*(T^{\phi}(\lambda)) = T^{\theta^*(\phi)}(\theta^*(\lambda)). \]
	This claim follows from the simplicial structure on $\RelCatdiff$ described in Remark~\ref{rem:relcatdiff} and commutativity of the following diagram for $0 \leq i < j \leq m$
	\[ \begin{diagram}
	\node{\mathsf{P}_{i,j}} \arrow[2]{s,l}{\theta(-)} \arrow{se,t}{\widetilde{\theta^*(\phi)}_{i,j}} \\
	\node[2]{\Fun_{\sset}(\finsset^{n_{\theta(j)}},\finsset^{n_{\theta(i)}})} \\
	\node{\mathsf{P}_{\theta(i),\theta(j)}} \arrow{ne,b}{\tilde{\phi}_{\theta(i),\theta(j)}}
	\end{diagram} \]
	That in turn follows from (\ref{eq:tilde-phi}) and the fact that the $i$th edge in the $m$-simplex $\theta^*(\phi)$ is precisely the edge in $\phi$ given by the composite
	\[ \phi_{\theta(i)} \phi_{\theta(i)-1} \cdots \phi_{\theta(i-1)+1}. \]
	
	It now follows that Definition~\ref{def:Tn} determines a map of simplicial sets
	\[ T: \Weilinfty \times \RelCatdiff \to \RelCatdiff. \]
	To show that $T$ corresponds to a strict action of the simplicial monoid $\Weilinfty$ on $\RelCatdiff$, consider first the case of $0$-simplexes. We have to check that
	\[ T^{A'}T^A(\cat{C},\cat{W}) = T^{A' \otimes A}(\cat{C},\cat{W}). \]
	Recall that we have chosen our functor $\infty$-categories so that there is an \emph{equality} (not just an isomorphism)
	\[ \Fun(\finbased^{n'},\Fun(\finbased^{n},\cat{C})) = \Fun(\finbased^{n+n'},\cat{C}) \]
	so it is sufficient to show that the subcategories $P_{A'}(P_A\cat{W})$ and $P_{A' \otimes A}\cat{W}$ are equal.
	
	In other words, let $f: L \to L',$ be a natural transformation between functors $\finbased^{n+n'} \to \cat{C}$. Then $f$ is in $P_{A'}(P_A\cat{W})$ if and only if $P_{A'}L \to P_{A'}L'$ takes values in $P_A\cat{W}$, i.e.\ if
	\[ P_AP_{A'}L \to P_AP_{A'}L' \]
	takes values in $\cat{W}$. Since $P_AP_{A'} \homeq P_{A' \otimes A}$, this is the case if and only if $f$ is in $P_{A' \otimes A}\cat{W}$.
	
	Now let us turn to higher degree simplexes. We have to show that, for $n$-simplexes $\phi,\phi'$ in $\Weilinfty$, and an $n$-simplex $\lambda$ in $\RelCatdiff$, we have
	\[ T^{\phi' \otimes \phi}(\lambda) = T^{\phi'}T^{\phi}(\lambda)). \]
	That is, we have to show that two functors
	\[ \mathsf{P}_{i,j} \to \Exinf \Fun(\Fun(\finbased^{n_i'+n_i},\cat{C}_i),\Fun(\finbased^{n_j'+n_j},\cat{C}_j)) \]
	are equal. This is another large diagram-chase; the key fact is that the following diagram of ordinary functors strictly commutes:
	\[ \begin{diagram}
		\node[2]{\Fun_{\sset}(\finsset^{n'_j},\finsset^{n'_i}) \times \Fun_{\sset}(\finsset^{n_j},\finsset^{n_i})} \arrow[2]{s,r}{\times} \\
		\node{\mathsf{P}_{i,j}} \arrow{se,b}{\widetilde{(\phi' \otimes \phi)}_{i,j}} \arrow{ne,t}{\langle \tilde{\phi}'_{i,j}, \tilde{\phi}_{i,j} \rangle} \\
		\node[2]{\Fun_{\sset}(\finsset^{n'_j+n_j},\finsset^{n'_i+n_i}).}
	\end{diagram} \]
	This claim follows from the fact that in Definition~\ref{def:phi} we have that if
	\[ \tilde{\phi'}(Y'_1, \dots, Y'_{n'_2}) = (X'_1,\dots,X'_{n'_1}) \quad \text{and} \quad \tilde{\phi}(Y_1,\dots,Y_{n_2}) = (X_1,\dots,X_{n_1}), \]
	then
	\[ \widetilde{(\phi' \otimes \phi)}(Y'_1,\dots,Y'_{n'_2},Y_1,\dots,Y_{n_2}) = (X'_1,\dots,X'_{n'_1},X_1,\dots,X_{n_1}). \]
	\end{proof}

We now transfer the $\Weilinfty$-action on $\RelCatdiff$ along the equivalence of $\infty$-categories $M: \Catdiff \weq \RelCatdiff$ of Corollary~\ref{cor:DK-relcatdiff}, using Lemma~\ref{lem:transfer}.

\begin{definition} \label{def:T-Catdiff}
	Let $T^\bullet: \Weilinfty \to \Fun(\Catdiff,\Catdiff)$ be the monoidal functor given by composing the transpose of the action map of Proposition~\ref{prop:T} with the monoidal equivalence
	\[ \Fun(\RelCatdiff,\RelCatdiff) \homeq \Fun(\Catdiff,\Catdiff) \]
	associated to the equivalence $M$ of Corollary~\ref{cor:DK-relcatdiff}, as given by Lemma~\ref{lem:transfer}.
\end{definition}

We finally have the following result (completing the definition of the Goodwillie tangent structure and the proof of Theorem~\ref{thm:main}).

\begin{theorem} \label{thm:T}
	The monoidal functor
	\[ T^\bullet: \Weilinfty \to \Fun(\Catdiff,\Catdiff) \]
	of Definition~\ref{def:T-Catdiff} determines a cartesian tangent structure on the $\infty$-category $\Catdiff$, in the sense of Definition~\ref{def:tangent-structure-infty}, whose underlying endofunctor and projection are, up to equivalence, as described in~\ref{def:tan-cat} and following. 
\end{theorem}
\begin{proof}
	To see that $T^\bullet$ is a tangent structure on $\Catdiff$, we have to show that it preserves the foundational and vertical lift pullbacks. This claim follows from Lemma~\ref{lem:foundational} and Proposition~\ref{prop:vertical-lift} with one proviso; we must show that the homotopy pullbacks in the Joyal model structure appearing in those results are pullbacks in the $\infty$-category $\Catdiff$.
	
	We prove that claim by applying a result of Riehl and Verity~\cite[6.3.12]{riehl/verity:2022} with $\mathcal{K} = \mathsf{CAT}_{\infty}$ the $\infty$-cosmos of $\infty$-categories. (See~\cite[Ch.\ 1]{riehl/verity:2022} for an introduction to the theory of $\infty$-cosmoses.) We deduce from that result that there is a `cosmologically-embedded' sub-$\infty$-cosmos $\mathsf{CAT}^{\N}_{\infty} \subseteq \mathsf{CAT}_\infty$ whose objects are the $\infty$-categories that admit sequential colimits, and whose $1$-morphisms are the functors that preserve sequential colimits. 
	
	The claim that $\mathsf{CAT}^{\N}_{\infty}$ is cosmologically-embedded~\cite[6.3.3]{riehl/verity:2022} implies that any square diagram in $\mathsf{CAT}^{\N}_{\infty}$ that is a pullback along a fibration in $\mathsf{CAT}_{\infty}$ is also a pullback along a fibration in $\mathsf{CAT}^{\N}_{\infty}$. The pullbacks of Lemma~\ref{lem:foundational} and Proposition~\ref{prop:vertical-lift} fit that bill, and so they determine pullbacks in the corresponding $\infty$-category $\mathbb{C}\mathrm{at}^{\N}_{\infty}$ (of $\infty$-categories that admit sequential colimits and functors that preserve them), and hence also in the full subcategory $\Catdiff \subseteq \mathbb{C}\mathrm{at}^{\N}_{\infty}$.
	
	Finally, it follows from Lemma~\ref{lem:cartesian} that $\Catdiff$ has finite products which are preserved by the tangent bundle functor, so the tangent structure $T$ is cartesian.
\end{proof}

\chapter{Differential Objects are Stable $\infty$-Categories} \label{sec:catdiff-diff}

Having constructed the Goodwillie tangent structure, we now turn to its initial study, and in this chapter we look at its differential objects. Since the objects of $\Catdiff$ are \emph{differentiable} $\infty$-categories, there is a serious danger of confusing the words `differential' and `differentiable' in this chapter.

Recall from~\cite[1.1.3.4]{lurie:2017} that an $\infty$-category $\cat{C}$ is \emph{stable} if it is pointed, admits finite limits and colimits, and a square in $\cat{C}$ is a pushout if and only if it is a pullback. (In particular, a stable $\infty$-category admits biproducts which we denote by $\oplus$.) Also recall from~\cite[6.1.1.7]{lurie:2017} that a stable $\infty$-category is differentiable if and only if it admits countable coproducts. We then have the following simple characterization.

\begin{theorem} \label{thm:diff-st}
	A differentiable $\infty$-category $\cat{C}$ admits a differential structure within the Goodwillie tangent structure if and only if $\cat{C}$ is a stable $\infty$-category.
\end{theorem}
\begin{proof}
	We apply Corollary~\ref{cor:tangent-space-differential}, so it is sufficient to show that the tangent spaces in $\Catdiff$ are the stable $\infty$-categories. For any differentiable $\infty$-category $\cat{C}$ and object $X \in \cat{C}$, we can identify $T_X\cat{C}$ with the $\infty$-category of excisive functors $\finbased \to \cat{C}$ that map $*$ to $X$. We therefore have
	\[ T_X\cat{C} \homeq \Exc_*(\finbased,\cat{C}_{/X}) \]
	where the right-hand side is the $\infty$-category of \emph{reduced} excisive functors from $\finbased$ to the slice $\infty$-category $\cat{C}_{/X}$. Thus $T_X\cat{C}$ is stable by~\cite[1.4.2.16]{lurie:2017}.
	
	Conversely, if $\cat{C}$ is stable, then
	\[ T_*\cat{C} \homeq \Exc_*(\finbased,\cat{C}) \]
	which is equivalent to $\cat{C}$ by~\cite[1.4.2.21]{lurie:2017}. Therefore, $\cat{C}$ is equivalent to a tangent space and so admits a differential structure by~\ref{cor:tangent-space-differential}.
\end{proof}

The last part of this proof provides a canonical identification of any stable (differentiable) $\infty$-category $\cat{C}$ with a tangent space, and hence, by Proposition~\ref{prop:tangent-space-differential}, a canonical differential structure on each such $\cat{C}$. This observation allows us to define a cartesian differential structure (see~\cite{blute/cockett/seely:2009}) on the homotopy category of stable (differentiable) $\infty$-categories.

\begin{theorem} \label{thm:stable-cartesian-differential-category}
	Let $\Catst$ be the full subcategory of $\Catdiff$ whose objects are the stable differentiable $\infty$-categories. Then the homotopy category $h\Catst$ has a cartesian differential structure in which the monoid structure on an object $\cat{C}$ is given by the biproduct functor $\oplus: \cat{C} \times \cat{C} \to \cat{C}$, and the derivative of a morphism $F: \cat{C} \to \cat{D}$ is the `directional derivative'
	\[ \nabla(F): \cat{C} \times \cat{C} \to \cat{D} \]
	given by
	\[ \nabla(F)(V,X) := D_1(F(X \oplus -))(V). \]
	This formula denotes Goodwillie's linear approximation $D_1$ applied to the functor $F(X \oplus -): \cat{C} \to \cat{D}$ and evaluated at object $V \in \cat{C}$.
\end{theorem}
\begin{proof}
	It is not hard to verify directly that the axioms in~\cite[2.1.1]{blute/cockett/seely:2009} hold for $h\Catst$, but we deduce this theorem from the results of Chapter~\ref{sec:differential} in order to illustrate how the cartesian differential structure on $h\Catst$ is related to the Goodwillie tangent structure on $\Catdiff$.
	
	Recall that Theorem~\ref{thm:cartesian-differential-structure} determines a cartesian differential structure on the category $\widehat{h\Diff}(\Catdiff)$ whose objects are the differential objects in $\Catdiff$, and whose morphisms are maps in $h\Catdiff$ between underlying objects.
	
	We define a functor
	\[ T_*: h\Catst \to \widehat{h\Diff}(\Catdiff) \]
	by sending the stable differentiable $\infty$-category $\cat{C}$ to the differential object in $\Catdiff$ with underlying object the tangent space
	\[ T_*\cat{C} = \Exc_*(\finbased,\cat{C}) = \spectra(\cat{C}) \]
	and differential structure determined by Proposition~\ref{prop:tangent-space-differential}. For a functor $F: \cat{C} \to \cat{D}$, we define $T_*(F)$ to be the morphism in $h\Catdiff$ given by the composite
	\[ \begin{diagram} \node{T_*\cat{C}} \arrow{e,tb}{\Omega^\infty}{\sim} \node{\cat{C}} \arrow{e,t}{F} \node{\cat{D}} \node{T_*\cat{D}} \arrow{w,tb}{\Omega^\infty}{\sim} \end{diagram} \]
	where $\Omega^\infty: T_*\cat{C} = \spectra(\cat{C}) \weq \cat{C}$ denotes evaluation of a reduced excisive functor at $S^0$; see~\cite[1.4.2.21]{lurie:2017}. This definition makes $T_*$ into a fully faithful embedding since morphisms in both categories are taken from $h\Catdiff$.
	
	Given an object of $\widehat{h\Diff}(\Catdiff)$, i.e.\ a differential object $\cat{D}$ in the Goodwillie tangent structure, we know from Theorem~\ref{thm:diff-st} that the underlying $\infty$-category of $\cat{D}$ is stable. The equivalence
	\[ \Omega^\infty: T_*\cat{D} \weq \cat{D} \]
	is an isomorphism in $\widehat{h\Diff}(\Catdiff)$, so $T_*$ is essentially surjective on objects. Hence $T_*$ is an equivalence of categories, and we can transfer the cartesian differential structure from Theorem~\ref{thm:cartesian-differential-structure} along $T_*$ to $h\Catst$.
	
	It remains to show that this inherited cartesian differential structure is as claimed in the statement of the theorem. To calculate that structure, we first examine the differential structure on the tangent space $T_*\cat{C}$.
	
	Working through the construction of the functor $\cat{T}_\bullet$ in Proposition~\ref{prop:tangent-space-differential}, we see that the $E_\infty$-monoid structure on $T_*\cat{C}$ is given by restricting the `addition' in the tangent bundle $T\cat{C}$, i.e.\ is the objectwise product
	\[ +: T_*\cat{C} \times T_*\cat{C} \to T_*\cat{C}; \quad (L_1,L_2) \mapsto L_1(-) \times L_2(-). \]
	Since this map commutes with evaluation at $S^0$, we deduce that the corresponding monoid structure on $\cat{C}$ is also the product (and hence biproduct)
	\[ \oplus: \cat{C} \times \cat{C} \to \cat{C}. \]
	The map $\hat{p}: T(T_*\cat{C}) \to T_*\cat{C}$ associated to the differential structure on $T_*\cat{C}$ is precisely the map $g$ appearing in the proof of Proposition~\ref{prop:vertical-lift}, that is,
	\[ \hat{p}(M)(X) = \hofib[M(X,S^0) \to M(*,S^0)] \]
	where $M: \finbased^2 \to \cat{C}$ is excisive in both variables and reduced in its second variable. Transferring back to $\cat{C}$ along $\Omega^\infty$, we deduce that $\hat{p}: T(\cat{C}) \to \cat{C}$ is given by
	\[ \hat{p}(L) := \hofib[L(S^0) \to L(*)]. \]
	The derivative $\nabla(F)$ of a morphism $F: \cat{C} \to \cat{D}$ in $h\Catst$ is the composite
	\[ \begin{diagram} \node{\cat{C} \times \cat{C}} \node{T(\cat{C})} \arrow{w,tb}{\langle p,\hat{p} \rangle}{\sim} \arrow{e,t}{T(F)} \node{T(\cat{D})} \arrow{e,t}{\hat{p}} \node{\cat{D}.} \end{diagram} \]
	To evaluate this composite at $(X,V) \in \cat{C} \times \cat{C}$ we have to identify an excisive functor $L: \finbased \to \cat{C}$ such that $L(*) \homeq X$ and $\hofib[L(S^0) \to L(*)] \homeq V$. We can express this functor as
	\[ L(-) = X \oplus (- \otimes V) \]
	where $\otimes$ denotes the canonical tensoring of a pointed $\infty$-category $\cat{C}$ with finite colimits over $\finbased$. (The functor $\otimes: \finbased \times \cat{C} \to \cat{C}$ can be constructed using the characterization in~\cite[1.4.2.6]{lurie:2017} of the $\infty$-category $\spaces_{\mathrm{fin}}$ of finite spaces.) It follows that
	\[ \nabla(F)(X,V) = \hofib[P_1(F(X \oplus (- \otimes V)))(S^0) \to F(X)]. \]
	Since $- \otimes V$ commutes with colimits, it also commutes with the construction of $P_1$, so that we have
	\[ \nabla(F)(X,V) = \hofib[P_1(F(X \oplus -))(S^0 \otimes V) \to F(X)] \]
	which is precisely $D_1(F(X \oplus -))(V)$ as claimed.
\end{proof}

Theorem~\ref{thm:stable-cartesian-differential-category} is closely related to~\cite[Cor. 6.6]{bauer/johnson/osborne/riehl/tebbe:2018}, which is the result that first inspired this paper. Let us briefly discuss the connection.

\begin{definition} \label{def:hCatAb}
	Let $h\Catab$ be the category in which:
	\begin{itemize}
		\item an object is an abelian category;
		\item a morphism from $\cat{A}$ to $\cat{B}$ is a pointwise-chain-homotopy class of functors
		\[ F: \cat{A} \to \mathrm{Ch}_+(\cat{B}) \]
		where the target is the category of non-negatively-graded chain complexes of objects in $\cat{B}$, and two such functors $F,G$ are \emph{pointwise-chain-homotopy equivalent} if for each object $A \in \cat{A}$ there is a chain-homotopy equivalence $F(A) \homeq G(A)$.
	\end{itemize}
	Composition of morphisms is achieved via `Dold-Kan prolongation' of such an $F$ to a functor $\mathrm{Ch}_+(\cat{A}) \to \mathrm{Ch}_+(\cat{B})$, see~\cite[3.2]{bauer/johnson/osborne/riehl/tebbe:2018}.
\end{definition}

\begin{theorem}[{Bauer-Johnson-Osborne-Riehl-Tebbe~\cite[6.6]{bauer/johnson/osborne/riehl/tebbe:2018}}] \label{thm:bjort}
	The category $h\Catab$ has a cartesian differential structure in which the monoid structure on an object $\cat{A}$ is given by the biproduct $\oplus: \cat{A} \times \cat{A} \to \cat{A}$ and the derivative of a morphism $F: \cat{A} \to \mathrm{Ch}_+(\cat{B})$ is the `directional derivative' $\nabla(F): \cat{A} \times \cat{A} \to \mathrm{Ch}_+(\cat{B})$ given by
	\[ \nabla(F)(X,V) := D_1(F(X \oplus -))(V) \]
	where $D_1$ denotes the \emph{linearization} of a chain-complex-valued functor in the sense of Johnson and McCarthy~\cite{johnson/mccarthy:2004}.
\end{theorem}

Despite the close similarity, there are slight differences between the context of Theorems~\ref{thm:bjort} and~\ref{thm:stable-cartesian-differential-category} that prohibit a direct comparison. In particular, note that $h\Catab$ is defined using \emph{pointwise} chain-homotopy equivalence, rather than natural equivalence (though we suspect that \cite{bauer/johnson/osborne/riehl/tebbe:2018} could have been written entirely in terms of natural chain-homotopy equivalence instead).

Modulo that distinction, we speculate that there is an equivalence of cartesian differential categories between a subcategory of $h\Catab$ (say, on those abelian categories $\cat{A}$ that admit countable coproducts and suitably continuous functors) and a subcategory of $h\Catst$ (say, on the corresponding stable $\infty$-categories $N_{dg}\mathrm{Ch}(\cat{A})$ given by the differential graded nerves~\cite[1.3.2.10]{lurie:2017} of the categories of chain complexes on such $\cat{A}$). We do not pursue a precise equivalence of this form here.

\chapter{Jets and $n$-Excisive Functors} \label{sec:jet}

Goodwillie's notion of excisive functor played a central role in the construction of what we have called the Goodwillie tangent structure on the $\infty$-category $\Catdiff$ of differentiable $\infty$-categories. Our goal in this chapter is to show that the notions of $n$-excisive functor, for $n > 1$, are implicit in that tangent structure, so that the entirety of Goodwillie's theory can be recovered from it.

We actually describe how the notion of $n$-excisive \emph{equivalence}, i.e.\ the condition that a natural transformation determines an equivalence between $n$-excisive approximations, relates to the Goodwillie tangent structure. The notion from ordinary differential geometry that corresponds to $n$-excisive equivalence is that of `$n$-jet'. Recall that we say two smooth maps $f,g: M \to N$ between smooth manifolds \emph{agree to order $n$ at $x \in M$} if $f(x) = g(x)$, and the (multivariable) Taylor series of $f$ and $g$ in local coordinates at $x$ agree up to degree $n$. The \emph{$n$-jet at $x$} of the map $f$ is its equivalence class under the relation of agreeing to order $n$ at $x$.

We can interpret the Taylor series condition in terms of the standard tangent structure on the category $\Mfld$. Let $T^n_x(M)$ denote the \emph{$n$-fold tangent space to $M$ at $x$}, that is, the fibre of the projection map $T^n(M) \to M$ over the point $x$, where $T^n(M)$ is the $n$-fold iterate of the tangent bundle functor $T$. A smooth map $f:M \to N$ then induces a smooth map
\[ T^n_x(f): T^n_x(M) \to T^n_{f(x)}(N) \subseteq T^n(N) \]
i.e.\ the restriction of $T^n(f)$ to $T^n_x(M)$.

\begin{lemma} \label{lem:jet}
Let $f,g:M \to N$ be smooth maps. Then $f$ and $g$ have Taylor series at $x$ that agree to degree $n$ if and only if $T^n_x(f) = T^n_x(g)$.
\end{lemma}

The main result of this chapter is an analogue of Lemma~\ref{lem:jet} that connects the higher excisive approximations in Goodwillie calculus to the Goodwillie tangent structure on $\Catdiff$ constructed in Chapter~\ref{sec:tangent-structure-catdiff}. We refer the reader to \cite{goodwillie:2003} for the original theory of $n$-excisive approximation, and to~\cite[6.1.1]{lurie:2017} for the generalization of that theory to (differentiable) $\infty$-categories.

\begin{definition} \label{def:PXn}
	Let $F: \cat{C} \to \cat{D}$ be a sequential-colimit-preserving functor between differentiable $\infty$-categories, and suppose that $\cat{C}$ admits finite colimits. For an object $X \in \cat{C}$, the \emph{$n$-excisive approximation to $F$ at $X$} is
	\[ P^X_nF := P_n(F_{/X}) : \cat{C}_{/X} \to \cat{D} \]
	that is, the $n$-excisive approximation of the restriction of $F$ to the slice $\infty$-category $\cat{C}_{/X}$ of objects over $X$.
\end{definition}

\begin{theorem} \label{thm:jet}
	Let $F,G: \cat{C} \to \cat{D}$ be morphisms in $\Catdiff$, where $\cat{C}$ admits finite colimits. Let $\alpha: F \to G$ be a natural transformation, and let $X$ be an object of $\cat{C}$. Then $\alpha$ induces an equivalence
	\[ P^X_n\alpha: P^X_nF \weq P^X_nG \]
	of $n$-excisive approximations at $X$ if and only if $\alpha$ induces an equivalence
	\[ T^n(\alpha)\iota_X: T^n(F)\iota_X \weq T^n(G)\iota_X \]
	in the functor $\infty$-category
	\[ \Fun(T^n_X(\cat{C}),T^n(\cat{D})), \]
	where $T^n_X(\cat{C})$ is the fibre over $X$ of the projection $T^n(\cat{C}) \to \cat{C}$, and
	\[ \iota_X: T^n_X(\cat{C}) \to T^n(\cat{C}) \]
	is the inclusion of that fibre.
\end{theorem}
\begin{proof}
	We start by noting that each of the conditions in question implies that $\alpha_X: F(X) \weq G(X)$ is an equivalence: if $P^X_n\alpha$ is an equivalence, then so is
	\[ P^X_0\alpha \homeq \alpha_X, \]
	and if $T^n(\alpha)\iota_X$ is an equivalence, then so is
	\[ p^n_{\cat{D}}T^n(\alpha)\iota_X \homeq \alpha_X. \]
	Replacing $\cat{C}$ with $\cat{C}_{/X}$ and $\cat{D}$ with $\cat{D}_{/G(X)}$, we can now reduce to the case that $X$ is a terminal object in $\cat{C}$, and $F,G$ are reduced.
	
	In this case, we can identify the objects of $T^n_X(\cat{C})$ with the functors $L: \finbased^n \to \cat{C}$ that are excisive in each variable separately and satisfy $L(*,\dots,*) \homeq *$. Then
	\[ T^n_X(F): T^n_X(\cat{C}) \to T^n(\cat{D}) \]
	corresponds to the map
	\[ L \mapsto P_{1,\dots,1}(FL). \]
	Our goal is therefore to show that $\alpha: F \to G$ induces an equivalence
	\[ P_n\alpha: P_nF \to P_nG \]
	if and only if it induces an equivalence
	\[ P_{1,\dots,1}(\alpha L): P_{1,\dots,1}(FL) \to P_{1,\dots,1}(GL) \]
	for every reduced and $(1,\dots,1)$-excisive functor $L: \finbased^n \to \cat{C}$.
	
	Suppose first that $P_n\alpha$ is an equivalence, and consider the commutative diagram
	\begin{equation} \label{eq:Pn-diag} \begin{diagram}
			\node{P_{1,\dots,1}(FL)} \arrow{e} \arrow{s} \node{P_{1,\dots,1}((P_nF)L)} \arrow{s} \\
			\node{P_{1,\dots,1}(P_n(FL))} \arrow{e} \node{P_{1,\dots,1}(P_n((P_nF)L))}
	\end{diagram} \end{equation}
	The vertical maps are given by $n$-excisive approximation and are equivalences since being $(1,\dots,1)$-excisive is a stronger condition than being $n$-excisive, by~\cite[6.1.3.4]{lurie:2017}. The bottom horizontal map is an equivalence by Lemma~\ref{lem:PnGF}, so the top map is too. From the assumption that $P_n\alpha: P_nF \to P_nG$ is an equivalence, it therefore follows that $P_{1,\dots,1}(\alpha L)$ is an equivalence too.
	
	Conversely, suppose that $P_{1,\dots,1}(\alpha L)$ is an equivalence for any reduced functor $L: \finbased^n \to \cat{C}$ that is excisive in each variable. Note that this same condition then holds for any reduced $L$ regardless of it being excisive; that claim follows from the equivalences
	\[ P_{1,\dots,1}(FL) \weq P_{1,\dots,1}(F(P_{1,\dots,1}L)) \]
	given by (\ref{eq:PAPA}).
	
	Our strategy for showing that $P_n\alpha$ is an equivalence is to use induction on the Taylor tower. Consider the pullback diagram
	\begin{equation} \label{eq:hom} \begin{diagram}
			\node{P_kF} \arrow{e} \arrow{s} \node{P_{k-1}F} \arrow{s} \\
			\node{*} \arrow{e} \node{R_kF}
	\end{diagram} \end{equation}
	of~\cite[6.1.2.4]{lurie:2017} (Goodwillie's delooping theorem for homogeneous functors), in which $R_kF: \cat{C} \to \cat{D}$ is $k$-homogeneous. By induction on $k$, it suffices to show that $\alpha$ induces an equivalence $R_k\alpha: R_kF \to R_kG$ for each $1 \leq k \leq n$.
	
	It is clear from (\ref{eq:hom}) that the construction $R_k$ naturally takes values in the $\infty$-category $\cat{D}_*$ of pointed objects in $\cat{D}$. By~\cite[6.1.2.11]{lurie:2017}, it is sufficient to show that $R_k\alpha$ is an equivalence of $k$-homogeneous functors $\cat{C}_* \to \cat{D}_*$ between pointed $\infty$-categories. Such functors naturally factor via $\spectra(\cat{D}_*)$ by~\cite[6.1.2.9]{lurie:2017}, so it is sufficient to show that $\Omega R_k\alpha$ is an equivalence, i.e.\ that $\alpha$ induces an equivalence
	\[ D_k\alpha: D_kF \to D_kG \]
	between the $k$-th layers of the Taylor towers of functors $\cat{C}_* \to \cat{D}_*$.
	
	To show that $D_k\alpha$ is an equivalence, we show that $\alpha$ induces an equivalence on multilinearized cross-effects. Since $\cat{C}_*$ is pointed and admits finite colimits, it has a canonical tensoring over $\finbased$, i.e.\ a functor
	\[ \otimes: \finbased \times \cat{C}_* \to \cat{C}_* \]
	that preserves finite colimits in each variable, and such that $S^0 \otimes Y \homeq Y$ for each $Y \in \cat{C}_*$. The functor $\otimes$ can be constructed from the characterization in~\cite[1.4.2.6]{lurie:2017} of the $\infty$-category $\spaces_{\mathrm{fin}}$ of finite spaces.
	
	Take objects $A_1,\dots,A_n \in \cat{C}_*$ and consider the functor
	\[ L: \finbased^n \to \cat{C}_*; \quad (X_1,\dots,X_n) \mapsto (X_1 \otimes A_1) \wdge \dots \wdge (X_n \otimes A_n) \]
	where $\wdge$ is the coproduct in $\cat{C}_*$. Our hypothesis on $\alpha$ implies that $P_{1,\dots,1}(\alpha L)$ is an equivalence.
	
	Since the functor
	\[ \finbased^n \to \cat{C}_*^n; \quad (X_1,\dots,X_n) \mapsto (X_1 \otimes A_1, \dots, X_n \otimes A_n) \]
	preserves colimits (including the null object) in each variable, it commutes with $P_{1,\dots,1}$, by~\cite[6.1.1.30]{lurie:2017}. It follows that $\alpha$ induces an equivalence
	\[ P_{1,\dots,1}(\alpha(- \wdge \dots \wdge -))(- \otimes A_1,\dots,- \otimes A_n). \]
	Evaluating at $(S^0,\dots,S^0)$, we deduce that $\alpha$ induces an equivalence (of functors $\cat{C}_*^n \to \cat{D}_*$):
	\begin{equation} \label{eq:creff} P_{1,\dots,1}(F(- \wdge \dots \wdge -)) \weq P_{1,\dots,1}(G(- \wdge \dots \wdge -)). \end{equation}
	For any $1 \leq k \leq n$, the $k$-th cross-effect of $F$, see~\cite[6.1.3.20]{lurie:2017}, is the total homotopy fibre of a $k$-cube whose entries are functors of the form
	\[ F(- \wdge \dots \wdge -) \]
	with some subset of its arguments replaced by the null object $*$ in $\cat{C}_*$. Since $P_{1,\dots,1}$ commutes with the construction of that total homotopy fibre, it follows from (\ref{eq:creff}) that $\alpha$ induces an equivalence
	\[ P_{1,\dots,1}(\creff_kF) \weq P_{1,\dots,1}(\creff_kG) \]
	for all $1 \leq k \leq n$. It follows by~\cite[6.1.3.23]{lurie:2017} that $\alpha$ then induces equivalences
	\[ \creff_k(D_kF) \to \creff_k(D_kG) \]
	and hence also, by~\cite[6.1.4.7]{lurie:2017}, equivalences
	\[ D_kF \to D_kG \]
	for $1 \leq k \leq n$, as required. This completes the proof that $P_n\alpha: P_nF \weq P_nG$ is an equivalence.
\end{proof}

\begin{remark}
	Theorem~\ref{thm:jet} explains how the notion of $n$-excisive approximation is related to the Goodwillie tangent structure on $\Catdiff$. However, it is not quite true to say that this notion is fully encoded in that tangent structure, since the statement of Theorem~\ref{thm:jet} relies on natural transformations that are not equivalences, and hence are not part of the $\infty$-category $\Catdiff$. We can include those natural transformations by replacing $\Catdiff$ with a corresponding $\infty$-bicategory $\CATdiff$, which is the subject of the next chapter.
\end{remark}

\chapter{The Tangent $(\infty,2)$-Category of Differentiable $\infty$-Categories} \label{sec:infty2-catdiff}

The goal of this chapter is to show that the Goodwillie tangent structure on $\Catdiff$ extends to a tangent structure, in the sense of Example~\ref{ex:tangent-infty-bicategories}, on an $\infty$-bicategory $\CATdiff$ of differentiable $\infty$-categories. We start by defining that object.

\begin{definition} \label{def:CATinf}
	Let $\CATdiff$ be the nerve of the simplicial category whose objects are the differentiable $\infty$-categories, with simplicial mapping objects given by the $\infty$-categories
	\[ \Hom_{\CATdiff}(\cat{C},\cat{D}) := \Fun_{\N}(\cat{C},\cat{D}) \]
	of sequential-colimit-preserving functors. Since each mapping object is an $\infty$-category, $\CATdiff$ is an $\infty$-bicategory.
\end{definition}

Our construction of a tangent structure on the $\infty$-bicategory $\CATdiff$ follows a similar path to that on the $\infty$-category $\Catdiff$. We again start with relative differentiable $\infty$-categories (Definition~\ref{def:diffrel}).

\begin{definition}
	Let $\RelzCATdiff$ be the simplicial category in which:
	\begin{itemize}
		\item objects are the differentiable relative $\infty$-categories $(\cat{C},\cat{W})$;
		\item mapping simplicial sets are the full subcategories
		\[ \Fun_{\N}((\cat{C}_0,\cat{W}_0),(\cat{C}_1,\cat{W}_1)) \subseteq \Fun(\cat{C}_0,\cat{C}_1) \]
		consisting of the differentiable relative functors, i.e.\ functors $G:\nolinebreak \cat{C}_0 \to\nolinebreak \cat{C}_1$ which preserve sequential colimits and for which $G(\cat{W}_0) \subseteq \cat{W}_1$.
	\end{itemize}
\end{definition}

The nerve of $\RelzCATdiff$ is an $\infty$-bicategory, but to get an accurate model for $\CATdiff$ we still need to invert the relative equivalences in the simplicial mapping objects, just as we did in constructing the $\infty$-category $\RelCatdiff$ in Definition~\ref{def:relzcatdiff}. 

To accomplish this inversion, we now work with the category $\scset$ of \emph{scaled simplicial sets} and the scaled model structure on $\scset$; see Definition~\ref{def:scaled}. Recall that $\infty$-bicategories can be identified with fibrant objects in that model structure.

\begin{definition} \label{def:pushout}
	Let $\ReloCATdiff$ be the scaled simplicial set given by the pushout (in the category of scaled simplicial sets):
	\begin{equation} \label{eq:pushout} \begin{diagram}
			\node{\RelzCatdiff} \arrow{e,t}{r} \arrow{s} \node{\RelCatdiff} \arrow{s} \\
			\node{\RelzCATdiff} \arrow{e} \node{\ReloCATdiff}
	\end{diagram} \end{equation}
	where the top horizontal map $r$ is described in Definition~\ref{def:relcatdiff}, and the left-hand vertical map is determined by the inclusions
	\[ \Fun_{\N}^{\homeq}((\cat{C}_0,\cat{W}_0),(\cat{C}_1,\cat{W}_1)) \subseteq \Fun_{\N}((\cat{C}_0,\cat{W}_0),(\cat{C}_1,\cat{W}_1)). \]
	Both maps are monomorphisms of $\infty$-bicategories, hence cofibrations in the scaled model structure on $\scset$, so the square is a homotopy pushout in that model structure. The scaled simplicial set $\ReloCATdiff$ is not itself an $\infty$-bicategory, but we will eventually take a fibrant replacement of it in the scaled model structure to obtain an $\infty$-bicategory $\RelCATdiff$ on which to define the Goodwillie tangent structure.
\end{definition}

Before that, however, we show that the scaled simplicial sets in (\ref{eq:pushout}) admit strict actions by the simplicial monoid $\Weilinfty$ that are compatible with the inclusions, and hence determine a strict $\Weilinfty$-action on the pushout $\ReloCATdiff$.

\begin{definition} \label{def:Weil-RelzCATdiff}
	We define a map of simplicial sets
	\[ T: \Weilinfty \times \RelzCATdiff \to \RelzCATdiff \]
	following a similar, but simpler, pattern to that in Definition~\ref{def:Tn}. Recall from there that an $n$-simplex $\phi$ in $\Weilinfty$, with underlying labelled Weil-algebra morphisms
	\[ \dgTEXTARROWLENGTH=2em A_0 \arrow{e,t}{\phi_1} A_1 \arrow{e,t}{\phi_2} \cdots \arrow{e,t}{\phi_n} A_n, \]
	determines a collection of functors
	\[ \tilde{\phi}_{i,j}: \mathsf{P}_{i,j} \to \Fun(\finbased^{n_j},\finbased^{n_i}). \]
	An $n$-simplex $\lambda$ in $\RelzCATdiff$ consists of a sequence of differentiable relative $\infty$-categories
	\[ (\cat{C}_0,\cat{W}_0), \dots, (\cat{C}_n,\cat{W}_n) \]
	and a collection of functors
	\[ \lambda_{i,j}: \mathsf{P}_{i,j} \to \Fun(\cat{C}_i,\cat{C}_j) \]
	which take values in the subcategories of differentiable relative functors
	\[ \Fun_{\N}((\cat{C}_i,\cat{W}_i),(\cat{C}_j,\cat{W}_j)) \subseteq \Fun(\cat{C}_i,\cat{C}_j). \]
	We define $T^{\phi}(\lambda)$ to consist of the sequence
	\[ T^{A_0}(\cat{C}_0,\cat{W}_0), \dots, T^{A_n}(\cat{C}_n,\cat{W}_n) \]
	together with the functors $T^{\phi}(\lambda)_{i,j}$ given by the composite
	\begin{equation} \label{eq:Weil-relcatdiff} \begin{split} 
			\mathsf{P}_{i,j} & \dgTEXTARROWLENGTH=4em \arrow{e,t}{\langle \tilde{\phi}_{i,j}, \lambda_{i,j} \rangle} \Fun(\finbased^{n_j},\finbased^{n_i}) \times \Fun(\cat{C}_i,\cat{C}_j) \\
			& \dgTEXTARROWLENGTH=4em \arrow{e,t}{c} \Fun(\Fun(\finbased^{n_i},\cat{C}_i),\Fun(\finbased^{n_j},\cat{C}_j)) \end{split} \end{equation}
	where
	\[ c: \Fun(\finbased^{n_j},\finbased^{n_i}) \times \Fun(\cat{C}_i,\cat{C}_j) \to \Fun(\Fun(\finbased^{n_i},\cat{C}_i),\Fun(\finbased^{n_j},\cat{C}_j)) \]
	is adjoint to the composition map.
\end{definition}

\begin{proposition} \label{prop:Weil-RelCATdiff}
	The $\Weilinfty$-action map $T$ of Definition~\ref{def:Weil-RelzCATdiff} agrees with that of Definition~\ref{def:Tn} on the simplicial subset $\Weilinfty \times \RelzCatdiff$, and so determines an action
	\[ T: \Weilinfty \times \ReloCATdiff \to \ReloCATdiff \]
	of the scaled simplicial monoid $\Weilinfty$ on the scaled simplicial set $\ReloCATdiff$.
\end{proposition}
\begin{proof}
	A similar argument to that of Lemma~\ref{lem:Tn}, but with a simpler argument for part (2) since the relevant edge is no longer a zigzag, implies that Definition~\ref{def:Weil-RelzCATdiff} defines a scaled morphism
	\[ T: \Weilinfty \times \RelzCATdiff \to \RelzCATdiff. \]
	A similar argument to that of Proposition~\ref{prop:T} implies that $T$ is an action of $\Weilinfty$ on $\RelzCATdiff$.
	
	The compatibility of Definitions~\ref{def:Tn} and~\ref{def:Weil-RelzCATdiff} with respect to the inclusions $\Fun(\cat{C}_i,\cat{C}_j) \subseteq \Exinf \Fun(\cat{C}_i,\cat{C}_j)$ implies that the $\Weilinfty$-actions on $\RelzCATdiff$ and $\RelCatdiff$ agree on their common simplicial subset $\RelzCatdiff$. It follows that these actions determine a $\Weilinfty$-action on the pushout $\ReloCATdiff$ as claimed.
\end{proof}

Having established an action of $\Weilinfty$ on $\ReloCATdiff$, we now take a fibrant replacement to obtain a corresponding action on an $\infty$-bicategory $\RelCATdiff$. For that purpose, we introduce the following model structure on the category of scaled simplicial sets with a strict action of $\Weilinfty$.

\begin{proposition} \label{prop:scaled-Weil-module}
	Let $\ModscWeil$ be the category of scaled $\Weilinfty$-modules, i.e.\ the category of modules over $\Weilinfty$, viewed as a monoid in $\scset$ with its maximal scaling where every $2$-simplex is thin. Then $\ModscWeil$ has a model structure in which a morphism is a weak equivalence (or fibration) if and only if the underlying map of scaled simplicial sets is a weak equivalence (or fibration).
\end{proposition}
\begin{proof}
	We apply Schwede-Shipley's result~\cite[4.1]{schwede/shipley:2000} to the scaled simplicial monoid $\Weilinfty$. By \cite[4.2]{schwede/shipley:2000} we must verify that the scaled model structure on $\scset$ is monoidal with respect to the cartesian product, which is done in~\cite[2.1.21]{devalapurkar:2016}.
\end{proof}

\begin{definition} \label{def:RELCATdiff}
	Let $\RelCATdiff$ be the scaled $\Weilinfty$-module given by a fibrant replacement, in the model structure of Proposition~\ref{prop:scaled-Weil-module}, of the $\Weilinfty$-action on $\ReloCATdiff$ described in Proposition~\ref{prop:Weil-RelCATdiff}.
	
	The scaled simplicial set $\RelCATdiff$ is an $\infty$-bicategory, and we can choose the fibrant replacement so that the comparison map
	\[ \ReloCATdiff \weq \RelCATdiff \]
	is a cofibration, hence a monomorphism, of scaled simplicial sets. Altogether we have produced the following diagram of inclusions of $\infty$-bicategories, with compatible $\Weilinfty$-actions, which is also a homotopy pushout of $(\infty,2)$-categories:
	\begin{equation} \label{eq:pushout3} \begin{diagram}
			\node{\RelzCatdiff} \arrow{e,t,V}{r} \arrow{s,V} \node{\RelCatdiff} \arrow{s,V} \\
			\node{\RelzCATdiff} \arrow{e,V} \node{\RelCATdiff}
	\end{diagram} \end{equation}
\end{definition}

We now show that $\RelCATdiff$ is a model for the $(\infty,2)$-category $\CATdiff$ of differentiable $\infty$-categories described in Definition~\ref{def:CATinf}.

\begin{proposition} \label{prop:CAT-RelCAT}
	The composite map
	\[ \dgTEXTARROWLENGTH=2em M: \CATdiff \arrow{e,t}{M_0} \RelzCATdiff \arrow{e,V} \RelCATdiff \]
	is an equivalence of $\infty$-bicategories, where $M_0: \CATdiff \to \RelzCATdiff$ is the $\qCat$-enriched functor given by $\cat{C} \mapsto (\cat{C},\cat{E}_{\cat{C}})$.  
\end{proposition}
\begin{proof}
	Our strategy is to model the inclusion $\RelzCATdiff \to\nolinebreak \RelCATdiff$ in the context of marked simplicial categories by translating the homotopy pushout (\ref{eq:pushout3}) back into that world. Consider the diagram of marked simplicial categories
	\begin{equation} \label{eq:pushout2} \begin{diagram}
			\node{(\RelzCatdiff,\homeq)} \arrow{s} \arrow{e} \node{(\RelzCatdiff,\RelzCatdiff)} \arrow{s} \\
			\node{(\RelzCATdiff,\homeq)} \arrow{e} \node{(\RelzCATdiff,\RelzCatdiff)}
	\end{diagram} \end{equation}
	where a pair $(\bcat{C},\bcat{M})$ denotes the simplicial category $\bcat{C}$ with markings given by those edges in the subcategory $\bcat{M}$. We use the notation $\homeq$ to denote the natural marking of a $\qCat$-category: the subcategory consisting of all the equivalences in the mapping objects. The horizontal functors in (\ref{eq:pushout2}) are given by the identity map, and the vertical functors by the inclusion $\RelzCatdiff \subseteq \RelzCATdiff$.
	
	We claim that (\ref{eq:pushout2}) is a homotopy pushout diagram in the model structure on marked simplicial categories described by Lurie in~\cite[A.3.2]{lurie:2009}. The top horizontal map is a cofibration because it has the left lifting property with respect to acyclic fibrations, and the square is a strict pushout of marked simplicial categories. Since the model structure on marked simplicial categories is left proper by~\cite[A.3.2.4]{lurie:2009}, it follows that (\ref{eq:pushout2}) is a homotopy pushout in that model structure.
	
	Consider the top-right corner of (\ref{eq:pushout2}): there is a marked simplicial functor of maximally marked simplicial categories
	\[ (\RelzCatdiff,\RelzCatdiff) \arrow{e,t}{r} (\RelCatdiff,\RelCatdiff) = (\RelCatdiff,\homeq) \]
	given on mapping objects by the map $r$ of Definition~\ref{def:relcatdiff}. To see that this functor is an equivalence of marked simplicial categories, we note that the maximal marking functor takes an acyclic cofibration of simplicial sets (in the Quillen model structure), such as each $r_Y$, to an equivalence in the marked model structure. This fact can be checked directly from the definition of marked (cartesian) equivalence in~\cite[3.1.3.3]{lurie:2009}.
	
	We have now done enough to establish that the homotopy pushout square (\ref{eq:pushout2}) corresponds, under the Quillen equivalence of~\cite[4.2.7]{lurie:2009a}, to the homotopy pushout square (\ref{eq:pushout3}), and hence that the map $\RelzCATdiff \to \RelCATdiff$ can be modelled by the marked simplicial functor
	\[ (\RelzCATdiff,\homeq) \to (\RelzCATdiff,\RelzCatdiff) \]
	that is the identity on the underlying simplicial category. The desired proposition is therefore reduced to showing that the functor
	\[ M: (\CATdiff,\homeq) \to (\RelzCATdiff,\RelzCatdiff); \quad \cat{C} \mapsto (\cat{C},\cat{E}_{\cat{C}}) \]
	is an equivalence of marked simplicial categories. Given two differentiable $\infty$-categories $\cat{C}_0,\cat{C}_1$, we have
	\[ \Fun_{\N}(\cat{C}_0,\cat{C}_1) = \Fun_{\N}((\cat{C}_0,\cat{E}_{\cat{C}_0}),(\cat{C}_1,\cat{E}_{\cat{C}_1}) \]
	so $M$ is fully faithful. The proof that $M$ is essentially surjective on objects follows from the construction in the proof of Proposition~\ref{prop:DK-relcatdiff} with no changes.
\end{proof}

We now transfer the $\Weilinfty$-action on the $\infty$-bicategory $\RelCATdiff$ along the equivalence $M: \CATdiff \weq \RelCATdiff$ of Proposition~\ref{prop:CAT-RelCAT}. This transfer requires an $\infty$-bicategorical version of Lemma~\ref{lem:transfer}.

\begin{lemma} \label{lem:transfer2}
	Let $i: \bcat{X} \weq \bcat{Y}$ be an equivalence of $\infty$-bicategories. Then there is an equivalence of monoidal $\infty$-categories
	\[ \Fun(\bcat{X},\bcat{X})^{\sim} \homeq \Fun(\bcat{Y},\bcat{Y})^{\sim} \]
	whose underlying functor is equivalent to $i(-)i^{-1}$.
\end{lemma}
\begin{proof}
	The method of proof for Lemma~\ref{lem:transfer} applies in exactly the same way, using the fact that the construction $\Fun(-,-)^{\homeq}$ takes an equivalence of $\infty$-bicategories (in either of its variables) to an equivalence of $\infty$-categories.
\end{proof}

\begin{definition} \label{def:CATdiff-T}
	Let
	\[ T: \Weilinfty \to \Fun(\CATdiff,\CATdiff)^{\sim} \]
	be the strict monoidal functor obtained by composing the action map
	\[ \Weilinfty \to \Fun(\RelCATdiff,\RelCATdiff)^{\sim} \]
	associated to Definition~\ref{def:RELCATdiff} with the equivalence of monoidal $\infty$-categories induced, via Lemma~\ref{lem:transfer2}, by the equivalence $\CATdiff \weq \RelCATdiff$ of Proposition~\ref{prop:CAT-RelCAT}.
\end{definition}

\begin{theorem} \label{thm:CATdiff-T}
	The map $T$ of Definition~\ref{def:CATdiff-T} is a tangent structure on the $\infty$-bicategory $\CATdiff$ which (up to equivalence) extends that of Theorem~\ref{thm:T} on the $\infty$-category $\Catdiff$.	
\end{theorem}
\begin{proof}
	To show that $T$ is a tangent structure, we apply Proposition~\ref{prop:2-pullback}. The only thing remaining to show is that the foundational and vertical lift pullbacks in $\Weilinfty$ determine homotopy $2$-pullbacks in $\CATdiff$ for each object $\cat{C} \in \CATdiff$. In the proof of Theorem~\ref{thm:T}, we showed each such square is a pullback along a fibration in the $\infty$-cosmos $\mathsf{CAT}_{\infty}^{\N}$, which immediately implies that claim.
\end{proof}

\subsection*{Goodwillie calculus in an $\infty$-bicategory}

In this final section of the paper, we show how to use Theorem~\ref{thm:jet} to define a notion of $P_n$-equivalence, and hence Taylor tower, in an arbitrary tangent $\infty$-bicategory $\bcat{X}$ which, when $\bcat{X} = \CATdiff$, recovers Goodwillie's theory.

\begin{definition} \label{def:higher-tangent-spaces}
	Let $(\bcat{X},T)$ be a tangent $\infty$-bicategory that admits a terminal object $*$, and let $x: * \to \cat{C}$ be a $1$-morphism in $\bcat{X}$, i.e.\ a generalized object in $\cat{C}$. We say that $\cat{C}$ \emph{admits higher tangent spaces} at $x$ if, for each $n$, there is a homotopy $2$-pullback in $\bcat{X}$ of the form
	\[ \begin{diagram}
		\node{T^n_x\cat{C}} \arrow{e,t}{\iota^n_x} \arrow{s} \node{T^n(\cat{C})} \arrow{s,r}{p^n} \\
		\node{*} \arrow{e,t}{x} \node{\cat{C}}
	\end{diagram} \]
	Here $p^n$ denotes the natural transformation associated to the labelled Weil-algebra morphism given by the augmentation $W^{\otimes n} \to \N$.
\end{definition}

\begin{example}
	When $\bcat{X} = \CATdiff$, a morphism $x: * \to \cat{C}$ as in Definition~\ref{def:higher-tangent-spaces} is an actual object of a differentiable $\infty$-category $\cat{C}$ which admits the higher tangent spaces $T^n_x\cat{C}$ as described in Theorem~\ref{thm:jet}.
\end{example}

\begin{definition} \label{def:Pn-equiv}
	Let $(\bcat{X},T)$ be a tangent $\infty$-bicategory, and suppose the object $\cat{C}$ in $\bcat{X}$ admits higher tangent spaces at $x$. For any $\cat{D} \in \bcat{X}$ and $n \geq 0$, we define the \emph{subcategory of $P^x_n$-equivalences}
	\[ \cat{P}^x_n\Hom_{\bcat{X}}(\cat{C},\cat{D}) \subseteq \Hom_{\bcat{X}}(\cat{C},\cat{D}) \]
	to be the subcategory of morphisms that map to equivalences under the functor
	\[ T^n\iota^n_x: \Hom_{\bcat{X}}(\cat{C},\cat{D}) \to \Hom_{\bcat{X}}(T^n_x\cat{C},T^n\cat{D}); \quad F \mapsto T^n(F)\iota^n_x. \]
\end{definition}

\begin{lemma} \label{lem:Pn-equiv}
	Let $(\bcat{X},T)$ and $x: * \to \cat{C}$ be as in Definition~\ref{def:Pn-equiv}. Then for all $n \geq 1$:
	\[ \mathsf{P}_{n}^x\Hom_{\bcat{X}}(\cat{C},\cat{D}) \subseteq \mathsf{P}_{n-1}^x\Hom_{\bcat{X}}(\cat{C},\cat{D}). \]
\end{lemma}
\begin{proof}
	It is sufficient to show that, up to natural equivalence, we have a factorization
	\begin{equation} \label{eq:Pn-equiv} \begin{diagram}
			\node{\Hom_{\bcat{X}}(\cat{C},\cat{D})} \arrow{e,t}{T^n\iota^n_x} \arrow{se,b}{T^{n-1}\iota^{n-1}_x} \node{\Hom_{\bcat{X}}(T^n_x\cat{C},T^n\cat{D})} \arrow{s} \\
			\node[2]{\Hom_{\bcat{X}}(T^{n-1}_x\cat{C},T^{n-1}\cat{D})}
	\end{diagram} \end{equation}
	where the vertical map is given by postcomposition with $p_{T^{n-1}\cat{D}}: T^n\cat{D} \to T^{n-1}\cat{D}$ and precomposition with the map $0_x^n: T^{n-1}_x\cat{C} \to T^n_x\cat{C}$ induced by the bottom homotopy $2$-pullback square in the following diagram
	\[ \begin{diagram}
		\node{T^{n-1}_x\cat{C}} \arrow{se,t,..}{0_x^n} \arrow{sse} \arrow{e} \node{T^{n-1}\cat{C}} \arrow{se,t}{0_{T^{n-1}\cat{C}}} \\
		\node[2]{T^n_x\cat{C}} \arrow{s} \arrow{e} \node{T^n\cat{C}} \arrow{s,r}{p^n_{\cat{C}}} \\
		\node[2]{*} \arrow{e,t}{x} \node{\cat{C}}
	\end{diagram} \]
	Note that $p^{n-1}_{\cat{C}} \homeq p^n_{\cat{C}}0_{T^{n-1}\cat{C}}$ by uniqueness of the augmentation map $W^{\otimes(n-1)} \to \N$. Finally, the diagram (\ref{eq:Pn-equiv}) commutes since
	\[ p_{T^{n-1}\cat{D}} T^nF 0_{T^{n-1}\cat{C}} \homeq T^{n-1}F \]
	by that same calculation, together with the naturality of $T$.
\end{proof}

\begin{example}
	Let $\bcat{X} = \CATdiff$, and let $x$ be an object in a differentiable $\infty$-category $\cat{C}$ which admits finite colimits. Then Theorem~\ref{thm:jet} tells us that
	\[ \cat{P}_n^x\Hom_{\CATdiff}(\cat{C},\cat{D}) \subseteq \Fun_{\N}(\cat{C},\cat{D}) \]
	is precisely the subcategory of $P^x_n$-equivalences between functors $\cat{C} \to \cat{D}$ which preserve sequential colimits.
	
	We recover Goodwillie's notion of $n$-excisive functor $\cat{C} \to \cat{D}$, and hence the Taylor tower, by observing that in $\CATdiff$ the subcategory of $P^x_n$-equivalences is associated with a left exact localization of $\Fun_{\N}(\cat{C},\cat{D})$. The local objects for that localization are the $n$-excisive functors. We generalize this observation to give a definition of Taylor tower in an arbitrary tangent $\infty$-bicategory.
\end{example}

\begin{definition} \label{def:Taylor-bicat}
	Let $\bcat{X}$ be a tangent $\infty$-bicategory, and suppose that an object $\cat{C}$ in $\bcat{X}$ admits higher tangent spaces at a generalized object $x: * \to \cat{C}$. We say that $\bcat{X}$ \emph{admits Taylor towers expanded at $x$} if, for each $\cat{D} \in \bcat{X}$ and each $n \geq 0$, there is a full subcategory
	\[ \Jet^n_{x}(\cat{C},\cat{D}) \subseteq \Hom_{\bcat{X}}(\cat{C},\cat{D}), \]
	for which the inclusion admits a left adjoint $P^x_n$ such that the subcategory of $P^x_n$-equivalences, in the sense of~\ref{def:Pn-equiv}, is equal to the subcategory of morphisms that are mapped to equivalences by $P^x_n$. We can refer to $\Jet^n_x(\cat{C},\cat{D})$ as \emph{the $\infty$-category of $n$-jets at $x$} for morphisms $\cat{C} \to \cat{D}$ in $\bcat{X}$.
	
	In that case, by Lemma~\ref{lem:Pn-equiv}, we necessarily have
	\[ \Jet^{n-1}_x(\cat{C},\cat{D}) \subseteq \Jet^n_x(\cat{C},\cat{D}) \]
	and for each $1$-morphism $F: \cat{C} \to \cat{D}$ in $\bcat{X}$, there is a sequence of morphisms in $\Hom_{\bcat{X}}(\cat{C},\cat{D})$ of the form
	\[ F \to \dots \to P^x_nF \to P^x_{n-1}F \to \dots \to P^x_0F \]
	which we can call the \emph{Taylor tower} of $F$ at $x$. Taking $\bcat{X}$ to be the Goodwillie tangent structure on the $\infty$-bicategory $\CATdiff$, we recover Goodwillie's notion of Taylor tower for functors between differentiable $\infty$-categories.
	
\end{definition}

\chapter*{Proposals for Future Work} \label{sec:future}

The work in this paper is intended to open up various avenues for further research, and we conclude by describing ideas for projects that build on the concepts developed here. Some of those ideas were mentioned in the introduction (under \hyperref[intro:con]{Connections and Conjectures}), and we now expand a little on some of those suggestions.

\subsection*{Vector bundles in Goodwillie calculus}

A central role in differential geometry is, of course, played by vector bundles. The corresponding notion in an abstract tangent category is a \emph{differential bundle} introduced by Cockett and Cruttwell and explored in detail in~\cite{cockett/cruttwell:2018}. MacAdam~\cite{macadam:2021} has proved that this abstract definition recovers the standard notion of vector bundle in the tangent category of smooth manifolds.

We proved in Chapter~\ref{sec:catdiff-diff} that the analogues of vector \emph{spaces} in Goodwillie calculus are the stable $\infty$-categories. It is clear, therefore, that a differential bundle in the tangent $\infty$-category $\Catdiff$ should consist of a functor
\[ q: \cat{E} \to \cat{M} \]
for which the fibre $\cat{E}_X$ over any object $X \in \cat{M}$ is a stable $\infty$-category. However, it is less clear what additional conditions the functor $q$ should be required to satisfy.

One challenge here is that Cockett and Cruttwell's original definition of differential bundle~\cite[Def. 2.3]{cockett/cruttwell:2018} does not easily extend to the $\infty$-categorical case. In order to describe the differential bundles in $\Catdiff$ we first need a characterization in terms of Weil-algebras akin to the description of differential objects given in Proposition~\ref{prop:diff}. 

Some guidance to finding that characterization might be given by the approach in~\cite[Sec. 3]{cockett/cruttwell:2018}, in which differential bundles can in some cases be identified with differential objects in a \emph{slice} tangent $\infty$-category, or in the work of MacAdam~\cite[2.2]{macadam:2021}, in which there is an alternative presentation of differential bundles which could be more amenable to translation into the $\infty$-categorical setting.

\subsection*{Connections, curvature, etc.}

In the context of abstract tangent categories, notions of connection were introduced by Cockett and Cruttwell in~\cite{cockett/cruttwell:2017} and further developed by Lucyshyn-Wright in~\cite{lucyshyn-wright:2018}. Roughly speaking, a \emph{connection} on an object $M$ in a tangent category $\mathbb{X}$ consists of a morphism
\[ K: T^2M \to TM \]
such that the triple
\[ \langle Tp, K, pT \rangle : T^2M \to TM \times_M TM \times_M TM \]
is an isomorphism (additional linearity conditions are also required). This decomposition of the double tangent bundle allows for constructions such as parallel transport~\cite[5.20]{cockett/cruttwell:2017} to be realized in an abstract tangent category.

In the tangent $\infty$-category $\Catdiff$, there is a natural candidate for such a connection on any differentiable $\infty$-category $\cat{C}$; the functor
\[ K: T^2\cat{C} \to T\cat{C} \]
is given by a form of multilinearization. That construction appears to be a connection under some limited circumstances, for example, if $\cat{C}$ is a stable $\infty$-category, but not in general. Nonetheless, $K$ does appear to be a \emph{vertical} connection in the sense of~\cite[3.2]{cockett/cruttwell:2017}, and any differentiable $\infty$-category $\cat{C}$ seems to also admit a \emph{horizontal connection} $H$, though typically $H$ and $K$ are not exactly compatible in the sense described in~\cite[5.2]{cockett/cruttwell:2017}. 

Curvature is defined in~\cite[3.17]{cockett/cruttwell:2017} for an object in a tangent category equipped with a (vertical) connection. With that definition, the vertical connection $K$ described above always appears to be `flat' (i.e.\ has zero curvature), so perhaps this version of curvature is not the right concept to focus on in the context of functor calculus.

As mentioned in the introduction, there are many other concepts from differential geometry that have been translated into the abstract setting, including affine spaces, differential forms, and Lie algebroids. We do not have specific ideas about how these concepts might manifest in the Goodwillie tangent structure, but each would be worthy of study in order to understand how they might be related to ideas from homotopy theory.

\subsection*{Other brands of functor calculus}

This paper is concerned with what is sometimes referred to as Goodwillie's `homotopy' calculus, but there are other versions of functor calculus we could consider.

Goodwillie and Weiss~\cite{weiss:1999, goodwillie/weiss:1999} developed a `manifold' calculus and used it to investigate spaces of embeddings. That theory focuses on presheaves (with values in some $\infty$-category) on a fixed smooth manifold $M$, and describes how global sections can be recovered from higher-order local information.\footnote{The premise of this paper is that Goodwillie's homotopy calculus admits a close analogy with the theory of smooth manifolds. The manifold calculus, on the other hand, actually concerns manifolds themselves, not via analogy. The imaginative reader might therefore consider how to generalize the manifold calculus to an abstract tangent $\infty$-category, and hence develop a `Goodwillie calculus' calculus, in which ideas from Goodwillie and Weiss are used to study presheaves on a site associated to a fixed differentiable $\infty$-category.} As in the homotopy calculus, there are Taylor towers whose terms play the role of `polynomial' approximations, and under suitable circumstances, that tower `converges' to the global sections of the presheaf of interest.

We are curious if there is a tangent $\infty$-category that encodes the Goodwillie-Weiss manifold calculus in the same way the Goodwillie tangent structure of this paper describes homotopy calculus. For example, the tangent bundle on an $\infty$-category $\cat{C}$ might be given by the $\infty$-category of presheaves on $M$, with values in $\cat{C}$, that are \emph{degree $\leq 1$} in the sense described in~\cite[2.2]{weiss:1999}. Preliminary calculations suggest that this construction of a tangent bundle does \emph{not} satisfy the full vertical lift axiom (see~\ref{prop:vertical-lift}) for a tangent category, though perhaps a weaker notion could still apply. A similar approach could be taken with the `orthogonal' calculus of Weiss~\cite{weiss:1995}, in which the objects of interest are functors on a certain $\infty$-category of real inner product spaces, or the related `unitary' calculus of Taggart~\cite{taggart:2022}. 

The homotopy calculus itself also has an equivariant version developed by Dotto~\cite{dotto:2016,dotto:2017} based on ideas of Blumberg~\cite{blumberg:2006}. In this context, excisive functors are classified by genuine equivariant spectra. Perhaps there is a suitable tangent $\infty$-category based on these constructions which bears a closer connection to modern equivariant homotopy theory than the version developed in this paper.

\subsection*{Goodwillie tangent structure on an $\infty$-cosmos}

In a series of papers starting with~\cite{riehl/verity:2017}, Riehl and Verity have developed the notion of an \emph{$\infty$-cosmos} which captures some of the features of the collection of $\infty$-categories and is intended as a model-independent context for $\infty$-category theory. We used that work in the proof of Theorem~\ref{thm:T} by considering the $\infty$-cosmos of $\infty$-categories that admit sequential colimits.

It seems reasonable to expect that the existence of the Goodwillie tangent structure could be extended to a wider collection of $\infty$-cosmoses. In particular, any $\infty$-cosmos $\mathbb{K}$ has cotensors by quasicategories, so that it makes sense to consider the object $\Fun(\finbased,\cat{C})$ for an object $\cat{C}$ in $\mathbb{K}$, and appropriate limits so that an analogue of the tangent bundle $T\cat{C} = \Exc(\finbased,\cat{C})$ can also be constructed. One might hope to identify the conditions on the $\infty$-cosmos $\mathbb{K}$ in which that tangent bundle underlies a tangent structure in the sense of this paper.

\subsection*{Operads and tangent $(\infty,2)$-categories}

In many of the proposals above, we have suggested how to apply the theory of abstract tangent categories to Goodwillie calculus and homotopy theory, but we can also look for applications in the reverse direction, by seeing if concepts from Goodwillie calculus can be generalized to the abstract setting.

Greg Arone and the third author examined in~\cite{arone/ching:2011} the role of operads in Goodwillie calculus, formulating a chain rule for derivatives in that context. It is reasonable to ask whether these operad structures reflect a more general phenomenon. We have seen that the natural setting for Goodwillie calculus is a tangent $(\infty,2)$-category, and that might be the appropriate place to identify such a generalization. Clues to the right approach here may be in work of Lemay~\cite{lemay:2018}, which relates the Fa\`{a} di Bruno formula (the chain rule for higher-order derivatives in ordinary calculus) to tangent categories, building on work of Cockett and Seely~\cite{cockett/seely:2011}.

\subsection*{Higher approximations to $\infty$-categories}

Heuts~\cite[1.7]{heuts:2021} has introduced a theory of Goodwillie towers for pointed compactly-generated $\infty$-\emph{categories}, instead of functors. That construction provides each such $\infty$-category $\cat{C}$ with a sequence of approximations $\cat{P}_n\cat{C}$. The first approximation $\cat{P}_1\cat{C}$ is the stabilization $\spectra(\cat{C})$, i.e.\ equal to the tangent space $T_*\cat{C}$ to $\cat{C}$ at the null object. It would be interesting to see if these higher approximations also admit a description in terms of the Goodwillie tangent structure, presumably somehow related to the theory of $n$-jets in Chapter~\ref{sec:jet}. Alternatively, one might find a connection between Heuts's theory and `higher' tangent structures based on Weil-algebras whose generating relations are higher degree monomials. For example, do Heuts's $n$-excisive approximations correspond to Weil-algebras of the form $\N[x]/(x^{n+1})$?

\subsection*{Cartesian differential $\infty$-categories}

In Chapter~\ref{sec:differential} we showed that a certain homotopy category of differential objects in a cartesian tangent $\infty$-category forms a cartesian differential category in the sense of Blute, Cockett and Seely~\cite{blute/cockett/seely:2009}. We are curious whether there is a sensible notion of cartesian differential $\infty$-category that refines this construction (as well as the corresponding construction of the first author and others~\cite{bauer/johnson/osborne/riehl/tebbe:2018} in the context of abelian categories) without the need to take a homotopy category. As usual, the original definition of cartesian differential category does not easily generalize to the $\infty$-category context. Based on Theorem~\ref{thm:stable-cartesian-differential-category}, we would expect the $\infty$-category $\Catst$ to be an example of a cartesian differential $\infty$-category, but we do not have a precise definition of such a structure.

\subsection*{Derived manifolds and synthetic differential geometry}

This paper has been primarily focused on the tangent structure that encodes functor calculus, but we did provide another example of a tangent $\infty$-category in Proposition~\ref{prop:dMfld}: the $\infty$-category of derived manifolds. We have not explored any features of that tangent structure in this paper, so much remains to be done. For example, what are the differential objects or differential bundles for derived manifolds, and how do they relate to ordinary vector spaces and vector bundles? We might also hope for there to be a close connection between derived manifolds and an $\infty$-categorical version of synthetic differential geometry~\cite{kock:1981}.

\subsection*{Tangent structures in spectral algebraic geometry}

We prove in Example~\ref{ex:ring-spec-op} that the opposite of the $\infty$-category of $E_\infty$-ring spectra admits a tangent structure which extends the ordinary tangent structures on the opposite of the category of commutative rings, or equivalently on the category of affine schemes. That ordinary tangent structure extends to a tangent structure on all schemes; see, for example,~\cite{cruttwell/lemay:2023}. It is therefore natural to ask whether the tangent structure on $\E\spectra^{op}$ extends to a tangent structure on any of the $\infty$-categories that appear in spectral algebraic geometry, e.g. Lurie's spectral Deligne-Mumford stacks~\cite[1.4.4.2]{lurie:2018}. The third author's work on tangent $\infty$-categories of $\infty$-toposes~\cite{ching:2021}, which are closely related to the Goodwillie tangent structure, may also have some role to play there.

\part*{Appendices}

\appendix

\renewcommand{\thesection}{\Alph{chapter}.\arabic{section}}

\chapter{$\infty$-Categories and $\infty$-Bicategories} \label{sec:infty}

In this section, we will collect some of the basic facts about $\infty$-categories and $\infty$-bicategories, which we use throughout the rest of this paper. None of this material is new, and it is largely based on Lurie's Higher Topos Theory~\cite{lurie:2009}, though the material on $\infty$-bicategories comes from the preprint~\cite{lurie:2009a} and from the Kerodon website~\cite{lurie:kerodon}. By default, therefore, we are using quasi-categories as our model for $(\infty,1)$-categories, though we expect that the general theory of tangent $\infty$-categories presented here could be written in any other model, or, for example, within the synthetic framework being developed by Cisinski, Cnossen, Nguyen, and Walde~\cite{cisinski/cnossen/nguyen/walde:2026}, without too many changes, though we do not attempt to prove that.  Our construction of the Goodwillie tangent structure in Part 2 of the paper relies more explicitly on the quasi-category model, though again we expect a more conceptual presentation could be made.

\section{$\infty$-Categories} \label{sec:infty-cat}

The intuitive notion of $(\infty,1)$-category is intended to describe a higher category in which all $k$-morphisms are invertible for $k \geq 2$. Many different concrete versions of this notion have been developed, and all are equivalent and interchangeable in a suitable sense. For definiteness, we primarily use the following model, introduced originally by Boardman and Vogt~\cite{boardman/vogt:1973} under the name of \emph{weak Kan complexes}, studied in detail by Joyal as \emph{quasi-categories}, and popularized by Lurie under the term \emph{$\infty$-categories}.

\begin{definition}
An \emph{$\infty$-category}\index{cat000@$\infty$-category} is a simplicial set $\bcat{X}$ for which every inner horn has a filler, i.e.\ every map of simplicial sets
\[ \Lambda^n_i \to \bcat{X} \]
for $0 < i < n$, extends to an $n$-simplex $\Delta^n \to \bcat{X}$. A \emph{functor}\index{functor!between $\infty$-categories} $F: \bcat{X} \to \bcat{Y}$ between $\infty$-categories is any map of simplicial sets. We denote by $\Fun(\bcat{X},\bcat{Y})$\index{Fun@$\Fun(\bcat{X},\bcat{Y})$, the $\infty$-category of functors between two $\infty$-categories} the usual simplicial mapping object between two $\infty$-categories. Then $\Fun(\bcat{X},\bcat{Y})$ is an $\infty$-category whose objects are the functors from $\bcat{X}$ to $\bcat{Y}$. The morphisms in $\Fun(\bcat{X},\bcat{Y})$ are maps of simplicial sets $\Delta^1 \times \bcat{X} \to \bcat{Y}$ and we refer to these as \emph{natural transformations}\index{natural transformation!between functors between $\infty$-categories}.
\end{definition}

\begin{example}
Let $\mathsf{C}$ be a category. Then the \emph{nerve}\index{nerve!of an ordinary category} of $\mathsf{C}$ is the simplicial set $N\mathsf{C}$\index{NC@$N\mathsf{C}$!nerve of a category $\mathsf{C}$} in which an $n$-simplex is a sequence
\[ \dgTEXTARROWLENGTH=2.5em C_0 \arrow{e,t}{f_1} C_1 \arrow{e,t}{f_2} \cdots \arrow{e,t}{f_n} C_n \]
of $n$ composable morphisms in $\mathsf{C}$. Face maps are given by composing adjacent morphisms (or deleting those on the end), and degeneracy maps insert an identity morphism from $\mathsf{C}$ at some point in the sequence.

The simplicial set $N\mathsf{C}$ is an $\infty$-category with the additional property that every inner horn has a \emph{unique} filler. Conversely, a simplicial set with unique horn fillers is isomorphic to the nerve of a category, which is unique up to unique isomorphism. We take the view that a category \emph{is} an $\infty$-category with this extra property, and we usually suppress notation for the nerve construction, writing $\mathsf{C}$ also for the corresponding $\infty$-category.

A functor $F: \mathsf{C} \to \mathsf{D}$ between ordinary categories determines a map of simplicial sets between their nerves, i.e.\ a functor between the corresponding $\infty$-categories. Conversely, any map of simplicial sets between the nerves of ordinary categories is the nerve of an ordinary functor. Thus, there is no ambiguity between the notions of functor between ordinary categories and between the corresponding $\infty$-categories.
\end{example}

\begin{remark}
The nerve functor from categories to simplicial sets has a left adjoint, which assigns to each simplicial set its \emph{homotopy category}\index{homotopy category of an $\infty$-category}. For an $\infty$-category $\cat{C}$, the homotopy category $h\cat{C}$ can (see~\cite[1.2.3]{lurie:2009}) be described as follows: $h\cat{C}$ has the same objects as $\cat{C}$, and the morphisms from $X$ to $Y$ in $h\cat{C}$ are equivalence classes of morphisms in $\cat{C}$, where $f,g: X \to Y$ are equivalent if there is a $2$-simplex in $\cat{C}$ of the form
\[ \begin{diagram}
	\node{X} \arrow[2]{e,t}{1_X} \arrow{se,b}{f} \node[2]{X} \arrow{sw,b}{g} \\
	\node[2]{Y}
\end{diagram} \]
\end{remark}

\begin{example} \label{ex:Kan}
A \emph{Kan complex}\index{Kan complex}, or \emph{$\infty$-groupoid}\index{groupoid@$\infty$-groupoid}, is a simplicial set $X$ in which every horn $\Lambda^n_i \to X$ has a filler for $0 \leq i \leq n$, that is, both inner and outer horns. Every Kan complex is an $\infty$-category. For an $\infty$-category $\cat{C}$, the \emph{core}\index{core $\infty$-groupoid!of an $\infty$-category} of $\cat{C}$ is the Kan complex $\cat{C}^{\simeq}$ consisting of those $n$-simplexes in $\cat{C}$ for which every edge is an equivalence.
\end{example}

\subsection*{Nerves of simplicial categories}

Another useful model for the notion of $(\infty,1)$-category is given by a \emph{simplicial category}, i.e.\ a category enriched in simplicial sets. That model is related to $\infty$-categories by the `simplicial nerve' construction, which we now recall.

\begin{definition} \label{def:simp-nerve}
Let $\mathsf{C}$ be a simplicial category. The \emph{simplicial nerve}\index{nerve!of a simplicial category} of $\mathsf{C}$ is the simplicial set $N\mathsf{C}$\index{NC@$N\mathsf{C}$, nerve of a simplicial category $\mathsf{C}$}\index{simplicial nerve} in which an $n$-simplex consists of:
\begin{itemize}
	\item a sequence $C_0,\dots,C_n$ of $n+1$ objects in $\mathsf{C}$;
	\item for each $0 \leq i < j \leq n$, a map of simplicial sets
	\[ f_{i,j}: \mathsf{P}_{i,j} \to \Hom_{\mathsf{C}}(C_i,C_j) \]
	where $\mathsf{P}_{i,j}$\index{Pij@$\mathsf{P}_{i,j}$, poset of subsets of $[i,j]$} denotes (the nerve of) the poset of subsets of $[n]:= \{0,\dots,n\}$ for which $i$ is the minimum element and $j$ is the maximum element, ordered by inclusion;
\end{itemize}
such that
\begin{itemize}
	\item for each $0 \leq i \leq n$, the map
	\[ \mathsf{P}_{i,i} \to \Hom_{\mathsf{C}}(C_i,C_i) \]
	sends the unique object of $\mathsf{P}_{i,i}$ to the identity element in $\Hom_{\mathsf{C}}(C_i,C_i)$;
	\item for each $0 \leq i \leq j \leq k \leq n$, the following diagram commutes
	\[ \begin{diagram}
		\node{\mathsf{P}_{i,j} \times \mathsf{P}_{j,k}} \arrow{s,l}{\cup} \arrow{e,t}{f_{i,j} \times f_{j,k}} \node{\Hom_{\mathsf{C}}(C_i,C_j) \times \Hom_{\mathsf{C}}(C_j,C_k)} \arrow{s,r}{\circ} \\
		\node{\mathsf{P}_{i,k}} \arrow{e,t}{f_{i,k}} \node{\Hom_{\mathsf{C}}(C_i,C_k).}
	\end{diagram} \]
\end{itemize}
\end{definition}

\begin{remark} \label{rem:simp-nerve}
Unpacking Definition~\ref{def:simp-nerve}, we obtain the following descriptions of low-dimensional simplexes in $N\mathsf{C}$:
\begin{enumerate} \setcounter{enumi}{-1}
	\item a $0$-simplex is an object $C_0$ of $\mathsf{C}$;
	\item a $1$-simplex is a morphism $f_{0,1}: C_0 \to C_1$, i.e.\ a vertex in the simplicial set $\Hom_{\mathsf{C}}(C_0,C_1)$;
	\item a $2$-simplex consists of three morphisms in $\mathsf{C}$ of the form
	\[ \begin{diagram}
		\node{C_0} \arrow{se,b}{f_{0,1}} \arrow[2]{e,t}{f_{0,2}} \node[2]{C_2} \\
		\node[2]{C_1} \arrow{ne,b}{f_{1,2}}
	\end{diagram} \]
	together with a $1$-simplex
	\[ \alpha_{0,1,2}: f_{0,2} \to f_{1,2} f_{0,1} \]
	in the simplicial set $\Hom_{\mathsf{C}}(C_0,C_2)$;
	\item a $3$-simplex consists of morphisms $f_{i,j}: C_i \to C_j$ for $0 \leq i < j \leq 3$, $1$-simplexes $\alpha_{i,j,k}: f_{i,k} \to f_{j,k}f_{i,j}$ for $0 \leq i < j < k \leq 3$, and a diagram in the simplicial set $\Hom_{\mathsf{C}}(C_0,C_3)$ of the following form:
	\[ \begin{diagram}
		\node{f_{0,3}} \arrow{e,t}{\alpha_{0,1,3}} \arrow{s,l}{\alpha_{0,2,3}} \node{f_{1,3}f_{0,1}} \arrow{s,r}{\alpha_{1,2,3}f_{0,1}} \\
		\node{f_{2,3}f_{0,2}} \arrow{e,t}{f_{2,3}\alpha_{0,1,2}} \node{f_{2,3}f_{1,2}f_{0,1}}
	\end{diagram} \]
	\item more generally, an $n$-simplex consists of an $n-1$-dimensional cubical diagram in the simplicial set $\Hom_{\mathsf{C}}(C_0,C_n)$ whose faces are built from lower-dimensional cubical diagrams in the simplicial sets $\Hom_{\mathsf{C}}(C_i,C_j)$ for $0 \leq i < j \leq n$.
\end{enumerate}
\end{remark}

\begin{proposition} \label{prop:simp-nerve}
Let $\mathsf{C}$ be a simplicial category for which the mapping objects $\Hom_{\mathsf{C}}(C,C')$ are Kan complexes. Then $N\mathsf{C}$ is an $\infty$-category.
\end{proposition}

\begin{example} \label{ex:ordinary-nerve}
Let $\mathsf{C}$ be an ordinary category. Then we can view $\mathsf{C}$ as a simplicial category in which each mapping object is a discrete simplicial set. Then the simplicial nerve of that simplicial category is isomorphic to the ordinary nerve of $\mathsf{C}$.
\end{example}

Several important examples of $\infty$-categories arise as the simplicial nerve of a suitable simplicial category.

\begin{example} \label{ex:spaces}
Let $\mathsf{S}$ denote the (large) simplicial category in which
\begin{itemize}
	\item an object is a (small) Kan complex;
	\item the mapping objects are given by the ordinary simplicial mapping objects
	\[ \Hom_{\mathsf{S}}(X,Y) := \Hom_{\sset}(X,Y) \]
	in which an $n$-simplex is a map of simplicial sets of the form $\Delta^n \times X \to Y$.
\end{itemize}
Then the simplicial nerve of $\mathsf{S}$ is an $\infty$-category $\spaces$\index{S@$\spaces$, the $\infty$-category of spaces}, which we refer to as the \emph{$\infty$-category of spaces}. 
\end{example}

\begin{example} \label{ex:bcatinf}
Let $\mathsf{Cat_\infty}$ denote the (very large)\footnote{Recall from the Introduction that we are implicitly using a hierarchy of Grothendieck universes to handle the different sizes of $\infty$-category involved in our theory.} simplicial category in which
\begin{itemize}
	\item an object is a (large) $\infty$-category;
	\item the mapping objects are given by
	\[ \Hom_{\mathsf{Cat_\infty}}(\bcat{X},\bcat{Y}) := \Fun(\bcat{X},\bcat{Y})^{\homeq}, \]
	i.e.\ the core of the $\infty$-category of functors, which is itself the ordinary mapping simplicial set.
\end{itemize}
Then the simplicial nerve of $\mathsf{Cat_\infty}$ is a (very large) $\infty$-category $\bCatinf$, which we refer to as the \emph{$\infty$-category of $\infty$-categories}\index{CATinf1@$\bCatinf$, the $\infty$-category of (large) $\infty$-categories}. Based on Remark~\ref{rem:simp-nerve}, we can identify the low-dimensional simplexes in $\bCatinf$ as follows:
\begin{itemize}
	\item a $0$-simplex is an $\infty$-category $\bcat{X}$;
	\item a $1$-simplex is a functor $F: \bcat{X}_0 \to \bcat{X}_1$;
	\item a $2$-simplex consists of functors
	\[ \begin{diagram}
		\node{\bcat{X}_0} \arrow[2]{e,t}{H} \arrow{se,b}{F} \node[2]{\bcat{X}_2} \\
		\node[2]{\bcat{X}_1} \arrow{ne,b}{G}
	\end{diagram} \]
	together with a natural equivalence $\alpha: H \weq GF$;
	\item an $n$-simplex consists of $\infty$-categories $\bcat{X}_0,\dots,\bcat{X}_n$ together with various cubical diagrams of functors and natural equivalences connecting them.
\end{itemize}
\end{example}

\begin{example}	 \label{ex:2,1-nerve}
Let $\mathsf{C}$ be a strict $(2,1)$-category\index{cat000@$(2,1)$-category}, i.e.\ a category enriched in groupoids. Applying the ordinary nerve construction to each mapping groupoid, we obtain a simplicial category in which each mapping space is a Kan complex. The simplicial nerve of that simplicial category is an $\infty$-category which we refer to simply as the \emph{simplicial nerve}\index{nerve!of a $(2,1)$-category} of $\mathsf{C}$.
\end{example}

\begin{example} \label{ex:cat}
Let $\mathsf{Cat}$ denote the strict $(2,1)$-category of (ordinary, large) categories, functors, and natural isomorphisms. Let $\mathbf{Cat}$ be the simplicial nerve of $\mathsf{Cat}$ in the sense of Example~\ref{ex:2,1-nerve}. We refer to $\mathbf{CAT}$\index{CAT1@$\mathbf{Cat}$, the $\infty$-category of ordinary categories} as the \emph{$\infty$-category of categories}. The nerve construction determines a fully faithful functor
\[ N: \mathbf{Cat} \to \bCatinf. \]
That functor has a right adjoint (in the $\infty$-categorical sense) given by the homotopy category construction
\[ h: \bCatinf \to \mathbf{Cat}. \]
\end{example}

\subsection*{The nerve of a bicategory in which every $2$-morphism is invertible}

Duskin described in~\cite{duskin:2001}, using a construction of Street, a nerve for bicategories, and proved that if every $2$-morphism is invertible in the bicategory, then the resulting nerve is an $\infty$-category. We refer to~\cite{johnson/yau:2021} for background on ordinary bicategories and their pseudofunctors.

\begin{definition}[Duskin nerve] \label{def:duskin-nerve}
Let $\mathsf{B}$ be a bicategory. Then there is an $\infty$-category $N\mathsf{B}$, the \index{Duskin nerve of a bicategory}\index{nerve!of a bicategory}\emph{nerve} of $\mathsf{B}$, in which an $n$-simplex is a strictly unital pseudofunctor
\[ [n] \to \mathsf{B}, \]
and the simplicial structure is determined by precomposition with the functors $[m] \to [n]$. Unpacking this definition, we can make the following observations about $N\mathsf{B}$:
\begin{itemize}
	\item an object is an object in $\mathsf{B}$;
	\item a morphism is a morphism in $\mathsf{B}$;
	\item a $2$-simplex consists of three morphisms $f_{0,1},f_{1,2},f_{0,2}$ as in Remark~\ref{rem:simp-nerve}(2), together with a $2$-morphism
	\[ \alpha_{0,1,2}: f_{0,2} \arrow{e} f_{1,2}f_{0,1}; \]
	this $2$-simplex is degenerate if $\alpha_{0,1,2}$ is an identity $2$-morphism;
	\item a $3$-simplex consists of $f_{i,j}$ and $\alpha_{i,j,k}$ as in Remark~\ref{rem:simp-nerve}(3), such that the following diagram commutes (in the relevant mapping groupoid of $\mathsf{B}$):
	\[ \begin{diagram}
		\node[2]{f_{03}} \arrow{sw,t}{\alpha_{0,1,3}} \arrow{se,t}{\alpha_{0,2,3}} \\
		\node{f_{1,3}f_{0,1}}  \arrow{s,l}{\alpha_{1,2,3}f_{0,1}} \node[2]{f_{2,3}f_{0,2}} \arrow{s,r}{f_{2,3}\alpha_{0,1,2}} \\
		\node{(f_{2,3}f_{1,2})f_{0,1}} \arrow[2]{e,t}{\isom} \node[2]{f_{2,3}(f_{1,2}f_{0,1})}
	\end{diagram} \]
	\item the simplicial set $N\mathsf{B}$ is $3$-coskeletal\index{coskeletal@$3$-coskeletal simplicial set}: every diagram indexed by the $3$-skeleton of an $n$-simplex extends uniquely to an $n$-simplex.
\end{itemize}
It follows from~\cite[8.6]{duskin:2001} that if every $2$-morphism in $\mathsf{B}$ is invertible, then $N\mathsf{B}$ is an $\infty$-category. When $\mathsf{B}$ is a strict $(2,1)$-category, $N\mathsf{B}$ agrees with the simplicial nerve of $\mathsf{B}$ described in Example~\ref{ex:2,1-nerve}.
\end{definition}

\subsection*{$\infty$-Categories associated to model categories}

Quillen~\cite{quillen:1967} introduced model categories as an abstraction of various classical ideas from homotopy theory. We will not review model categories here, referring the reader to~\cite{hovey:1999} for an extensive exposition of their theory. We do recall how an $\infty$-category can be recovered from a \emph{simplicial} model category, i.e. one in which there is a simplicial enrichment which is compatible with the model structure.

\begin{definition} \label{def:model-infty}
Let $\mathsf{M}$ be a simplicial model category. Let $\mathsf{M}^{\circ}$ be the full subcategory whose objects are the fibrant-cofibrant objects in $\mathsf{M}$. Then the simplicial nerve $N(\mathsf{M}^{\circ})$ is an $\infty$-category, which we refer to as \emph{the $\infty$-category associated to $\mathsf{M}$}\index{cat000@$\infty$-category!associated to a simplicial model category}.
\end{definition}

The $\infty$-category $\bCatinf$ described in Example~\ref{ex:bcatinf} above can itself be realized as the $\infty$-category associated to a simplicial model category: of \emph{marked} simplicial sets.

\begin{definition} \label{def:msset}
A \emph{marked simplicial set}\index{marked simplicial set} consists of a simplicial set $X$ together with a set $E$ of edges in $X$ which includes all degenerate edges. We refer to the elements of $E$ as the \emph{marked edges}\index{marked edges in a marked simplicial set}. There is a simplicial category $\msset$\index{sset+@$\msset$, the category of marked simplicial sets} whose objects are marked simplicial sets and whose mapping objects are the maximal Kan complexes
\[ \Hom_{\msset}((X,E),(X',E'))^{\homeq} \subseteq \Hom_{\sset}(X,X') \]
whose vertices take marked edges to marked edges.

Lurie constructs in~\cite[3.1.3.7]{lurie:2009} a simplicial model structure on $\msset$ with the following properties:
\begin{itemize}
	\item a map $(X,E) \to (X',E')$ is a cofibration if and only if the underlying map of simplicial sets $X \to X'$ is a monomorphism, i.e. injective in each dimension; in particular, every object is cofibrant;
	\item an object $(X,E)$ is fibrant if and only if $X$ is an $\infty$-category and $E$ is the set of \emph{equivalences} in $X$.
\end{itemize}
We refer to this model structure as the \emph{marked model structure}\index{marked model structure on $\msset$}.
\end{definition}

\begin{proposition} \label{prop:msset}
There is an isomorphism of simplicial categories
\[ \mathsf{Cat_\infty} \isom \msset^{\circ} \]
and hence $\bCatinf$ can be identified with the $\infty$-category associated to the marked model structure on $\msset$.
\end{proposition}
\begin{proof}
Both simplicial categories have as their objects the $\infty$-categories, so it is sufficient to show that there are isomorphisms of mapping objects
\[ \Fun(\bcat{X},\bcat{Y})^{\homeq} \isom \Hom_{\msset}((\bcat{X},E_{\bcat{X}}),(\bcat{Y},E_{\bcat{Y}}))^{\homeq}. \]
In fact, these two simplicial subsets of $\Fun(\bcat{X},\bcat{Y}) = \Hom_{\sset}(\bcat{X},\bcat{Y})$ are equal because any functor between $\infty$-categories preserves equivalences.
\end{proof}

\subsection*{Mapping spaces in an $\infty$-category}

An $\infty$-category $\bcat{X}$ can be viewed as a category enriched (in a weak sense) in spaces. That is, for any pair of objects in $\bcat{X}$ we have a mapping space, and these admit composition operations that are unital and associative up to coherent homotopy.

\begin{definition} \label{def:hom}
Let $\bcat{X}$ be an $\infty$-category, and let $A,B$ be objects in $\bcat{X}$. We write\index{Hom1@$\Hom_{\bcat{X}}(A,B)$, the space of maps between objects $A,B$ in an $\infty$-category $\bcat{X}$}
\[ \Hom_{\bcat{X}}(A,B), \quad \text{or} \quad \bcat{X}(A,B) \]
for the space (or $\infty$-groupoid) of maps from $A$ to $B$ in $\bcat{X}$. This object is defined precisely as in~\cite[1.2.2]{lurie:2009} (where it is denoted $\mathrm{Map}_{\bcat{X}}(A,B)$).
\end{definition}

\section{$\infty$-Bicategories} \label{sec:infty-bicat}

The theory of $(\infty,2)$-categories is less well-developed than that of $(\infty,1)$-categories, although recently several different models have been produced and proved to be equivalent; see~\cite{moser/ozornova/rovelli:2022}. We will maintain as close a connection as possible between our approaches to $(\infty,2)$- and $(\infty,1)$-categories by relying again on work of Lurie, in which an $(\infty,2)$-category is modelled as a type of simplicial set. The drawback of this approach is that the conditions on that simplicial set are rather more involved than in the case of $\infty$-categories, though in this paper we have no need to work with those conditions explicitly.

\begin{definition} \label{def:inf-bicat}
An \index{bicategory@$\infty$-bicategory}\emph{$\infty$-bicategory} is a simplicial set $\bcat{B}$ which satisfies the conditions of~\cite[\href{https://kerodon.net/tag/01W9}{01W9}]{lurie:kerodon}. We will not write out those conditions in detail here. We note that certain $2$-simplexes in an $\infty$-bicategory play a special role in the theory: those which are \emph{thin}\index{thin $2$-simplex} in the sense of~\cite[\href{https://kerodon.net/tag/00AD}{00AD}]{lurie:kerodon}. Intuitively, a $2$-simplex in an $\infty$-bicategory $\bcat{B}$ can be viewed as a diagram of the form
\[ \begin{diagram}
	\node{B_0} \arrow[4]{e,t}{H} \arrow[2]{se,b}{F} \node[2]{} \node[2]{B_2} \\
	\node[3]{\Downarrow} \\
	\node[3]{B_1} \arrow[2]{ne,b}{G}
\end{diagram} \]
where $H \Rightarrow GF$ is a $2$-morphism. That $2$-simplex is thin if the corresponding $2$-morphism is invertible.
\end{definition}

\begin{example}
A simplicial set is an $\infty$-category if and only if it is an $\infty$-bicategory in which all $2$-simplexes are thin. 
\end{example}

\begin{definition} \label{def:fun-bicat}
A \emph{functor}\index{functor!between $\infty$-bicategories} $F: \bcat{B} \to \bcat{B}'$ between $\infty$-bicategories is a map of simplicial sets which takes thin $2$-simplexes to thin $2$-simplexes, and we write
\[ \Fun(\bcat{B},\bcat{B}') \subseteq \Hom_{\sset}(\bcat{B},\bcat{B}') \]
for the simplicial subset whose $n$-simplexes are those maps of simplicial sets
\[ \alpha: \Delta^n \times \bcat{B} \to \bcat{B}' \]
such that if $\sigma$ is a degenerate $2$-simplex in $\Delta^n$, and $\tau$ is a thin $2$-simplex in $\bcat{B}$, then $\alpha(\sigma,\tau)$ is a thin $2$-simplex in $\bcat{B}'$.

The simplicial set $\Fun(\bcat{B},\bcat{B}')$ is an $\infty$-bicategory, and a $2$-simplex $\alpha: \Delta^2 \times \bcat{B} \to \bcat{B}'$ is thin if and only if $\alpha(\sigma,\tau)$ is thin for \emph{any} $2$-simplex $\sigma$ in $\Delta^2$ and any thin $2$-simplex $\tau$ in $\bcat{B}$.
\end{definition}

\begin{example}
If $\bcat{C}$ is an $\infty$-category, i.e.\ an $\infty$-bicategory in which all $2$-simplexes are thin, then $\Fun(\bcat{B},\bcat{C})$ is the entire simplicial mapping object, which is then an $\infty$-category.
\end{example}

\begin{definition} \label{def:core}
The \emph{core $\infty$-category}\index{core $\infty$-category of an $\infty$-bicategory} of an $\infty$-bicategory $\bcat{B}$ is the $\infty$-category $\bcat{B}^{\sim}$ consisting of those simplexes in $\bcat{B}$ whose $2$-faces are all thin. Intuitively, this corresponds to the process of forgetting all non-invertible $2$-morphisms. The core satisfies the universal property that for any $\infty$-category $\bcat{C}$ the restriction map
\[ \Fun(\bcat{B},\bcat{C}) \to \Fun(\bcat{B}^{\sim},\bcat{C}) \]
is an equivalence of $\infty$-categories.

The \emph{core $\infty$-groupoid}\index{core $\infty$-groupoid!of an $\infty$-bicategory} of $\bcat{B}$ is the core (in the sense of Example~\ref{ex:Kan}) of the $\infty$-category $\bcat{B}^{\sim}$. That is the simplicial subset $\bcat{B}^{\simeq} \subseteq \bcat{B}$ consisting of those simplexes for which all $2$-faces are thin, and all edges are equivalences in $\bcat{B}^{\sim}$.
\end{definition}

\subsection*{Nerves and examples of $\infty$-bicategories}

\begin{proposition} \label{prop:simp-nerve-bicat}
Let $\mathsf{C}$ be a simplicial category in which each mapping object $\Hom_{\mathsf{C}}(C,C')$ is an $\infty$-category. Then the simplicial nerve $N\mathsf{C}$ is an $\infty$-bicategory.
\end{proposition}

\begin{example} \label{ex:bicat-cat}
Let $\mathsf{CAT}_\infty$ be the (very large) simplicial category in which
\begin{itemize}
	\item an object is a (large) $\infty$-category;
	\item the mapping objects are the ordinary simplicial mapping objects
	\[ \Hom_{\mathsf{CAT}_\infty}(\bcat{X},\bcat{Y}) := \Fun(\bcat{X},\bcat{Y}). \]
\end{itemize}
Then the simplicial nerve $N\mathsf{CAT}_\infty$ is an $\infty$-bicategory $\mathbf{CAT}_\infty$ which we refer to as \emph{the $\infty$-bicategory of $\infty$-categories}\index{CATinf12@$\mathbf{CAT}_\infty$, the $\infty$-bicategory of $\infty$-categories}. There is an isomorphism of $\infty$-categories
\[ \mathbf{CAT}_\infty^{\sim} \isom \bCatinf \]
where $\bCatinf$ is the $\infty$-category of $\infty$-categories from Example~\ref{ex:bcatinf}.
\end{example}

\begin{example} \label{ex:CAT2}
Let $\mathsf{CAT}_{(\infty,2)}$ be the (very large) simplicial category in which
\begin{itemize}
	\item an object is a (large) $\infty$-bicategory;
	\item the mapping objects are the $\infty$-categories
	\[ \Hom_{\mathsf{CAT}_{(\infty,2)}}(\bcat{B},\bcat{B'}) := \Fun(\bcat{B},\bcat{B'})^{\sim} \]
	given by taking the core $\infty$-category of the $\infty$-bicategory of functors. 
\end{itemize}
Then the simplicial nerve $N\mathsf{CAT}_{(\infty,2)}$ is an $\infty$-bicategory $\mathbf{CAT}_{(\infty,2)}$ which we refer to as \emph{the $\infty$-bicategory of $\infty$-bicategories}\index{CATinf22@$\mathbf{CAT}_{(\infty,2)}$, the $\infty$-bicategory of $\infty$-bicategories}. Treating each $\infty$-category as an $\infty$-bicategory determines a fully faithful functor (of $\infty$-bicategories)
\[ \mathbf{CAT}_\infty \to \mathbf{CAT}_{(\infty,2)}. \]
We also write
\[\mathbf{Cat}_{(\infty,2)} := \mathbf{CAT}_{(\infty,2)}^{\sim} \]
for the core $\infty$-category of the $\infty$-bicategory of $\infty$-bicategories, and refer to this object as \emph{the $\infty$-category of $\infty$-bicategories}\index{CATinf21@$\mathbf{Cat}_{(\infty,2)}$, the $\infty$-category of $\infty$-bicategories}. We can identify $\mathbf{Cat}_{(\infty,2)}$ with the simplicial nerve of the simplicial category whose objects are $\infty$-bicategories, and whose mapping spaces are the core $\infty$-groupoids $\Fun(\bcat{B},\bcat{B}')^{\simeq}$.
\end{example}

We turn to the relationship between ordinary bicategories and $\infty$-bicategories. We use~\cite{johnson/yau:2021} as a general reference for the theory of bicategories. Recall from Definition~\ref{def:duskin-nerve} that Duskin described in~\cite{duskin:2001} a nerve construction for bicategories; see~\cite[5.4.7]{johnson/yau:2021} for details. The following result summarizes the properties of this construction.

\begin{proposition} \label{prop:bicat-nerve}
For a bicategory $\mathsf{B}$, the nerve $N(\mathsf{B})$ is an $\infty$-bicategory. The $0$- and $1$-simplexes in $N\mathsf{B}$ are the objects and $1$-morphisms of $\mathsf{B}$, respectively. A $2$-simplex in $N\mathsf{B}$ can be identified with a diagram in $\mathsf{B}$ of the form displayed in Definition~\ref{def:inf-bicat}, comprising a $2$-morphism $\alpha: H \to GF$. That $2$-simplex is thin if and only if $\alpha$ is invertible. If all $2$-morphisms in $\mathsf{B}$ are invertible, then $N\mathsf{B}$ is an $\infty$-category.

A strictly unitary pseudofunctor $F: \mathsf{B} \to \mathsf{B}'$ induces a map of simplicial sets $NF: N\mathsf{B} \to N\mathsf{B}'$, which preserves thin $2$-simplexes, i.e.\ a functor of $\infty$-bicategories. Conversely, any functor of $\infty$-bicategories $N\mathsf{B} \to N\mathsf{B}'$ arises from a strictly unitary pseudofunctor in this way. Moreover, there is an equivalence of $\infty$-bicategories
\[ \Fun(N\mathsf{B},N\mathsf{B}') \homeq N\mathsf{Bicat}(\mathsf{B},\mathsf{B}') \]
where the right-hand side is the nerve of the bicategory whose objects are strictly unitary pseudofunctors $\mathsf{B} \to \mathsf{B}'$, whose $1$-morphisms are (strong) natural transformations, and whose $2$-morphisms are modifications. 
\end{proposition}
\begin{proof}
See~\cite[\href{https://kerodon.net/tag/01WD}{01WD}]{lurie:kerodon}.
\end{proof}

\begin{definition} \label{def:CAT2}
Let $\mathbf{CAT}_2$\index{CAT2@$\mathbf{CAT}_2$, the $\infty$-bicategory of ordinary bicategories} denote the full subcategory of the $\infty$-bicategory $\mathbf{CAT}_{(\infty,2)}$ whose objects are the nerves of ordinary bicategories. The inclusion of that full subcategory admits a left $2$-adjoint
\[ h_2: \mathbf{CAT}_{(\infty,2)} \to \mathbf{CAT}_2, \]
which takes an $\infty$-bicategory $\bcat{B}$ to its \emph{homotopy bicategory}\index{homotopy bicategory of an $\infty$-bicategory} $h_2\bcat{B}$\index{h2B@$h_2\bcat{B}$, the homotopy bicategory of an $\infty$-bicategory $\bcat{B}$}.

Various versions of the homotopy bicategory have been constructed in different contexts; see, for example,~\cite{campbell:2020} and~\cite{romo:2023}. Essentially, $h_2$ is given by applying the ordinary homotopy category construction to the mapping $\infty$-categories of an $\infty$-bicategory. In particular, when $\bcat{B}$ is the simplicial nerve of a category enriched in $\infty$-categories, $h_2\bcat{B}$ can be identified with the strict $2$-category obtained by taking homotopy categories of those mapping $\infty$-categories. We will not give a more precise construction here. For us, the key fact of $h_2$ is the following universal characterization, which follows from~\cite[9.3]{campbell:2020}.
\end{definition}

\begin{proposition} \label{prop:homotopy-bicat}
For an $\infty$-bicategory $\bcat{B}$ and ordinary bicategory $\mathsf{B}'$, there are equivalences of $\infty$-bicategories
\[ \Fun(\bcat{B},N\mathsf{B}') \homeq N\mathsf{Bicat}(h_2\bcat{B},\mathsf{B}'). \]
\end{proposition}

\subsection*{Scaled simplicial sets and $\infty$-bicategories}

Recall that $\infty$-categories can be identified with the fibrant-cofibrant objects in a model structure on the category of marked simplicial sets. Similarly, we can identify $\infty$-bicategories with the fibrant-cofibrant objects in a suitable model category.

\begin{definition} \label{def:scaled}
A \index{scaled simplicial set}\emph{scaled simplicial set} consists of a simplicial set $X$ and a set $T$ of $2$-simplexes in $X$ which includes all degenerate $2$-simplexes. We refer to the elements of $T$ as the \emph{thin} $2$-simplexes in the scaled simplicial set $(X,T)$. In particular, each $\infty$-bicategory determines a scaled simplicial set in which $T$ is the set of $2$-simplexes which are thin in the sense of~\cite[\href{https://kerodon.net/tag/00AD}{00AD}]{lurie:kerodon}.

There is a simplicial category \index{ssetsc@$\scset$, the category of scaled simplicial sets}$\scset$ whose objects are the scaled simplicial sets, with mapping objects
\[ \Hom_{\scset}((X,T),(X',T')) \subseteq \Hom_{\sset}(X,X') \]
the simplicial subsets consisting of those $n$-simplexes $\Delta^n \times X \to X'$ which take elements of $T$ to elements of $T'$. 

Lurie constructs in~\cite{lurie:2009a} a simplicial model structure on $\scset$ in which
\begin{itemize}
	\item a map $(X,T) \to (X',T')$ is a cofibration if and only if the underlying map $X \to X'$ is a monomorphism of simplicial sets;
	\item an object $(X,T)$ is fibrant if and only if $X$ is an $\infty$-bicategory and $T$ is the set of thin $2$-simplexes in that $\infty$-bicategory.
\end{itemize}
We refer to this as the \index{scaled model structure on $\scset$}\emph{scaled model structure} on $\scset$. The $\infty$-bicategories are the fibrant-cofibrant objects of $\scset$, and the $\infty$-category $N(\scset^{\circ})$ associated to the scaled model structure is isomorphic to $\mathbf{CAT}_{(\infty,2)}$, the $\infty$-category of $\infty$-bicategories.
\end{definition}

\begin{remark}
The model categories of marked simplicial sets and scaled simplicial sets are closely related to Verity's more general framework of \emph{stratified simplicial sets} and various model structures whose fibrant objects, the \emph{$n$-complicial sets}, model $(\infty,n)$-categories for $0 \leq n \leq \infty$. 
\end{remark}

\subsection*{Mapping $\infty$-categories in an $\infty$-bicategory}

We can view an $\infty$-bicategory as a category (weakly) enriched in $\infty$-categories. In particular, for any two objects $X,Y$ in an $\infty$-bicategory $\bcat{B}$, there is an $\infty$-category of maps from $X$ to $Y$.

\begin{definition} \label{def:hom-bicat}
Let $X,Y$ be objects in an $\infty$-bicategory $\bcat{B}$. We write\index{Hom2@$\Hom_{\bcat{B}}(X,Y)$, the $\infty$-category of maps between objects $X,Y$ in an $\infty$-bicategory $\bcat{B}$}
\[ \Hom_{\bcat{B}}(X,Y), \quad \text{or} \quad \bcat{B}(X,Y) \]
for the $\infty$-category of maps from $X$ to $Y$ in $\bcat{B}$. A precise definition of that $\infty$-category is given in~\cite[\href{https://kerodon.net/tag/01KY}{01KY}]{lurie:kerodon}. The morphisms in this $\infty$-category provide a notion of `$2$-morphism' between $1$-morphisms in $\bcat{B}$. 
\end{definition}

The core $\infty$-groupoid of the mapping $\infty$-category $\Hom_{\bcat{B}}(X,Y)$ is equivalent to the mapping space in the core $\infty$-category of $\bcat{B}$. That is:
\[ \Hom_{\bcat{B}}(X,Y)^{\simeq} \simeq \Hom_{\bcat{B}^{\sim}}(X,Y). \]

\chapter{Monoidal and Module $\infty$-categories} \label{sec:monoidal}

The goal of this part of the appendix is to collect information about monoidal $\infty$-categories and module $\infty$-categories over them, which play a central role in the theory of tangent structures elsewhere in this paper. Most of this material is not new, and our aim is to give sufficient references for the reader to discover what they need to know from other sources.

In~\ref{sec:mon} we define an $\infty$-category $\MonCat$ of monoidal $\infty$-categories and (strong) monoidal functors and explain how both ordinary monoidal categories and the strict monoidal $\infty$-categories of Definition~\ref{def:mon-quasi} give rise to objects in $\MonCat$.

In~\ref{sec:module} we discuss module $\infty$-categories over a monoidal $\infty$-category $\bcat{W}$ and construct an $\infty$-bicategory $\Mod_{\bcat{W}}$ of $\bcat{W}$-modules, $\bcat{W}$-functors, and $\bcat{W}$-natural transformations. That definition forms the basis of the $\infty$-bicategory $\Tan_\infty$ of tangent $\infty$-categories, tangent functors, and tangent transformations; see Definition~\ref{def:tangent-functor}.

\section{Monoidal $\infty$-categories} \label{sec:mon}

As with $(\infty,1)$-categories themselves, there is no single `correct' notion of monoidal $\infty$-category. Instead, there are different models, each of which determines a different, but equivalent, $\infty$-category. We can view any of these models as an equally good incarnation of `the' $\infty$-category of monoidal $\infty$-categories, and different models have different benefits. In this section, we describe three different approaches to the notion of monoidal $\infty$-category:
\begin{itemize}
	\item monoids in the $\infty$-category of $\infty$-categories;
	\item marked simplicial monoids;
	\item $\infty$-bicategories with a single object. 
\end{itemize}
In each of these cases, we can construct an $\infty$-category $\MonCat$ whose objects are a suitable version of monoidal $\infty$-categories and whose morphisms correspond to (strong) monoidal functors. The reader who prefers a single concrete definition of monoidal $\infty$-category (and monoidal functor) can choose any one of these three approaches as giving that definition.

\subsection{Monoids in the $\infty$-category of $\infty$-categories} \label{sec:monoid-catinf}

Perhaps the most intuitive way to define monoidal $\infty$-categories is as monoids (in the $\infty$-categorical sense, which we recall below) in the $\infty$-category of $\infty$-categories. Monoids in this context are automatically associative and unital only up to coherent equivalences, and the associator and unit isomorphisms in an ordinary monoidal category are naturally part of that structure. The following definition comes from~\cite[4.1.2]{lurie:2017}.

\begin{definition} \label{def:monoid}
Let $\mathbf{C}$ be an $\infty$-category which admits finite products, and let $\bDelta$\index{simplicial indexing category@$\bDelta$, simplicial indexing category} denote the simplicial indexing category (whose objects are the finite ordered sets $[n] = \{0,\dots,n\}$ and whose morphisms are the order-preserving functions). Then a \emph{monoid}\index{monoid in an $\infty$-category} in $\mathbf{C}$ is a functor
\[ M: \bDelta^{op} \to \mathbf{C} \]
such that for each $n$ there is an equivalence
\[ M([n]) \weq M([1])^n \]
induced by the $n$ maps in $\bDelta$ of the form
\[ \rho_i: [1] \to [n]; \quad 0 \mapsto i-1, \quad 1 \mapsto i \]
for $i = 1,\dots,n$. The \emph{$\infty$-category of monoids}\index{MonC@$\mathbf{Mon}(\mathbf{C})$, the $\infty$-category of monoids in an $\infty$-category $\mathbf{C}$} in $\mathbf{C}$ is the full subcategory
\[ \mathbf{Mon}(\mathbf{C}) \subseteq \Fun(\bDelta^{op},\mathbf{C}) \]
whose objects are the monoids in the above sense. 
\end{definition}

\begin{remark}
When $\mathbf{C}$ is an ordinary category with finite products, the $\infty$-category $\mathbf{Mon}(\mathbf{C})$ is equivalent to the ordinary category of monoids in $\mathbf{C}$ in the usual sense: an object $M$ together with morphisms $M \times M \to M$ and $* \to M$ satisfying associativity and unit conditions.
\end{remark}

The $\infty$-category $\bCatinf$ of Example~\ref{ex:bcatinf} admits finite products, which are given by the ordinary cartesian product of simplicial sets. We can therefore make the following definition.

\begin{definition} \label{def:moncat1}
The $\infty$-category $\MonCat$\index{MonCat@$\MonCat$, the $\infty$-category of monoidal $\infty$-categories} of monoidal $\infty$-categories\index{monoidal $\infty$-category} and monoidal functors is
\[ \mathbf{Mon}(\bCatinf), \]
the $\infty$-category of monoids in $\bCatinf$. Thus, from this perspective, a monoidal $\infty$-category is a certain kind of functor $\bcat{M}: \bDelta^{op} \to \bCatinf$, and a monoidal functor is a certain kind of functor $\bDelta^{op} \times \Delta^1 \to \bCatinf$. 
\end{definition}

In Definition~\ref{def:mon-quasi} we introduced a notion of \emph{strict} monoidal $\infty$-category; that is, a simplicial monoid whose underlying simplicial set is an $\infty$-category. Here is how such an object determines a monoidal $\infty$-category in the sense of Definition~\ref{def:moncat1}.

\begin{definition} \label{def:strict-mon}
Let $\bcat{W}$ be a strict monoidal $\infty$-category. We define an ordinary functor
\[ W: \bDelta^{op} \to \sset \]
by
\[ W([n]) := \bcat{W}^n, \]
and for an order-preserving map $\alpha: [m] \to [n]$, we define $W(\alpha): W([n]) \to W([m])$ by
\[ (A_n,\dots,A_1) \mapsto (A_{\alpha(m)} \otimes \dots \otimes A_{\alpha(m-1)+1}, \dots, A_{\alpha(1)} \otimes \dots \otimes A_{\alpha(0)+1}). \]
Note that if $\alpha(j) = \alpha(j-1)$, then the relevant term in this expression is the unit element for the simplicial monoid $\bcat{W}^{\otimes}$. The strict associativity and unit conditions for $\bcat{W}^{\otimes}$ imply that $W$ is a functor as claimed.

We can view $W$ as a simplicially-enriched functor from $\bDelta^{op}$ (interpreted as a simplicial category with discrete mapping objects) to the simplicial category $\mathsf{Cat_\infty}$ whose nerve is $\bCatinf$ (from Example~\ref{ex:bcatinf} above). Therefore, applying the simplicial nerve construction to $W$, we obtain a functor
\[ \bcat{W}^\bullet: \bDelta^{op} \to \bCatinf; \quad [n] \mapsto \bcat{W}^n. \]
(Recall that we use the same notation for the ordinary category $\bDelta^{op}$ and for its nerve.) The functor $\bcat{W}^\bullet$\index{W@$\bcat{W}^\bullet$, the monoidal $\infty$-category associated to a strict monoidal $\infty$-category $\bcat{W}$} is a monoid in the sense of Definition~\ref{def:monoid}, and we take it as the object of $\MonCat$ which represents the strict monoidal $\infty$-category $(\bcat{W},\otimes)$.
\end{definition}

Ordinary monoidal categories also fit into the framework described above for monoidal $\infty$-categories. Let us write $\mathbf{Cat}$ for simplicial nerve, in the sense of Example~\ref{ex:2,1-nerve}, of the strict $(2,1)$-category of ordinary categories, functors, and natural isomorphisms. Then $\mathbf{Cat}$ admits finite products (given by the usual product of categories), and so we have an $\infty$-category\index{MonCat@$\mathbf{MonCat}$, the $\infty$-category of ordinary monoidal categories}
\[ \mathbf{MonCat} := \mathbf{Mon}(\mathbf{Cat}) \]
of monoids in $\mathbf{Cat}$. 

\begin{proposition} \label{prop:mon-nerve}
The $\infty$-category $\mathbf{MonCat}$ is equivalent to the nerve of the strict $(2,1)$-category $\mathsf{MonCat}$ of ordinary monoidal categories, strong monoidal functors, and monoidal natural isomorphisms. Moreover, the nerve and homotopy category constructions determine an adjunction of $\infty$-categories
\[ h: \MonCat \rightleftarrows \mathbf{MonCat}: N \]
where $N$ is fully faithful. The image of $N$ is the full subcategory of $\mathbf{MonCat}$ consisting of those monoidal $\infty$-categories for which the underlying $\infty$-category is (the nerve of) an ordinary category.
\end{proposition}
\begin{proof}
We assume that the first part is well-known, but we were unable to find a complete proof in the literature, so we sketch one here. We define a map of $\infty$-categories
\[ \iota: N\mathsf{MonCat} \to \Fun(\bDelta^{op},\mathbf{Cat}) \]
as follows. For a monoidal category $\bcat{W}$, we define a functor $\iota_{\bcat{W}}: \bDelta^{op} \to \mathbf{Cat}$ by
\[ \iota_{\bcat{W}}([n]) := \bcat{W}^n \]
with simplicial structure determined by the monoidal unit and multiplication maps by the same pattern described in Definition~\ref{def:strict-mon}. In this case, that construction does not determine an ordinary functor $\bDelta^{op} \to \sset$, but it does still define a functor of $\infty$-categories. The necessary higher coherences are given by the associativity and unit isomorphisms for the monoidal structure. 

A strong monoidal functor $F: \bcat{W} \to \bcat{V}$ determines a natural transformation $\iota_F: \iota_{\bcat{W}} \to \iota_{\bcat{V}}$, i.e.\ a map $\Delta^1 \times \bDelta^{op} \to \mathbf{Cat}$ with components the induced functors
\[ \iota_F([n]) := F^n: \bcat{W}^n \to \bcat{V}^n \]
and naturality data determined by the monoidal structure on $F$.

A $2$-simplex in $N\mathsf{MonCat}$ consists of a monoidal natural isomorphism $\alpha: H \isom GF$ for some triple of compatible strong monoidal functors $F,G,H$. The corresponding $2$-simplex in $\Fun(\bDelta^{op},\mathbf{Cat})$ is the corresponding natural isomorphism $\iota_\alpha: \iota_H \to \iota_G\iota_F$ with components
\[ \iota_{\alpha}([n]) := \alpha^n: H^n \to (GF)^n = G^nF^n. \]
Higher-dimensional simplexes on both sides are larger diagrams of compatible natural isomorphisms, which are preserved by the construction of $\iota$ so far. Thus, we obtain the required map
\[ \iota: N\mathsf{MonCat} \to \Fun(\bDelta^{op},\mathbf{Cat}). \]
We can then show that $\iota$ is fully faithful with essential image $\mathbf{MonCat}$, by analysing the effect of $\iota$ on mapping spaces.

The adjunction $h: \bCatinf \rightleftarrows \mathbf{Cat} : N$ determines an adjunction between the $\infty$-categories of functors
\[ h: \Fun(\bDelta^{op},\bCatinf) \rightleftarrows \Fun(\bDelta^{op},\mathbf{Cat}) : N \]
in which $N$ is still fully faithful. Since both $h$ and $N$ preserve finite products, this adjunction restricts to the desired adjunction between $\infty$-categories of monoids. For the last part, we note that a monoidal $\infty$-category $\bcat{W}$ is in the image of $N$ if and only if $\bcat{W}^n$ is an ordinary category for each $n$. That is the case if and only if $\bcat{W}$ itself is an ordinary category. 
\end{proof}

\subsection{Marked simplicial monoids} \label{sec:msmon}

Our second model for the notion of monoidal $\infty$-category is described by Lurie in~\cite[4.1.8.7]{lurie:2017}, and is based on an $\infty$-categorical version of MacLane's coherence theorem: that any monoidal category is monoidally equivalent to a strict monoidal category. Following this idea, we can take as objects of $\MonCat$ the strict monoidal $\infty$-categories of Definition~\ref{def:mon-quasi}, i.e.\ the simplicial monoids for which the underlying simplicial set is an $\infty$-category.

In order to describe the correct notion of monoidal functor between strict monoidal $\infty$-categories, we describe a model structure on the category of marked simplicial monoids, in which the strict monoidal $\infty$-categories are the fibrant objects. We now review that category and its model structure.

\begin{definition} \label{def:msmon}
A \emph{marked simplicial monoid}\index{marked simplicial monoid} is a monoid in the category $\msset$ of marked simplicial sets of Definition~\ref{def:msset}, or equivalently, a simplicial monoid together with a marking of the underlying simplicial set such that a product of marked edges is marked. The category $\msmon$\index{smon@$\msmon$, category of marked simplicial monoids} of marked simplicial monoids has as its morphisms the maps of simplicial monoids which preserve markings.

The category $\msmon$ inherits a simplicial model structure (the \emph{projective model structure}\index{projective model structure on marked simplicial monoids}) from $\msset$ in which a morphism is a weak equivalence (or, respectively, a fibration) if the underlying map of marked simplicial sets is a weak equivalence (or, respectively, a fibration) in the marked model structure. In particular, the fibrant objects in $\msmon$ are the strict monoidal $\infty$-categories (in the sense of Definition~\ref{def:mon-quasi}) with marked edges given by those morphisms which are equivalences in the underlying $\infty$-category.
\end{definition}

\begin{remark} \label{rem:msmon}
We have already described, in Definition~\ref{def:strict-mon}, how a strict monoidal $\infty$-category determines a monoid in $\bCatinf$. That construction underlies an equivalence of $\infty$-categories
\[ N(\msmon^{\circ}) \weq \mathbf{MonCat} \]
outlined in~\cite[4.1.8.7]{lurie:2017}, and based on~\cite[4.1.8.4]{lurie:2017}, together with~\cite[2.4.2.5]{lurie:2017}.

In particular, any monoidal $\infty$-category is equivalent, in $\MonCat$, to one arising from a strict monoidal $\infty$-category by the construction of Definition~\ref{def:strict-mon}. Monoidal functors $\bcat{W} \to \bcat{X}$ between strict monoidal $\infty$-categories are given by morphisms of marked simplicial monoids $\tilde{\bcat{W}} \to \bcat{X}$, where $\tilde{\bcat{W}}$ denotes a cofibrant replacement for $\bcat{W}$ in the projective model structure on marked simplicial monoids.
\end{remark}

\subsection{One-object $\infty$-bicategories}

Ordinary monoidal categories can be characterized as bicategories with a single object. That description of monoidal structures also extends to the $\infty$-categorical setting, though with a caveat. The objects in an $\infty$-bicategory naturally form a \emph{space} (the core $\infty$-groupoid), and monoidal $\infty$-categories correspond to those $\infty$-bicategories for which that space is \emph{nonempty and connected}, but not necessarily contractible. In fact, this distinction is there for ordinary monoidal categories too, where the core groupoid of the corresponding one-object bicategory is the Picard group of the monoidal category (the group of objects invertible under the monoidal structure).

\begin{proposition} \label{prop:moncat-bicat}
There is an equivalence of $\infty$-categories
\[ \MonCat \homeq \mathbf{Cat}_{(\infty,2)}^0 \]
where \index{CATinf220@$\mathbf{Cat}_{(\infty,2)}^0$, the $\infty$-category of $\infty$-bicategories with a single object}$\mathbf{Cat}_{(\infty,2)}^0 \subseteq \mathbf{Cat}_{(\infty,2)}$ denotes the full subcategory consisting of those $\infty$-bicategories with a single $0$-simplex.
\end{proposition}
\begin{proof}
Lurie proves in~\cite[4.2.7]{lurie:2009a} that the `marked nerve' construction is a right Quillen equivalence from the category of `marked simplicial categories', i.e.\ categories enriched in marked simplicial sets, and the scaled model structure on $\scset$. We can also identify the projective model structure on marked simplicial monoids $\msmon$ with the full subcategory of one-object marked simplicial categories. Combining these two constructions with the equivalence of Remark~\ref{rem:msmon}, we obtain an embedding
\[ \MonCat \to \mathbf{Cat}_{(\infty,2)} \]
whose image is the full subcategory $\mathbf{Cat}_{(\infty,2)}^0$.
\end{proof}

\begin{remark} \label{rem:icon}
Proposition~\ref{prop:moncat-bicat} suggests that we might extend $\MonCat$ to an $\infty$-bicategory by considering the full subcategory $\mathbf{CAT}_{(\infty,2)}^0 \subseteq \mathbf{CAT}_{(\infty,2)}$ of the $\infty$-bicategory of $\infty$-bicategories (and functors and natural transformations) of Example~\ref{ex:CAT2}. However, this does \emph{not} recover the standard notion of monoidal natural transformation between monoidal functors. Those instead correspond to certain \emph{oplax} natural transformations, or more specifically \emph{icons}; see~\cite{lack:2010}.
\end{remark}

We can now describe explicitly the $\infty$-bicategory which represents a strict monoidal $\infty$-category.

\begin{example} \label{ex:moncat-bicat}
Let $\bcat{W}$ be a strict monoidal $\infty$-category. Then there is a simplicial category with a single object and mapping space equal to $\bcat{W}$, where composition is given by the monoid structure. Applying the simplicial nerve to this simplicial category, we get a simplicial set $\un{\bcat{W}}$ in which
\begin{itemize}
	\item there is a single vertex which we denote $\bullet$;
	\item there is an edge for each object of $\bcat{W}$;
	\item a $2$-simplex of the form
	\[ \begin{diagram}
		\node{\bullet} \arrow[2]{s,l}{A_{01}} \arrow{ese,t}{A_{02}} \\ \node[3]{\bullet} \\
		\node{\bullet} \arrow{ene,b}{A_{12}}
	\end{diagram} \]
	consists of a morphism in $\bcat{W}$ of the form $f_{012}: A_{02} \to A_{01} \otimes A_{12}$;
	\item a $3$-simplex consists of a square diagram (i.e. two $2$-simplexes) in $\bcat{W}$ of the form
	\[ \begin{diagram}
		\node{A_{03}} \arrow{se} \arrow{e,t}{f_{013}} \arrow{s,l}{f_{023}} \node{A_{01} \otimes A_{13}} \arrow{s,r}{A_{01} \otimes f_{123}} \\
		\node{A_{02} \otimes A_{23}} \arrow{e,b}{f_{012} \otimes A_{23}} \node{A_{01} \otimes A_{12} \otimes A_{23}}
	\end{diagram} \]
	\item an $n$-simplex consists of an $(n-2)$-cubical diagram in $\bcat{W}$ extending the pattern above.
\end{itemize}
It follows from~\cite[4.2.7]{lurie:2009a} that $\un{\bcat{W}}$\index{W@$\un{\bcat{W}}$, the one-object $\infty$-bicategory associated to a strict monoidal $\infty$-category $\bcat{W}$} is an $\infty$-bicategory in which a $2$-simplex is thin if and only if the corresponding morphism $A'' \to A \otimes A'$ is an equivalence in $\bcat{W}$.
\end{example}

\section{Module $\infty$-categories} \label{sec:module}

We now turn to module $\infty$-categories over a monoidal $\infty$-category $\bcat{W}$. Intuitively, such a \emph{$\bcat{W}$-module} consists of an $\infty$-category $\bcat{X}$ together with a functor
\[ \bcat{W} \times \bcat{X} \to \bcat{X}, \]
which satisfies the usual axioms of an action up to (coherent) equivalence. For $\bcat{W}$-modules $\bcat{X},\bcat{Y}$, we also define \emph{$\bcat{W}$-functors} $\bcat{X} \to \bcat{Y}$ to be functors that commute (up to coherent equivalence) with the $\bcat{W}$-actions. Our main goal in this section is to define an $\infty$-bicategory $\Mod_{\bcat{W}}$ of $\bcat{W}$-modules, $\bcat{W}$-functors, and $\bcat{W}$-natural transformations. We will consider $\bcat{W}$-modules from each of the three perspectives used above to present a monoidal $\infty$-category.

\subsection{Modules over a monoid in the $\infty$-category of $\infty$-categories} \label{sec:mod1}

Recall that we have defined a monoidal $\infty$-category to be a monoid (in the $\infty$-categorical sense) in the $\infty$-category $\bCatinf$. Associated to any monoid in an $\infty$-category is a notion of module, which we now recall from~\cite[4.2.2.2]{lurie:2017}.

\begin{definition} \label{def:mod1}
Let $M: \bDelta^{op} \to \mathbf{C}$ be a monoid in the $\infty$-category $\mathbf{C}$ in the sense of Definition~\ref{def:monoid}. Then \emph{an $M$-module in $\mathbf{C}$}\index{module over a monoid in an $\infty$-category} consists of another functor
\[ X : \bDelta^{op} \to \mathbf{C} \]
and a natural transformation
\[ p: X \to M \]
such that for each $n \geq 0$, the following naturality square for $p$ is a pullback in $\bCatinf$,
\begin{equation} \label{eq:mod1} \begin{diagram}
	\node{X([n])} \arrow{s,l}{X(\rho)} \arrow{e,t}{p_{[n]}} \node{M([n])} \arrow{s,r}{M(\rho)} \\
	\node{X([0])} \arrow{e,t}{p_{[0]}} \node{M([0])}
\end{diagram} \end{equation}
where $\rho: [0] \to [n]; 0 \mapsto 0$. The \emph{underlying object} of an $M$-module is $X([0])$, and the pullback condition above, together with the product condition in Definition~\ref{def:monoid}, implies that
\[X([n]) \weq M([1])^n \times X([0]). \]
We often abuse notation by writing $X([0])$ simply as $X$, and by writing $M([1])$ as $M$. 

The module structure map for an $M$-module $X$ can be viewed as the composite
\[ \dgTEXTARROWLENGTH=2.5em M \times X \lweq X([1]) \arrow{e,t}{X(\mu)} X([0]) = X \]
where the second map is induced by $\mu: [0] \to [1]; 0 \mapsto 1$.

The \emph{$\infty$-category of $M$-modules in $\mathbf{C}$}\index{ModVC@$\mathbf{Mod}_{M}(\mathbf{C})$, the $\infty$-category of modules over a monoid $M$ in an $\infty$-category $\mathbf{C}$} is the full subcategory
\[ \mathbf{Mod}_{M}(\mathbf{C}) \subseteq \Fun(\bDelta^{op},\mathbf{C})_{/M} \]
of the slice $\infty$-category, whose objects are the $M$-modules. In particular, an \emph{$M$-module morphism}\index{module morphism} from $X$ to $Y$ is a diagram
\[ \begin{diagram}
	\node{X} \arrow[2]{e,t}{f} \arrow{se,b}{p} \node[2]{Y} \arrow{sw,b}{p'} \\
	\node[2]{M}
\end{diagram} \]
in the $\infty$-category of functors $\bDelta^{op} \to \mathbf{C}$. This definition amounts to a morphism from $X$ to $Y$ in $\mathbf{C}$ that commutes with the $M$-module structure maps up to coherent equivalence. 
\end{definition}

Now let $\bcat{W}$ be a monoidal $\infty$-category, i.e.\ a monoid in $\bCatinf$. Then we have the following notion.

\begin{definition} \label{def:mod-mon}
Let $\bcat{W}$ be a monoidal $\infty$-category. Then the $\infty$-category of $\bcat{W}$-module $\infty$-categories\index{module $\infty$-category!over a monoidal $\infty$-category} is
\[ \Mod_{\bcat{W}}(\bCatinf) \subseteq \Fun(\bDelta^{op},\bCatinf)_{/\bcat{W}}, \]
the $\infty$-category of modules\index{ModW1@$\Mod_{\bcat{W}}(\bCatinf)$, the $\infty$-category of module $\infty$-categories over a monoidal $\infty$-category $\bcat{W}$} in $\bCatinf$ over the monoid $\bcat{W}$. Given $\bcat{W}$-module $\infty$-categories $\bcat{X}$ and $\bcat{Y}$, a \emph{$\bcat{W}$-functor}\index{functor@$\bcat{W}$-functor!between module $\infty$-categories over a monoidal $\infty$-category $\bcat{W}$} from $\bcat{X}$ to $\bcat{Y}$ is a morphism from $\bcat{X}$ to $\bcat{Y}$ in $\Mod_{\bcat{W}}(\bCatinf)$.
\end{definition}

\begin{definition} \label{def:mon-mod-functor}
Given a monoidal functor $F: \bcat{W} \to \bcat{W}'$ between monoidal $\infty$-categories, composition with $F$ determines a map between the slice $\infty$-categories, and hence a functor
\[ F^*: \Mod_{\bcat{W}'}(\bCatinf) \to \Mod_{\bcat{W}}(\bCatinf) \]
which is the identity on underlying $\infty$-categories. This construction describes the usual pullback of module structure along $F$.
\end{definition}

\begin{example} \label{ex:mon-self-mod}
Let $\bcat{W}: \bDelta^{op} \to \bCatinf$ be a monoidal $\infty$-category in the sense of Definition~\ref{def:moncat1}. Define $\bcat{W}' : \bDelta^{op} \to \bCatinf$ by 
\[ \bcat{W}' := \bcat{W} \circ s \]
where $s: \bDelta^{op} \to \bDelta^{op}$ is given by
\[ s([n]) := [n+1] \]
and for $\alpha: [m] \to [n]$, $s(\alpha): [m+1] \to [n+1]$ is given by
\[ s(\alpha)(j) := \begin{cases} s(j) & \text{if $j \leq m$}; \\ n+1 & \text{if $j = m+1$}. \end{cases} \]
There is a natural transformation $\pi: s \to 1_{\bDelta^{op}}$ with components
\[ \pi_{[n]}: [n+1] \to [n] \]
given by $j \mapsto \min(j,n)$. That induces a natural transformation $p := \bcat{W} \pi : \bcat{W}' \to \bcat{W}$, which is a $\bcat{W}$-module $\infty$-category, whose underlying $\infty$-category is $\bcat{W}'([0]) := \bcat{W}([1])$. This module encodes the canonical action\index{monoidal $\infty$-category!action on itself} of a monoidal $\infty$-category on itself; see~\cite[4.2.2.4]{lurie:2017}.
\end{example}

\begin{remark}
The construction of Definition~\ref{def:mod-mon} produces the $\infty$-\emph{category} of $\bcat{W}$-modules and $\bcat{W}$-functors, encoding only the invertible $\bcat{W}$-natural transformations between those functors. In Definition~\ref{def:mod-infty-bicat} below we extend the $\infty$-category $\Mod_{\bcat{W}}(\bCatinf)$ to an $\infty$-bicategory $\Mod_{\bcat{W}}(\mathbf{CAT}_\infty)$, which incorporates also non-invertible $\bcat{W}$-natural transformations.
\end{remark}

\subsection{Modules over a marked simplicial monoid}

We noted in Section~\ref{sec:msmon} that the $\infty$-category $\MonCat$ is presented by the projective model structure on marked simplicial monoids, whose fibrant objects are the \emph{strict} monoidal $\infty$-categories of Definition~\ref{def:mon-quasi}. We show now that module $\infty$-categories over a strict monoidal $\infty$-category $\bcat{W}$ are also presented by a model structure on \emph{strict} $\bcat{W}$-modules, that is (marked) simplicial sets with a strictly associative and unital action by the simplicial monoid $\bcat{W}$.

\begin{definition} \label{def:mod3}
Let $\bcat{W}$ be a strict monoidal $\infty$-category. A \emph{strict $\bcat{W}$-module $\infty$-category}\index{module $\infty$-category!strict}\index{strict module $\infty$-category} is an $\infty$-category $\bcat{X}$ together with a strictly associative and unital action of the simplicial monoid $\bcat{W}$ on the simplicial set $\bcat{X}$.

More generally, let $\mathsf{Mod}_{\bcat{W}}^+$\index{ModZW+@$\mathsf{Mod}_{\bcat{W}}^+$, the category of marked modules over a marked simplicial monoid $\bcat{W}$} denote the (ordinary) category of modules over the marked simplicial monoid $\bcat{W}$. That is, an object of $\mathsf{Mod}_{\bcat{W}}^+$ consists of a marked simplicial set $\bcat{X}$ together with a (left) $\bcat{W}$-action
\[ \bcat{W} \times \bcat{X} \to \bcat{X} \]
satisfying the usual (strict) axioms for associativity and unitality. A morphism is a map of marked simplicial sets which commutes (strictly) with the $\bcat{W}$-action maps.

We endow $\mathsf{Mod}_{\bcat{W}}^+$ with the projective model structure in which a morphism is a weak equivalence or fibration if and only if it is so in the marked model structure on $\msset$ of Definition~\ref{def:msset}. In particular, a fibrant object in $\mathsf{Mod}_{\bcat{W}}^+$ is a strict $\bcat{W}$-module $\infty$-category, with marked edges given by the equivalences in the underlying $\infty$-category.
\end{definition}

Our goal in this section is to show that the $\infty$-category associated to the projective model structure on $\mathsf{Mod}_{\bcat{W}}^+$ is equivalent to $\Mod_{\bcat{W}}(\bCatinf)$, and hence gives another model for the $\infty$-category of $\bcat{W}$-modules and their functors. We start by describing how a strict $\bcat{W}$-module $\infty$-category gives rise to a module over the monoidal $\infty$-category corresponding to $\bcat{W}$.

\begin{definition} \label{def:marked-module-monoid}
Let $\bcat{W}$ be a strict monoidal $\infty$-category, and let $\bcat{X}$ be a strict $\bcat{W}$-module $\infty$-category with action map $\odot: \bcat{W} \times \bcat{X} \to \bcat{X}$. Then we define an ordinary functor $X: \bDelta^{op} \to \sset$ as follows\index{module $\infty$-category!associated to a strict module}.

We set $X([n]) := \bcat{W}^n \times \bcat{X}$. For $\alpha: [m] \to [n]$, we define $X([n]) \to X([m])$ to be the map which sends $(A_n,\dots,A_1,C)$ to
\[ (A_{\alpha(m)} \otimes \dots \otimes A_{\alpha(m-1)+1}, \dots, A_{\alpha(1)} \otimes \dots \otimes A_{\alpha(0)+1}, A_{\alpha(0)} \otimes \dots \otimes A_1 \odot C). \]
The strictness of the $\bcat{W}$-action on $\bcat{X}$ ensures that this construction does indeed define a functor $X: \bDelta^{op} \to \sset$. By construction, there is a simplicial natural transformation $p: X \to W$ given by the projection maps $\bcat{W}^n \times \bcat{X} \to \bcat{W}^n$. Applying the simplicial nerve, we get a natural transformation between functors $\bDelta^{op} \to \bCatinf$, and by construction, this natural transformation satisfies the condition of Definition~\ref{def:mod1}.
\end{definition}

We now use a result of Lurie to deduce that the projective model structure on $\mathsf{Mod}_{\bcat{W}}^+$ presents the $\infty$-category of $\bcat{W}$-modules.

\begin{proposition} \label{prop:mod3}
Let $\bcat{W}$ be a strict monoidal $\infty$-category, and let $W: \bDelta^{op} \to \bCatinf$ be the functor given by $W([n]) := \bcat{W}^n$ as in Definition~\ref{def:strict-mon}. Then there is an equivalence of $\infty$-categories
\[ (\mathsf{Mod}_{\bcat{W}}^+)^{\circ} \weq \Mod_{W}(\bCatinf) \]
given on objects by the construction of Definition~\ref{def:marked-module-monoid}.
\end{proposition}
\begin{proof}
This result is an application of~\cite[4.3.3.17]{lurie:2017} in the case where (using the notation therein) $\bcat{A}$ is the category of marked simplicial sets (with the marked model structure), $A = *$ is the trivial object, and $B = \bcat{W}$.
\end{proof}

\begin{corollary}
Let $\bcat{W}$ be a strict monoidal $\infty$-category. Then any $\bcat{W}$-module $\infty$-category $\bcat{X}$ (in the sense of Definition~\ref{def:mod1}) is equivalent to a strict $\bcat{W}$-module $\infty$-category (in the sense of Definition~\ref{def:mod3}).
\end{corollary}

Proposition~\ref{prop:mod3} also gives us a way to understand $\bcat{W}$-functors between strict $\bcat{W}$-module $\infty$-categories via an explicit cofibration replacement in the model structure on $\mathsf{Mod}^+_{\bcat{W}}$.

\begin{definition} \label{def:bar}
Let $\bcat{W}$ be a strict monoidal $\infty$-category, and let $\bcat{X}$ be a strict $\bcat{W}$-module $\infty$-category. We define the bar construction\index{bar construction}\index{BX@$B\bcat{X}$, the bar construction on a strict module $\infty$-category $\bcat{X}$}
\[ B\bcat{X} = B(\bcat{W},\bcat{W},\bcat{X}) \]
to be the geometric realization of the simplicial object $B_\bullet\bcat{X}$ in $\mathsf{Mod}_{\bcat{W}}^+$ given by
\[ B_k\bcat{X} := \bcat{W}^{k+1} \times \bcat{X} \]
with free left $\bcat{W}$-action, with face maps given by the monoid structure on $\bcat{W}$ and the $\bcat{W}$-action on $\bcat{X}$, and degeneracy maps given by the unit map $* \to \bcat{W}$ to the identity object of $\bcat{W}$. There is a canonical map of marked $\bcat{W}$-modules
\[ \epsilon: B\bcat{X} \to \bcat{X} \]
given by applying the $\bcat{W}$-action of $\bcat{X}$.
\end{definition}

\begin{remark}
Explicitly, the geometric realization can be described as follows. The underlying simplicial set of $B\bcat{X}$ is the diagonal of the bisimplicial set $B_\bullet\bcat{X}$. That is, we have
\[ B_k\bcat{X} = \bcat{W}^{k+1}_k \times \bcat{X}_k \]
with free left $\bcat{W}_k$-action, and face and degeneracy maps given by combining those for $B_\bullet\bcat{X}$ and for the $\bcat{W}$ and $\bcat{X}$ separately. An edge in $B\bcat{X}$ is a triple of $1$-morphisms $(\alpha_1,\alpha_0,\beta) \in \bcat{W}_1 \times \bcat{W}_1 \times \bcat{X}_1$, and such an edge is marked if $\alpha_1,\alpha_0$ are equivalences in the $\infty$-category $\bcat{W}$ and $\beta$ is an equivalence in $\bcat{X}$.
\end{remark}

\begin{lemma} \label{lem:bar-cof}
The map $\epsilon: B\bcat{X} \to \bcat{X}$ exhibits $B\bcat{X}$ as a cofibrant replacement of $\bcat{X}$ in the projective model structure on $\mathsf{Mod}_{\bcat{W}}^+$.
\end{lemma}
\begin{proof}
Cofibrations in $\msset$ are the monomorphisms, so the simplicial bar construction $B_\bullet\bcat{X}$ is levelwise cofibrant, and the degeneracies are cofibrations. Hence, this simplicial bar construction is Reedy cofibrant, and its geometric realization is cofibrant in $\mathsf{Mod}_{\bcat{W}}^+$. The map $\epsilon$ is a weak equivalence in $\msset$ (and hence in the projective model structure) since the underlying simplicial object in $\msset$ of the bar construction has an extra degeneracy.
\end{proof}

\begin{remark} \label{rem:bar}
It follows from Lemma~\ref{lem:bar-cof} that any $\bcat{W}$-functor\index{functor@$\bcat{W}$-functor!between strict $\bcat{W}$-module $\infty$-categories} between strict $\bcat{W}$-module $\infty$-categories $\bcat{X},\bcat{Y}$ can be represented by a morphism of marked $\bcat{W}$-modules
\[ F_\bullet: B\bcat{X} \to \bcat{Y}, \]
that is, a map of marked simplicial sets which commutes strictly with the $\bcat{W}$-actions. We can unpack the structure of such a map as follows, using the description of $B\bcat{X}$ as a coend in the category of marked $\bcat{W}$-modules.

A $\bcat{W}$-functor from $\bcat{X}$ to $\bcat{Y}$ comprises, for each $n \geq 0$, an $n$-simplex $F_n$ in the simplicial mapping space
\[ \Hom_{\mathsf{Mod}^+_{\bcat{W}}}(\bcat{W}^{n+1} \times \bcat{X},\bcat{Y}), \]
which, by the free-forgetful adjunction between marked simplicial sets and marked $\bcat{W}$-modules, can be identified with the simplicial mapping space
\[ \Hom_{\msset}(\bcat{W}^n \times \bcat{X},\bcat{Y}), \]
which, in turn, is the ordinary simplicial mapping space between these two $\infty$-categories, that is
\[ \Fun(\bcat{W}^n \times \bcat{X},\bcat{Y}). \]
In particular, we can view $F_n$ as a functor
\[ F_n: \Delta^n \times \bcat{W}^n \times \bcat{X} \to \bcat{Y}, \]
or equivalently as a functor
\[ \bcat{W}^n \times \bcat{X} \to \Fun(\Delta^n,\bcat{Y}). \]
When written this way, we can express the required compatibility between the maps $F_n$ in terms of face and degeneracy maps between simplexes in $\bcat{Y}$. We have the following conditions: for the face maps
\[ d^i(F_n(A_n,\dots,A_1,X)) = \begin{cases} F_{n-1}(A_n,\dots,A_2,A_1 \odot_{\bcat{X}} X) & \text{if $i = 0$}; \\
	F_{n-1}(A_n,\dots,A_{i+1} \otimes A_i,\dots, A_1, X) & \text{if $0 < i < n$}; \\
	A_n \odot_{\bcat{Y}} F_{n-1}(A_{n-1},\dots,A_1,X) & \text{if $i = n$}; \end{cases} \]
and for the degeneracies
\[ s^j(F_n(A_n,\dots,A_1,X)) = F_{n+1}(A_n,\dots,A_{j+1},I,A_j,\dots,A_1,X) \quad \text{for $0 \leq j \leq n$}, \]
where $I$ denotes the unit object in the simplicial monoid $\bcat{W}$. 

Moreover, the condition that $F_\bullet$ is a map of \emph{marked} $\bcat{W}$-modules is equivalent to the requirement that all the edges in the $n$-simplexes $F_n(A_n,\dots,A_1,X)$ are equivalences in the $\infty$-category $\bcat{Y}$. Equivalently, using the face map compatibilities above, it is sufficient to require that each $1$-simplex $F_1(A,X): A \odot_{\bcat{Y}} F_0(X) \to F_0(A \odot_{\bcat{X}} X)$ is an equivalence in $\bcat{Y}$.

We can unpack the condition above further in low degrees. A $\bcat{W}$-functor from $\bcat{X}$ to $\bcat{Y}$ includes:
\begin{itemize}
\item a functor $F: \bcat{X} \to \bcat{Y}$
\item a natural equivalence $\alpha_{A,X}: A \odot_{\bcat{Y}} F(X) \isom F(A \odot_{\bcat{X}} X)$ between functors $\bcat{W} \times \bcat{X} \to \bcat{Y}$, such that $\alpha_{I,X} = 1_{F(X)}$;
\item a (natural) collection of $2$-simplexes $\alpha_{A_2,A_1,X}$ in $\bcat{Y}$ of the form
\[ \begin{diagram}
\node{A_2 \otimes A_1 \odot_{\bcat{Y}} F(X)} \arrow[2]{e,t}{\alpha_{A_2 \otimes A_1,X}} \arrow{se,b}{A_2 \otimes \alpha_{A_1,X}} \node[2]{F(A_2 \otimes A_1 \odot_{\bcat{X}} X)} \\
\node[2]{A_2 \odot_{\bcat{Y}} F(A_1 \odot_{\bcat{X}} X)} \arrow{ne,b}{\alpha_{A_2,A_1 \odot_{\bcat{X}} X}}
\end{diagram} \]
which are degenerate when either $A_2$ or $A_1$ is the unit object of $\bcat{W}$;
\item similar natural collections of $n$-simplexes $\alpha_{A_n,\dots,A_1,X}$ for higher $n$ with boundary prescribed by the lower degree information, and that are suitably degenerate when one of $A_n,\dots,A_1$ equals $I$.
\end{itemize}
In this presentation, we can most easily see the connection to the notion of functor between actegories over an ordinary monoidal category; see~\cite[3.3.2]{capucci/gavranovic:2022}. 
\end{remark}

\subsection{Functors between $\infty$-bicategories}

A monoidal $\infty$-category $\bcat{W}$ can be represented by a one-object $\infty$-bicategory, which we denote $\un{\bcat{W}}$. Then the notion of a module $\infty$-category over that monoidal $\infty$-category is simple to describe.

\begin{definition} \label{def:mod4}
A \emph{$\bcat{W}$-module $\infty$-category}\index{module $\infty$-category!as a functor of $\infty$-bicategories} is a functor of $\infty$-bicategories
\[ \bcat{X}: \un{\bcat{W}} \to \mathbf{CAT}_\infty, \]
where $\mathbf{CAT}_\infty$ is the $\infty$-bicategory of $\infty$-categories of Example~\ref{ex:bicat-cat}. That is, $\bcat{X}$ is a map of simplicial sets which sends thin $2$-simplexes in $\un{\bcat{W}}$ to thin $2$-simplexes in $\mathbf{CAT}_\infty$.

We can then define the \emph{$\infty$-bicategory of $\bcat{W}$-module $\infty$-categories}\index{ModW2@$\Mod_{\bcat{W}}(\mathbf{CAT}_\infty)$, the $\infty$-bicategory of module $\infty$-categories over the monoidal $\infty$-category $\bcat{W}$} to be 
\[ \Mod_{\bcat{W}}(\mathbf{CAT}_\infty) := \Fun(\un{\bcat{W}},\mathbf{CAT}_\infty). \]
For $\bcat{W}$-module $\infty$-categories $\bcat{X},\bcat{Y}$, we write
\[ \Mod_{\bcat{W}}(\bcat{X},\bcat{Y}) := \Hom_{\Mod{\bcat{W}}(\mathbf{CAT}_\infty)}(\bcat{X},\bcat{Y}) \]
for the $\infty$-category\index{ModWXY@$\Mod_{\bcat{W}}(\bcat{X},\bcat{Y})$, the $\infty$-category of $\bcat{W}$-functors between $\bcat{W}$-module $\infty$-categories $\bcat{X}$ and $\bcat{Y}$} of $\bcat{W}$-functors from $\bcat{X}$ to $\bcat{Y}$. The morphisms in this $\infty$-category are the \emph{$\bcat{W}$-natural transformations}\index{natural transformation!between functors between module $\infty$-categories} between such functors.

Informally, given $\bcat{W}$-functors $F,G: \bcat{X}^{\odot} \to \bcat{Y}^{\odot}$, a $\bcat{W}$-natural transformation $\alpha: F \Rightarrow G$ is one such that, for each $X \in \bcat{X}$ and $A \in \bcat{W}$ the following diagram commutes up to higher equivalence:
\[ \begin{diagram}
	\node{A \odot F(X)} \arrow{s,l}{A \odot \alpha_X} \arrow{e,t}{\sim} \node{F(A \odot X)} \arrow{s,r}{\alpha_{A \odot X}} \\
	\node{A \odot G(X)} \arrow{e,t}{\sim} \node{G(A \odot X)}
\end{diagram} \]
where the horizontal maps express the compatibility of $F$ and $G$ with the $\bcat{W}$-actions on $\bcat{X}$ and $\bcat{Y}$.
\end{definition}

To see the connection between the $\infty$-bicategorical approach to module $\infty$-categories and those before, recall that $\infty$-bicategories are, roughly speaking, a model for categories enriched in $\infty$-categories. Then a functor of $\infty$-bicategories $\bcat{X}: \un{\bcat{W}} \to \mathbf{CAT}_\infty$ comprises a chosen $\infty$-category $\bcat{X} := \bcat{X}(\bullet)$ together with a functor
\[ \Hom_{\un{\bcat{W}}}(\bullet,\bullet) \to \Hom_{\mathbf{CAT}_\infty}(\bcat{X}(\bullet),\bcat{X}(\bullet)) \]
which commutes, up to equivalence, with composition and preserves identities. In other words, $\bcat{X}$ is given by a map of monoidal $\infty$-categories
\[ \bcat{W} \to \Fun(\bcat{X},\bcat{X}), \]
which corresponds to the action of $\bcat{W}$ on $\bcat{X}$.

Here is a more precise construction, from a \emph{strict} monoidal $\infty$-category $\bcat{W}$ and \emph{strict} $\bcat{W}$-module $\infty$-category $\bcat{X}$, of a functor of $\infty$-bicategories of the form in Definition~\ref{def:mod4}.

\begin{definition} \label{def:mon-mod-bicat}
Let $\bcat{W}$ be a strict monoidal $\infty$-category, and recall that we can view $\bcat{W}$ as a one-object simplicial category with $\bcat{W}$ as the simplicial mapping space and the monoid structure on $\bcat{W}$ giving the identity morphism and composition. Denote that simplicial category by $\mathsf{B}\bcat{W}$.\index{BW@$\mathsf{B}\bcat{W}$, the one-object simplicial category associated to a strict monoidal $\infty$-category $\bcat{W}$}

Now let $\bcat{X}$ be a strict $\bcat{W}$-module $\infty$-category. Then $\bcat{X}$ defines a simplicially-enriched functor $\mathsf{B}\bcat{W} \to \mathsf{CAT_\infty}$, which sends the one object to $\bcat{X}$ and is given on mapping spaces by the $\bcat{W}$-action map
\[ \bcat{W} \to \Fun(\bcat{X},\bcat{X}) = \Hom_{\mathsf{CAT_\infty}}(\bcat{X},\bcat{X}). \]
Applying the simplicial nerve to this functor, we get a map of $\infty$-bicategories of the form
\[ \un{\bcat{X}}: \un{\bcat{W}} \to \mathbf{CAT}_\infty. \]
\end{definition}

We can now show that the $\infty$-bicategory $\Mod_{\bcat{W}}(\mathbf{CAT}_\infty)$ described in Definition~\ref{def:mod4} agrees, after passing to the core $\infty$-category, with the approaches described in the previous parts of this section.

\begin{proposition} \label{prop:mon-mod-bicat}
Let $\bcat{W}$ be a strict monoidal $\infty$-category, and let $\un{\bcat{W}}$ be the corresponding $\infty$-bicategory with one object described in Example~\ref{ex:moncat-bicat}. Then there is an equivalence of $\infty$-categories
\[ (\mathsf{Mod}^+_{\bcat{W}})^{\circ} \weq \Fun(\un{\bcat{W}},\mathbf{Cat}_\infty) \isom \Fun(\un{\bcat{W}},\mathbf{CAT}_\infty)^{\sim}, \]
where the right-hand side is the core $\infty$-category of the $\infty$-bicategory of functors of $\infty$-bicategories $\un{\bcat{W}} \to \mathbf{CAT}_\infty$; see Definition~\ref{def:inf-bicat}.
\end{proposition}
\begin{proof}
The construction of~\ref{def:mon-mod-bicat} determines an equivalence (of ordinary categories)
\[ \mathsf{Mod}^+_{\bcat{W}} \homeq \Fun_{\msset}(\mathsf{B}\bcat{W},\msset), \]
where the right-hand side is the category of $\msset$-enriched functors $F: \mathsf{B}\bcat{W} \to \msset$. The inverse equivalence sends such a functor $F$ to the marked simplicial set $F(\bullet)$ with $\bcat{W}$-action determined by the enrichment. The projective model structure on $\mathsf{Mod}^+_{\bcat{W}}$ then corresponds to the projective model structure on this enriched functor category; in both the weak equivalences and fibrations are detected in the underlying category $\msset$.

It is therefore sufficient to show that the simplicial nerve determines an equivalence of $\infty$-categories
\[ \Fun_{\msset}(\mathsf{B}\bcat{W},\msset)^{\circ} \weq \Fun(\un{\bcat{W}},\mathbf{Cat}_\infty). \]
That claim follows from~\cite[A.3.4.14]{lurie:2009}.
\end{proof}

\begin{remark}
The $2$-simplexes in $\Mod_{\bcat{W}}$ determine a notion of $\bcat{W}$-natural transformation between $\bcat{W}$-functors of $\bcat{W}$-module $\infty$-categories. Informally, given two such functors $F,G: \bcat{X}^{\odot} \to \bcat{Y}^{\odot}$, a $\bcat{W}$-natural transformation $\alpha: F \Rightarrow G$ is one such that, for each $X \in \bcat{X}$ and $A \in \bcat{W}$ the following diagram commutes up to higher equivalence:
\[ \begin{diagram}
	\node{A \odot F(X)} \arrow{s,l}{A \odot \alpha_X} \arrow{e,t}{\sim} \node{F(A \odot X)} \arrow{s,r}{\alpha_{A \odot X}} \\
	\node{A \odot G(X)} \arrow{e,t}{\sim} \node{G(A \odot X)}
\end{diagram} \]
where the horizontal maps express the compatibility of $F$ and $G$ with the $\bcat{W}$-actions on $\bcat{X}$ and $\bcat{Y}$.
\end{remark}

Definition~\ref{def:mod4} generalizes easily to a notion of $\bcat{W}$-module in an arbitrary $\infty$-bicategory $\mathbf{C}$, which we use in Chapter~\ref{sec:infty2} to define tangent objects in $\mathbf{C}$.

\begin{definition} \label{def:mod-infty-bicat}
Let $\mathbf{C}$ be an $\infty$-bicategory, and let $\bcat{W}$ be a monoidal $\infty$-category. Then a \emph{$\bcat{W}$-module in $\mathbf{C}$}\index{module over a monoidal $\infty$-category $\bcat{W}$ in an $\infty$-bicategory $\mathbf{C}$} is a functor of $\infty$-bicategories
\[ \un{\bcat{W}} \to \mathbf{C} \]
and we write
\[ \Mod_{\bcat{W}}(\mathbf{C}) := \Fun(\un{\bcat{W}},\mathbf{C}) \]
for the $\infty$-bicategory of $\bcat{W}$-modules in $\mathbf{C}$\index{ModWC@$\Mod_{\bcat{W}}(\mathbf{C})$, the $\infty$-bicategory of $\bcat{W}$-modules in the $\infty$-bicategory $\mathbf{C}$}
\end{definition}

We also note that monoidal functors induce functors between $\infty$-bicategories of modules.

\begin{definition} \label{def:n-module}
Let $n: \bcat{V} \to \bcat{W}$ be a monoidal functor between monoidal $\infty$-categories. Then $n$ determines a map of $\infty$-bicategories
\[ \un{n}: \un{\bcat{V}} \to \un{\bcat{W}}. \]
Precomposition with $\un{n}$ then determines a functor of $\infty$-bicategories
\[ n^*: \Mod_{\bcat{W}}(\mathbf{C}) = \Fun(\un{\bcat{W}},\mathbf{C}) \to \Fun(\un{\bcat{V}},\mathbf{C}) = \Mod_{\bcat{V}}(\mathbf{C}). \]
\end{definition}

We conclude this section by examining the relationship between modules over a monoidal $\infty$-category and over its homotopy category. First we show that when $\bcat{W}$ is an ordinary monoidal category, Definition~\ref{def:mod-infty-bicat} recovers the notion of $\bcat{W}$-module, or $\bcat{W}$-actegory; see~\cite{capucci/gavranovic:2022}.

\begin{proposition} \label{prop:nerve-mod-bicat}
Let $\bcat{W}$ be an ordinary monoidal category, viewed also as a monoidal $\infty$-category via the equivalence of Proposition~\ref{prop:mon-nerve}. Let $\mathbf{CAT}$ be the $\infty$-bicategory given by the nerve of the ordinary bicategory of categories, functors, and natural transformations. Then
\[ \Mod_{\bcat{W}}(\mathbf{CAT}) \]
is equivalent to the nerve of the ordinary bicategory of $\bcat{W}$-modules, $\bcat{W}$-functors, and $\bcat{W}$-natural transformations. 
\end{proposition}
\begin{proof}
A monoidal equivalence $\bcat{W} \weq \bcat{W}'$ induces an equivalence of $\infty$-bicategories $\un{\bcat{W}} \weq \un{\bcat{W}'}$ and hence an equivalence
\[ \Mod_{\bcat{W}'} \weq \Mod_{\bcat{W}}. \]
That monoidal equivalence also induces an equivalence between the ordinary bicategories of $\bcat{W}$-modules and of $\bcat{W}'$-modules. Thus, it is sufficient to consider the case that $\bcat{W}$ is strict monoidal.

In this case, the one-object $\infty$-bicategory $\un{\bcat{W}}$ described in Example~\ref{ex:moncat-bicat} is isomorphic to the simplicial nerve of the one-object $2$-category $\mathsf{B}\bcat{W}$ corresponding to the strict monoidal category $\bcat{W}$. It therefore follows from Proposition~\ref{prop:bicat-nerve} that
\[ \Mod_{\bcat{W}}(\mathbf{CAT}) = \Fun(\un{\bcat{W}},\mathbf{CAT}) \homeq N\mathsf{Bicat}(\mathsf{B}\bcat{W},\mathsf{CAT}). \]
It is then sufficient to identify pseudofunctors $\mathsf{B}\bcat{W} \to \mathsf{CAT}$ with $\bcat{W}$-actegories, natural transformations between such pseudofunctors with $\bcat{W}$-functors, and modifications with $\bcat{W}$-natural transformations. That identification follows directly by comparing the definitions in, say,~\cite{johnson/yau:2021} with those in~\cite{capucci/gavranovic:2022}.
\end{proof}

Now let $\bcat{W}$ be a monoidal $\infty$-category, and let $h\bcat{W}$ be its homotopy category with the induced monoidal structure; see Proposition~\ref{prop:mon-nerve}. We next show that $\bcat{W}$-module structures on an ordinary category $\bcat{X}$ are completely determined by an action of $h\bcat{W}$. 

\begin{proposition} \label{prop:homotopy-monoidal-module}
Let $\bcat{W}$ be a monoidal $\infty$-category. Then the monoidal functor $n: \bcat{W} \to h\bcat{W}$ induces an equivalence of $\infty$-bicategories
\[ \Mod_{h\bcat{W}}(\mathbf{CAT}) \weq \Mod_{\bcat{W}}(\mathbf{CAT}). \]
\end{proposition}
\begin{proof}
For the same reason as in the proof of~\ref{prop:nerve-mod-bicat}, it is sufficient to consider the case where $\bcat{W}$ is a strict monoidal $\infty$-category. The associated one-object $\infty$-bicategory is the simplicial nerve of the one-object simplicial category with mapping object $\bcat{W}$. The homotopy $2$-category of an $\infty$-bicategory is given by applying the homotopy category functor to mapping objects. We therefore have
\[ h_2\un{\bcat{W}} \isom \un{h\bcat{W}} \]
where the right-hand side is the one-object $2$-category whose mapping object is $h\bcat{W}$ with its induced strict monoidal structure. Moreover, we can identify the functor of $\infty$-bicategories
\[ \un{\bcat{W}} \to \un{h\bcat{W}}, \]
induced by $n$ with the canonical map
\[ \un{\bcat{W}} \to h_2\un{\bcat{W}} \]
from an $\infty$-bicategory to (the nerve of) its homotopy bicategory. Since $\mathbf{CAT}$ is (the nerve of) a $2$-category, it follows from Proposition~\ref{prop:homotopy-bicat} that the induced map
\[ \Fun(\un{h\bcat{W}},\mathbf{CAT}) \to \Fun(\un{\bcat{W}},\mathbf{CAT}) \]
is an equivalence, which yields the desired result.
\end{proof}

It follows from Proposition~\ref{prop:homotopy-monoidal-module} that if $\bcat{X},\bcat{Y}$ are ordinary categories with $\bcat{W}$-module structure, then the $\infty$-category of $\bcat{W}$-functors $\bcat{X} \to \bcat{Y}$ is equivalent to the ordinary category of $h\bcat{W}$-functors $\bcat{X} \to \bcat{Y}$ and $h\bcat{W}$-natural transformations. In fact, a similar result holds when we assume only that $\bcat{Y}$ is an ordinary category.

\begin{proposition} \label{prop:Y-module-functor}
Let $\bcat{X},\bcat{Y}$ be $\bcat{W}$-module $\infty$-categories such that $\bcat{Y}$ is an ordinary category. Then, taking homotopy categories determines an equivalence of $\infty$-categories
\[ \Mod_{\bcat{W}}(\bcat{X},\bcat{Y}) \homeq \Mod_{h\bcat{W}}(h\bcat{X},\bcat{Y}). \]
\end{proposition}
\begin{proof}
Without loss of generality we can assume that $\bcat{W}$ is strict monoidal, and $\bcat{X},\bcat{Y}$ are strict $\bcat{W}$-modules. We can then apply the description of $\bcat{W}$-functors in Remark~\ref{rem:bar}. The $n$-simplexes in the $\infty$-category $\Mod_{\bcat{W}}(\bcat{X},\bcat{Y})$ can be identified with $\bcat{W}$-functors
\[ \bcat{X} \to \Fun(\Delta^n,\bcat{Y}) \]
where $\Fun(\Delta^n,\bcat{Y})$ is a strict $\bcat{W}$-module with the pointwise $\bcat{W}$-action on $\bcat{Y}$. By Remark~\ref{rem:bar}, these $\bcat{W}$-functors are given by sequences of functors
\[ \Delta^m \times \bcat{W}^m \times \bcat{X} \to \Fun(\Delta^n,\bcat{Y}) \]
satisfying certain conditions. Since $\bcat{Y}$ is an ordinary category, each such functor of $\infty$-categories corresponds to an ordinary functor
\[ \Delta^m \times h\bcat{W}^m \times h\bcat{X} \isom h(\Delta^m \times \bcat{W}^m \times \bcat{X})\to \Fun(\Delta^n,\bcat{Y}) \]
satisfying corresponding conditions. By Remark~\ref{rem:bar} again, this data corresponds to an $n$-simplex in the $\infty$-category $\Mod_{h\bcat{W}}(h\bcat{X},\bcat{Y})$.
\end{proof}

\chapter{Weil-Algebras and $E_\infty$-Semirings} \label{sec:e-infty}

The goal of this section is to prove Theorem~\ref{thm:comparison} by showing that the monoidal $\infty$-category $\Weilinfty$ of Proposition~\ref{prop:weil-tensor} is (monoidally) equivalent to the full subcategory of $\ERig$ on the `$E_\infty$-Weil-algebras' of Definition~\ref{def:infty-Weil}, thus establishing the equivalence between our two descriptions of the monoidal $\infty$-category of Weil-algebras.

In this section, we construct a fully faithful symmetric monoidal functor
\[ \sigma: \Weilinfty \to \ERig. \]
For a Weil-algebra $A$, the object $\sigma(A)$ is the $E_\infty$-semiring introduced in Example~\ref{ex:e-infty} and there denoted $\N_\infty^{M_A}$, where $M_A$ is the partial commutative monoid of nonzero monomials in $A$. In the course of defining the functor $\sigma$, we give a precise construction of the $E_\infty$-semiring $\N_\infty^M$ for any finite partial commutative monoid $M$.

To construct $\sigma$, recall the characterization of the $\infty$-category $\ERig$ given by Gepner, Groth, and Nikolaus in~\cite{gepner/groth/nikolaus:2015}. They show that the $\infty$-category $\EMon$ of $E_\infty$-monoids has a canonical symmetric monoidal structure $\otimes$, and they define\index{ERig@$\ERig$, the $\infty$-category of $E_\infty$-semirings}
\[ \ERig := \Alg_{\E}(\EMon) \]
to be the $\infty$-category of $E_\infty$-algebras for that symmetric monoidal $\infty$-category.

Cranch~\cite{cranch:2010} showed that the $\infty$-category $\EMon$ admits an equivalence
\[ \EMon \homeq \Fun^{\times}(\SpanFin^{op},\spaces) \]
where $\SpanFin$\index{SpanFin@$\SpanFin$, the $\infty$-category of spans of finite sets (also denoted $\E$)} is the $\infty$-category of spans of finite sets. That $\infty$-category already appeared in Chapter~\ref{sec:differential} of this paper, where we denote it $\E$; see Definition~\ref{def:E} for a precise construction. We use different notation in this section to reflect the different role that this $\infty$-category is playing.

Cnossen, Haugseng, Lenz, and Linskens prove in~\cite[3.3.5]{cnossen/haugseng/lenz/linskens:2024} that the Gepner-Groth-Nikolaus tensor product on $\EMon$ is equivalent to (a localization of) the Day convolution symmetric monoidal structure on the functor $\infty$-category $\Fun(\SpanFin,\spaces)$ with respect to the cartesian product on $\spaces$ and the symmetric monoidal structure $\otimes$ on $\SpanFin$ given by cartesian product of finite sets. Note that this symmetric monoidal structure is not the categorical product in $\SpanFin$, which is given by the disjoint union of finite sets.

It then follows, for example, by work of Glasman~\cite{glasman:2016}, that we have a fully faithful symmetric monoidal embedding
\[ \SpanFin \to \Fun^{\times}(\SpanFin^{op},\spaces) \homeq \EMon \]
given by the Yoneda embedding for $\SpanFin$, and hence a fully-faithful symmetric monoidal embedding
\[ \Alg_{\E}(\SpanFin) \to \ERig. \]
The image of this embedding consists of those $E_\infty$-semirings for which the underlying $E_\infty$-monoid is free on a finite set, which includes the $E_\infty$-Weil-algebras, which are in the image of our proposed functor $\sigma: \Weilinfty \to \ERig$. 

It is therefore sufficient for us to construct a fully faithful (and symmetric monoidal) embedding
\[ \sigma: \Weilinfty \to \Alg_{\E}(\SpanFin), \]
and the rest of this section is devoted to that goal.

Our first step is to describe an explicit model for $\Alg_{\E}(\SpanFin)$. In Lurie's approach, a symmetric monoidal structure on an $\infty$-category $\cat{C}$ is represented as a certain type of cocartesian fibration of simplicial sets
\[ q: \cat{C}^{\otimes} \to \Fin_*, \]
over the category of finite pointed sets, whose fibre over the object $\langle 1 \rangle = \{1\}_+$ is equivalent to $\cat{C}$; see~\cite[2.1.2.18]{lurie:2017}. The $\infty$-category of $E_\infty$-algebras in $\cat{C}$ is then given by the simplicial set
\[ \Alg_{\E}(\cat{C}^{\otimes}) := \Fun^{\mathrm{inert}}_{\Fin_*}(\Fin_*,\cat{C}^{\otimes}) \]
of sections of $q$ which map inert morphisms in $\Fin_*$ to $q$-cocartesian morphisms in $\cat{C}^{\otimes}$. To construct $\sigma$ we therefore have to give a map of $\infty$-categories of the form
\[ \sigma: \Weilinfty \to \Fun^{\mathrm{inert}}_{\Fin_*}(\Fin_*,\mathrm{SpanFin}^{\otimes}). \]
We start by giving an explicit model, due to Barwick, Glasman, and Shah~\cite[2.6]{barwick/glasman/shah:2020} for the cocartesian fibration $q: \mathrm{SpanFin}^{\otimes} \to \Fin_*$ that represents the symmetric monoidal structure on $\mathrm{SpanFin}$ given by the cartesian product of finite sets.

\begin{definition} \label{def:SpanFin}
Let $\SpanFin^{\otimes}$\index{SpanFino@$\SpanFin^{\otimes}$} be the $\infty$-category given by the Duskin nerve of the bicategory with the following data.
\begin{itemize} \itemsep=10pt
\item An object of $\SpanFin^{\otimes}$ is a finite set $I$ and a collection $(X_i)_{i \in I}$ of finite sets $X_i$ indexed by $I$.
\item A $1$-morphism of $\SpanFin^{\otimes}$ from $(X_i)_{i \in I}$ to $(Y_j)_{j \in J}$ is a morphism $f: I_+ \to J_+$ in $\Fin_*$ together with, for each $j \in J$, a span of finite sets of the form
\[ \begin{diagram} \dgARROWLENGTH=1em
\node[2]{Z_j} \arrow{sw} \arrow{se} \\
\node {\prod_{i \in f^{-1}(j)} X_i} \node[2]{Y_j}
\end{diagram} \]
\item A $2$-morphism exists only between $1$-morphisms with the same underlying map $f: S_+ \to T_+$ and consists of a collection of isomorphisms of spans. Composition and identities for $2$-morphisms are the usual composition and identities for isomorphisms of spans.
\item The identity $1$-morphism on $(X_i)_{i \in I}$ consists of the identity function on $I_+$ together with the identity spans on each $X_i$.
\item The composition of $1$-morphisms in $\SpanFin^{\otimes}$ is given by combining the cartesian product and pullbacks of finite sets in the following way. Given $1$-morphisms $(X_i)_{i \in I} \to (Y_j)_{j \in J}$ and $(Y_j)_{j \in J} \to (Z_k)_{k \in K}$, we compose the underlying functions in $\Fin_*$ and then form the spans
\[ \begin{diagram} \dgARROWLENGTH=1em
\node[3]{L''_k} \arrow{sw} \arrow{se} \\
\node[2]{\prod_{j \in g^{-1}(k)} L_j} \arrow{sw} \arrow{se} \node[2]{L'_k} \arrow{sw} \arrow{se} \\
\node{\prod_{i \in f^{-1}g^{-1}(k)} X_i} \node[2]{\prod_{j \in g^{-1}(k)} Y_j} \node[2]{Z_k}
\end{diagram} \]
for $k \in K$, where the central square is a pullback of finite sets, and the bottom-left part of the diagram is a product, over $j \in g^{-1}(k)$, of the spans appearing in the morphism from $(X_i)_{i \in I}$ to $(Y_j)_{j \in J}$.
\item The associator and unit isomorphisms are the unique maps induced by the universal property of the pullbacks, which define composition.
\end{itemize}
\end{definition}

\begin{lemma} \label{lem:SpanFin}
The forgetful functor $q: \SpanFin^{\otimes} \to \Fin_*$ is a cocartesian fibration, in which a morphism is $q$-cocartesian if and only if each of the legs of the spans it comprises is a bijection. The map $q$ encodes a symmetric monoidal structure on $\SpanFin$ given by the cartesian product of finite sets.
\end{lemma}
\begin{proof}
For any $\infty$-category $C$ with finite products, Barwick, Glasman, and Shah construct a cocartesian fibration~\cite[2.14]{barwick/glasman/shah:2020} that represents the induced symmetric monoidal structure on $\mathrm{Span}(C)$ (which they denote as $A^{\mathrm{eff}}(C)^{\otimes}$). In our case, $C$ is the ordinary category $\Fin$, and their construction then reduces to that of Definition~\ref{def:SpanFin}. 
\end{proof}

The $\infty$-category $\SpanFin^{\otimes}$ is the nerve of a bicategory and so, by Proposition~\ref{prop:bicat-nerve}, the $\infty$-category
\[ \Fun^{\mathrm{inert}}_{\Fin_*}(\Fin_*,\SpanFin^{\otimes}) \]
is the nerve of the bicategory of strictly unitary pseudofunctors $\Fin_* \to \SpanFin^{\otimes}$. Moreover, since $\Weilinfty$ is also the nerve of a bicategory, the desired map $\sigma$ also corresponds to a strictly unitary pseudofunctor, and we will define it in those terms. See~\cite[4.1.2]{johnson/yau:2021} for the definition of pseudofunctor we use here.

We start by defining $\sigma$ on objects. Recall that an object of $\Weilinfty$ is the partial monoid of nonzero monomials in a Weil-algebra. In fact, we can define $\sigma$ on arbitrary finite partial commutative monoids.

\begin{definition} \label{def:sigma(M)}
Let $M$ be a finite partial commutative monoid. We define a pseudofunctor $\sigma(M): \Fin_* \to \SpanFin^{\otimes}$\index{sigmaM@$\sigma(M)$, the pseudofunctor $\Fin_* \to \SpanFin^{\otimes}$ associated to a finite partial commutative monoid $M$} as follows:
\begin{itemize} \itemsep=10pt
\item For each object $I_+ \in \Fin_*$ we take
\[ \sigma(M)(I_+) := (M)_{i \in I}, \]
the constant $I$-tuple of finite sets all equal to $M$.

\item For $f: I_+ \to J_+$, we set $\sigma(M)(f)$ to be the morphism in $\SpanFin^{\otimes}$ from $(M)_{i \in I}$ to $(M)_{j \in J}$ consisting of the spans of finite sets
\[ \begin{diagram}
\node[2]{M_{f^{-1}(j)}} \arrow{sw,V} \arrow{se,t}{\mu_M} \\
\node{M^{f^{-1}(j)}} \node[2]{M}
\end{diagram} \]
for $j \in J$, where for any finite set $K$ we write $M_{K}$ for the subset of the $K$-fold product $M^K$ consisting of those tuples $(m_k)_{k \in K}$ such that the product $\prod_{k \in K} m_k$ is defined in the partial monoid structure on $M$. We write $\mu_M$ for the function $M_K \to M$ that forms that product. Note that $\sigma(M)(1_{I_+})$ is precisely the identity morphism on the object $(M)_{i \in I}$ in $\SpanFin^{\otimes}$, so that $\sigma(M)$ is strictly unital.

\item It remains to give the lax functoriality constraint for $\sigma(M)$. For morphisms
\[ I_+ \arrow{e,t}{f} J_+ \arrow{e,t}{g} K_+ \]
in $\Fin_*$, we have to define a $2$-isomorphism in $\SpanFin^{\otimes}$ of the form
\[ \sigma(M)(gf) \isom \sigma(M)(g) \circ \sigma(M)(f). \]
That $2$-isomorphism amounts to a collection of bijections of the following form, for each $k \in K$
\[  M_{f^{-1}g^{-1}(k)} \isom \left( \prod_{j \in g^{-1}(k)} M_{f^{-1}(j)} \right) \times_{ M^{g^{-1}(k)}} M_{g^{-1}(k)}, \]
which commute with the maps to $M^{f^{-1}g^{-1}(k)}$ and $M$. Equivalently, we must choose dashed maps making the following diagram commute and the central square a pullback:
\[\begin{tikzcd}
	{M^{f^{-1}g^{-1}(k)}} && {M_{f^{-1}g^{-1}(k)}} && M \\
	\\
	& {\prod_{j \in g^{-1}(k)}M_{f^{-1}(j)}} && {M_{g^{-1}(k)}} \\
	\\
	&& {M^{g^{-1}(k)}}
	\arrow[from=1-3, to=1-1]
	\arrow[from=1-3, to=1-5]
	\arrow["{(\sigma(M)_{g,f,j})_{j \in g^{-1}(k)}}"{description}, dashed, from=1-3, to=3-2]
	\arrow["{\sigma(M)_{g,f,k}}"{description}, dashed, from=1-3, to=3-4]
	\arrow[from=3-2, to=1-1]
	\arrow[from=3-2, to=5-3]
	\arrow[from=3-4, to=1-5]
	\arrow[from=3-4, to=5-3]
\end{tikzcd} \]
Notice that the three spans on the edges of this triangle are formed from the components of $\sigma(M)_f$, $\sigma(M)_g$, and $\sigma(M)_{gf}$.

The desired functions are given by
\[ \sigma(M)_{g,f,j}: M_{f^{-1}g^{-1}(k)} \to M_{f^{-1}(j)}; \quad (m_i)_{i \in f^{-1}g^{-1}(k)} \mapsto (m_i)_{i \in f^{-1}(j)} \]
and
\[ \sigma(M)_{g,f,k}: M_{f^{-1}g^{-1}(k)} \to M_{g^{-1}(k)}; \quad (m_i)_{i \in f^{-1}g^{-1}(k)} \mapsto \left( \prod_{i \in f^{-1}(j)} m_i \right)_{j \in g^{-1}(k)}. \]
That these choices make the central square in the diagram commute corresponds to the following property of a partial commutative monoid $M$: for any $I$-tuple $(m_i)_{i \in I}$ of elements in $M$, and any partition of $I$ into subsets $(I_j)_{j \in J}$, the product $\prod_{i \in I} m_i$ is defined in $M$ if and only if each of the products $\prod_{i \in I_j} m_i$ is defined for $j \in J$ \emph{and} the product
\[ \prod_{j \in J} \left( \prod_{i \in I_j} m_i \right) \]
is defined in $M$.
\end{itemize}

To ensure that $\sigma(M)$ defines a pseudofunctor $\Fin_* \to \SpanFin^{\otimes}$, we have to check compatibility of the data above with respect to the associativity constraints. We do that as part of the proof of the following lemma.
\end{definition}

\begin{lemma} \label{lem:sigma(M)}
The data in Definition~\ref{def:sigma(M)} form an object
\[ \sigma(M) \in \Fun^{\mathrm{inert}}_{\Fin_*}(\Fin_*,\SpanFin^{\otimes}). \]
\end{lemma}
\begin{proof}
We have to check the lax unit and lax associativity conditions of~\cite[4.1.2]{johnson/yau:2021}. The unit conditions say that when we take either $f$ or $g$ to be an identity morphism in $\Fin_*$, then the chosen isomorphism
\[ \sigma(M)_{gf} \isom \sigma(M)_g \circ \sigma(M)_f \]
reduces to the canonical bijections of the form
\[ M_{I} \isom M_{I} \times_M M, \quad \text{and} \quad M_{J} \isom M^{J} \times_{M^{J}} M_{J}, \]
which are straightforward to check.

The lax associativity condition says that for any maps
\[ I_+ \arrow{e,t}{f} J_+ \arrow{e,t}{g} K_+ \arrow{e,t}{h} L_+ \]
in $\Fin_*$ and for any $\ell \in L$, the following diagrams commute:
\[\begin{tikzcd}
	{\displaystyle M_{f^{-1}g^{-1}h^{-1}(\ell)}} && {\displaystyle \prod_{k \in h^{-1}(\ell)} M_{f^{-1}g^{-1}(k)}} \\
	\\
	& {\displaystyle \prod_{j \in g^{-1}h^{-1}(\ell)} M_{f^{-1}(j)}}
	\arrow["{(\sigma(M)_{h,gf,k})_{k \in h^{-1}(\ell)}}", from=1-1, to=1-3]
	\arrow["{(\sigma(M)_{hg,f,j})_{j \in g^{-1}h^{-1}(\ell)}}"', from=1-1, to=3-2]
	\arrow["{\prod_{k \in h^{-1}(\ell)}(\sigma(M)_{g,f,j})_{j \in g^{-1}(k)}}", from=1-3, to=3-2]
\end{tikzcd} \]
and
\[\begin{tikzcd}
	{M_{f^{-1}g^{-1}h^{-1}(\ell)}} && {M_{g^{-1}h^{-1}(\ell)}} \\
	& {M_{h^{-1}(\ell)}}
	\arrow["{\sigma(M)_{hg,f,\ell}}", from=1-1, to=1-3]
	\arrow["{\sigma(M)_{h,gf,\ell}}"', from=1-1, to=2-2]
	\arrow["{\sigma(M)_{h,g,\ell}}", from=1-3, to=2-2]
\end{tikzcd}\]
For the first diagram, all the maps involved can be identified with inclusions of subsets of $M^{f^{-1}g^{-1}h^{-1}(\ell)}$, and the second commutes by associativity of the partial multiplication on $M$. We therefore have a pseudofunctor $\sigma(M): \Fin_* \to \SpanFin^{\otimes}$ as claimed.

It is immediate from the definition that $\sigma(M)$ is a section of the cocartesian fibration $q$, and if $f: I_+ \to J_+$ is an inert map, then every span of finite sets in the morphism $\sigma(M)(f)$ is the identity span on $M$, so $\sigma(M)(f)$ is $q$-cocartesian. We therefore have $\sigma(M) \in \Fun^{\mathrm{inert}}_{\Fin_*}(\Fin_*,\SpanFin^{\otimes})$.
\end{proof}

\begin{remark} \label{rem:sigma(M)}
The object $\sigma(M)$ is an $E_\infty$-algebra in $\SpanFin^{\otimes}$. The underlying object of that algebra is the finite set $M$, and the multiplication operation is the span of finite sets
\[ \begin{diagram}
\node[2]{M_2} \arrow{sw,V} \arrow{se,t}{\mu_M} \\
\node{M \times M} \node[2]{M}
\end{diagram} \]
where $M_2 = \{(m_1,m_2) \in M \times M \; | \; m_1 \cdot m_2 \in M\}$ is the set of pairs of elements in $M$ for which the product is defined in the partial monoid structure.

An equivalent construction of $\sigma(M)$ appears in work of Contreras, Mehta, and Stern in~\cite{contreras/mehta/stern:2025}. They show that commutative pseudomonoids in the symmetric monoidal bicategory of spans of finite sets (which are essentially the same as $E_\infty$-algebras in the symmetric monoidal $\infty$-category $\SpanFin$) correspond to $2$-Segal $\Gamma$-sets, i.e.\ functors $\Fin_* \to \Fin$ satisfying $2$-Segal conditions analogous to those of Dyckerhoff and Kapranov~\cite{dyckerhoff/kapranov:2019}. It is therefore sufficient to show that a finite partial commutative monoid $M$ gives rise to a $2$-Segal $\Gamma$-set in this sense, which Contreras et al.\ do in~\cite[5.10]{contreras/mehta/stern:2025}.
\end{remark}

Now we turn to the definition of $\sigma$ on morphisms.

\begin{definition} \label{def:sigma(phi)}
Take a morphism $\phi: M \to M'$ in $\Weilinfty$, i.e.\ a span of finite partial commutative monoids
\begin{equation} \label{eq:sigma(phi)} \begin{diagram}
\node[2]{L} \arrow{sw,t}{s} \arrow{se,t}{t} \\
\node{M} \node[2]{M'}
\end{diagram} \end{equation}
satisfying the conditions of Definition~\ref{def:weil-mor}. We define a natural transformation $\sigma(\phi): \sigma(M) \to \sigma(M')$\index{sigmaphi@$\sigma(\phi)$, the natural transformation associated to a span $\phi$ of finite partial commutative monoids} between these functors $\Fin_* \to \SpanFin^{\otimes}$ as follows.
\begin{itemize} \itemsep=10pt
\item For $I_+ \in \Fin_*$, the component $\sigma(\phi)_{I_+}$ of that natural transformation is the morphism $(M)_{i \in I} \to (M')_{i \in I}$ in $\SpanFin^{\otimes}$ consisting of copies of the diagram (\ref{eq:sigma(phi)}) indexed by $I$.
\item For $f: I_+ \to J_+$, we have a naturality constraint for $\sigma(\phi)$ which is a square diagram in $\SpanFin^{\otimes}$ whose commutativity is witnessed by a $2$-isomorphism. That square consists of a collection of diagrams of finite sets of the form
\[\begin{tikzcd}
	{M^{I'}} && {M_{I'}} & M \\
	&& {M_{I'} \times _M L} & L \\
	{L^{I'}} & {L^{I'} \times_{(M')^{I'}} M'_{I'}} \\
	{(M')^{I'}} & {M'_{I'}} && {M'}
	\arrow[tail, from=1-3, to=1-1]
	\arrow["{\mu_M}", from=1-3, to=1-4]
	\arrow[from=2-3, to=1-3]
	\arrow["\lrcorner"{anchor=center, pos=0.125, rotate=90}, draw=none, from=2-3, to=1-4]
	\arrow[from=2-3, to=2-4]
	\arrow["s"', from=2-4, to=1-4]
	\arrow["t", from=2-4, to=4-4]
	\arrow["{s^{I'}}", from=3-1, to=1-1]
	\arrow["{t^{I'}}"', from=3-1, to=4-1]
	\arrow["{\sigma(\phi)_{f,j}}"{description}, from=3-2, to=2-3]
	\arrow[tail, from=3-2, to=3-1]
	\arrow["\lrcorner"{anchor=center, pos=0.125, rotate=-90}, draw=none, from=3-2, to=4-1]
	\arrow[from=3-2, to=4-2]
	\arrow[tail, from=4-2, to=4-1]
	\arrow["{\mu_{M'}}"', from=4-2, to=4-4]
\end{tikzcd}\]
for $j \in J$, where $I' := f^{-1}(j)$. The bottom-left and top-right squares are the pullbacks of finite sets which provide the composite spans in the square diagram, and the map labelled $\sigma(\phi)_{f,j}$ is the required isomorphism between those composite spans, which forms the $2$-morphism in $\SpanFin$ that witnesses the commutativity.

The pullback of the bottom-left square can be identified with the set of $I'$-tuples in $L$ which, after applying $t$, give an $I'$-tuple in $M'$ for which the product exists. Because the partial monoid homomorphism $t:L \to M'$ satisfies condition (2) of Definition~\ref{def:weil-mor}, that set is equal to $L_{I'}$. (That condition covers explicitly the case where $|I'| = 2$, but other cases are either trivial (for $|I'| = 0,1$) or follow by induction.)

We then define $\sigma(\phi)_{f,j}$ to consist of the map $L_{I'} \to M_{I'}$ given by applying $s$ elementwise and the map $L_{I'} \to L$ given by the partial monoid structure on $L$. The resulting diagram commutes because $s$ and $t$ are both partial monoid homomorphisms. Finally, the map $\sigma(\phi)_{f,j}$ is a bijection because $s: L \to M$ satisfies condition (1) of Definition~\ref{def:weil-mor}. 
\end{itemize}
\end{definition}

\begin{lemma} \label{lem:sigma(phi)}
The data of Definition~\ref{def:sigma(phi)} form a natural transformation $\sigma(\phi): \sigma(M) \to \sigma(M')$ between these functors $\Fin_* \to \SpanFin^{\otimes}$.
\end{lemma}
\begin{proof}
For each pair of morphisms $I_+ \arrow{e,t}{f} J_+ \arrow{e,t}{g} K_+$ in $\Fin_*$, we have to check the lax naturality constraint~\cite[4.2.3]{johnson/yau:2021} relating the $2$-isomorphisms $\sigma(\phi)_f$, $\sigma(\phi)_g$, and $\sigma(\phi)_{gf}$. To simplify the notation a bit, fix $k \in K$, and write $J' = g^{-1}(k)$ and $I' = f^{-1}(J')$. Then that constraint amounts to commutativity of the following diagram of finite sets
\[\begin{tikzcd}
	{L^{I'} \times_{(M')^{I'}} M'_{I'}} && {L^{I'} \times_{(M')^{I'}} (\prod_{j \in J'} M'_{f^{-1}(j)}) \times_{(M')^{J'}} M'_{J'}} \\
	\\
	&& {(\prod_{j \in J'} M_{f^{-1}(j)}) \times_{M^{J'}} L^{J'} \times_{(M')^{J'}} M'_{J'}} \\
	\\
	{M_{I'} \times_{M} L} && {(\prod_{j \in J'} M_{f^{-1}(j)}) \times_{M^{J'}} M_{J'} \times_M L}
	\arrow["\cong", from=1-1, to=1-3]
	\arrow["{\sigma(\phi)_{gf,k}}"', from=1-1, to=5-1]
	\arrow["{(\prod_j \sigma(\phi)_{f,j}) \times_{(M')^{J'}} M'_{J'}}", from=1-3, to=3-3]
	\arrow["{(\prod_j M_{f^{-1}(j)}) \times_{M^{J'}} \sigma(\phi)_{g,1}}", from=3-3, to=5-3]
	\arrow["\cong", from=5-1, to=5-3]
\end{tikzcd}\]
This claim holds via a lengthy but easy diagram chase using associativity of the partial monoid structure on $L$. 
\end{proof}

One convenient feature of the construction of Definition~\ref{def:sigma(phi)} is that $\sigma$ strictly preserves both composition and identity $1$-morphisms.

\begin{lemma} \label{lem:sigma(alpha)}
For an identity $1$-morphism $1_M$ in $\Weilinfty$, $\sigma(1_M)$ is the identity natural transformation on $\sigma(M)$. For $1$-morphisms $\phi:M \to M'$ and $\phi':M' \to M''$ in $\Weilinfty$, there is an equality of natural transformations
\[ \sigma(\phi'\phi) = \sigma(\phi')\sigma(\phi). \]
\end{lemma}
\begin{proof}
The $1$-components $\sigma(1_M)_{I_+}$ are the identity on $\sigma(M)_{I_+} = (M)_{i \in I}$, and the $2$-components $\sigma(1_M)_f$ are made up of the obvious isomorphisms
\[ \sigma(1_M)_{f,j}: M^{I'} \times_{M^{I'}} M_{I'} \isom M_{I'} \isom M_{I'} \times_M M, \]
where $I' := f^{-1}(j)$. These choices agree with the components of the identity natural transformation on $\sigma(M)$ defined in~\cite[4.2.12]{johnson/yau:2021}. 

The composite of morphisms $\phi: M \to M'$ and $\phi': M' \to M''$ in $\Weilinfty$ is represented by a diagram of finite partial commutative monoids
\begin{equation} \label{eq:sigma(alpha)} \begin{diagram}
\node[3]{L''} \arrow{sw,t}{s''} \arrow{se,t}{t''} \\
\node[2]{L} \arrow{sw,t}{s} \arrow{se,t}{t} \node[2]{L'} \arrow{sw,t}{s'} \arrow{se,t}{t'} \\
\node{M} \node[2]{M'} \node[2]{M''}
\end{diagram} \end{equation}
where $L''$ is (a fixed choice of) the pullback $L \times_{M'} L'$. Composition of natural transformations $\sigma(M) \to \sigma(M'')$ is defined in~\cite[4.2.15]{johnson/yau:2021}, and equality between the $1$-components $\sigma(\phi'\phi)_{I_+}$ and $(\sigma(\phi')\sigma(\phi))_{I_+}$ amounts to the pullback of finite partial commutative monoids $L''$ using the same underlying pullback of finite sets as in the definition of composition in $\SpanFin$. Equality of $2$-components then amounts to commutativity of the diagrams, for each $f: I_+ \to J_+$ and each $j \in J$
\[\begin{tikzcd}
	{(L'')^{I'} \times_{(M'')^{I'}} M''_{I'}} & {L^{I'} \times_{(M')^{I'}} (L')^{I'} \times_{(M'')^{I'}} M''_{I'}} \\
	& {L^{I'} \times_{(M')^{I'}} M'_{I'} \times_M L'} \\
	{M_{I'} \times_M L''} & {M_{I'} \times_{M'} L \times_{M'} L'}
	\arrow[from=1-1, to=1-2]
	\arrow["{\sigma(\phi'\phi)_{f,j}}"', from=1-1, to=3-1]
	\arrow["{L^{I'} \times_{(M')^{I'}}\sigma(\phi')_{f,j}}", from=1-2, to=2-2]
	\arrow["{\sigma(\phi)_{f,j} \times_{M'} L'}", from=2-2, to=3-2]
	\arrow[from=3-1, to=3-2]
\end{tikzcd}\]
where $I' := f^{-1}(j)$, and the horizontal maps are induced by the equation $L'' = L \times_M L'$. Identifying an element of the top-left object as an $I'$-tuple in $L''$ whose product is defined in the partial monoid structure, it is a straightforward diagram chase to show that this diagram commutes.
\end{proof}

The final part of the definition of a pseudofunctor
\[ \sigma: \Weilinfty \to \Fun_{\Fin_*}(\Fin_*,\SpanFin^{\otimes}) \]
is to give the value of $\sigma$ on a $2$-morphism in the bicategory $\Weilinfty$. 

\begin{definition} \label{def:sigma(alpha)}
A $2$-isomorphism in the bicategory $\Weilinfty$ is an isomorphism $\alpha: \phi \to \phi'$ between spans of finite partial commutative monoids of the form
\[ \begin{diagram}
\node[2]{L} \arrow{sw,t}{s} \arrow{se,t}{t} \arrow[2]{s,lr}{\alpha}{\isom} \\
\node{M} \node[2]{M'} \\
\node[2]{L'} \arrow{nw,b}{s'} \arrow{ne,b}{t'}
\end{diagram} \]
We define a modification (see~\cite[4.4.1]{johnson/yau:2021}) $\sigma(\alpha): \sigma(\phi) \to \sigma(\phi')$ between the corresponding natural transformations. That modification consists of the following data:
\begin{itemize}
\item For each finite set $I_+$, a $2$-isomorphism in $\SpanFin^{\otimes}$ from $\sigma(\phi)_{I_+}$ to $\sigma(\phi')_{I_+}$ consisting of copies of $\alpha$ indexed by $I$.
\end{itemize}
The modification axioms require the commutativity, for each $f: I_+ \to J_+$ and each $j \in J$, of the following diagrams
\[ \begin{diagram}
\node{L^{I'} \times_{(M')^{I'}} M'_{I'}} \arrow{s,l}{\alpha^{I'} \times_{(M')^{I'}} M'_{I'}} \arrow{e,t}{\sigma(\phi)_{f,j}} \node{M_{I'} \times_M L} \arrow{s,r}{M_{I'} \times_M \alpha} \\
\node{(L')^{I'} \times_{(M')^{I'}} M'_{I'}}\arrow{e,t}{\sigma(\phi')_{f,j}} \node{M_{I'} \times_M L'}
\end{diagram} \]
which follows from the fact that $\alpha$ is a homomorphism of partial monoids.

That completes the definition of the $2$-morphism $\sigma(\alpha)$. It is straightforward to check that $\sigma$ preserves composition and identity $2$-morphisms. 
\end{definition}

We now have the following result.

\begin{proposition} \label{prop:weil-calg}
The constructions of Definitions~\ref{def:sigma(M)},~\ref{def:sigma(phi)}, and~\ref{def:sigma(alpha)} form a strict functor of bicategories\index{sigma@$\sigma$, the functor of bicategories $\Weilinfty \to \Alg_{E_\infty}(\SpanFin^{\otimes})$}
\[ \sigma: \Weilinfty \to \Fun^{\mathrm{inert}}_{\Fin_*}(\Fin_*,\SpanFin^{\otimes}). \]
\end{proposition}
\begin{proof}
All that remains is to check the lax associativity and lax unity conditions of~\cite[4.1.3, 4.1.4]{johnson/yau:2021}, which are significantly simplified by Lemma~\ref{lem:sigma(alpha)}. For the associativity condition take morphisms $\phi_i: M_i \to M_{i+1}$ in $\Weilinfty$ for $i = 0,1,2$. Let $L_{ij}$ denote the apex of the span of finite partial commutative monoids that describes the composite morphism from $M_i$ to $M_j$ for $0 \leq i < j \leq 3$. The required condition is the equality of two modifications
\[\begin{tikzcd}
	{(\sigma(\phi_2)\sigma(\phi_1))\sigma(\phi_0)} && {\sigma(\phi_2)(\sigma(\phi_1)\sigma(\phi_0))} \\
	{\sigma((\phi_2\phi_1)\phi_0)} && {\sigma(\phi_2(\phi_1\phi_0))}
	\arrow["{a_{\sigma(\phi_2),\sigma(\phi_1),\sigma(\phi_0)}}", curve={height=-18pt}, from=1-1, to=1-3]
	\arrow["{=}"', no head, from=1-1, to=2-1]
	\arrow["{=}", no head, from=1-3, to=2-3]
	\arrow["{\sigma(a_{\phi_2,\phi_1,\phi_0})}"', curve={height=18pt}, from=2-1, to=2-3]
\end{tikzcd}\]
given by the associator isomorphisms in $\Weilinfty$ and $\Fun_{\Fin_*}(\Fin_*,\SpanFin^{\otimes})$, respectively. Each of these maps has components at $I_+ \in \Fin_*$ given by copies, indexed by $I$, of the canonical isomorphism
\[ L_{01} \times_{M_1} (L_{12} \times_{M_2} L_{23}) \isom L_{01} \times_{M_1} (L_{12} \times_{M_2} L_{23}). \]
Similarly, the modifications appearing in the lax unity conditions are given in each case by copies of the canonical isomorphisms
\[ M \times_M L \isom L \quad \text{and} \quad L \times_{M'} M' \isom L, \]
so they too are equal.
\end{proof}

Having constructed the desired map $\sigma$, we now complete the proof of Theorem~\ref{thm:comparison}, which is a consequence of the following result.

\begin{proposition} \label{prop:sigma-ff}
The functor $\sigma$ of Proposition~\ref{prop:weil-calg} is fully faithful.
\end{proposition}
\begin{proof}
Fix objects $M,M'$ in $\Weilinfty$. We wish to show that the induced map on mapping groupoids
\[ \sigma_{M,M'}: \Hom_{\Weilinfty}(M,M') \to \Hom_{\Alg_{\E}(\SpanFin)}(\sigma(M),\sigma(M')) \]
is an equivalence. We start by proving that $\sigma_{M,M'}$ is fully faithful. 

So let $\phi,\phi': M \to M'$ be morphisms in $\Weilinfty$, given by spans of finite partial commutative monoids
\[ \begin{diagram}
\node[2]{L} \arrow{sw,t}{s} \arrow{se,t}{t} \\
\node{M} \node[2]{M'} \\
\node[2]{L'} \arrow{nw,b}{s'} \arrow{ne,b}{t'}
\end{diagram} \]
We have associated morphisms
\[ \sigma(\phi),\sigma(\phi'): \sigma(M) \to \sigma(M') \]
in the bicategory of pseudofunctors $\Fin_* \to \SpanFin^{\otimes}$, and suppose given a $2$-morphism $\gamma: \sigma(\phi) \to \sigma(\phi')$, that is, a modification between these natural transformations. Thus $\gamma$ consists of isomorphisms of spans
\[ \gamma_{I_+}: \sigma(\phi)_{I_+} \to \sigma(\phi')_{I_+} \]
for each $I_+ \in \Fin_*$. So $\gamma_{I_+}$ consists of a collection of isomorphisms of spans of finite sets of the form
\[ \begin{diagram}
\node[2]{L} \arrow{sw,t}{s} \arrow{se,t}{t} \arrow[2]{s,lr}{\isom}{\gamma_i} \\
\node{M} \node[2]{M'} \\
\node[2]{L'} \arrow{nw,b}{s'} \arrow{ne,b}{t'}
\end{diagram} \]
indexed by $i \in I$. These choices must satisfy modification axioms with respect to $f: I_+ \to J_+$, which amount to commutative diagrams of finite sets
\[ \begin{diagram} \dgARROWLENGTH=4em
\node{L^{I'} \times_{(M')^{I'}} M'_{I'}} \arrow{s,l}{(\gamma_i)_{i \in I'}} \arrow{e,t}{\sigma(\phi)_{f,j}} \node{M_{I'} \times_M L} \arrow{s,r}{\gamma_j} \\
\node{(L')^{I'} \times_{(M')^{I'}} M'_{I'}} \arrow{e,t}{\sigma(\phi')_{f,j}} \node{M_{I'} \times_M L'}
\end{diagram} \]
for all $j \in J$, where $I' := f^{-1}(j)$. Applying this property to the various inert maps $I_+ \to \langle 1 \rangle$, we deduce that all the $\gamma_i$ are necessarily equal to a single map $\alpha := \gamma_1: L \to L'$. Applying the property to the active maps $\langle 2 \rangle \to \langle 1 \rangle$ and $\langle 0 \rangle \to \langle 1 \rangle$, we deduce that $\alpha$ is an isomorphism of partial monoids, i.e.\ $\alpha$ is a $2$-morphism in the bicategory $\Weilinfty$. It is then straightforward to see that the assignment $\gamma \mapsto \alpha$ is inverse to the map induced by $\sigma_{M,M'}$ on those $2$-morphisms $\phi \to \phi'$ in $\Weilinfty$. So $\sigma_{M,M'}$ is fully faithful.

It remains to show that $\sigma_{M,M'}$ is essentially surjective on objects, so take an arbitrary natural transformation $\psi: \sigma(M) \to \sigma(M')$ between these functors $\Fin_* \to \SpanFin^{\otimes}$. We can unpack the structure of $\psi$ as follows.
\begin{itemize} \itemsep=10pt
\item For each $I_+$ in $\Fin_*$, we have a collection of spans of finite sets of the form
\begin{equation} \label{eq:Psi} \begin{diagram}
\node[2]{\Psi_i} \arrow{sw,t}{s_i} \arrow{se,t}{t_i} \\
\node{M} \node[2]{M'}
\end{diagram} \end{equation}
which make up the component $\psi_{I_+}$ of the natural transformation.
\item For each $f: I_+ \to J_+$, we have a naturality square witnessed by a collection of diagrams of finite sets of the form
\[\begin{tikzcd}
	{M^{f^{-1}(j)}} && {M_{f^{-1}(j)}} & M \\
	&& \bullet & {\Psi_j} \\
	{\prod_{i \in f^{-1}(j)} \Psi_i} & \bullet \\
	{(M')^{f^{-1}(j)}} & {M'_{f^{-1}(j)}} && {M'}
	\arrow[tail, from=1-3, to=1-1]
	\arrow["{\mu_M}", from=1-3, to=1-4]
	\arrow[from=2-3, to=1-3]
	\arrow["\lrcorner"{anchor=center, pos=0.125, rotate=90}, draw=none, from=2-3, to=1-4]
	\arrow[from=2-3, to=2-4]
	\arrow["{s_j}"', from=2-4, to=1-4]
	\arrow["{t_j}", from=2-4, to=4-4]
	\arrow["{\prod_i s_i}", from=3-1, to=1-1]
	\arrow["{\prod_i t_i}"', from=3-1, to=4-1]
	\arrow["{\psi_{f,j}}", from=3-2, to=2-3]
	\arrow[tail, from=3-2, to=3-1]
	\arrow["\lrcorner"{anchor=center, pos=0.125, rotate=-90}, draw=none, from=3-2, to=4-1]
	\arrow[from=3-2, to=4-2]
	\arrow[tail, from=4-2, to=4-1]
	\arrow["{\mu_{M'}}"', from=4-2, to=4-4]
\end{tikzcd}\]
for each $j \in J$, where the top-right and bottom-left squares are pullbacks. That is, we must have bijections
\[ \psi_{f,j}: (\prod_{i \in f^{-1}(j)} \Psi_i) \times_{(M')^{f^{-1}(j)}} M'_{f^{-1}(j)} \isom M_{f^{-1}(j)} \times_M \Psi_j \]
which make these diagrams commute. As for the maps $\sigma(\phi)_{f,j}$ in the proof of Lemma~\ref{lem:sigma(phi)}, these bijections must satisfy additional conditions related to composition in $\Fin_*$, and we will make use of those conditions later in the proof.
\end{itemize}

First, though, we prove that all the spans (\ref{eq:Psi}) are isomorphic. For any $I_+$ and $i \in I$, let $f_i: I_+ \to \langle 1 \rangle$ be the inert morphism with $f_i(i) = 1$. Then $\psi_{f_i,1}$ is a bijection
\[ \Psi_i \times_{M'} M' \isom M \times_M \Psi_1 \]
which gives us an isomorphism of spans $\psi_i: \Psi_i \isom \Psi_1$. Up to isomorphism we can then assume $\psi$ has the property that all the spans (\ref{eq:Psi}) making up $\psi_{I_+}$ are equal to a single span
\begin{equation} \label{eq:psi} \begin{diagram}
\node[2]{L} \arrow{sw,t}{s} \arrow{se,t}{t} \\
\node{M} \node[2]{M'}
\end{diagram} \end{equation}

We next show that (\ref{eq:psi}) underlies a morphism $\phi$ in $\Weilinfty$. First, we show that $L$ admits the structure of a partial commutative monoid. Consider the active map $\mu_2: \langle 2 \rangle \to \langle 1 \rangle$, for which we have an associated commutative diagram
\\
\begin{equation} \label{eq:cdot-lambda} \begin{tikzcd}
	{M \times M} && {M_2} & M \\
	&& {M_2 \times_M L} & L \\
	{L \times L} & {L_2} \\
	{M' \times M'} & {M'_2} && {M'}
	\arrow[tail, from=1-3, to=1-1]
	\arrow["{\mu_M}", from=1-3, to=1-4]
	\arrow[from=2-3, to=1-3]
	\arrow["\lrcorner"{anchor=center, pos=0.125, rotate=90}, draw=none, from=2-3, to=1-4]
	\arrow["{\pi_L}"', from=2-3, to=2-4]
	\arrow["s"', from=2-4, to=1-4]
	\arrow["t", from=2-4, to=4-4]
	\arrow["{s \times s}", from=3-1, to=1-1]
	\arrow["{t \times t}"', from=3-1, to=4-1]
	\arrow["{\psi_{\mu_2,1}}", from=3-2, to=2-3]
	\arrow[from=3-2, to=3-1]
	\arrow["\lrcorner"{anchor=center, pos=0.125, rotate=-90}, draw=none, from=3-2, to=4-1]
	\arrow[from=3-2, to=4-2]
	\arrow[tail, from=4-2, to=4-1]
	\arrow["{\mu_{M'}}"', from=4-2, to=4-4]
\end{tikzcd}\end{equation}
where the bottom-left and top-right squares are pullbacks, and the top and bottom rows are the spans forming the partial monoid operation for $M$ and $M'$, respectively.

Since monomorphisms are stable under pullback, we can assume that $L_2$ is a subset of $L \times L$, and we define a partial operation on the finite set $L$ by
\begin{equation} \label{eq:cdot-L} \ell_1 \cdot \ell_2 := \begin{cases} \pi_L \psi_{\mu_2,1}(\ell_1,\ell_2) & \text{if $(\ell_1,\ell_2) \in L_2$}; \\ \text{undefined} & \text{otherwise}. \end{cases} \end{equation}
We will show below that this operation makes $L$ into a finite partial commutative monoid. Assuming that fact for now, it follows from the commutativity of the diagram (\ref{eq:cdot-lambda}) that $s$ and $t$ are partial monoid homomorphisms. Since the bottom-left square is a pullback, we know that $\ell_1 \cdot \ell_2$ is defined in $L$ if and only if $t(\ell_1) \cdot t(\ell_2)$ is defined in $M'$, so $t$ satisfies condition (2) of Definition~\ref{def:weil-mor}.

For condition (1) of~\ref{def:weil-mor}, the top-right square of (\ref{eq:cdot-lambda}) being a pullback (and the map $\lambda$ being a bijection) gives us the required unique factorization property in the case $r = 2$. By induction, we get that property for all $r \geq 2$, and $r = 1$ is trivial. The case $r = 0$ remains, which is equivalent to showing that $s^{-1}(1_M) = \{1_L\}$. To prove that claim, consider the naturality square for $\psi$ corresponding to the unique map $\mu_0: \langle 0 \rangle \to \langle 1 \rangle$. This square takes the form
\[\begin{tikzcd}
	{*} && {*} & M \\
	& {} & {s^{-1}(1_M)} & L \\
	{*} & {L_0} \\
	{*} & {*} && {M'}
	\arrow[from=1-3, to=1-1]
	\arrow["{1_M}", from=1-3, to=1-4]
	\arrow[from=2-3, to=1-3]
	\arrow["\lrcorner"{anchor=center, pos=0.125, rotate=90}, draw=none, from=2-3, to=1-4]
	\arrow[tail, from=2-3, to=2-4]
	\arrow["s"', from=2-4, to=1-4]
	\arrow["t", from=2-4, to=4-4]
	\arrow[from=3-1, to=1-1]
	\arrow[from=3-1, to=4-1]
	\arrow["\cong", from=3-2, to=2-3]
	\arrow[from=3-2, to=3-1]
	\arrow["\lrcorner"{anchor=center, pos=0.125, rotate=-90}, draw=none, from=3-2, to=4-1]
	\arrow[from=3-2, to=4-2]
	\arrow[from=4-2, to=4-1]
	\arrow["{1_{M'}}"', from=4-2, to=4-4]
\end{tikzcd}\]
Since the bottom-left square is a pullback, $L_0$ has one element, and hence so does $s^{-1}(1_M)$, which completes the proof that $s$ satisfies condition (1) of Definition~\ref{def:weil-mor}.

We now return to the claim that our chosen partial operation (\ref{eq:cdot-L}) makes $L$ into a partial commutative monoid. Consider the compatibility conditions for the components $\psi_{f,j}$ with respect to pairs of morphisms
\[ I_+ \arrow{e,t}{f} J_+ \arrow{e,t}{g} K_+ \]
in $\Fin_*$. That compatibility requires that each of the following diagrams commute, where we fix $k \in K$ and write $J' := g^{-1}(k)$ and $I' := f^{-1}(J')$.
\begin{equation} \label{eq:psi-2} \begin{tikzcd}
	{L^{I'} \times_{(M')^{I'}} M'_{I'}} && {L^{I'} \times_{(M')^{I'}} (\prod_{j \in J'} M'_{f^{-1}(j)}) \times_{(M')^{J'}} M'_{J'}} \\
	\\
	&& {(\prod_{j \in J'} M_{f^{-1}(j)}) \times_{M^{J'}} L^{J'} \times_{(M')^{J'}} M'_{J'}} \\
	\\
	{M_{I'} \times_{M} L} && {(\prod_{j \in J'} M_{f^{-1}(j)}) \times_{M^{J'}} M_{J'} \times_M L}
	\arrow["\cong", from=1-1, to=1-3]
	\arrow["{\psi_{gf,k}}"', from=1-1, to=5-1]
	\arrow["{(\prod_j \psi_{f,j}) \times_{(M')^{J'}} M'_{J'}}", from=1-3, to=3-3]
	\arrow["{(\prod_j M_{f^{-1}(j)}) \times_{M^{J'}} \psi_{g,k}}", from=3-3, to=5-3]
	\arrow["\cong", from=5-1, to=5-3]
\end{tikzcd}\end{equation}

To see that the operation is commutative, we look at this diagram for the sequence
\[ \langle 2 \rangle \arrow{e,t}{\tau} \langle 2 \rangle \arrow{e,t}{\mu_2} \langle 1 \rangle \]
where $\tau$ permutes the non-basepoints and the composite is again $\mu_2$. In that case, (\ref{eq:psi-2}) reduces to
\[\begin{tikzcd}
	{L^2 \times_{(M')^2} M'_2} && \\
	&& {L^2 \times_{(M')^2} M'_2} \\
	{M_2 \times_M L}
	\arrow[from=1-1, to=2-3]
	\arrow["{\psi_{\mu_2,1}}"', from=1-1, to=3-1]
	\arrow["{\psi_{\mu_2,1}}", from=2-3, to=3-1]
\end{tikzcd}\]
where the top map permutes the factors in each of the terms $L^2$, $(M')^2$, $M'_2$. It follows that the binary operation for $L$, which is derived from $\psi_{\mu_2,1}$, is commutative.

For associativity the story is similar, but based on the diagram in $\Fin_*$ of the form
\[ \begin{diagram}
\node{\langle 3 \rangle} \arrow{e,t}{1(23)} \arrow{s,l}{(12)3} \arrow{se} \node{\langle 2 \rangle} \arrow{s,r}{\mu_2} \\
\node{\langle 2 \rangle} \arrow{e,t}{\mu_2} \node{\langle 1 \rangle}
\end{diagram} \]
where $1(23)$ and $(12)3$ denote active maps $\langle 3 \rangle \to \langle 2 \rangle$ based on the given pattern. For the unit property, we use the diagram
\[ \begin{diagram}
\node{\langle 1 \rangle} \arrow{e,t}{i_1} \arrow{se,b}{=} \node{\langle 2 \rangle} \arrow{s,r}{\mu_2} \\
\node[2]{\langle 1 \rangle}
\end{diagram} \]
where $i_1(1) = 1$, together with the diagram
\[ \begin{diagram}
\node{\langle 1 \rangle} \arrow{s} \arrow{e,t}{i_1} \node{\langle 2 \rangle} \arrow{s,r}{p_1} \\
\node{\langle 0 \rangle} \arrow{e,t}{\mu_0} \node{\langle 1 \rangle} 
\end{diagram} \]
that relates $i_1$ by inert morphisms to the map $\mu_0$ used previously to identify the identity element $1_L$. 

We have now shown that the span of finite sets (\ref{eq:psi}) underlies a morphism $\phi: M \to M'$ in $\Weilinfty$. We complete the proof of essential surjectivity of $\sigma_{M,M'}$ by showing that $\sigma(\phi)$ is equal to $\psi$. (Recall that we already replaced $\psi$ up to isomorphism in order to simplify our construction of $\phi$.)

We have already arranged for $\psi_{I_+} = \sigma(\phi)_{I_+}$, so it remains to show that for each $f: I_+ \to J_+$ and each $j \in J$, the two bijections
\[ \sigma(\phi)_{f,j}, \psi_{f,j}: L_{f^{-1}(j)} = L^{f^{-1}(j)} \times_{(M')^{f^{-1}(j)}} M'_{f^{-1}(j)} \isom M_{f^{-1}(j)} \times_M L \]
are equal. Let $P_{f,j}$ denote the statement that $\sigma(\phi)_{f,j} = \psi_{f,j}$.

By construction, $P_{f,j}$ holds when $f = \mu_2$ (each map gives the partial monoid operation on $L$), when $f$ is inert (when each map is essentially the identity on $L$), and when $f = \mu_0$ (when each map gives the identity element in $L$). To deduce $P_{f,j}$ for all other $f,j$, recall that the maps $\sigma(\phi)_{f,j}$ and $\psi_{f,j}$ both satisfy the compatibility conditions expressed in diagrams of the form (\ref{eq:psi-2}). We can deduce from those diagrams the following facts:
\begin{enumerate}
\item If $P_{g,k}$ and $P_{f,j}$ for all $j \in g^{-1}(k)$, then $P_{gf,k}$.
\item If $g$ is inert, then $P_{gf,f(j)}$ if and only if $P_{f,j}$.
\end{enumerate}
Since every morphism in $\Fin_*$ can be factored as an active morphism followed by an inert morphism, by (2) it suffices to prove $P_{f,j}$ when $f$ is active. Every active morphism can be factored into active morphisms for which each fibre has at most $2$ elements, so by (1) we need only prove $P_{f,j}$ when $|f^{-1}(j)| = r \leq 2$. In that case, we have a commutative diagram in $\Fin_*$ of the form
\[ \begin{diagram}
\node{I_+} \arrow{e,t}{f} \arrow{s} \node{J_+} \arrow{s,r}{j \mapsto 1} \\
\node{\langle r \rangle} \arrow{e,t}{\mu_r} \node{\langle 1 \rangle} 
\end{diagram} \]
in which the vertical maps are inert. Since we already know that $P_{\mu_r,1}$ holds, (1) and (2) now imply $P_{f,j}$, completing the proof that $\sigma(\phi) = \psi$.

Thus $\sigma_{M,M'}$ is essentially surjective, and so is an equivalence. That completes the proof that $\sigma$ is fully faithful. Since the symmetric monoidal structures on both $\Weilinfty$ and $\Alg_{\E}(\SpanFin)$ are given by the cartesian product of finite sets, $\sigma$ is also symmetric monoidal.
\end{proof}

\begin{remark} \label{rem:2Segal}
We noted in Remark~\ref{rem:sigma(M)} that Contreras, Mehta, and Stern have identified $E_\infty$-algebras in $\SpanFin$ with $2$-Segal $\Gamma$-sets, and also shown how a finite partial commutative monoid $M$ gives rise to such an object; see~\cite{contreras/mehta/stern:2025}. Proposition~\ref{prop:sigma-ff} extends their result to morphisms, showing that maps between these $E_\infty$-algebras correspond to certain spans of finite partial commutative monoids, namely those described by Definition~\ref{def:weil-mor}. We suspect that our proof extends to $2$-Segal $\Gamma$-sets that do not arise from partial commutative monoids, so that $\Alg_{\E}(\SpanFin)$ is equivalent to a certain $\infty$-category of spans of $2$-Segal $\Gamma$-sets whose component maps satisfy conditions analogous to those of~\ref{def:weil-mor}.

In fact, G\"{o}dicke~\cite{godicke:2024} has proved exactly this result, but for $A_\infty$-algebras (that is, those which are associative, but not necessarily commutative, up to higher coherent homotopy) in $\SpanFin$, which correspond to $2$-Segal (simplicial) sets. He shows that $\Alg_{A_\infty}(\SpanFin)$ is equivalent to an $\infty$-category of spans of $2$-Segal sets in which the component maps satisfy conditions introduced by G\'{a}lvez-Carrillo, Kock, and Tonks~\cite{galvezcarrillo/kock/tonks:2018}; see~\cite[2.2]{godicke:2024}. In fact, G\"{o}dicke goes much further, showing that for any $\infty$-category $\cat{C}$ with finite limits, $A_\infty$-algebras in $\mathrm{Span}(\cat{C})$ correspond to $2$-Segal objects in $\cat{C}$ and their spans. We conjecture that G\"{o}dicke's result extends to $E_\infty$-algebras, and perhaps even to algebras over other $\infty$-operads. 
\end{remark}

\backmatter


\bibliographystyle{amsalpha}
\bibliography{mcching}

\printindex

\end{document}